\documentclass[opre,nonblindrev]{informs4}
\OneAndAHalfSpacedXI
\renewcommand{\arraystretch}{0.8}
\usepackage{makecell}
\usepackage{bm}
\usepackage{tikz}
\usetikzlibrary{arrows.meta,decorations.pathreplacing}
\usepackage{endnotes}
\usepackage{booktabs}
\usepackage{bbm}
\usepackage{setspace}
\usepackage{float}
\usepackage{xcolor}
\let\footnote=\endnote

\usepackage{hyperref}
\makeatletter
\DeclareFontShape{OT1}{ptm}{m}{scit}{<->ssub * ptm/m/sc}{}
\makeatother

\usepackage{natbib}
 \bibpunct[, ]{(}{)}{,}{a}{}{,}%
 \def\bibfont{\small}%
\TheoremsNumberedThrough     
\ECRepeatTheorems

\EquationsNumberedThrough    

\usepackage{multirow}
\usepackage{enumitem}
\usepackage{fix-cm}
\usepackage{etoolbox}  
\makeatletter
\newcommand{\appendixtableofcontents}{%
  \section*{Appendices}
  \@starttoc{app}%
}

\let\originaladdcontentsline\addcontentsline

\newcommand{\appendixaddtocline}[3]{%
  \originaladdcontentsline{#1}{#2}{#3}
  \ifstrequal{#1}{toc}{%
    \ifstrequal{#2}{section}{%
      \phantomsection
      \originaladdcontentsline{app}{section}{#3}
    }{%
      \ifstrequal{#2}{subsection}{%
        \phantomsection
        \originaladdcontentsline{app}{subsection}{#3}%
      }{}%
    }%
  }{}%
}

\AtBeginEnvironment{APPENDICES}{%
  \let\addcontentsline\appendixaddtocline
}

\hypersetup{bookmarksdepth=subsection}
\makeatother
\MANUSCRIPTNO{}
\usepackage[ruled,vlined, noend,linesnumbered]{algorithm2e}

\newcommand{\NN}{\mathbf{N}} 
\newcommand{\MM}{\mathbf{M}} 
\newcommand{\vecx}{\mathbf{x}}

\newcommand{\vecc}{\mathbf{c}} 
\newcommand{\vece}{\mathbf{e}}

\newcommand{\vecv}{\mathbf{v}}

\newcommand{\vecq}{\mathbf{q}}

\newcommand{\1}{\mathbf{1}} 
\newcommand{\matc}{\mathbf{C}} 
\newcommand{\mate}{\mathbf{E}} 
\newcommand{\matz}{\mathbf{Z}} 
 
\newcommand{\Prb}{\mathbb{P}}

\newcommand{\EE}{\mathbb{E}} 
\newcommand{\RR}{\mathbb{R}}

\newcommand{\CFP}{\mathbf{x}^{\mathbf{c}}}

\newcommand{\cFPi}{x_i^{\vecc}}
\newcommand{\CijFPs}{\mathbf{x}^{\mathbf{c}^{+ij}}}
\newcommand{\CjFPs}{\mathbf{x}^{\mathbf{c}^{+j}}}
\newcommand{\CiFPs}{\mathbf{x}^{\mathbf{c}^{+i}}}

\newcommand{\fLP}{f^{\operatorname{LP}}} 
\newcommand{\fLPM}{f^{\operatorname{LP}}_{\MM}}
\newcommand{\fFP}{f^{\operatorname{FP}}}
\newcommand{\fDA}{f}
\newcommand{\fpi}{f^{\pi}} 

\newcommand{\fIND}{f^{\operatorname{IND}}}
\newcommand{\ZFP}{Z^{\operatorname{FP}}} 
\newcommand{\ZLP}{Z^{\operatorname{LP}}} 
\newcommand{\ZDA}{Z} 
\newcommand{\Zpi}{Z^{\pi}} 
 
\makeatletter
\AtBeginDocument{%
  \setlength{\abovedisplayskip}{4pt plus 2pt minus 2pt}%
  \setlength{\belowdisplayskip}{4pt plus 2pt minus 2pt}%
  \setlength{\abovedisplayshortskip}{2pt plus 1pt}%
  \setlength{\belowdisplayshortskip}{3pt plus 1pt}%
  \setlength{\textfloatsep}{8pt plus 2pt minus 2pt}%
  \setlength{\intextsep}{8pt plus 2pt minus 2pt}%
  \setlength{\floatsep}{8pt plus 2pt minus 2pt}%
}
\AtBeginDocument{\renewcommand\normalsize{\@setfontsize\normalsizeXI{11}{12.66}%
  \baselineskip18pt plus 1pt minus 1pt \@bls=18pt \@em=11pt\relax}%
  \let\AppendixFontSize\normalsize
  \normalsize}
\makeatother
\begin{document}




\TITLE{A Unified Algorithmic Framework for Dynamic Assortment Optimization under MNL Choice}

\ARTICLEAUTHORS{%
\AUTHOR{Shuo Sun}
\AFF{Department of Industrial and Systems Engineering, University of Wisconsin--Madison, Madison, WI 53706}
\AUTHOR{Rajan Udwani}
\AFF{Department of Industrial Engineering and Operations Research, UC Berkeley, Berkeley, CA 94720}
\AUTHOR{Zuo-Jun Max Shen}
\AFF{Faculty of Engineering and Faculty of Business and Economics, University of Hong Kong, Hong Kong
}
\AFF{Department of Industrial Engineering and Operations Research, UC Berkeley, Berkeley, CA 94720}

} 

\ABSTRACT{%
We consider assortment and inventory planning problems with dynamic stockout-based substitution effects, and without replenishment, in two different settings: (1) Customers can see all available products when they arrive, a typical scenario in physical stores. (2) The seller can choose to offer a subset of available products to each customer, which is more common on online platforms. Both settings are known to be computationally challenging, and the current approximation algorithms for the two settings are quite different. We develop a unified algorithmic framework under the Multinomial Logit (MNL) choice model for both settings. Our algorithms improve on the state-of-the-art algorithms in terms of approximation guarantees, runtime, and the ability to manage uncertainty in the total number of customers and handle more complex constraints. In the process, we establish various novel properties of dynamic assortment planning (for the MNL choice model) that may be useful more broadly.
}%


\maketitle

%

\section{Introduction}

Choosing which products to offer and how much inventory to stock for each product are central operational challenges in many industries, including retail and airlines. Extensive research has focused on \emph{static} assortment optimization, where an optimal assortment is selected for a single customer. This setting assumes that we can maintain a sufficiently large inventory for each product, which is suitable for scenarios where inventory is not a primary concern, such as for digital goods or services. In many settings, however, such as retail and the sharing economy, inventory is finite. As customers arrive sequentially, products may stock out, potentially resulting in substantial revenue losses.  
Consequently, beyond the selection of products to offer to customers, the seller must also determine the inventory level of each product. This joint decision problem is known as \emph{dynamic} assortment optimization. The problem is notoriously challenging because of stockout-driven substitution: customers choose among the products currently available, while their purchases determine which products remain available to future customers \citep{mahajan2001stocking, aouad2018greedy}. More specifically, demand for each product evolves as a \emph{stochastic consumption process} jointly determined by the customer arrival sequence, customer preferences as specified by a choice model, and the occurrence of stockouts. In fact, even efficiently evaluating the expected value of the total revenue given a starting inventory (and assortment policy) remains challenging for many choice models, including the Multinomial Logit (MNL) model. 

In this work, we study dynamic assortment optimization under cardinality and budget constraints. 
We consider two operating regimes. In the first, which is natural for brick-and-mortar retail, the seller chooses the initial inventory, and every arriving customer sees all products that remain in stock. The assortment therefore changes only through stockouts \citep{aouad2018greedy, aouad2019approximation}. We refer to this setting as dynamic assortment optimization without personalization (\textsc{da}). In the second regime, which is natural for online platforms, customers of different types arrive according to a possibly nonstationary distribution. The seller observes the type of each arriving customer and may offer a personalized subset of the products in stock \citep{bai2022coordinated}. We refer to this setting as dynamic assortment optimization with personalization (\textsc{dap}).

Existing algorithmic approaches for \textsc{da} and \textsc{dap} are fundamentally different. For \textsc{dap}, a natural linear programming (LP) relaxation provides a useful fluid approximation to the expected revenue of the stochastic consumption process when the total number of customers is known. \citet{bai2022coordinated} optimize a submodular surrogate of the LP objective, at the cost of a constant-factor loss in the approximation guarantee. However, this LP relaxation implicitly allows the seller to offer a subset of available products, making it too loose for \textsc{da}. Moreover, to the best of our knowledge, no analogous relaxation was known to provide a constant-factor approximation for \textsc{da}. Instead, existing algorithms for \textsc{da} work directly with the underlying stochastic consumption process: they use Monte Carlo sampling and dynamic programming to estimate total expected revenue, and combine these estimates with computationally intensive enumeration procedures to identify high-quality solutions under cardinality constraints \citep{aouad2018greedy, aouad2022stability}.

A further source of difficulty is uncertainty in the total number of customers, denoted by $T$. On modern e-commerce and other online platforms, the distribution of $T$ may exhibit high variability \citep{szpektor2011improving, aouad2022nonparametric,bai2023fluid} and need not satisfy the increasing failure rate (IFR) property. Designing algorithms that accommodate uncertain customer volumes, particularly non-IFR distributions of $T$, has remained an open problem. \begingroup Together, these challenges have led to a fragmented landscape in which the best-known guarantee depends on the operating regime, the distribution of $T$, and the type of inventory constraint. This raises a natural question: \emph{Can a unified and computationally efficient algorithmic framework provide improved guarantees for both operating regimes, including when $T$ is uncertain?} \endgroup

\begingroup
Our main contribution is a unified relax-and-transform approach for both \textsc{da} and \textsc{dap}. The approach provides several advantages: it is computationally efficient, achieves improved approximation guarantees, accommodates uncertainty in the total number of customers, and extends to more general constraints. It consists of two key steps:
(i) construct LP-based fluid objectives over the feasible integral initial inventories and compute an inventory vector that approximately maximizes this relaxed objective; and
(ii) transform the first-step inventory vector into an implementable inventory vector for \textsc{da}, or pair it with an online assortment policy for \textsc{dap}, such that the resulting solutions perform well under the original stochastic consumption process. Beyond the performance guarantees summarized in Section~\ref{sec:results}, our work develops a collection of technical results and new ideas that may be of broader interest beyond dynamic assortment optimization. We discuss these contributions in more detail in Section~\ref{sec:tech}. Finally, we conduct extensive numerical experiments to benchmark our algorithms against existing approximation methods. The numerical results demonstrate that the proposed algorithms are computationally efficient, scalable, and achieve higher expected revenue across a wide range of settings.
\endgroup

\subsection{Main Results}\label{sec:results}

\begingroup
We summarize our main approximation guarantees in Table~\ref{tab:main-results}. 
The table focuses on cardinality constraints; our results for budget constraints are discussed in the following paragraphs.

\begin{table}[htp]
    \centering
    \small
    \setcellgapes{3.6pt}
    \makegapedcells
    \begin{tabular}{|c|c|c|}
        \hline
        & \textsc{da} & \textsc{dap} \\
        \hline
        Deterministic
        & \makecell{$0.158-\epsilon \rightarrow \bm{0.474-\epsilon}$\\
                    {\scriptsize (Thm.~\ref{thm:thr-dp})}}
        & \makecell{$\frac{1}{4}\!\left(1-\frac{1}{e}\right) \rightarrow
                    \bm{\frac{1}{2}\!\left(1-\frac{1}{e}\right)-\epsilon}$\\
                    {\scriptsize (Thm.~\ref{thm:gua-DAP})}}
        \\
        \hline
        IFR
        & \makecell{$0.122-\epsilon \rightarrow \bm{0.268-\epsilon}$\\
                    {\scriptsize (Thm.~\ref{thm:thr-dp})}}
        & \multirow{2}{*}{\makecell{$?\rightarrow \bm{0.25-\epsilon}$\\
                    {\scriptsize (Thm.~\ref{thm:gua-DAP-stoch})}}}
        \\
        \cline{1-2}
        General Distribution
        & $?$
        &
        \\
        \hline
    \end{tabular}
    \caption{Comparison of best-known approximation ratio guarantees in relevant settings. Entries in bold represent our results, with the corresponding theorem given in parentheses; an arrow points from the previous best-known guarantee to our new guarantee. A ``?'' indicates that no result was known prior to our work; a ``?'' not followed by an arrow indicates that the setting remains open. The previous best-known guarantees are due to \citet{aouad2018greedy} for \textsc{da} and \citet{bai2022coordinated} for \textsc{dap}. Note that $\epsilon>0$ is a parameter that we can choose and the runtime of the algorithms is an increasing function of $\frac{1}{\epsilon}$.}
    \label{tab:main-results}
\end{table}

\renewcommand{\arraystretch}{1.}
\begin{table}[htp]
    \centering
    \small
    \setcellgapes{3.6pt}
    \makegapedcells
    \begin{tabular}{|c|c|c|}
        \hline
        & \textsc{da} & \textsc{dap} \\
        \hline
        Deterministic
        & \makecell{$0.158-\epsilon \rightarrow \bm{0.474-\epsilon}$\\
                    {\scriptsize (Thm.~\ref{thm:thr-dp})}}
        & \makecell{$\frac{1}{4}\!\left(1-\frac{1}{e}\right) \rightarrow
                    \bm{\frac{1}{2}\!\left(1-\frac{1}{e}\right)-\epsilon}$\\
                    {\scriptsize (Thm.~\ref{thm:gua-DAP})}}
        \\
        \hline
        IFR
        & \makecell{$0.122-\epsilon \rightarrow \bm{0.268-\epsilon}$\\
                    {\scriptsize (Thm.~\ref{thm:thr-dp})}}
        & \multirow{2}{*}{\makecell{$?\rightarrow \bm{0.25-\epsilon}$\\
                    {\scriptsize (Thm.~\ref{thm:gua-DAP-stoch})}}}
        \\
        \cline{1-2}
        General Distribution
        & $?$
        &
        \\
        \hline
    \end{tabular}
    \caption{Comparison of best-known approximation ratio guarantees in relevant settings. Entries in bold represent our results, with the corresponding theorem given in parentheses; an arrow points from the previous best-known guarantee to our new guarantee. A ``?'' indicates that no result was known prior to our work; a ``?'' not followed by an arrow indicates that the setting remains open. The previous best-known guarantees are due to \citet{aouad2018greedy} for \textsc{da} and \citet{bai2022coordinated} for \textsc{dap}. Note that $\epsilon>0$ is a parameter that we can choose and the runtime of the algorithms is an increasing function of $\frac{1}{\epsilon}$.}
    \label{tab:main-results}
\end{table}

We first discuss the results for \textsc{da}. For this setting, we show that the fluid relaxation based on the choice-based deterministic linear program (CDLP) is in fact a good proxy for the original problem: an approximately optimal inventory vector for the fluid relaxation can be transformed into an effective initial inventory vector for \textsc{da}. Thus, we avoid the Monte Carlo sampling and enumeration procedures used in prior approximation algorithms for \textsc{da} \citep{aouad2018greedy}, thereby significantly improving the computational efficiency. In particular, our algorithm has a runtime of $O(\frac{n} {\epsilon^{2}} \log^2(\frac{K}{\epsilon})\log K\log n)$, where $K$ is the cardinality constraint on the total number of inventory units and $n$ is the total number of products. Moreover, this approach leads to improved approximation guarantees for both deterministic $T$ and stochastic $T$ following an IFR distribution:  $0.319-\epsilon$ and $0.188-\epsilon$, respectively\footnote{Assuming the existence of an efficient oracle, the algorithm in \cite{aouad2018greedy} provides a $(0.139-\epsilon)$-approximate solution with a probability of at least $1-\delta$ for any $\epsilon\in (0,1/4)$ and $\delta >0$. 
To the best of our knowledge, the existence of the oracle remains an open problem. Note that our algorithms have better approximation guarantees than this result as well. Also, \cite{aouad2022stability} give a near-optimal algorithm for \textsc{da}, but this algorithm is not polynomial-time (see Section~\ref{sec-lit}).}. Inspired by \cite{chen2022approximation}, we also propose a dynamic-programming-based fully polynomial-time approximation scheme (FPTAS) for the first step, which leads to improved factors of $0.474-\epsilon$ (deterministic $T$) and $0.268-\epsilon$ (IFR $T$), at the cost of a larger runtime.
 Moreover, our work resolves the open problem posed in \cite{aouad2018greedy} of finding a constant-factor approximation algorithm for \textsc{da} under a budget constraint. Our threshold-based algorithm achieves the same approximation guarantees under a budget constraint as under a cardinality constraint: $0.319-\epsilon$ for deterministic $T$ and $0.188-\epsilon$ for IFR $T$.
\endgroup

We next present the results for \textsc{dap}. 
For deterministic (known) $T$, we work directly with the CDLP with multiple customer types (multi-type CDLP) and improve the previous approximation guarantee from 
\(\frac{1}{4}(1-\frac{1}{e})\) in \citet{bai2022coordinated} to
\(\frac{1}{2}(1-\frac{1}{e})-\epsilon\). Note that for the special case of $T=1$ (single-customer setting), \cite{el2023joint} show that, unless P$=$NP, no polynomial-time algorithm can achieve an approximation ratio better than $1-\frac{1}{e}$. Perhaps more importantly, our framework extends to the setting where $T$ is uncertain and follows an arbitrary (known) distribution, including high-variance distributions, where no prior approximation results were known. In this more general setting, we establish a lower bound of $0.25-\epsilon$ on the approximation ratio of our algorithm. The same framework extends to budget constraints, where the first-stage approximation factor changes from \(\frac{1}{2}\) to \(\frac{1}{3}\) and the overall constant factor changes, resulting in a \((\frac{1}{3}(1-\frac{1}{e})-\epsilon)\)-approximation for deterministic \(T\) and a \((\frac{1}{6}-\epsilon)\)-approximation when \(T\) follows an arbitrary distribution (Theorem~\ref{thm:gua-budget-dap}).

\vspace{-5pt}

\vspace{-5pt}
\subsection{Key Algorithmic Ideas} \label{sec:tech}

\begingroup
We now describe the two-step algorithmic framework and the key technical ideas. In the first step, we approximately optimize the appropriate fluid objectives over the feasible set of initial inventories. 
Specifically, let \(F(\vecc)\) denote the fluid objective, corresponding to variants of the CDLP for \textsc{da} and \textsc{dap}; its exact form in each setting is summarized in Table~\ref{tab:fr} of Section~\ref{sec:fluid_relaxations}. The first step approximately solves
\(
\max_{\vecc\in\mathcal F} F(\vecc),
\)
where \(\mathcal F\) denotes the feasible set of initial inventories. Throughout, we refer to this optimization problem as the \emph{fluid relaxation}. Suppose the algorithm returns \(\widehat{\vecc}\in\mathcal F\) satisfying
\(
F(\widehat{\vecc})
\geq
\alpha \max_{\vecc\in\mathcal F}F(\vecc)
\)
for some $\alpha\in(0,1]$. One can show that the fluid objective provides an upper bound on the optimal value of the original stochastic problem, denoted by $\mathrm{OPT}$, i.e., \(\max_{\vecc\in\mathcal F}F(\vecc)\ge \mathrm{OPT}\).

The second step uses the integral inventory vector \(\widehat{\vecc}\) obtained in the first step to construct an implementable solution with good performance under the original stochastic consumption process. For \textsc{da}, this means constructing an integral initial inventory vector, denoted by \(\vecc^{\operatorname{ALG}}\), that can be used in the show-all stochastic consumption process. For \textsc{dap}, this means using the selected inventory together with an online assortment policy. If the transformation step guarantees that the expected revenue of the resulting policy is at least a \(\beta\)-fraction of \(F(\widehat{\vecc})\), then the overall approximation factor is \(\alpha\beta\). More precisely, if
\(
\operatorname{Rev}(\widehat{\vecc})
\geq
\beta F(\widehat{\vecc}),
\)
where \(\operatorname{Rev}(\widehat{\vecc})\) denotes the expected revenue achieved in the original stochastic problem by the decision obtained from \(\widehat{\vecc}\): the inventory \(\vecc^{\operatorname{ALG}}\) for \textsc{da}, and the inventory together with the online policy for \textsc{dap}, then
\[
\operatorname{Rev}(\widehat{\vecc})
\geq
\alpha\beta
\max_{\vecc\in\mathcal F} F(\vecc)
\geq
\alpha\beta\cdot \operatorname{OPT},
\]
where \(\operatorname{OPT}\) is the optimal value of the corresponding original problem.

The unification in this paper is algorithmic: the same relax-and-transform framework applies to both \textsc{da} and \textsc{dap}. Moreover, both analyses rely on two key structural properties of the customers' fluid consumption process under MNL, as we will discuss in the following paragraphs. 
\endgroup

\textbf{Technical ideas for \textsc{da}.} \begingroup
A central challenge for \textsc{da} is that it is unknown whether there exists a good fluid proxy for the stochastic process in which customers can see all the available products. We show that, under the MNL choice model, the CDLP provides such a fluid proxy. The gap between the CDLP objective and the expected revenue of the corresponding \textsc{da} policy has two sources. First, there is a \emph{personalization gap}: the CDLP can effectively choose which subsets of products to offer to an arriving customer, whereas in \textsc{da}, every stocked product remains available until it sells out. Consequently, lower-revenue products may cannibalize demand from higher-revenue products. Second, there is a \emph{stochasticity gap}: the CDLP optimizes a deterministic fluid consumption process (see Section~\ref{sec:fluid_relaxations} for more details), whereas the actual consumption process is stochastic. Under general choice models, the overall gaps can be arbitrarily large.

Perhaps surprisingly, we show that, under the MNL choice model, a small but careful modification of the inventory vector is sufficient to control the personalization gap. Our first observation is that, given an optimal CDLP solution, one can construct a possibly fractional initial inventory vector that achieves the same revenue in a non-personalized fluid consumption process. Therefore, the personalization gap reduces to bounding the integrality gap of the initial inventory vector in this fluid process. We then establish two structural properties that control this integrality gap. First, this inventory vector has a generalized revenue-ordered structure: products above a revenue threshold are fully stocked, products below the threshold are excluded, and at most one threshold product has fractional inventory. Second, we prove a \emph{separability} property, which characterizes the effect of changing  the inventory of a single product as replacing a terminal segment of the fluid consumption process. Together, these properties imply that the loss from converting the fractional CDLP solution into an integral inventory vector for \textsc{da} is localized to rounding the inventory of a single product in a fluid process. Importantly, the same fluid-process representation and structural properties will be reused to establish the diminishing-return structure underlying our algorithm for \textsc{dap}. These properties also lead to a nearly linear-time algorithm for solving the single-type CDLP. Given the widespread application of CDLP in revenue management, this fast algorithm may be of independent interest. 

We then bound the stochasticity gap. When $T$ is deterministic, we use the separability property again and structural properties of the stochastic consumption process to reduce the analysis to a single-product problem. After this reduction, the standard $1-1/e$ bound for a single item with stochastic customer arrivals applies \citep{alaei2012online}. When $T$ is stochastic, the analysis is more delicate because the fluid objective is evaluated at the mean horizon $\mu_T$, while the stochastic process is evaluated over different realizations of $T$. Under the IFR assumption, we establish a new single-product comparison showing that expected sales under the random horizon are at least a $1/e$ fraction of the corresponding fluid sales at the mean horizon. To prove this bound, we formulate a factor-revealing program over IFR distributions and show that the worst case is attained by shifted geometric distributions, for which the bound approaches $1/e$. 

It remains to solve the first step, and we give two algorithms with different runtime and guarantees. The first one uses the sales-based formulation of the CDLP, known as the SBLP \citep{gallego2004managing,liu2008choice}. Fixing the number of no-purchase customers $y_0$ yields an objective that is monotone and
diminishing-returns submodular in the inventory vector. Motivated by this observation, we search over a geometric grid
of $y_0$ values and apply a threshold-based submodular maximization routine \cite{soma2018maximizing}, which yields
a \emph{fast} $(1-\tfrac{1}{e}-\epsilon)$-approximation. We also propose a dynamic-programming-based fully polynomial-time approximation scheme (FPTAS) that gives a first-stage guarantee of $1-\epsilon$ at the cost of a significantly larger runtime.
\endgroup

\textbf{Technical ideas for \textsc{dap}.} \begingroup Since \textsc{dap} allows the seller to personalize assortments based on customer types, a generalization of the single-type CDLP that incorporates heterogeneous customers is a natural fluid objective. Unlike in \textsc{da}, there is no personalization gap. Moreover, the stochasticity gap can be controlled by online assortment policies from the literature \citep{golrezaei2014real, ma2021dynamic}, see the end of this subsection. The main challenge therefore lies in the first stage: optimizing the multi-type CDLP over the initial inventory vector.

Our main contribution is an improved approximation algorithm for this relaxation for deterministic $T$ and its generalization to arbitrary distributions of $T$. For deterministic $T$, \cite{bai2022coordinated} optimize a further surrogate LP whose objective is monotone submodular. This additional relaxation, however, incurs a constant-factor loss. In contrast, we work directly with the multi-type CDLP and develop an approximation algorithm for it. Prior work also leaves the case of stochastic $T$ open.  The mean-horizon relaxation used for \textsc{da}, which replaces $T$ by $\mu_T$, can be arbitrarily loose for general distributions of $T$ \citep{bai2023fluid}. We instead optimize the distribution-aware relaxation obtained by taking the expectation of the multi-type CDLP value over the distribution of $T$. To introduce the algorithmic idea, we start from single-type CDLP.   \endgroup

\noindent \emph{Single-type CDLP:} We show that decreasing price order is a diminishing-return submodular order (DR-SO) for the CDLP objective. This property generalizes the submodular order (SO) property \citep{sof}, which has been used in static assortment optimization. Proving the DR-SO property for the relaxation is technically challenging because small changes in inventory can cause significant changes in the LP optimal solution. Here we reuse the generalized revenue-ordered structure and the separability property developed above for the MNL non-personalized fluid process to analyze the influence of changing the starting inventory on the revenue in the fluid process. 
\smallskip

\begingroup
\noindent \emph{Multi-type CDLP:} Perhaps surprisingly, the decreasing price order is not a DR-SO for the multi-type CDLP. The reason is that the multi-type CDLP determines not only the total sales of each product but also how its inventory is allocated across customer types. Adding inventory may change this allocation and thereby violate the diminishing-return property. To address this issue, we reformulate the multi-type CDLP as an inventory-allocation problem. In this formulation, the algorithm chooses both the total inventory and a fractional allocation of each product across customer types. This reformulation allows us to exploit the DR-SO structure of the single-type CDLP within a customized threshold-augmentation algorithm. The resulting algorithm achieves essentially the same constant-factor guarantee that would follow if the multi-type CDLP itself possessed a DR-SO structure. This threshold-augmentation approach may also be useful for other resource allocation problems with submodular-like objectives.

 Finally, once the initial inventory is decided in \textsc{dap}, we use standard algorithmic techniques from the online personalized assortment literature to decide the personalized assortments \citep{golrezaei2014real, ma2021dynamic}. For deterministic $T$, we use a randomized CDLP-based sampling policy and bound its loss relative to the multi-type CDLP. We show an improved competitive ratio of $1-1/e$ compared with \cite{bai2022coordinated}, by extending a result from \cite{ma2021dynamic}. For stochastic $T$, we use the online greedy policy, which, despite being oblivious to $T$, attains at least half of the value of the distribution-aware relaxation described above \citep{chan2009stochastic,golrezaei2014real}. Combining these online policies with our first-stage approximation algorithms yields the overall guarantees stated in Table~\ref{tab:main-results}.

\endgroup

\subsection{Related Work} \label{sec-lit}
There is a large body of work on the single-customer (static) assortment optimization problem. For the MNL choice model, optimal or approximately optimal algorithms are known for a variety of constraints. In particular, \cite{talluri2004revenue} show the optimal unconstrained assortment follows a revenue-ordered structure, where only the highest-priced products are offered. \cite{rusmevichientong2010dynamic} propose a polynomial-time algorithm to solve the cardinality-constrained assortment optimization problem. \cite{desir2014near} show that the problem is NP-hard under a budget constraint and give a fully polynomial-time approximation scheme (FPTAS) for the problem. Various other variants of this problem have been considered; for further details, see \cite{rusmevichientong2012robust, sumida2021revenue, el2023joint, sof}. 
 For a comprehensive review of static assortment optimization under other choice models, we refer the reader to \cite{kok2015assortment} and \cite{strauss2018review}.
 
 The dynamic setting poses greater challenges due to the stochastic nature of customer behavior. Several works explore the theoretical difficulties of this problem and develop heuristics based on continuous relaxation or probabilistic assumptions on demand distributions \citep{smith2000management, mahajan2001stocking, gaur2006assortment, honhon2010assortment}. 
Recent research has developed provably good approximation algorithms in a variety of settings. We first review the existing literature for \textsc{da}. \cite{aouad2018greedy} propose the first constant-factor approximation algorithm for \textsc{da} under the MNL choice model, assuming the distribution of the number of customers satisfies the IFR property. \cite{aouad2022stability} further explore the structural properties of the stochastic consumption process under the MNL choice model and develop a Monte Carlo algorithm to compute the starting inventory within an arbitrary degree of accuracy $\epsilon$ for general distributions of $T$. However, their algorithm has a runtime of $O(\frac{1}{\delta}(nK)^{O(\frac{\log ({v_{\text{max}}/v_{\text{min}}})}{\epsilon^3})})$, which is exponential in $\log({v_{\text{max}}/v_{\text{min}}})$, where $v_{\text{max}}$ and $v_{\text{min}}$ denote the maximum and minimum weights of the products, respectively, in the MNL choice model, and the guarantee holds with probability at least $1-\delta$. Observe that this algorithm is not a polynomial-time algorithm and thus does not directly compare with our polynomial-time results.
 In addition, other works such as \cite{goyal2016near}, \cite{segev2015assortment}, \cite{aouad2019approximation} and \cite{GoldsteinFeder2023} propose approximation algorithms for dynamic assortment planning in more specialized settings, such as products consumed in a particular order or horizontally differentiated products.
These specialized structures enable innovative algorithm designs, in some cases for more general choice models (for example, \citet{GoldsteinFeder2023}), but the techniques are not directly applicable to our problem setting and vice versa. \begingroup Relatedly, \cite{elmachtoub2022revenue} study product retirement and customer selection under the MNL choice model, where the seller may permanently retire products and choose which customers to serve. They similarly use a sales-based LP relaxation and exploit a nested-assortment structure in the corresponding fluid dynamics. \endgroup

Another line of work in dynamic assortment optimization focuses on the setting where the assortment offered to each customer is decided in real time, depending on the remaining inventory and customer type. The inventory may be fixed exogenously or may be jointly decided with the assortment. This stream of work is related to the \textsc{dap} setting. \begingroup The work that is most closely related to ours is by \cite{bai2022coordinated}. They consider the \textsc{dap} problem through the lens of a surrogate to the CDLP  
and give a $\frac{1}{4}(1-\frac{1}{e})$-approximate solution for the MNL choice model with deterministic $T$ and a $(1-(\sqrt{2}+1) \sqrt[3]{\frac{n}{K}})$-approximate solution for a general choice model. We improve the approximation guarantee for MNL to $\frac{1}{2}(1-\frac{1}{e})-\epsilon$ for deterministic $T$ and show an approximation guarantee of $\frac{1}{4}-\epsilon$ for stochastic $T$, which is the first non-trivial approximation guarantee for stochastic $T$. When the inventory is fixed exogenously, the \textsc{dap} problem reduces to an online assortment optimization problem. Of particular interest to us is the horizon-oblivious online Greedy algorithm, which offers the revenue maximizing assortment to every customer and does not require knowledge of the number of customers $T$. \cite{chan2009stochastic, golrezaei2014real} show that the total revenue of this algorithm is always at least half of the optimal value of the clairvoyant and CDLP, respectively. For deterministic $T$, \cite{bai2022coordinated} show that a simple randomized rounding algorithm based on the CDLP also achieves the same guarantee as Greedy. A second policy of interest is the independent sampling and rounding algorithm, which samples an assortment from an optimal CDLP solution and offers the available subset. Unlike Greedy, this policy uses knowledge of the horizon since the CDLP encodes its value. \cite{bai2022coordinated} show that this policy achieves a competitive ratio of $\frac{1}{2}$. \cite{ma2021dynamic} show that for stationary customer arrivals, the randomized rounding policy yields a competitive ratio of $(1-\frac{1}{e})$. We extend the analysis to establish the $1-\frac{1}{e}$ guarantee for the nonstationary setting (Lemma~\ref{lem:tran-d}). Relatedly, \cite{huang2024basic} consider the joint inventory and online assortment problem under the MNL choice model with reusable resources. For further discussion of online assortment optimization problems, we refer the reader to \cite{ma2020algorithms}, \cite{rusmevichientong2020dynamic}, and \cite{gong2022online}. 
\endgroup

Our work also closely relates to the literature on joint inventory and resource allocation problems in a two-stage approximation framework. \cite{chen2022approximation} study a joint inventory selection and online matching problem. While their setting focuses on matching rather than assortment optimization, they use a similar two-step framework based on a fluid relaxation. \cite{epstein2024optimizing} consider an inventory placement problem with different dynamic fulfillment policies. 
For a comprehensive review of this area, we refer readers to the literature reviews in these two papers. 

\begingroup
Another related line of research studies joint assortment and inventory problems with a focus on achieving regret-type guarantees. The models in this stream price inventory directly, through ordering costs or leftover penalties, whereas our objective is gross revenue and excess inventory is controlled through the cardinality or budget constraint, as in \cite{goyal2016near, aouad2018greedy, aouad2019approximation, bai2022coordinated}. The guarantees are also of a different type: the regret-type work establishes additive asymptotic optimality, whereas we establish worst-case multiplicative guarantees. Thus, the analyses of these two streams of work do not apply directly to each other. \cite{mouchtaki2021joint} propose a sample average approximation-based algorithm for the unconstrained problem with ordering cost under a Markov chain choice model that achieves sublinear regret. \citet{liang2021assortment} consider the joint assortment and inventory problem under MNL choice with nonzero initial inventory, ordering costs, and leftover penalties. They show that the fluid relaxation can be solved by sequentially solving a series of simpler linear programs and that the optimal solutions to these LPs exhibit a ``quasi-gain-ordered'' property. This property is conceptually similar to our ``generalized revenue-ordered'' property for the optimal initial inventory of the non-personalized fluid process, as both arise from the knapsack-like structure of the SBLP under MNL choice. Nevertheless, the two models consider different settings and yield distinct structural properties. In our model, the generalized revenue-ordered property characterizes the structure of the optimal initial inventory for the non-personalized fluid process; see Section~\ref{sec:structural} for details.
\cite{guo2023bilevel} and \cite{zhang2024leveraging} show improved performance guarantees for this problem under a general choice model. Additionally, some work considers the scenario where the choice model is not known. \cite{liang2023online} study the dynamic assortment planning problem under the MNL choice model with unknown preference weights, focusing on learning the weights and using an existing solution algorithm for the dynamic assortment optimization problem as an oracle.  
\endgroup

\subsection{Outline}
Section~\ref{sec:2} formally defines the \textsc{da} and \textsc{dap} problems, introduces the fluid relaxations used in our algorithms, and develops structural properties of the MNL fluid process that underlie the analyses of both settings. \begingroup Section~\ref{sec:da} develops the algorithms and performance guarantees for \textsc{da}, and Section~\ref{sec:dap} does the same for \textsc{dap}. Section~\ref{sec:extensions} presents extensions, including budget constraints and a dynamic-programming-based variant with stronger guarantees. We then report numerical experiments in Section~\ref{sec:exp}. Finally, Section~\ref{sec:conclusion} concludes the paper and discusses open problems. \endgroup

\section{Models, Fluid Relaxations and Structural Properties}\label{sec:2}

\subsection{Problem Formulations}\label{sec:2-dadap}

We start with the common setup for both \textsc{dap} and \textsc{da}. Let $\NN=\{1,2,\ldots, n\}$ denote the ground set of products, ordered by descending (fixed) nonnegative per-unit prices $r_1 \geq r_2 \geq\cdots \geq r_n\geq 0$. We use \(\mathbb{Z}_+^n\) and \(\mathbb{R}_+^n\) to denote the sets of \(n\)-dimensional nonnegative integer and real vectors, respectively. Before the selling horizon begins, the seller chooses an integral initial inventory vector $\vecc= (c_i)_{i \in \NN}\in\mathbb Z_+^n$, where $c_i$ is the number of units of product $i$ stocked. Throughout, $\vece_i\in\{0,1\}^n$ denotes the $i$th unit vector, and $\1\{\cdot\}$ denotes the indicator function. 
We write $\mathbf{0}$ for the zero vector, with its dimension indicated by a superscript when needed (e.g., $\mathbf{0}^n$). \begingroup
We consider two types of inventory constraints. Under a cardinality constraint, $\sum_{i\in \NN} c_i\leq K$, where $K\in \mathbb{Z}_+ $. Under a budget constraint, $\sum_{i\in \NN} b_i c_i\leq B$, where $B>0$ and $b_i>0$ is the inventory cost or resource consumption associated with one unit of product $i$. We use $\mathcal F$ to denote the feasible set of initial inventories. \endgroup

The selling horizon consists of $T$ customers, which can alternatively be viewed as $T$ periods with one customer arriving in each period. We use $t\in [T]$  (where $[T]$ denotes the set $\{1,2,\ldots,T\}$) to refer interchangeably to period \(t\) and the customer arriving in that period. The value of $T$ may be deterministic or random. When $T$ is random, its distribution $D_T$ is known before the inventory decision is made, but its realization is not, and $T$ is independent of customer choices. We use $\mu_T=\EE[T]$ to denote its mean. We assume $0<\mu_T<\infty$. Throughout the paper, customer choices follow the MNL model. For the \textsc{da} setting, let \(v_i\) be the attraction parameter of product \(i\), and normalize the attraction parameter of the outside option to \(v_0=1\). We assume \(v_i>0\) for every \(i\in\NN\). When assortment \(S\subseteq\NN\) is displayed, the choice probabilities are 
$\phi(i,S)=\frac{v_i}{v_0+\sum_{j\in S}v_j}$ for $i\in S\cup\{0\}$.
The expected revenue from displaying \(S\) to one customer is  $R(S)=\sum_{i\in S}r_i\phi(i,S)$.

 \subsubsection{The DA Problem.} \begingroup
We first introduce \textsc{da}, the dynamic assortment problem without personalization. After choosing the initial inventory vector \(\vecc\), the seller makes no further assortment decisions. Every arriving customer observes all products that remain in stock.

For a deterministic horizon \(T\), let
\(\ZDA_i(\vecc,T)\) denote the expected number of units of product \(i\) sold over the horizon, starting from inventory \(\vecc\). The corresponding expected total revenue is \(\fDA(\vecc, T)=\sum_{i \in \NN } r_i \ZDA_i(\vecc,T).
\) Thus, our goal is to identify a feasible starting inventory $\vecc \in \mathcal{F}$ that maximizes the expected total revenue: 
\(
\max_{\vecc\in \mathcal{F}} \fDA(\vecc,T). 
\) \endgroup
For stochastic $T$, the optimization problem becomes 
\(
\max_{\vecc\in \mathcal{F}}\EE_T[\fDA(\vecc,T)]=\max_{\vecc\in \mathcal{F}}\sum_{i \in \NN } r_i \EE_T[\ZDA_i(\vecc,T)].
\)

Similar to previous work, we assume that $D_T$ satisfies the IFR property. This means that $\frac{\Prb(T=k)}{\Prb(T\geq k)}$ is non-decreasing for any $k\in \mathbb{Z}_+$. The IFR property holds for various distributions commonly used in the literature, such as Poisson and geometric distributions \citep{goyal2016near, aouad2018greedy}.

 \subsubsection{The DAP Problem.}     
We next define the dynamic assortment problem with personalization, denoted by \textsc{dap}. In this setting, there are $m$ customer types, indexed by $\MM=\{1,2,\ldots,m\}$. \begingroup  A type $k\in \MM$ customer follows an MNL model with attraction parameters $(v_{ki})_{i \in \NN }$, where $v_{ki}\geq 0$, with outside-option parameter \(v_{k0}=1\); a type-$k$ customer offered assortment $S$ purchases product $i\in S$ with probability $\phi_k(i,S)=\frac{v_{ki}}{1+\sum_{j\in S}v_{kj}}$. Let $\lambda_{k,t}$ denote the probability that customer $t$ is of type $k$, where $\sum_{k\in \MM} \lambda_{k,t}=1$. Customer types are drawn independently, and the distribution
\((\lambda_{k,t})_{k\in\MM}\) is known before the inventory decision is made. When \((\lambda_{k,t})_{k\in\MM}\) is constant over $t$, we say that the arrivals are stationary. When $T$ is random, it is independent of the customer-type sequence and of customer choices. Unlike the \textsc{da} problem where customers see all available products, in \textsc{dap}, when a customer arrives at time $t$, the seller may choose a subset of available products (assortment) to offer. Let \(N_t\subseteq\NN\) denote the set of products with positive remaining inventory immediately before customer \(t\) arrives. A feasible online policy \(\pi\) chooses an assortment \(S_t\subseteq N_t\) after observing the current customer's type. This choice may depend on the purchase history of the first \(t-1\) customers, the current inventory state, the current customer's type, and the known distributions of future customer types, but not on future type realizations. The \textsc{da} model corresponds to the fixed policy \(S_t=N_t\) in every period. Let $\Zpi_i(\vecc,T)$ denote the expected consumption of product \(i\) under policy \(\pi\), starting from inventory \(\vecc\), and let \(\fpi(\vecc,T)
=
\sum_{i\in\NN}r_i\Zpi_i(\vecc,T)\) denote the corresponding expected total revenue.
When $T$ is deterministic, the joint inventory and assortment optimization problem is
\[
\max_{\vecc\in\mathcal F}\;
\sup_{\pi\in\Pi_T}\fpi(\vecc,T),
\]
where \(\Pi_T\) denotes the set of feasible nonanticipative policies over a horizon of length \(T\). This formulation coincides with the setting studied by \cite{bai2022coordinated}.

When $T$ is stochastic, the initial inventory vector is chosen using only the distribution \(D_T\). The objective is \(
\max_{\vecc\in \mathcal{F}}\sup_{\pi} \EE_{T\sim D_T} [ \fpi (\vecc,T)].\) Unlike \textsc{da}, here we allow \(D_T\) to be an arbitrary distribution. The deterministic $T$ formulation is recovered when $D_T$ is degenerate at a single value. Here the supremum is over nonanticipative policies that do not observe the realized $T$, and the inventory decision uses only $D_T$. 
\endgroup

\subsection{Fluid Objectives} \label{sec:fluid_relaxations}

We next introduce the CDLP-based fluid objectives used throughout this paper. For a fixed inventory vector, these CDLP values provide fluid upper bounds on the expected revenue of the corresponding stochastic consumption process. 

For \textsc{da}, we use the single-type CDLP as the fluid objective. 
Let $y(S)$ denote the fraction of time assortment $S$ is offered. Then define:
\begin{equation}\label{cdlp}
\fLP(\vecc,T)=\max_{y(S)\geq 0,  \forall S\subseteq \NN } \left\{ T \sum_{S\subseteq \NN} y(S)R(S) :  T \sum_{S\subseteq \NN} y(S)\phi(i,S)\leq c_i, \forall i \in \NN ; \sum_{S\subseteq \NN } y(S)= 1 \right\}.
\end{equation}
$\fLP(\vecc,T)$ relaxes the original problem in two ways: (i) Customers arrive continuously from time $0$ to $T$. The customer arriving at time $t$ consumes a deterministic infinitesimal amount $\phi(i, S(t))\;\mathrm{d}t$ of product $i$ where $S(t)$ is the assortment offered at $t$. (ii) We can decide the mix of assortments offered to each customer, i.e., personalization. For stochastic $T$ following an IFR distribution, we use $\fLP(\vecc,\mu_T)$ as the fluid approximation for the expected revenue. One can show that $\fLP(\vecc,\mu_T)$ is an upper bound for $\EE_T[\fDA(\vecc,T)]$; see Proposition~\ref{prop:ub} in Appendix~\ref{app:prop-ub}. 

Under the MNL model, \cite{gallego2004managing} introduce the following more compact \emph{sales-based formulation} (SBLP). In this formulation, each $y_i$ represents the sales of product $i$, and $y_0$ denotes the number of customers making no purchases. 
\begin{equation}
\label{sblp}
\fLP(\vecc,T)=\max_{\mathbf{y} \in \RR^{n+1}_+} \Biggl\{ \sum_{i \in \NN } r_i y_i :\sum_{i \in \NN } y_i + y_0 = T; y_i \leq c_i,  y_i \leq  y_0 v_i, \forall i \in \NN  \Biggr\}.
\end{equation}
 The constraints \(y_i\le c_i\) enforce the inventory limits, and the constraints \(y_i\le v_i y_0\) encode the MNL choice structure.

\begingroup
For \textsc{dap}, we use a natural generalization of the single-type CDLP that accounts for multiple customer types, defined as multi-type CDLP. Let $R_k(S)=\sum_{i\in S} r_i\phi_k(i,S)$ denote the expected revenue of assortment $S$ for customer type $k$. The formulation is as follows: 
\begin{equation}
\label{eq:multi-cdlp}
\begin{aligned}
  \vspace{-5pt}
\fLPM(\vecc,T)
=
\max_{\mathbf{y}_k(S)\geq 0,  \forall S\subseteq \NN, \forall k\in \MM}&
\Biggl\{
 \sum_{t=1}^T \sum_{S\subseteq \NN}\sum_{k \in \MM}\lambda_{k,t}  y_k(S)R_k(S):\\& 
\sum_{t=1}^T\sum_{k \in \MM} \sum_{S\subseteq \NN} y_k(S)\lambda_{k,t} \phi_k(i,S)\leq c_i, \forall i \in \NN\;;\sum_{S\subseteq \NN} y_k(S)= 1 , \forall k\in \MM
\Biggr\}.
\end{aligned}
\end{equation} 
Notably, for any starting inventory vector $\vecc$, $\fLPM(\vecc,T)$ is an upper bound for the expected revenue of the optimal assortment personalization policy (Lemma 1 in \cite{golrezaei2014real}, and we restate it in Appendix~\ref{app:multi-ub}). The multi-type CDLP also has a polynomial-size equivalent sales-based formulation, which we defer to Appendix~\ref{app:multi-sblp}.
\endgroup

\begingroup
When \(T\) is stochastic, we use the distribution-aware fluid objective $\EE_{T\sim D_T}\bigl[\fLPM(\vecc,T)\bigr]$. This choice is important when \(D_T\) has high variance. A relaxation that simply replaces \(T\) by its mean can be arbitrarily loose: the ratio between the revenue of the original stochastic problem and the corresponding mean-horizon fluid value may approach zero
\citep{bai2023fluid, aouad2022nonparametric}. By contrast, taking the expectation of \(\fLPM(\vecc,T)\) over \(T\) preserves the uncertainty in the horizon. It relaxes the stochasticity of customer types and choices, but does not replace the random horizon by its mean. Indeed, the pointwise upper bound for deterministic $T$ implies \(\EE_{T\sim D_T}
\left[
\sup_{\pi\in\Pi_T}\fpi(\vecc,T)
\right]
\le
\EE_{T\sim D_T}
\bigl[
\fLPM(\vecc,T)
\bigr].\) In fact, $\EE_T[\fLPM(\vecc,T)]$ upper bounds the revenue under the optimal personalization policy that knows the full sequence of customer types \citep{golrezaei2014real}, whereas the formulations in \cite{aouad2022nonparametric} upper bound the revenue under the optimal \emph{online} policy\footnote{Our result does not contradict the impossibility results in \cite{aouad2022nonparametric} and \cite{bai2023fluid}. First, the impossibility results of \citet{bai2023fluid} and the \textsc{Indep} model of \cite{aouad2022nonparametric} are against the mean-horizon fluid benchmark, which aggregates horizon uncertainty through the expected demand \(\mu_T\) under stationary arrivals. Second, the impossibility result (no constant-factor competitive ratio) of the \textsc{Correl} model of \cite{aouad2022nonparametric} arises from a different mechanism in correlated arrivals with arrival-dependent rewards while in our setting, products have fixed prices. }. 

\endgroup
Table~\ref{tab:fr} summarizes the fluid objective used in each setting. In every case, the full fluid relaxation is obtained by maximizing the objective over \(\vecc\in\mathcal F\).
\begin{table}[htp]
\centering
\renewcommand{\arraystretch}{1.1}  
\setlength{\tabcolsep}{6pt}        
\begin{tabular}{|c|c|c|}
\hline
& \textsc{da} & \textsc{dap}\\
\hline
Deterministic $T$ & $\fLP (\vecc,T)$ & $\fLPM (\vecc,T)$ \\
\hline
Stochastic $T$ & $\fLP(\vecc,\mu_T)$ & $\EE_T[\fLPM (\vecc,T)]$  \\
\hline
\end{tabular}
\caption{Fluid objective values used in different settings. }\label{tab:fr}
\vspace{-5pt}
\end{table}
\vspace{-5pt}

\subsection{The MNL Fluid Process and Key Structural Properties}     \label{sec:structural} \begingroup
We now introduce a deterministic MNL fluid process that is closely related to the CDLP, but does not personalize the displayed assortment and instead offers every product with positive remaining inventory. This fluid process will serve as a central analytical tool
throughout the paper. Specifically, let \textsc{fp}$(\vecx,T)$ denote the fluid process that starts with inventory $\vecx$ and has $T$ periods. In this process, customers are represented 
by a continuum that arrives at unit rate over $[0,T]$. At each time $t$, all 
products with positive remaining inventory are offered, so the assortment
changes only when a product stocks out. \endgroup Let $S(t)$ denote the set of products
remaining in stock at time $t$. Each product $i\in S(t)$ is consumed
deterministically at rate $\phi(i,S(t))$. We write
$\ZFP_i(\vecx,T)$ for the total consumption of product $i$ and
\(
\fFP(\vecx,T)
=
\sum_{i\in\NN}r_i\ZFP_i(\vecx,T)
\)
for the resulting total revenue. For example, consider two products with initial inventories $x_1=1$ and
$x_2=2$, attraction parameters $v_1=v_2=v_0=1$, and horizon $T=4$.
Initially, both products are available and are consumed at rate $1/3$.
Product~1 therefore stocks out at time $t=3$, after which only product~2
remains available. The complete sequence of assortments and stockout times can be obtained straightforwardly by tracking the order in which products stock out; see the corresponding Sequencing algorithm in Appendix~\ref{app:alg-seq}.

We next establish two structural properties of this fluid process that play a central role throughout the paper.

\paragraph{Generalized revenue-ordered inventories.}
When $T=1$, the problem $\max_{\vecx\preceq\vecc}\fFP(\vecx,T)$ reduces to the unconstrained static assortment problem, whose optimal solution is revenue-ordered: there is a price threshold above which products are offered and below which they are not \citep{talluri2004revenue}. For general $T$, we establish an analogous structure at the level of \emph{inventory}.

\begingroup
Under the MNL model, it is known that the
CDLP value $\fLP(\vecc,T)$ can be attained by \textsc{fp} starting from a
suitable reduced inventory $\vecx\preceq\vecc$, where $\preceq$ denotes
componentwise inequality (Theorem~2 of \citet{goyal2022dynamic}). Conversely, for every $\vecx\preceq\vecc$, the
consumption quantities generated by \textsc{fp}$(\vecx,T)$ form a feasible
solution to the CDLP with inventory $\vecc$. Hence, \(
\max_{\vecx\preceq\vecc}\fFP(\vecx,T)=\fLP(\vecc,T)\). We next construct such a $\CFP$ explicitly with the \emph{generalized revenue-ordered property}.
\vspace{-3pt}
\begin{lemma}[Generalized revenue-ordered inventory]\label{lem:3-rev}
For any inventory vector $\vecc \in \mathbb{R}_+^n $, let $\{y^*_i\}_{i\in\NN\cup\{0\}}$ denote an optimal solution to the sales-based formulation of the single-type CDLP $\fLP(\vecc,T)$. Define $i_{\tau}$ as the highest index $i$ for which $\sum_{k=1}^{i-1}\min\{c_k, v_k y_0^*\}\leq T-y_0^*$ and construct an inventory vector $\CFP\in\RR^n$ by
\[
\cFPi=c_i \ \text{ for } 1\le i<i_\tau,\qquad
\cFPi=0 \ \text{ for } i>i_\tau,\qquad
x_{i_{\tau}}^{\vecc}=T-y_0^*-\sum_{k=1}^{i_{\tau}-1}\min\{c_k, v_k y_0^*\}.
\]
Then $\CFP\preceq\vecc$ and $\fFP(\CFP,T)=\max_{\vecx\preceq\vecc}\fFP(\vecx,T)=\fLP(\vecc,T)$.
\end{lemma}
The proof, which constructs $\CFP$ from the knapsack structure of the sales-based formulation \eqref{sblp}, is given in Appendix~\ref{app:prf-3-rev}. 
Moving forward, we refer to $\CFP$ as this revenue-ordered inventory vector and, when the dependence on $T$ is clear, suppress it from the notation. By construction, when $\vecc$ is integral, at most one product (the threshold product $i_\tau$), has fractional inventory; every higher-priced product is stocked at its full level $c_i$, and every lower-priced product is excluded. This is precisely the fractional component rounded by the transformation algorithm in Section~\ref{sec:3-bad}. 
\endgroup
\begingroup
\paragraph{A Separability Property.}
Our second structural result characterizes the effect of adding inventory
of any single product on the fluid process. At the level of cumulative consumption, the augmented
process can be represented as the original process run over a shorter
horizon, together with additional consumption of the augmented product.
\begin{lemma}[Separability]
\label{lem:sep}
Fix an inventory vector $\vecx\in\RR_+^n$, a horizon $T$, a product
$i\in\NN$, and an additional inventory level $\delta_i\ge0$. Consider
\textsc{fp}$(\vecx,T)$ and
\textsc{fp}$(\vecx+\delta_i\vece_i,T)$. Define
\(
\alpha_i
=
\ZFP_i(\vecx+\delta_i\vece_i,T)
-
\ZFP_i(\vecx,T),
\)
the increase in the total consumption of product $i$. Then:
\begin{enumerate}[label=(\roman*)]
\item
For product $i$,
\(
\ZFP_i(\vecx+\delta_i\vece_i,T)
=
\alpha_i+\ZFP_i(\vecx,T-\alpha_i),
\)
and, for every $j\neq i$,
\(
\ZFP_j(\vecx+\delta_i\vece_i,T)
=
\ZFP_j(\vecx,T-\alpha_i).
\)
\item
The total revenue satisfies
\(
\fFP(\vecx+\delta_i\vece_i,T)
=
\fFP(\vecx,T-\alpha_i)+\alpha_i r_i.
\)
\end{enumerate}
\end{lemma}
\endgroup

\begin{figure}[htp]
\centering
\includegraphics[width=.82\textwidth]{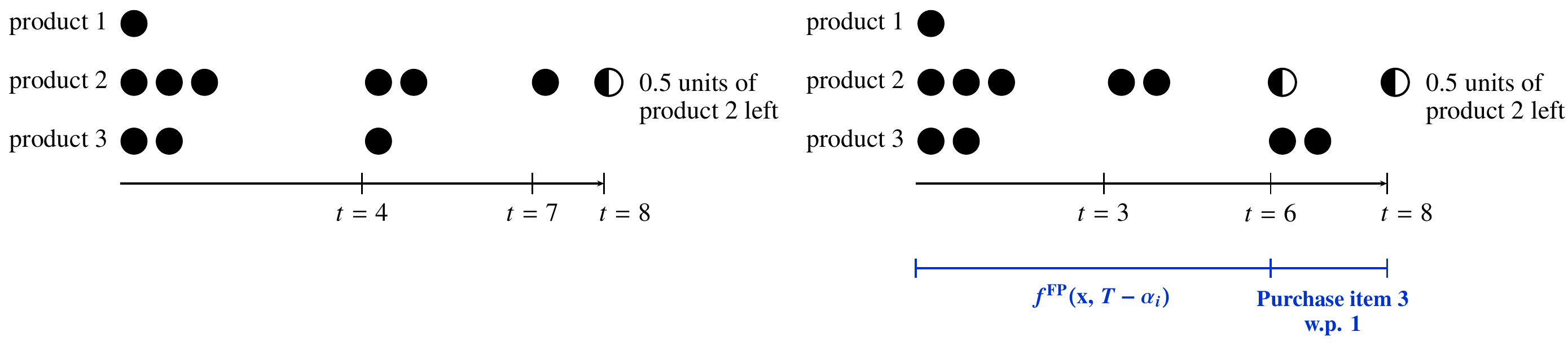}
\caption{Illustration of the separability property.
(a) The fluid process \textsc{fp}$(\vecx+\delta_i\vece_i,T)$.
(b) An equivalent representation consisting of
\textsc{fp}$(\vecx,T-\alpha_i)$ and an artificial terminal segment
accounting for $\alpha_i$ units of consumption of product $i$.}
 \vspace{-3pt}
\label{fig:sep-exp}
\end{figure}

Thus, in terms of cumulative product consumption,
\textsc{fp}$(\vecx+\delta_i\vece_i,T)$ is equivalent to first running
\textsc{fp}$(\vecx,T-\alpha_i)$ and then adding $\alpha_i$ units of
consumption of product $i$. This decomposition reproduces the total consumption of each product, but not necessarily the chronological assortment path of the augmented fluid process. To illustrate, consider three products with $x_1=1$, $x_2=3$, and $x_3=0$, attraction parameters $v_1=v_2=v_3=v_0=1$ and $T=8$. In \textsc{fp}$(\vecx,T)$, the initial assortment is $\{1,2\}$. By time $t=3$, product 1 is out of stock, and the assortment changes to $\{2\}$. By $t=7$, product 2 is also out of stock, and customers thereafter see no products. 
Now add two units of product~3. The actual assortment sequence of the
augmented process is shown in Figure~\ref{fig:sep-exp}(a). In this example,
$\alpha_3=2$. By Lemma~\ref{lem:sep}, the same cumulative consumption vector
is obtained by running the original process until
$T-\alpha_3=6$ and then adding two units of consumption of product~3, as
shown in Figure~\ref{fig:sep-exp}(b). Under both representations, $0.5$
units of product~2 remain at the end of the horizon.
 \vspace{-5pt}

The proof relies on the Independence of Irrelevant Alternatives (IIA) property of the MNL model. A key observation is that, with the outside option $0$ considered as a product that never stocks out, the purchase probability of any product $i \in \NN $ before it runs out of stock is $v_i$ times the purchase probability of option $0$. Moreover, the total consumption of all products, including option $0$, equals $T$. These two properties yield the matched consumption in the two processes. The detailed proof is in Appendix~\ref{prf:lem-sep}.
\begin{remark}
These two structural properties also lead to a new $O(n\log n)$ algorithm for solving the single-type CDLP; see Appendix~\ref{app:fast_cdlp} for more details. 
\end{remark}

\section{Dynamic Assortment without Personalization}
\label{sec:da}\begingroup
 In this section, we develop the algorithms and analysis for \textsc{da}, beginning with the case of deterministic $T$. Section~\ref{sec:da-fluid-opt} presents the first-stage algorithm for optimizing the fluid relaxation. In Section~\ref{sec:da-transform}, we first analyze the personalization and stochasticity gaps, then show how to use the CDLP solution to 
construct an integral starting inventory for \textsc{da} while controlling
both gaps. Section~\ref{sec:da-stochasticT} extends the analysis to a stochastic horizon $T$ satisfying the IFR property. Throughout this section, we focus on the cardinality constraint and discuss the extension to budget constraints in Section~\ref{sec:ext-budget}.
\endgroup
\subsection{Approximately Solving the Fluid Relaxation}      \label{sec:da-fluid-opt}

Recall that under a cardinality constraint, the fluid optimization problem is
\begin{equation}\label{eq:lpsstg}
\max_{\vecc\in\mathcal F}\fLP(\vecc,T).
\end{equation}
\begingroup
Here we present a threshold-based algorithm that approximates \eqref{eq:lpsstg} within a factor of $1-\frac{1}{e}-\epsilon$. Note that our two-step analysis is modular: the transformation step in Section~\ref{sec:da-transform} applies to any inventory vector, so any approximation algorithm for \eqref{eq:lpsstg} can be plugged into the first stage, and the analysis carries over with the corresponding first-stage factor (see Section~\ref{sec:ext-dp}). We focus on the threshold-based algorithm because it is fast and practical, and because it relies on threshold-based augmentation, the same idea underlying our algorithm for \textsc{dap} in Section~\ref{sec:dap}. 

The algorithm exploits the sales-based formulation~\eqref{sblp}. For a
fixed value of the outside-option fluid mass $y_0\in[0,T]$, define
\[
g(\vecc,y_0)=\max_{\mathbf{y} \in \RR^{n}_+} \left\{ \sum_{i \in \NN } r_i y_i :\sum_{i \in \NN } y_i \leq T-y_0; y_i \leq c_i,  y_i \leq  y_0 v_i, \forall i \in \NN  \right\},
\]
the SBLP with $y_0$ fixed. It is straightforward to show 
$\fLP(\vecc,T)=\max_{y_0\in[0,T]}g(\vecc,y_0)$. 
A key observation established by \cite{bai2022coordinated} is that
$g(\vecc,y_0)$ is monotone and diminishing-returns submodular
(DR-submodular) in the inventory vector $\vecc$. \endgroup Therefore, for each fixed
$y_0$, we can approximately maximize $g(\vecc,y_0)$ using a standard
algorithm for monotone DR-submodular maximization under a cardinality
constraint \citep{nemhauser1978analysis,sviridenko2004,soma2018maximizing}. \begingroup The optimal value of $y_0$ is not known in advance. We therefore consider a
geometric grid of candidate values. For each candidate, we approximately
maximize the corresponding fixed-$y_0$ objective, evaluate the resulting
inventory using the original CDLP objective $\fLP(\vecc,T)$, and return the 
best candidate. \endgroup

\begin{algorithm}[htp]
\SetAlgoNoLine
\begingroup
\caption{Threshold-Based Optimization for Cardinality-Constrained DA Fluid Relaxation}
\label{alg:da-thr}
\KwIn{Cardinality \(K\), horizon \(T\), accuracy parameter \(\epsilon\)}
\KwOut{Inventory vector \(\widehat{\vecc}\)}
Construct a geometric grid \(\mathcal G\) of candidate values for the outside-option fluid mass \(y_0\)\;
\For{\(y_0\in\mathcal G\)}{
Approximately maximize the fixed-\(y_0\) objective \(g(\vecc,y_0)\) over \(\{\vecc\in\mathbb Z_+^n:\|\vecc\|_1\le K\}\) using a cardinality-constrained DR-submodular maximization routine\;
}
Return the inventory vector \(\widehat{\vecc}\) among all candidates with the largest value of \(\fLP(\widehat{\vecc},T)\)\;
\endgroup
\end{algorithm}

\begingroup
\begin{theorem}[Approximating the \textsc{da} fluid relaxation]
\label{thm:da-fluid-opt}
For any \(\epsilon\in(0,1)\), Algorithm~\ref{alg:da-thr} returns an inventory vector \(\widehat{\vecc}\in\mathbb Z_+^n\) with \(\|\widehat{\vecc}\|_1\le K\) such that
\(
\fLP(\widehat{\vecc},T)
\geq
\left(1-\frac{1}{e}-\epsilon\right)
\max_{\|\vecc\|_1\le K}\fLP(\vecc,T).
\)
The running time is
\(
O\left(
\frac{n}{\epsilon^2}\log^2\left(\frac{K}{\epsilon}\right)\log K\log n
\right).
\)
\end{theorem}
The proof is given in Appendix~\ref{app:single-T}. The
$1/\epsilon^2$ dependence arises from the outer geometric grid over $y_0$
and the inner geometric threshold schedule used by the DR-submodular
maximization routine; the $\log K\log n$ factor arises from evaluating
marginal benefits of adding units of products by binary search. To simplify the complexity expressions below, throughout the paper we use
$\log K$ as shorthand for $\max\{1,\log K\}$ in running-time bounds. In Section~\ref{sec:extensions}, we
present a dynamic-programming-based variant that improves the first-stage
approximation ratio from \(1-\frac1e-\epsilon\) to \(1-\epsilon\) at the 
expense of a significantly higher running time. \endgroup
\vspace{-10pt}
\subsection{Transforming the Inventory Solution from the Fluid Relaxation}    \label{sec:da-transform}
\begingroup
Let \(\vecc\in\mathbb Z_+^n\) be any integral inventory vector, for example, the one returned by Algorithm~\ref{alg:da-thr}. In this section, we show how to transform it into an integral starting inventory whose expected revenue in the original stochastic \textsc{da} process is at least a constant factor of \(\fLP(\vecc,T)\). As discussed in Section~\ref{sec:fluid_relaxations}, the CDLP relaxes \textsc{da} by allowing the seller to personalize the assortments offered to customers, and it replaces stochastic purchases with deterministic fluid consumption. Consequently, the CDLP can be a loose proxy for the revenue obtained in \textsc{da} from the same inventory vector, even under an MNL model.  \endgroup To see the gap, consider two products and one customer, with $r_1=M$, $r_2=1$, $v_1=v_0=1$, and $v_2=M$, and starting inventory $\vecc=(1,1)$. In $\fLP$, the optimal solution offers only product $1$ to the customer, so $\fLP(\vecc,1)=\frac{M}{2}$. In \textsc{da}, both products must be displayed because both are initially in stock, so $\fDA(\vecc,1)=\frac{2M}{M+2}$. Consequently, $\frac{\fDA(\vecc,1)}{\fLP(\vecc,1)}=\frac{4}{M+2}$, which approaches zero as $M\to\infty$.

\begingroup
We exploit the structural properties of the non-personalized MNL fluid process developed in Section~\ref{sec:structural} to handle these gaps. Given \(\vecc\), let \(\CFP\) denote the (possibly fractional) revenue-ordered inventory constructed from \(\vecc\) (Lemma~\ref{lem:3-rev}). Recall that $\fFP(\CFP,T)=\fLP(\vecc,T)$. If \(\CFP\) is integral and both \textsc{da} and the fluid process start from this inventory, $\fFP(\CFP,T)$ and $\fDA(\CFP,T)$ differ only through the stochasticity gap.  In Section~\ref{sec:da-transform-alg} we show that, under the MNL model, this gap is at most the factor $1-\frac1e$ for \emph{every} integral starting inventory, even though it can be arbitrarily large for general choice models. 
However, \(\CFP\) may be fractional and, by Lemma~\ref{lem:3-rev}, at most one product in \(\CFP\) has a fractional inventory level, which we call the \emph{bad} product. In Section~\ref{sec:3-bad} we show that one can round the bad product's inventory while losing at most a $\frac{1}{T+1}$ fraction of the fluid revenue. \endgroup
\subsubsection{The Stochasticity Gap without a Bad Product}  \label{sec:da-transform-alg}
We now bound the stochasticity gap between the revenue of the fluid process and \textsc{da} for any integral starting inventory.

\textbf{Single-product stochasticity bound.}
We begin with a single product $i$ with initial inventory $c_i$ and a single customer type. In this case, both the CDLP and \textsc{da} offer the product to every customer until either the inventory is depleted or the selling horizon ends. Let $p$ denote the probability that $i$ is purchased by one customer before it stocks out, so the consumption in $\fLP$ is $\min\{p T,c_i\}$. The following standard result relates the expected consumption of product $i$ in \textsc{da} to that in $\fLP$.

\begin{lemma}[Lemma 5.3 in \cite{alaei2012online}]
Let $\ZDA_i(p,T,c_i)$ denote the expected number of units of product $i$ sold in \textsc{da} when the initial inventory is $c_i$ and each of the $T$ customers purchases with probability $p$. Then: $\ZDA_i(p,T,c_i)\geq (1-\frac{1}{e})\min\{p T,c_i\}$.
\label{lem:4-uni}
\end{lemma}

With multiple products, however, the stochasticity gap could be very large, due to the dynamic cannibalization effect: along some sample paths, a lower-revenue product may remain in stock in \textsc{da} after it has already depleted in the fluid process, thereby continuing to cannibalize demand from higher-revenue products. \begingroup For a general choice model, this effect can make the stochastic revenue arbitrarily smaller than the fluid revenue; see Example~\ref{exm:fp-da-ranking-gap} in Appendix~\ref{app:4-1}. Perhaps surprisingly, under the MNL model, we show the single-product factor $1-1/e$ continues to hold for every integral starting inventory. \endgroup
\begin{lemma}\label{lem:4-product}
For every $\vecx\in \mathbb{Z}_+^{n}$ and $i \in \NN $, we have $\ZDA_i(\vecx,T)\geq (1-\frac{1}{e}) \ZFP_i(\vecx,T)$. 
\end{lemma}
\begingroup
To illustrate the main idea, we consider a system with two products and integral initial inventories $x_1$ and $x_2$; the proof for a general inventory vector is deferred to Appendix~\ref{app:prf-thm-pro}. We consider two cases. \endgroup

\textbf{Case 1: Both products remain available in the fluid process until time $T$.} \begingroup
Without loss of generality, consider product $1$.  Because both products remain available throughout the horizon, product $1$ is consumed at the constant rate $\phi(1,\{1,2\})$ in the fluid process. 
In the stochastic \textsc{da} process, the purchase probability of product $1$ is initially also $\phi(1,\{1,2\})$. Moreover, if product $2$ stocks out before product $1$, then removing product $2$ can only increase the purchase probability of product $1$. Therefore, whenever product $1$ remains available, its purchase probability in \textsc{da} is at least $\phi(1,\{1,2\})$. Then we directly apply the single-item bound from Lemma~\ref{lem:4-uni} to obtain $\ZDA_1(\vecx,T)
\geq
\left(1-\frac{1}{e}\right)\ZFP_1(\vecx,T)$. \endgroup

\textbf{Case 2: One product depletes before the other in the fluid process.}
Without loss of generality, suppose that product $2$ depletes before product $1$ in the fluid process. Product $2$ can be handled as in Case~1: until its depletion time, it is consumed in the fluid process at rate $\phi(2,\{1,2\})$, and in the stochastic process its purchase probability is at least this value whenever it is still available, so Lemma~\ref{lem:4-uni} again gives the factor $1-\frac1e$.

\begingroup
The main difficulty concerns product $1$. In the fluid process, once product $2$ is depleted, product $1$ is no longer cannibalized by product $2$. In the stochastic process, however, product $2$ may remain available beyond its fluid depletion time and may therefore continue to reduce the purchase probability of product $1$; we must show that the loss caused by this delayed stockout can be controlled. The key observation is that product $2$ has only $x_2$ units of inventory. Hence, on every sample path, at most $x_2$ of the $T$ arriving customers can purchase product $2$, so at least $T-x_2$ customers do not. Focus on any $T-x_2$ such customers. Conditional on a customer not purchasing product $2$, the MNL IIA property implies that, whenever product $1$ remains available, the customer purchases it with probability $\phi(1,\{1\})$. We can therefore couple these customers with a single-product stochastic system that has $x_1$ units of product $1$, a horizon of $T-x_2$, and purchase probability $\phi(1,\{1\})$; under this coupling, the sales of product $1$ in the original \textsc{da} system dominate those in the single-product system, that is, $\ZDA_1(\vecx,T)\geq \ZDA_1(x_1\vece_1,T-x_2)$. For the fluid process, the Separability Lemma (Lemma~\ref{lem:sep}) shows that the consumption of product $1$ starting from inventory $\vecx$ is exactly the consumption of product $1$ in a single-product fluid process with $x_1$ units over a period of \(T-x_2\), that is, $\ZFP_1(\vecx,T) = \ZFP_1(x_1\vece_1, T-x_2)$. Applying the single-product bound (Lemma~\ref{lem:4-uni}) to this auxiliary system and chaining the three relations, we obtain
\[
\ZDA_1(\vecx,T)
\geq \ZDA_1(x_1\vece_1,T-x_2) \geq \left(1-\frac{1}{e}\right)
\ZFP_1(x_1\vece_1,T-x_2) =
\left(1-\frac{1}{e}\right)\ZFP_1(\vecx,T).
\] \endgroup
The same reduction extends to a general inventory vector. Fix a product \(i\), and let \(\operatorname{Out}(i)\) be the set of products that deplete before \(i\) in the fluid process. These products can be purchased by at most \(\sum_{j\in\operatorname{Out}(i)}x_j\) customers on any sample path. Deleting those customers reduces the consumption of \(i\) to that in a system with a uniform purchase probability and \(T-\sum_{j\in\operatorname{Out}(i)}x_j\) effective arrivals, and applying Lemma~\ref{lem:4-uni} then gives the desired \(1-\frac1e\) bound.

\subsubsection{Rounding the Bad Product} \label{sec:3-bad}  \begingroup
It remains to round a possibly fractional revenue-ordered inventory \(\CFP\) to an integral inventory vector. Recall that the bad product is the unique fractional component of \(\CFP\); by construction, it is the threshold product with the lowest price among all products included in \(\CFP\). We show that rounding its inventory either up or down, and choosing the better of the two resulting inventories, retains at least a $(1-\frac{1}{T+1})$ fraction of the fluid revenue. The overall transformation step, which computes \(\CFP\) from an inventory vector \(\vecc\) and then rounds the bad product, is formalized in Algorithm~\ref{alg:trans-T}. \endgroup

{\SetAlgoNoLine
\begin{algorithm}[htp]
\caption{Transformation for \textsc{da} with deterministic $T$}\label{alg:trans-T}
\begingroup
\textbf{Input:} inventory $\vecc$, horizon $T$. \textbf{Output:} an integral inventory vector\;
Compute the revenue-ordered inventory $\CFP$ (Lemma~\ref{lem:3-rev}) and let $i$ be the lowest-priced product in its support, i.e., the largest index with $x^{\vecc}_i>0$\;
\endgroup
Let $\bar{\vecx}^{\vecc}$ and $\underline{\vecx}^{\vecc}$ be obtained from $\CFP$ by replacing $x^{\vecc}_i$ with $\lceil x^{\vecc}_i\rceil$ and $\lfloor x^{\vecc}_i\rfloor$, respectively\;
\Return{$\arg\max_{\vecx\in\{\bar{\vecx}^{\vecc},\,\underline{\vecx}^{\vecc}\}}\fFP(\vecx,T)$}\;
\end{algorithm}}

\begingroup
The transformation algorithm changes only the inventory of the bad product and returns the better of the rounded-up and rounded-down candidates. Since \(\vecc\) is integral and \(\CFP\preceq\vecc\), both candidates are componentwise bounded by \(\vecc\), so the transformation preserves any cardinality or budget constraint satisfied by \(\vecc\). The algorithm is also efficient: the revenue-ordered inventory \(\CFP\) can be computed in \(O(n\log n)\) time by running the fast CDLP algorithm in Appendix~\ref{app:fast_cdlp} and applying the construction of Lemma~\ref{lem:3-rev} to the optimal SBLP solution it induces, and after preprocessing their common depletion order, the fluid revenues of the two rounded candidates can be evaluated and compared in \(O(n)\) time using the sequencing algorithm in Appendix~\ref{app:alg-seq}. Thus, the overall running time is \(O(n\log n)\).
\endgroup

\begin{lemma}\label{lem:4-per-gap}
Given any inventory vector \(\vecc\in \mathbb{Z}^n_+\) and any real horizon $T>0$, let \(\CFP\) be the revenue-ordered inventory constructed from \(\vecc\), and let \(\vecc^{\operatorname{ALG}}\in\mathbb Z_+^n\) be the output of Algorithm~\ref{alg:trans-T}. Then 
\(
\fFP(\vecc^{\operatorname{ALG}},T)
\geq
\left(1-\frac{1}{T+1}\right)\fFP(\CFP,T).
\)
\end{lemma}

To prove the lemma, it suffices to show $\max\{\fFP(\bar{\vecx}^{\vecc}, T), \fFP(\underline{\vecx}^{\vecc},T)\}\geq (1-\frac{1}{T+1})\fFP(\CFP,T)$. \begingroup The Separability Lemma allows us to isolate the effect of rounding the fractional inventory of the bad product while keeping the consumption of all other products unchanged. The problem therefore reduces to a factor-revealing program that compares the original fractional inventory with its two integral roundings. Rounding the bad product down removes only the terminal consumption associated with this lowest-priced product. Rounding it up, by contrast, can only extend the terminal portion of the fluid process during which this product is consumed. Comparing these two losses and choosing the better rounding yields the factor $(1-1/(T+1))$. \endgroup A key observation used is that the bad product has the lowest price. The complete proof is given in Appendix~\ref{app:4-prf-T1}.
\begingroup
\paragraph{Overall Performance Guarantee for Deterministic $T$.}
We can now state the complete algorithm for deterministic $T$. When $T$ is small, we use the static assortment to avoid the larger loss of $1/(T+1)$ incurred by the transformation procedure. Let \(S^{\operatorname{static}}\) be an optimal cardinality-constrained static assortment and let \(\vecc^{\operatorname{static}}\) be the corresponding inventory vector with $c^{\operatorname{static}}_i=1$ if $i\in S^{\operatorname{static}}$ and $c^{\operatorname{static}}_i=0$ otherwise. The complete procedure is given in Algorithm~\ref{alg:da-det}.
{\SetAlgoNoLine
\begin{algorithm}[htp]
\caption{Threshold-Based \textsc{da} Algorithm for Deterministic $T$}\label{alg:da-det}
\textbf{Input:} cardinality $K$, horizon $T$, accuracy parameter $\epsilon$. \textbf{Output:} an integral inventory vector\;
\lIf{$T\leq 3$}{compute $S^{\operatorname{static}}$ and \Return{$\vecc^{\operatorname{static}}$}}
Run Algorithm~\ref{alg:da-thr} with input $(K,T,\epsilon)$ to obtain $\widehat{\vecc}$\;
\Return{the output of Algorithm~\ref{alg:trans-T} with input $(\widehat{\vecc},T)$}\;
\end{algorithm}}
\endgroup
\begingroup
\begin{theorem}\label{thm:gua-DA-det}
For any \(\epsilon\in(0,1)\), under a cardinality constraint and \(T\) customers, the threshold-based \textsc{da} algorithm has approximation factor
\(
0.319-\epsilon.
\)
\end{theorem}
\endgroup
\proof{Proof of Theorem~\ref{thm:gua-DA-det}.}
\begingroup
For \(T\geq1\), Theorem~\ref{thm:da-fluid-opt}, Lemma~\ref{lem:4-product} and Lemma~\ref{lem:4-per-gap} give the factor
\(
\left(1-\frac1e\right)^2(1-\frac{1}{T+1})-\epsilon
\)
relative to the fluid upper bound. Moreover, the static assortment has a factor of at least \(1/T\) relative to the optimal objective value since
\(
\max_{\|\vecc\|_1\leq K}\fDA(\vecc,T)
\leq
T R(S^{\operatorname{static}}).
\) Indeed, since $\vecc^{\operatorname{static}}$ stocks one unit of each product in $S^{\operatorname{static}}$, the first customer is offered $S^{\operatorname{static}}$ in full, so $\fDA(\vecc^{\operatorname{static}},1)=R(S^{\operatorname{static}})$; by the monotonicity of $\fDA(\vecc,t)$ in $t$, it follows that $\fDA(\vecc^{\operatorname{static}},T)\geq \fDA(\vecc^{\operatorname{static}},1)\geq \frac{1}{T}\max_{\|\vecc\|_1\leq K}\fDA(\vecc,T)$. Minimizing $\max \{\left(1-\frac1e\right)^2(1-\frac{1}{T+1}),  \frac{1}{T} \} $ over $T$ yields the uniform guarantee of
\(
\frac{4}{5}\left(1-\frac1e\right)^2-\epsilon\geq 0.319-\epsilon,
\)
with the minimum attained at $T=4$.
\Halmos
\endgroup
\endproof

\begingroup
\paragraph{Overall Runtime.} The first-stage approximation algorithm runs in $O(\frac{n}{\epsilon^2}\log^2(\frac{K}{\epsilon})\log K\log n)$ time (Theorem~\ref{thm:da-fluid-opt}). The runtime of the second stage is dominated by solving a CDLP, which can be solved using our fast near-linear-time algorithm in $O(n\log n)$ time (Appendix~\ref{app:fast_cdlp}). So the overall runtime is $O(\frac{n}{\epsilon^2}\log^2(\frac{K}{\epsilon})\log K\log n)$. 
\endgroup
\begingroup
\subsection{Stochastic \texorpdfstring{$T$}{T} with IFR Distributions}\label{sec:da-stochasticT}
We now extend the deterministic-horizon algorithm to a stochastic horizon
\(T\) that follows an IFR distribution. In both stages, the algorithm uses
the mean horizon \(\mu_T\). The first-stage optimization carries over directly: by Theorem~\ref{thm:da-fluid-opt}, Algorithm~\ref{alg:da-thr}
returns an inventory vector \(\widehat{\vecc}\) satisfying
\[
\fLP(\widehat{\vecc},\mu_T)
\geq
\left(1-\frac{1}{e}-\epsilon\right)
\max_{\lVert\vecc\rVert_1\leq K}
\fLP(\vecc,\mu_T).
\] Thus the main challenge lies in analyzing the transformation step that bounds the stochasticity gap. For deterministic $T$, both the fluid and stochastic processes operate over the same horizon. For stochastic $T$, in contrast, \(\fLP(\vecx,\mu_T)\) replaces the random $T$ by its mean, whereas
\(\EE_T[\ZDA_i(\vecx,T)]\) averages the stochastic consumption over different realizations of $T$.

As before, we begin with the single-product gap between \textsc{da} and $\fLP$. For stochastic $T$, the fluid consumption is $\min\{\mu_T p, c_i\}$, and the following bound follows from \cite{aouad2018greedy}.
\endgroup
\begin{lemma}[Lemma 11 in \cite{aouad2018greedy}]
Given any inventory vector $\vecc \in \mathbb{Z}_+^n$, a product $i \in \NN$ with purchase probability $p$, and $T$ following an IFR distribution with mean $\mu_T$, it holds that $\EE_{T}[\ZDA_i(p, T, c_i)] \geq \frac{1}{2} \min\{p \mu_T, c_i\}$.
\label{lem:4-uni-DT}
\end{lemma}

We show this $\frac{1}{2}$ bound is tight: it is approached by geometric distributions as the purchase probability vanishes; see Appendix~\ref{app:geo-12}. For non-IFR distributions, however, the gap can be arbitrarily large (example in Appendix~\ref{app:exm-nonIFR}).

\begingroup
Recall the reduction for deterministic $T$ in the two-product example: when product $2$ stocks out in the fluid process after selling its $a:=x_2$ units, we delete the $a$ customers who may purchase it and reduce to a single-product process over the remaining $T-a$ customers. With stochastic $T$, however, the realized horizon may already have ended before the deletion leaves any customers, i.e., \(T\leq a\), so the deterministic argument no longer applies. Instead, we lower bound the expected sales over the remaining $T-a$ customers on the event $\{T>a\}$. It suffices to establish a lower bound for a single product with a single unit and uniform purchase probability $p$. On the event $\{T>a\}$, the unit is sold if at least one of the $T-a$ remaining customers purchases it, which occurs with probability $1-(1-p)^{T-a}$. Thus the expected sales are $\EE_T\!\left[(1-(1-p)^{T-a}) \1\{T>a\}\right]$. The corresponding single-unit fluid consumption is $\min\{p(\mu_T-a),1\}$. \endgroup We establish novel stochastic
inequalities to show the $1/e$ approximation factor in this single-item setting, which may provide new insights for other revenue
management problems with IFR distributions. Interestingly, the $1/e$ ratio frequently appears in the pricing literature for IFR distributions \citep{dhangwatnotai2010revenue,chen2021revenue}. However,
to the best of our understanding, this is just a coincidence because the ratio we analyze involves
substantially different quantities.

\begin{lemma}
Let $p\in (0,1)$ be a uniform purchase probability. If $T$ follows an IFR distribution with finite mean $\mu_T$, then for any integer $a\in \mathbb{Z}_+$ with $a<\mu_T$, $\EE_T [ (1-(1-p)^{(T-a)})\1\{T>a\} ] \geq \frac{1}{e}\min\{p(\mu_T-a),1\}$.
\label{lem:4-Stoch-DT}
\end{lemma}
\vspace{-5pt}
We introduce the high-level idea of the proof, with full details in Appendix~\ref{app:IFRe-2}. 
We first formulate the factor-revealing program. Define $q_k = \Prb(T\geq k|T\geq k-1)$, for $k=1,2,3,\cdots$; then the sequence $\vecq=(q_1,q_2,\cdots )$ uniquely characterizes the distribution of $T$, $D_T$. Given that $D_T$ is an IFR distribution, we have $1\geq q_1\geq q_2\geq q_3\geq \cdots\geq 0$. We can express $\Prb(T\geq a)=\prod_{k=1}^{a} q_k $. \begingroup Thus $\EE_T[(1-(1-p)^{T-a})\1\{T>a\}]=q_1q_2\cdots q_a\,\EE_T[1-(1-p)^{T-a}\,|\,T\geq a]$. When $(\mu_T-a)p\leq 1$, the target lower bound equals $\frac{1}{e}p(\mu_T-a)$; the complementary case $(\mu_T-a)p>1$ reduces to the boundary value $p=\frac{1}{\mu_T-a}$ by the monotonicity of the left-hand side in $p$. \endgroup Therefore, Lemma~\ref{lem:4-Stoch-DT} is equivalent to showing that $\frac{1}{e}$ is a lower bound for the following factor-revealing program:
 \begin{align}
 \label{OPT:mu-c}
  & \inf_{\mu_T, a, p,\vecq}  \frac{q_1q_2 \cdots q_a \EE_T[(1-(1-p)^{T-a})|T\geq a]}{(\mu_T-a)p} \\
   & \text{s.t.} \quad   1\geq q_1\geq q_2\geq \cdots \geq 0, \quad  (\mu_T-a)p\leq 1, \quad  0< p\leq 1, \quad \mu_T-a>0, \quad a\in \mathbb{Z}_+  \nonumber  \\
    &\quad 
    \quad \mu_T = q_1 +q_1q_2+\cdots +q_1q_2\cdots q_a(1+q_{a+1}+q_{a+1}q_{a+2}+\cdots). \nonumber   
\end{align}
The constraints define the feasible values of $a$, $p$, and $\vecq$, with the last constraint derived from the tail-sum formula for the expectation over $T$. The proof has three steps. First, we simplify the program by showing that the optimal solution $\vecq^*$ stabilizes to a constant $q_i= b$ for any $i>a$, which implies that $T-a$, conditioned on $T>a$, follows a geometric distribution. Next, we identify the structure of the distribution given by $\vecq^*$, referred to as a ``shifted geometric distribution'' (see Lemma~\ref{lem:shiftgeo} below). This step is the most technically challenging, and we introduce a novel ``step-removal'' approach to prove it. Finally, we show that the infimum of the program is $\frac{1}{e}$ when restricted to shifted geometric distributions. 

\begin{lemma}\label{lem:shiftgeo}
 Given any feasible $\mu_T$ and $a$ with respect to \eqref{OPT:mu-c}, an optimal solution $\vecq^*$ exists and has the following structure: There exists a constant $b \in (0,1)$ such that either (i) there is a threshold index $1\leq k \leq a$, where $q_i = 1$ for all $i\leq k$ and $q_i = b$ for $ i\geq k+1$, or (ii) $q_i = b$ for all $i\geq 1$.
\end{lemma}

\begingroup
With Lemma~\ref{lem:4-Stoch-DT} established, it is straightforward to apply the same multi-product reduction procedure as in deterministic $T$ to bound the stochasticity gap with stochastic $T$. The complete proof is in Appendix~\ref{app:prf-lem-IFRe}. \endgroup

\begin{lemma}
Suppose $T$ follows an IFR distribution. For any inventory vector $\vecx\in \mathbb{Z}_+^n$ and any product $i \in \NN$, $\EE_T[\ZDA_i(\vecx,T)]\geq \frac{1}{e}\ZFP_i(\vecx,\mu_T)$.
\label{lem:IFRe}
\end{lemma}

\begingroup
\paragraph{Overall performance guarantee.}
As in the deterministic $T$ case, we use the optimal static assortment to handle small $\mu_T$ values: we use $\vecc^{\operatorname{static}}$ when $\mu_T<4.31$ and $\vecc^{\operatorname{ALG}}$ otherwise.
\begin{theorem}\label{thm:gua-DA-ifr}
Consider any accuracy level $\epsilon\in (0,1)$. Under a cardinality constraint, when $T$ follows an IFR distribution, \textsc{da} can be approximated within a factor of $0.188-\epsilon$.
\end{theorem}
The proof is in Appendix~\ref{app:5-da}.

Our stochastic-horizon analysis for \textsc{da} relies on the IFR property. Obtaining a constant-factor approximation for \textsc{da} when $T$ follows a general distribution remains an open problem.
\endgroup

\begingroup
\section{Dynamic Assortment with Personalization}
\label{sec:dap}
We now develop the algorithm and analysis for \textsc{dap}. The relax-and-transform scheme is the same as for \textsc{da}, but the bottleneck moves. Since the seller can personalize assortments, the gap between \textsc{dap} and the CDLP is essentially the stochasticity gap and is easier to bound; the difficulty now lies in the first step, maximizing the multi-type CDLP $\fLPM(\vecc,T)$ over the initial inventory. 

We first consider deterministic \(T\). Our starting point is the
single-type CDLP: using the non-personalized fluid process and its
structural results in Section~\ref{sec:structural}, we show that decreasing 
price order is a diminishing-returns submodular order. Although this property does
not extend directly to the multi-type CDLP, we  reformulate the
multi-type relaxation as an inventory-allocation problem and establish an 
analogous ordered diminishing-returns structure for the reformulated
objective. This leads to a threshold-augmentation algorithm
that achieves a $(\frac{1}{2}-\epsilon)$-approximation for the first-stage problem. \endgroup
\subsection{A DR-Submodular Order for the Single-Type CDLP}
\label{sec:dr-so}
We begin with the single-type CDLP $\fLP(\vecc,T)$ used for \textsc{da}, viewed as a function of the inventory vector $\vecc$. We show that products ordered by decreasing price satisfy a
diminishing-returns property: the marginal value of adding a lower-priced
product can only decrease after additional inventory of a higher-priced 
product is added.

Given a set function $g:2^{\NN}\rightarrow \RR_+$ on the ground set $\NN$, an order $\sigma$ of $\NN$ is a (strong) \emph{submodular order} (SO) for $g$ if
\begin{equation}
 g(A\cup\{ i, j\})- g(A \cup \{j\})\leq g(A\cup\{i\})- g(A), \quad  \forall A\subseteq \NN \text{ and } i>r_{\sigma}(A\cup \{j\}),
 \label{eq:so}
\end{equation}
where $i>r_{\sigma}(A\cup \{j\})$ means that item $i$ lies to the right of $A\cup \{j\}$ in the order $\sigma$ \citep{sof}. The SO property is defined for set functions, i.e., functions on the Boolean hypercube $\{0,1\}^{n}$. For functions on the integer lattice, we propose a natural generalization of SO
called the diminishing returns submodular order property (DR-SO), inspired by the notion of diminishing-returns submodularity \citep{soma2018maximizing}. For maximizing a monotone function with the SO property, 
 \cite{sof} gives a $(\frac{1}{2}-\epsilon)$-approximation under a cardinality constraint and a $(\frac{1}{3}-\epsilon)$-approximation under a budget constraint\footnote{\cite{sof} also gives a $(\frac{1}{2}-\epsilon)$-approximation under a budget constraint using partial enumeration.}. Under a cardinality constraint, the algorithm idea is simple: process the ground set in the submodular order and add every element whose marginal value exceeds a threshold, subject to the constraint. 

\textbf{DR-SO property.} For a multiset $A$ over $\NN$, let $\vecx_A\in \mathbb{Z}_+^n$ denote its multiplicity vector, whose $i$th component is the number of copies of element $i$ in $A$. An order $\sigma$ over $\NN$ is a DR-SO for $g:\mathbb{Z}_+^n\rightarrow\RR_+$ if
\[
g(\vecx_A+ \vece_i+\vece_j)-g(\vecx_A+\vece_j)\leq g(\vecx_A+ \vece_i)-g(\vecx_A)
\]
for every multiset $A$ over $\NN$ and every $i\ge r_{\sigma}(A\cup \{j\})$, where $r_{\sigma}(A\cup \{j\})$ is evaluated on the set of distinct elements of $A\cup\{j\}$. The threshold-based algorithms and performance guarantees for SO functions can be naturally generalized to DR-SO functions. In Section \ref{sec:dap-alloc}, we generalize
these algorithms to obtain new algorithms with stronger performance guarantees for solving \textsc{dap}.

We now establish the DR-SO property for the single-type CDLP. For $T=1$, the single-type CDLP $\fLP(\vecc,T)$ reduces to static assortment optimization, for which the descending price order is known to be a DR-SO of the expected revenue (Lemma~10 in \citet{sof}). We show that this holds for $T>1$. 
\begin{theorem}\label{thm:sof} 
For any given $T$, sorting products in descending order of price, with ties broken arbitrarily, is a DR-SO for $\fLP(\vecc, T)$. Formally, suppose the products are sorted in non-increasing price order, such that $r_1 \geq r_2 \geq \cdots \geq r_n$. Then for any $\vecc \in \mathbb{Z}_+^n$, and for any $i, j \in \NN$ such that $i \geq j$ and $i \geq k$ for every $k \in \NN$ with $c_k > 0$, we have:
\begin{equation}
\label{eq:sof}
\fLP(\vecc + \vece_j + \vece_i, T) - \fLP(\vecc + \vece_j, T) \leq \fLP(\vecc + \vece_i, T) - \fLP(\vecc, T).
\end{equation}
\end{theorem}
Proving Theorem~\ref{thm:sof} is technically challenging because it requires comparing the optimal values of LPs whose solutions may change discontinuously when one unit of inventory is added. \begingroup Our key idea is to use the fluid-process representation from Section~\ref{sec:structural} to reduce this LP sensitivity problem to a one-dimensional comparison between terminal segments of two fluid processes. The proof relies on two structural properties established in Section~\ref{sec:structural}. First, the generalized revenue-ordered structure (Lemma~\ref{lem:3-rev}) allows us to represent each CDLP optimum by a fluid process starting from a highly structured inventory vector. Because product $i$ is the lowest-priced product involved in \eqref{eq:sof}, this structure localizes the marginal comparison to product $i$. Second, the Separability Lemma represents the marginal value of adding product $i$ as the gain from replacing a terminal segment of the original fluid trajectory by certain sales of product $i$. We then show that adding the higher-priced product $j$ has two effects: it weakly shortens this replacement segment and weakly increases the revenue of the trajectory being replaced. Both effects reduce the marginal value of product $i$. 

Fix an inventory vector $\vecc$, let $N(\vecc)=\{k\in\NN : c_k>0\}$ denote its support, and let product $i$ lie weakly to the right of every product in $N(\vecc)\cup\{j\}$ in the decreasing-price order, so that $i$ is the lowest-priced product in the comparison. Write
\(
\vecc^{+i}=\vecc+\vece_i\), \(\vecc^{+j}=\vecc+\vece_j\), \(\vecc^{+ij}=\vecc+\vece_i+\vece_j\).
By Lemma~\ref{lem:3-rev}, for every inventory vector $\mathbf{d}$, the CDLP value $\fLP(\mathbf{d},T)$ equals the revenue of a fluid process starting from a revenue-ordered inventory vector $\mathbf{x}^{\mathbf{d}}\preceq \mathbf{d}$. Recall that this inventory vector has a threshold structure: all products above the threshold are fully stocked, all products below the threshold are excluded, and at most one threshold product is fractionally stocked.  \endgroup

If
\(
\fLP(\vecc^{+ij},T)-\fLP(\vecc^{+j},T)=0,
\)
then the desired DR-SO inequality is immediate by monotonicity of $\fLP$. \begingroup Hence suppose this marginal value is positive. The revenue-ordered structure then implies a strong localization property: in the four relevant fluid processes corresponding to $\vecc,\vecc^{+i},\vecc^{+j},\vecc^{+ij}$, every product priced above product $i$ is fully included in the revenue-ordered inventory. Consequently, the marginal CDLP comparison can be reduced to the effect
of increasing the inventory of product \(i\) in the non-personalized 
fluid process \textsc{fp}. \endgroup

\begingroup
Next we use the Separability Lemma to convert this isolated inventory perturbation into
a comparison between terminal segments of the two fluid processes. \endgroup Let
\(
\alpha_1
=
\ZFP_i(\mathbf{x}^{\vecc^{+ij}},T)
-
\ZFP_i(\mathbf{x}^{\vecc^{+j}},T)
\)
denote the additional amount of product \(i\) sold when one unit of
\(i\) is added to \(\vecc^{+j}\). \begingroup By the Separability Lemma, in terms of cumulative consumption, the fluid
process associated with \(\vecc^{+ij}\) is equivalent to running the
process associated with \(\vecc^{+j}\) over a horizon of \(T-\alpha_1\)
and then selling \(\alpha_1\) additional units of product \(i\), which
earn revenue \(\alpha_1 r_i\). \endgroup In particular, \(\fLP(\vecc^{+ij},T)
=
\fFP(\mathbf{x}^{\vecc^{+j}},T-\alpha_1)
+
\alpha_1r_i\), whereas
\(
\fLP(\vecc^{+j},T)
=
\fFP(\mathbf{x}^{\vecc^{+j}},T-\alpha_1)
+
\int_{T-\alpha_1}^{T}R(S^{+j}(t))\,dt,
\)
where $S^{+j}(t)$ is the assortment available at time $t$ in the fluid process starting with inventory $\vecx^{\vecc^{+j}}$. Subtracting the common prefix gives
\begin{equation}\label{eq:dap-diff}
\fLP(\vecc^{+ij},T)-\fLP(\vecc^{+j},T)
=
\int_{T-\alpha_1}^{T}
\bigl(r_i-R(S^{+j}(t))\bigr)\,dt.
\end{equation}

Similarly, let 
\(
\alpha_2
=
\ZFP_i(\mathbf{x}^{\vecc^{+i}},T)
-
\ZFP_i(\mathbf{x}^{\vecc},T).
\)
Then
\(
\fLP(\vecc^{+i},T)-\fLP(\vecc,T)
=
\int_{T-\alpha_2}^{T}
\bigl(r_i-R(S(t))\bigr)\,dt,
\)
where $S(t)$ is the assortment in the fluid process starting with inventory $\vecx^{\vecc}$.
 
\begin{figure}[htp]
\centering
\includegraphics[width=0.78\textwidth]{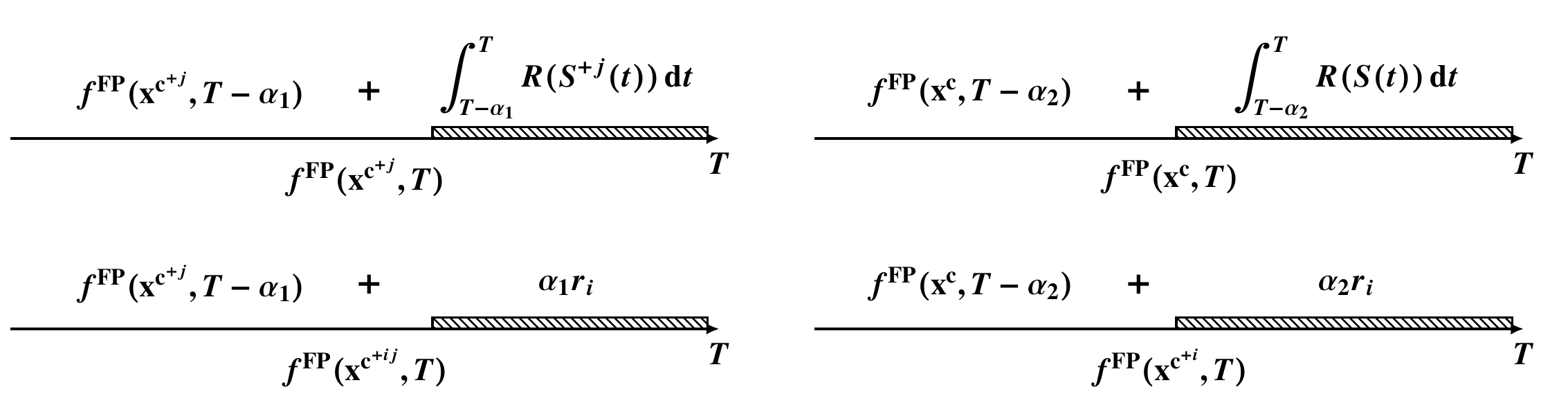}
\caption{Left: $\fFP(\CijFPs,T)-\fFP(\CjFPs,T)$. 
Right: $\fFP(\CiFPs,T)-\fFP(\CFP,T)$.}
\label{fig:mar-dif}
\end{figure}

 We illustrate this in Figure~\ref{fig:mar-dif}. The left
column compares adding product \(i\) to \(\vecc^{+j}\), while the right
column compares adding product \(i\) directly to \(\vecc\).

\begingroup
The DR-SO inequality now reduces to comparing these two terminal
segments. Adding the higher-priced product \(j\) before product \(i\)
reduces the marginal value of \(i\) in two ways. First, the corresponding 
terminal segment is weakly shorter:
\(
\alpha_1\leq \alpha_2.
\)
Intuitively, the extra unit of the higher-priced product $j$ competes with product $i$ in the MNL fluid process and can only reduce the amount of additional product $i$ that is profitably sold. Second, throughout the common terminal replacement interval $[T-\alpha_1,T]$, the assortment with the extra product $j$ has weakly larger expected revenue:
\(
R(S(t))\leq R(S^{+j}(t))\leq r_i.
\)
The upper bound \(R(S^{+j}(t))\leq r_i\) ensures that this
replacement has nonnegative marginal value.
\endgroup

\begingroup
Combining \eqref{eq:dap-diff} with $\alpha_1\leq \alpha_2$ and $R(S(t))\leq R(S^{+j}(t))\leq r_i$, we have: 
\[
\begin{aligned}
\fLP(\vecc^{+ij},T)-\fLP(\vecc^{+j},T)
&=
\int_{T-\alpha_1}^{T}
\bigl(r_i-R(S^{+j}(t))\bigr)\,dt \leq
\int_{T-\alpha_1}^{T}
\bigl(r_i-R(S(t))\bigr)\,dt\\
&\leq
\int_{T-\alpha_2}^{T}
\bigl(r_i-R(S(t))\bigr)\,dt=
\fLP(\vecc^{+i},T)-\fLP(\vecc,T).
\end{aligned}
\]
The first inequality follows from
\(R(S(t))\leq R(S^{+j}(t))\), while the second follows from
\(\alpha_1\leq\alpha_2\) and the nonnegativity of the integrand over the
additional interval. This proves \eqref{eq:sof} and hence the DR-SO
property. The complete proof is given in Appendix~\ref{app:sof}.
\endgroup

\subsection{An Allocation Reformulation and Threshold-based Augmentation Algorithm}
\label{sec:dap-alloc}
For \textsc{dap}, we optimize the multi-type CDLP $\fLPM(\vecc,T)$ in~\eqref{eq:multi-cdlp} of the first step. Perhaps surprisingly, the decreasing-price order is \emph{not} a DR-SO for $\fLPM$. The reason is that the multi-type LP determines not only how much of each product is
sold, but also how its inventory is allocated across customer types. Adding a unit of inventory can change this allocation, and the resulting reallocation breaks the diminishing-returns inequality. We give a concrete counterexample in Appendix~\ref{app:ce-SO-multitype}.

\begingroup
Our key idea is to make the inventory allocation an explicit decision 
rather than optimizing it implicitly within the multi-type CDLP. Let $\matc=(\matc_k)_{k\in\MM}$, where  $\matc_k\in\RR_+^n$ denotes the inventory allocated to type $k$ customers. For a fixed allocation \(\matc\), define 
\(
\hat{f}(\matc,T)=\sum_{k\in\MM}\fLP_k(\matc_k,\Lambda_k)\), where \( \Lambda_k=\sum_{t=1}^{T}\lambda_{k,t}
\)
and $\fLP_k$ is the single-type CDLP for type $k$ with its own attraction parameters $(v_{ki})_{i\in \NN}$. For a given inventory vector \(\vecc\), an allocation is feasible if \(\sum_{k\in\MM}C_{k,i}\leq c_i,
\forall i\in\NN.\) It is straightforward to show that optimizing $\hat f$ over all valid allocations recovers the multi-type CDLP.

\begin{lemma}\label{lem:alloc-equiv}
For any inventory $\vecc$ and horizon $T$,
\(
\fLPM(\vecc,T)=\max_{\matc}\ \Bigl\{\hat{f}(\matc,T):\ \textstyle\sum_{k\in\MM}C_{k,i}\le c_i\ \ \forall i\in\NN,\ \ \matc\ge 0\Bigr\}.
\)
\end{lemma}
The proof is in Appendix~\ref{app:hat}.
\endgroup

The advantage of this reformulation is that it separates the inventory-allocation decision from the marginal-value calculation. In $\fLPM(\vecc,T)$, the allocation across customer types is reoptimized whenever the
inventory vector changes. In $\hat{f}(\matc,T)$, by contrast, the allocation is specified explicitly by $\mathbf{C}$, so the marginal
effect of adding inventory is evaluated relative to a fixed current allocation. Moreover, the descending order of price is a DR-SO for $\hat{f}$. \begingroup Let \(\bar{\NN}:=\NN\times\MM\) be the expanded ground set,   where an element \((i,k)\in\bar{\NN}\)
corresponds to the coordinate \(C_{k,i}\), the amount of product \(i\)
allocated to customer type \(k\). Order the elements of \(\bar{\NN}\)
by nonincreasing product price \(r_i\), breaking ties across customer
types arbitrarily. By Theorem~\ref{thm:sof}, decreasing price order is a
DR-SO for each single-type objective
\(\fLP_k(\matc_k,\Lambda_k)\). Since
\(
\hat f(\matc,T)
=
\sum_{k\in\MM}\fLP_k(\matc_k,\Lambda_k)
\)
and the DR-SO property is preserved under addition, the decreasing-price order is also a DR-SO for
\(\hat f\) on the expanded ground set \(\bar{\NN}\); see Lemma~\ref{lem:sof-mint} in
Appendix~\ref{app:subhatf}. We therefore consider the following reformulation of the first-stage 
optimization problem:  \endgroup
\begin{align}\label{OPT-matc}
\vspace{-5pt}
\underset{\underset{k\in\MM}{\sum} \matc_k = \vecc,\,  \vecc\in \mathcal{F}}{\max}   \hat{f}( \matc,T).
\vspace{-4pt}
\end{align} The equality can be imposed without loss of generality by monotonicity. 
Note that the problem~\eqref{OPT-matc} jointly determines the quantity of each product
to stock and its allocation across customer types. 

\begingroup
Suppose we restrict $\matc$ to be integer-valued.
Recall that threshold-based algorithms provide constant-factor approximation guarantees for maximizing monotone DR-SO functions. In these algorithms, items are ordered in the submodular order and a unit of each item is added to the set if its marginal benefit exceeds a certain threshold value $\tau$. \endgroup In \eqref{OPT-matc}, however, we only require $\vecc$ to be integer-valued, not $\matc$.  
In general, the optimal solution $\matc^{\operatorname{OPT}}$ may be fractional and forcing $\matc$ to be integer-valued might not be a good approximation to \eqref{OPT-matc}. 
Consider, for example, one product and two
customer types. Suppose $\vecc^{\operatorname{OPT}}=(1)$, $C^{\operatorname{OPT}}_{1,1}=0.7$, $C^{\operatorname{OPT}}_{2,1}=0.3$. The optimal solution divides the single unit between the two types. An
integral threshold algorithm would instead consider only allocating the
entire unit to type \(1\) or the entire unit to type \(2\). It may happen
that neither integral allocation has marginal value at least \(\tau\),
even though the fractional allocation \((0.7,0.3)\) does. The integral
algorithm would then reject the product altogether.

To address this issue, we propose a novel threshold-based algorithm that allows for fractional allocation among customer types. \begingroup When examining each item, unlike the standard threshold-based augmentation algorithms for SO functions, which calculate the marginal benefit of adding whole units, our algorithm computes the optimal fractional allocation of each product across customer types that maximizes the marginal benefit relative to the current allocation.  
\endgroup

Since the horizon \(T\) remains fixed throughout, we suppress it below and write \(\hat f(\matc)\) for \(\hat f(\matc,T)\). The 
fixed-threshold routine is given below. The full algorithm tries a
geometric sequence of threshold values; see
Appendix~\ref{app:cdlp-m-T-alg}.

\begingroup
\begingroup
\SetAlgoNoEnd
\SetAlgoNoLine
\begin{algorithm}[htp]
    \caption{Cardinality-Constrained Threshold Add for DAP}
    \label{alg:cdlp-multi-type}
    \DontPrintSemicolon
    \KwIn{Cardinality $K$, threshold $\tau$, ground set $\NN$, number of customers $T$}
    \KwOut{Inventory vector $\vecc$, inventory allocation vector $\matc$}
    Initialize $\vecc = \mathbf{0}^n$, $\matc = \mathbf{0}^{mn}$\;
    \For{$i \in \NN$}{
        \While{$\|\vecc\|_1<K$}{
        Find the optimal allocation $\matz^*$ at the current $\matc$: $
        \matz^* = \argmax_{\matz \in \mathbb{R}^{mn}_+, \sum_{k \in \MM} Z_{k,i} = 1, \sum_{k \in \MM} Z_{k,j} = 0, \forall j\not= i} \hat{f}(\matc + \matz) - \hat{f}(\matc)
        $ \label{line:m-4}\;
        \lIf{$\hat{f}(\matc + \matz^*) - \hat{f}(\matc)<\tau$}{\textbf{break}}
        $\vecc \leftarrow \vecc + \mathbf{e}_i$, $\matc \leftarrow \matc + \matz^*$\;
        }
    }
    \Return{$\vecc,\matc$}\;
\end{algorithm}
\endgroup
At each iteration, line~\ref{line:m-4} computes the optimal allocation
of one additional unit of the current product \(i\), while holding the
existing allocation \(\matc\) fixed. If the resulting marginal value is
at least \(\tau\), the algorithm adds one unit of product \(i\) to
\(\vecc\) and updates \(\matc\) using the optimal fractional allocation
\(\matz^*\). The allocation subproblem is then recomputed because the
base allocation has changed. Once the marginal value falls below
\(\tau\), the algorithm moves to the next product. Consider the previous example where $\vecc^{\operatorname{OPT}}=(1)$, $C^{\operatorname{OPT}}_{1,1}=0.7$, and $C^{\operatorname{OPT}}_{2,1}=0.3$. In this new algorithm, we first calculate the optimal allocation for this item, which is $\matz^*=(0.7, 0.3)$. For the appropriate value of $\tau$ (the algorithm tries several different guesses), we have $\hat{f}(\matz^*)\geq \tau$ and we set $\vecc=(1)$ and $\matc=(0.7, 0.3)$. 
 While iterating over the expanded ground set may seem time-consuming (each product corresponds to $m$ elements of $\bar{\NN}$, each of which may receive up to $T$ units), we show in Appendix~\ref{app:acc} that it is sufficient to iterate over the $n$ distinct products instead.
\endgroup

We show that the full threshold-based augmentation algorithm (see Algorithm~\ref{alg:tau-m} in Appendix~\ref{app:cdlp-m-T-alg}) provides a $(\frac{1}{2}-\epsilon)$-approximate solution to \eqref{OPT-matc} under a cardinality constraint. The proof relies on the fact that the submodular-order structure of $\hat{f}$ continues to hold when the current allocation \(\matc\) is
fractional. Due to the equivalence between $\max_{\vecc\in \mathcal{F}}\fLPM(\vecc,T)$ and ${\max}_{\sum_{k\in\MM}, \matc_k = \vecc,\,  \vecc\in \mathcal{F}}   \hat{f}( \matc,T)$, the first-step optimization can be approximated within a factor of $1/2$. The complete proof is in Appendix~\ref{app:prf-3-ap-m-ra-T}. 

\begin{theorem}\label{thm:3-ap-m-ra-T}
For any $\epsilon\in(0,1)$, the problem $\max_{\vecc\in \mathcal{F}}\fLPM(\vecc,T)$ can be approximated within a factor of $\frac{1}{2}-\epsilon$ under a cardinality constraint.
\end{theorem}

\begingroup
We then note the difference between our approach and that of \citet{bai2022coordinated}. They 
introduce a further surrogate relaxation of the CDLP whose objective is
monotone submodular, and incur an additional constant-factor loss when
relating the surrogate value back to the CDLP. We instead work with the multi-type CDLP directly through the equivalent allocation reformulation \eqref{OPT-matc} and leverage the DR-SO property from Theorem~\ref{thm:sof}.
\endgroup
\vspace{-10pt}
\subsection{The Online Assortment Policy and Overall Performance Guarantee}
\label{sec:dap-online}

We use standard online assortment policies in the second step. The policy depends on whether the horizon $T$ is deterministic or stochastic. The results in this subsection apply to any choice model satisfying \emph{substitutability}: for any $S'\subseteq S$ and $i\in S'$, $\phi(i,S')\ge \phi(i,S)$. 

\begingroup
 For deterministic $T$, we use a randomized policy that samples assortments from the multi-type CDLP solution. This policy has been used in \textsc{dap} in \cite{bai2022coordinated}, where a $\frac{1}{2}$-competitive ratio is established. By extending the analysis of \cite{ma2021dynamic} from stationary customer arrivals ($\lambda_{k,t}$ is the same for all $t$) to nonstationary arrivals, we strengthen this guarantee to $1-1/e$.
\endgroup

\vspace{-5pt}
{\SetAlgoNoLine%
  \begin{algorithm}[htp]
     \caption{Assortment Policy for Deterministic $T$}\label{alg:per-T}
 Solve the multi-type CDLP and obtain an optimal solution $\{y_k(S)\}_{k\in \MM,S\subseteq \NN}$\;
 For each time $t\in [T]$, when a type $k$ customer arrives, sample an assortment $S_t$ from the distribution where $\Prb(S_t=S)=y_k(S)$ for all subsets $S\subseteq \NN$\;
 Show assortment $S_t\cap N_t$, where $N_t$ is the set of available products at time $t$.
  \end{algorithm}}
\vspace{-10pt}
\begin{lemma}\label{lem:tran-d}
Given any $\vecc\in \mathbb{Z}_+^n$, any number $T$ of customers, and any choice model with the substitutability property, let $\fIND(\vecc,T)$ denote the total expected revenue obtained by the assortment policy above. Then, $\fIND(\vecc,T)\geq (1-\frac{1}{e}) \fLPM(\vecc,T)$.
\end{lemma}
The proof is given in Appendix~\ref{app:thm-tran-d}.

\paragraph{Overall Performance Guarantees}
Theorem~\ref{thm:3-ap-m-ra-T} returns an inventory $\widehat{\vecc}$ with $\fLPM(\widehat{\vecc},T)\ge(\tfrac12-\epsilon)\max_{\vecc\in\mathcal F}\fLPM(\vecc,T)$, and the sampling policy of Lemma~\ref{lem:tran-d} earns at least $(1-\tfrac1e)\fLPM(\widehat{\vecc},T)$. Since $\fLPM(\vecc,T)$ upper bounds the expected revenue of the optimal personalization policy (Lemma~\ref{lem:multi-ub}), chaining the two factors gives
\(
\sup_\pi\fpi(\widehat{\vecc},T)\ \ge\ \Bigl(\tfrac12(1-\tfrac1e)-\epsilon\Bigr)\max_{\vecc\in\mathcal F}\sup_\pi\fpi(\vecc,T).
\)
\begin{theorem}\label{thm:gua-DAP}
Consider any accuracy level $\epsilon\in(0,1)$. Under a cardinality constraint, when $T$ is deterministic, \textsc{dap} can be approximated within a factor of $\tfrac12(1-\tfrac1e)-\epsilon$.
\end{theorem}
The proof is in Appendix~\ref{app:5}.
\begingroup
\paragraph{Runtime summary. }For \textsc{dap}, the first-stage algorithm iterates over the items of the expanded ground set and computes the maximum marginal benefit of each item. We improve the efficiency by iterating only over the $n$ distinct products rather than the expanded ground set, determining the inventory level of each product through binary search; this reduces the query complexity from $O(\frac{nmT}{\epsilon}\log K)$ to $O(\frac{n}{\epsilon}\log K)$. The total runtime is $O(mT+\frac{n}{\epsilon}\log K(mn+\log K)\log(mn))$, see Appendix~\ref{app:acc} for details. The second stage solves the multi-type CDLP once via its polynomial-size sales-based formulation \eqref{eq:multi-sblp} and converts each type's solution into an assortment mixture using the recovery and Sequencing procedures of Appendices~\ref{app:fast_cdlp} and \ref{app:alg-seq}, in $O(mn\log n)$ additional time plus one LP solve.
\endgroup
\vspace{-10pt}
\subsection{Stochastic \texorpdfstring{$T$}{T} with Arbitrary Distributions}
\vspace{-4pt}
\label{sec:dap-stochasticT} \begingroup
We now allow $T$ to follow an arbitrary known distribution $D_T$, including high-variance and non-IFR distributions. If we continue to use the mean-horizon fluid relaxation as in \textsc{da}, i.e., $\fLPM(\vecc,\mu_T)$, it can be an arbitrarily loose fluid upper bound \citep{bai2023fluid, aouad2022nonparametric}. Instead, we optimize the distribution-aware relaxation $\EE_T[\fLPM(\vecc,T)]$, which retains the full distribution of $T$ and relaxes only the choice stochasticity. The threshold-augmentation algorithm extends directly to this objective: the DR-SO properties hold for each realization of $T$ and are therefore preserved under expectation. Hence, Theorem~\ref{thm:3-ap-m-ra-T} carries over with the same $\tfrac12-\epsilon$ first-stage approximation factor (See Appendix \ref{app:multi-stochT} for a formal proof). For the online assortment decisions, we use the \emph{greedy} policy, which offers each arriving customer the revenue-maximizing assortment among the products currently in stock. The policy does not require knowledge of the realized horizon $T$ and achieves a $\frac{1}{2}$ guarantee relative to the distribution-aware relaxation.
\endgroup
\vspace{-10pt}
\begin{lemma}[\cite{golrezaei2014real}]\label{lem:4-DAP-stoc}
For any $\vecc\in \mathbb{Z}_+^n$ and any distribution of $T$, let $\EE_T[f^{\textsc{Greedy}}(\vecc,T)]$ denote the expected revenue obtained by the online greedy policy given inventory $\vecc$. Then, $\EE_T[f^{\textsc{Greedy}}(\vecc,T)]\geq \frac{1}{2}\EE_T[\fLPM(\vecc,T)]$.
\end{lemma}
\vspace{-4pt}
\begingroup
Combining this $\tfrac12$ factor with the $\tfrac12-\epsilon$ first-stage guarantee yields the following result.
\vspace{-4pt}
\begin{theorem}\label{thm:gua-DAP-stoch}
Consider any accuracy level $\epsilon\in(0,1)$. Under a cardinality constraint, when $T$ follows an arbitrary known distribution $D_T$, \textsc{dap} can be approximated within a factor of $\frac{1}{4}-\epsilon$.
\end{theorem}
\vspace{-4pt}
The proof is in Appendix~\ref{app:multi-stochT}. Here we assume the ability to calculate the expectation of $\hat{f}$, which is a standard assumption, as seen in \cite{bai2023fluid} and \cite{aouad2022nonparametric}. \endgroup
\begingroup
\section{Extensions}\label{sec:extensions}
In this section, we extend our results in two directions. First, we consider budget constraints in place of
cardinality constraints. Second, for \textsc{da}, we develop a dynamic-programming-based first-stage algorithm
that improves the approximation factors at the cost of a larger runtime. \endgroup\vspace{-10pt}
\subsection{Budget Constraints}\label{sec:ext-budget}
For \textsc{da}, we apply the standard algorithm for DR-submodular maximization under a budget constraint \citep{soma2018maximizing}, which yields the same first-stage approximation guarantee as under a cardinality constraint. The transformation step in Section~\ref{sec:da-transform} is independent of the form of the inventory constraint and preserves budget feasibility. Thus, the overall guarantees for \textsc{da} carry over to the budget-constrained setting.

For \textsc{dap}, we generalize the threshold-augmentation algorithm by scaling the threshold with the budget consumption $b_i$ of each product, and the first-stage factor becomes $\frac{1}{3}-\epsilon$ for both deterministic and stochastic $T$ (Theorems~\ref{thm:3-ap-m-bud-T} and \ref{thm:3-ap-m-bud-stoch-T} in Appendices~\ref{app:multi-budget} and \ref{app:multi-stochT}, respectively). Combining these first-stage guarantees with the same transformation and online policies as in Sections~\ref{sec:da} and \ref{sec:dap} yields the following guarantees. In particular, the \textsc{da} guarantees resolve the open problem posed by \cite{aouad2018greedy} of finding a constant-factor approximation algorithm for \textsc{da} under a budget constraint. The proofs of the following guarantees are in Appendix~\ref{app:5}.

\begingroup
\begin{theorem}\label{thm:gua-budget-da}
Consider any accuracy level $\epsilon\in (0,1)$ under a budget constraint.
(i) When $T$ is deterministic, \textsc{da} can be approximated within a factor of $0.319-\epsilon$.
(ii) When $T$ is stochastic and follows a distribution with the IFR property, \textsc{da} can be approximated within a factor of $0.188-\epsilon$. 
\end{theorem}
\begin{theorem}\label{thm:gua-budget-dap}
Consider any accuracy level $\epsilon\in (0,1)$ under a budget constraint.
(i) When $T$ is deterministic, \textsc{dap} can be approximated within a factor of $\frac{1}{3}(1-\frac{1}{e})-\epsilon$.
(ii) When $T$ is stochastic and follows an arbitrary known distribution, \textsc{dap} can be approximated within a factor of $\frac{1}{6}-\epsilon$.
\end{theorem}
\vspace{-5pt}
\subsection{A Dynamic-Programming-Based First-Stage Algorithm for \textsc{da}}\label{sec:ext-dp}
Recall that the threshold-based Algorithm~\ref{alg:da-thr} yields a factor of $1-\frac{1}{e}-\epsilon$ in the first stage. Inspired by \cite{chen2022approximation}\footnote{This result is included in an earlier version of the paper that was shared with us by Jake Feldman. Unfortunately, the final published manuscript does not include this result. We reproduce the result from that unpublished preprint version of \cite{chen2022approximation}, with permission from the authors.}, we propose a dynamic-programming-based algorithm that first guesses the quantity of no-purchase customers $y_0$ on a geometric grid as in the threshold-based algorithm. For each fixed \(y_0\), however, instead of applying threshold-based augmentation for DR-submodular functions, it directly solves the resulting problem using dynamic programming over discretized cumulative sales and discretized revenue. The algorithm achieves a \((1-\epsilon)\)-approximation and is a fully polynomial-time approximation scheme. We refer to this dynamic-programming-based algorithm as {\small DP}.
\vspace{-5pt}
\begin{theorem}\label{thm:dp-fluid-opt}
For any accuracy level $\epsilon\in (0,1)$, under the cardinality constraint $\mathcal F=\{\vecc:\|\vecc\|_1\leq K\}$, {\small DP} returns an inventory vector $\widehat{\vecc}\in\mathbb Z_+^n$ with $\|\widehat{\vecc}\|_1\leq K$ and
\(
\fLP(\widehat{\vecc},T)\geq (1-\epsilon)\max_{\vecc\in \mathcal{F}} \fLP(\vecc,T).
\)
Its arithmetic running time is $\tilde{O}\left(\frac{n^4}{\epsilon^4}\right)$, where $\tilde{O}(\cdot)$ suppresses factors polylogarithmic in the model parameters and in $\frac{1}{\epsilon}$, and its bit complexity is polynomial in the binary input size and $\frac{1}{\epsilon}$.
\end{theorem}
The algorithm, the proof, and the running-time derivation are given in Appendix~\ref{app:fptas}. Using {\small DP} as the first stage improves the overall constants for \textsc{da} at the cost of a larger runtime:

\begin{theorem}\label{thm:thr-dp}
Suppose we use the {\small DP} algorithm for the first-stage approximation. Consider any accuracy level $\epsilon\in (0,1)$.
When $T$ is deterministic, \textsc{da} can be approximated within a factor of $0.474-\epsilon$ under a cardinality constraint. When $T$ is stochastic and follows a distribution with the IFR property, \textsc{da} can be approximated within a factor of $0.268-\epsilon$ under a cardinality constraint.
\end{theorem}
The proof is in Appendix~\ref{app:5}.

For \textsc{dap} with a single customer type, the multi-type CDLP reduces to the single-type CDLP~\eqref{sblp}, so {\small DP} can also be used in its first stage. Combining it with the sampling policy of Lemma~\ref{lem:tran-d} approximates \textsc{dap} with deterministic $T$ within a factor of $1-\frac{1}{e}-\epsilon$, improving on the factor $\frac{1}{2}(1-\frac{1}{e})-\epsilon$ of Theorem~\ref{thm:gua-DAP}.
\vspace{-3pt}
\section{Numerical Experiments}\label{sec:exp}
We conduct numerical experiments to evaluate the empirical performance and computational efficiency of our algorithms against existing benchmarks for both \textsc{da} and \textsc{dap}. In this section, we focus on \textsc{da}, the setting in which our algorithms give the largest runtime improvement over prior work. Our goal is to test whether this improvement carries over to the implementation and whether it comes at a cost in solution quality, even though our approximation guarantees also improve on those of \citet{aouad2018greedy}. The corresponding experiments for \textsc{dap} are deferred to Appendix~\ref{app:experiments}.
\paragraph{Instance generation.}
Our instance generation procedures for both \textsc{da} and \textsc{dap} are unified, and largely follow those of \citet{aouad2018greedy} and \citet{bai2022coordinated}, previous work on \textsc{da} and \textsc{dap}, respectively. Product revenues are drawn independently from
$\operatorname{Unif}[0,10]$, and the MNL weights are drawn from $\operatorname{Unif}[1,10]$. The outside-option weight is calibrated so that the no-purchase probability is $0.05$ when all products are offered; equivalently, all attraction weights are divided by the calibrated outside-option weight, matching the normalization $v_0=1$ of Section~\ref{sec:2}.  We vary the number of products over $n\in\{50,100\}$ and the horizon length over $\{100,200,400,800,1600\}$, and we consider two horizon models: a deterministic horizon $T$, and a Poisson horizon with mean $\mu_T$ truncated at $3\mu_T$. Truncation preserves the IFR property, so the guarantees for IFR horizons apply. This gives the $20$ configurations. For the cardinality constraint, we draw $\eta\sim\operatorname{Unif}[0.2,0.5]$ and set $K=\lfloor\eta\,\EE[T]\,p_{\rm buy}(S^{\rm myopic})\rfloor$. Here $p_{\rm buy}(S)=\sum_{i\in S}\phi(i,S)$ is the probability
of purchasing any product from assortment $S$, $S^{\rm myopic}$ is the unconstrained revenue-maximizing static assortment and $\EE[T]$ is the nominal horizon mean ($T$ or $\mu_T$). For each of the $(n,\text{horizon length})$ pairs we generate $20$ instances and evaluate each instance under both horizon models, and then report the average results in Table~\ref{tab:exp-da-summary}.

\paragraph{Benchmarks.}
We evaluate both variants of our algorithm: the threshold-based first-stage algorithm ({\small THR}) and the dynamic-programming-based first-stage algorithm ({\small DP}). We compare them with the selective-greedy approximation algorithm of \citet{aouad2018greedy} ({\small SG}) and the exponential-time near-optimal algorithm of \citet{aouad2022stability} ({\small AS}). The full implementation details are in Appendix \ref{app:experiments}.

\paragraph{Performance measures and computing environment.} We report both average revenue and runtime. The revenue is computed by taking the average of $10{,}000$ sample paths.
Following the instance-level comparison used by \citet{aouad2018greedy},
since computing the exact optimum is intractable, we normalize, on each instance, the revenue
of each algorithm by the highest revenue achieved among the methods on that instance, and report the average of these ratios over the $20$ instances of each configuration. Each method has a $1{,}200$-second limit; if an algorithm reaches the limit, we record the best solution found so far. The experiments are run on AWS Linux x86-64 machines with 24 worker processes, 16 vCPUs, and 32 GiB of memory.

\begin{table}[htbp]
\centering
\scriptsize
\setlength{\tabcolsep}{3.2pt}
\renewcommand{\arraystretch}{0.92}
\resizebox{0.9\textwidth}{!}{%
\begin{tabular}{llrrrrrrrrr}
\toprule
& & & \multicolumn{4}{c}{Normalized revenue (\%)} &
\multicolumn{4}{c}{Mean runtime (seconds)}\\
\cmidrule(lr){4-7}\cmidrule(lr){8-11}
Horizon & $n$ & $T$ or $\mu_T$ &
THR & DP & SG & AS &
THR & DP & SG & AS\\
\midrule

\multirow{10}{*}{Poisson}
& 50
& 100
& 98.4 & 98.7 & 96.5 & 45.0
& 49.8 & 1{,}200.0 & 1{,}200.0 & 1{,}181.0\\
&
& 200
& 99.4 & 99.8 & 97.0 & -
& 60.1 & 1{,}200.0 & 1{,}200.0 & -\\
&
& 400
& 98.6 & 99.1 & 99.2 & -
& 66.5 & 1{,}200.0 & 460.0 & -\\
&
& 800
& 99.1 & 99.5 & 99.8 & -
& 71.1 & 1{,}200.0 & 484.6 & -\\
&
& 1600
& 99.5 & 99.9 & 96.6 & -
& 75.2 & 1{,}200.0 & 1{,}181.7 & -\\

\addlinespace
& 100
& 100
& 98.0 & 98.1 & 97.8 & 25.0
& 111.5 & 1{,}200.0 & 1{,}200.0 & 1{,}200.0\\
&
& 200
& 98.6 & 98.8 & 97.1 & -
& 132.7 & 1{,}200.0 & 1{,}200.0 & -\\
&
& 400
& 97.9 & 98.4 & 99.6 & -
& 155.6 & 1{,}200.0 & 576.7 & -\\
&
& 800
& 99.3 & 99.9 & 97.0 & -
& 169.8 & 1{,}200.0 & 1{,}196.0 & -\\
&
& 1600
& 99.1 & 100.0 & 96.3 & -
& 180.1 & 1{,}200.0 & 1{,}200.0 & -\\

\midrule

\multirow{10}{*}{Deterministic}
& 50
& 100
& 98.5 & 98.8 & 96.0 & 64.9
& 76.8 & 1{,}200.0 & 1{,}200.0 & 1{,}083.9\\
&
& 200
& 99.0 & 99.4 & 96.2 & -
& 90.6 & 1{,}200.0 & 1{,}200.0 & -\\
&
& 400
& 98.9 & 99.4 & 99.1 & -
& 100.6 & 1{,}200.0 & 481.5 & -\\
&
& 800
& 99.3 & 99.7 & 99.7 & -
& 109.3 & 1{,}200.0 & 601.2 & -\\
&
& 1600
& 99.6 & 100.0 & 96.5 & -
& 120.9 & 1{,}200.0 & 1{,}200.0 & -\\

\addlinespace
& 100
& 100
& 98.4 & 98.5 & 97.3 & 25.0
& 169.7 & 1{,}200.0 & 1{,}200.0 & 1{,}200.0\\
&
& 200
& 97.9 & 98.1 & 95.9 & -
& 198.7 & 1{,}200.0 & 1{,}200.0 & -\\
&
& 400
& 98.5 & 99.0 & 99.7 & -
& 227.3 & 1{,}200.0 & 621.4 & -\\
&
& 800
& 99.5 & 100.0 & 96.7 & -
& 254.7 & 1{,}200.0 & 1{,}200.0 & -\\
&
& 1600
& 99.2 & 100.0 & 96.2 & -
& 282.8 & 1{,}200.0 & 1{,}200.0 & -\\

\midrule
\multicolumn{3}{r}{\textbf{Average}}
& 98.8 & 99.3 & 97.5 & -
& 135.2 & 1{,}200.0 & 1{,}000.2 & -\\
\bottomrule
\end{tabular}}
\caption{\textsc{da} performance for Poisson and deterministic horizons. Each entry is averaged over $20$ instances. For {\small AS}, $-$ marks the configurations in which it returned no feasible inventory vector on any instance.}
\label{tab:exp-da-summary}
\end{table}

Table~\ref{tab:exp-da-summary} shows that {\small THR} offers the strongest
overall trade-off between revenue and runtime. Across all
configurations, {\small THR} achieves an average normalized revenue of $98.8\%$ while having the lowest mean runtime, $135.2$ seconds. In comparison, {\small DP} achieves the highest average normalized revenue, $99.3\%$, but exhausts the $1{,}200$-second budget in every configuration, and {\small SG} achieves $97.5\%$ normalized revenue with a mean runtime of $1{,}000.2$ seconds.
Moreover, {\small THR} has the smallest average runtime in every configuration, and exceeds {\small SG}'s normalized revenue in $14$ of the $20$ configurations; in the remaining six configurations it trails {\small SG} by at most $1.7$ percentage points. {\small AS} encounters a more severe computational bottleneck: it returns a feasible inventory vector within the $1{,}200$-second time limit for only a few instances when $T$ or $\mu_T$ is $100$, and for none when $T$ or $\mu_T\in \{200,400,800,1600\}$. On each instance where {\small AS} returns no feasible inventory vector within the time limit, its normalized revenue is recorded as zero. The low average normalized revenues reported for {\small AS} therefore reflect these failures rather than poor solution quality. In fact, on almost every instance where {\small AS} returns a feasible inventory vector, its revenue is the highest among all methods. We therefore do not report an overall average for {\small AS}.

We also note that {\small DP} attains the highest average normalized revenue. However, it pays for this accuracy in runtime under the common $1{,}200$-second time budget: the revenue it gains over {\small THR} is at most $0.9$ percentage points in any configuration, at roughly nine times the average runtime. Overall, {\small THR} combines the scalability of a fast first-stage algorithm with revenue that is consistently close to the best observed across all methods.

\section{Conclusion}\label{sec:conclusion}
In this work, we have developed a unified algorithmic framework for both \textsc{da} and \textsc{dap} under the MNL choice model. This framework outperforms existing polynomial-time approximation algorithms in the settings we study in terms of speed and approximation guarantees, while also accommodating a stochastic number of customers and both cardinality and budget constraints; in particular, it resolves the open problem of \cite{aouad2018greedy} on budget-constrained \textsc{da}. From a technical perspective, we believe that our algorithmic ideas and the novel structural properties of the MNL model can shed light on algorithm design for assortment planning and broader revenue management problems. Further research could consider extending this approximation framework to other choice models, such as the Nested Logit and Markov Chain choice models. 
In work subsequent to ours, \cite{huang2024basic} introduce a joint inventory and assortment planning problem for reusable resources. They extend techniques from \cite{bai2022coordinated}, and it would be interesting to see if our results can also be used to achieve better approximation ratios in this new setting. Technically, we note that settling the computational hardness of \textsc{da} and developing an approximation algorithm for \textsc{da} when $T$ follows a general distribution are also interesting open problems for future work. Finally, incorporating the usual cost components is a natural direction, and nonasymptotic approximation guarantees for cost-aware versions of \textsc{da} and \textsc{dap} remain open. 
\endgroup

%
%
%

\theendnotes



\bibliographystyle{informs2014} 
\bibliography{papers.bib} 


\newpage

\begin{APPENDICES}


\appendixtableofcontents
\newpage
\renewcommand{\thesubsection}{\Alph{section}.\arabic{subsection}}

\section{Additional Formulations, Algorithms and Proofs in Section \ref{sec:2}}\label{app:sec2}
This appendix collects the proof of fluid upper bounds, the CDLP formulations omitted from the body, and the Sequencing algorithm used to compute the fluid consumption process.
\hypersetup{bookmarksdepth=section}
\subsection{Proposition~\ref{prop:ub}}\label{app:prop-ub}
\begin{proposition}\label{prop:ub}
$\fLP(\vecc,\mu_T)$ serves as an upper bound on $\EE_T[\fDA(\vecc,T)]$.
\end{proposition}

\proof{Proof of Proposition~\ref{prop:ub}.}
Consider any assortment policy $\pi$ that, based on past customer choices, offers each customer an assortment of in-stock products, i.e., $S_t\subseteq N_t$ for every $t$. Recall from Section~\ref{sec:2-dadap} that the \textsc{da} process corresponds to the fixed policy $S_t=N_t$.

Let $y_\pi(S)$ denote the expected number of times assortment $S$ is offered under policy $\pi$. This is given by
$y_\pi(S) = \EE\left[ \sum_{t=1}^T \1\{ S_t = S \} \right]$,
where the expectation is taken over both $T$ and the purchase history.

For each product $i$, define the expected number of purchases as $x_i = \sum_{S \subseteq \NN} y_\pi(S) \phi(i, S)$. The expected revenue under policy $\pi$ is then
$\EE_T[\fpi(\vecc,T)] = \sum_{S \subseteq \NN} y_\pi(S) R(S)$.

Define $y(S)=\frac{1}{\mu_T}y_\pi(S)$ for all $S\subseteq \NN$. We now show that $\{y(S)\}_{S \subseteq \NN}$ is a feasible solution to $\fLP(\vecc,\mu_T)$ by verifying the following constraints: first, since $\pi$ only offers in-stock products, the sales of product $i$ never exceed its initial inventory $c_i$ on any sample path, so
$x_i = \sum_{S \subseteq \NN} y_\pi(S) \phi(i, S) \leq c_i$,
which implies $\mu_T \sum_{S\subseteq\NN}y(S)\phi(i,S)\leq c_i$. Second, the total expected number of assortments offered equals the expected number of customers: $\sum_{S \subseteq \NN} y_\pi(S) = \mu_T$,
which implies $\sum_{S \subseteq \NN} y(S) = 1$. Finally, $y(S) \geq0$ for all $S \subseteq \NN$ follows from $y_\pi(S) \geq0$.

These constraints confirm that $\{y(S)\}_{S \subseteq \NN}$ is feasible for $\fLP(\vecc,\mu_T)$. Therefore,
$\EE_T[\fpi(\vecc,T)] = \sum_{S \subseteq \NN} y_\pi(S)R(S) = \mu_T \sum_{S\subseteq \NN} y(S)R(S) \leq \fLP(\vecc,\mu_T)$.
Since this holds for the \textsc{da} policy $S_t=N_t$ in particular, we conclude that $\fLP(\vecc,\mu_T) \geq\EE_T[\fDA(\vecc,T)]$.
\Halmos
\endproof
For deterministic $T$, taking $D_T$ to be degenerate at $T$ gives $\fLP(\vecc,T)\geq \fDA(\vecc,T)$.
\subsection{Lemma~\ref{lem:multi-ub}}\label{app:multi-ub}
The following lemma is directly from Lemma 1 in \cite{golrezaei2014real}. 
\begin{lemma}\label{lem:multi-ub}
Fix an inventory vector $\vecc\in\mathbb{Z}_+^n$. For deterministic $T$, the expected revenue of any assortment policy, including a policy with knowledge of the full sequence of customer types, is at most $\fLPM(\vecc,T)$. When $T$ is random, the expected revenue of any such policy is at most $\EE_T[\fLPM(\vecc,T)]$.
\end{lemma}
\proof{Proof of Lemma~\ref{lem:multi-ub}.}
Condition on the realized horizon $T=t$ and on the realized type sequence, and let $N_k$ denote the number of type-$k$ customers among the $t$ arrivals. Conditional on this realization, the expected revenue of any policy that offers only in-stock products is at most the optimal value of the choice-based formulation \eqref{eq:multi-cdlp} with each type mass $\Lambda_k$ replaced by $N_k$ (Lemma 1 in \citealt{golrezaei2014real}); the bound there is established against the optimal clairvoyant solution for the realized arrival sequence, so it covers policies with foreknowledge of the sequence. We denote the objective value by $V(\vecc,N)$ for $N=(N_k)_{k\in\MM}$.

The masses $\Lambda_k$ enter \eqref{eq:multi-sblp} only through the right-hand sides of the constraints $y_{k0}+\sum_{i\in\NN}y_{ki}=\Lambda_k$, and the program is feasible for every $N\geq 0$ (take $y_{k0}=N_k$ and $y_{ki}=0$), so $V(\vecc,\cdot)$ is concave in the masses, by the standard concavity of the optimal value of a maximization LP in its right-hand side. Since customer types are drawn independently across periods and $T$ is independent of the type sequence, $\EE[N_k\mid T=t]=\sum_{s=1}^{t}\lambda_{k,s}$, which are exactly the type masses defining $\fLPM(\vecc,t)$. Jensen's inequality then gives
\(
\EE[V(\vecc,N)\mid T=t]\leq V(\vecc,\EE[N\mid T=t])=\fLPM(\vecc,t).
\)
Taking expectations over the type sequence and then over $T$ completes the proof. 
\Halmos
\endproof

\subsection{Sequencing Algorithm} \label{app:alg-seq}
The Sequencing algorithm determines the sequence of assortments that appear in the fluid process and their corresponding time durations \citep{goyal2022dynamic}. In the fluid process, the assortment changes only when a product stocks out. The algorithm iteratively identifies the next product to stock out and the time until that stockout. We note that \citet{goyal2022dynamic} formally consider the algorithm for general choice models, but the idea was already known for the MNL choice model from prior work such as \citet{topaloglu2013joint}.

\begin{algorithm}[htp]
\SetAlgoNoLine
\KwIn{products $\NN$, number of customers $T$, and inventory vector $\vecc$}
\KwOut{A sequence of assortments $\{S_l\}_{l=1}^L:S_L\subseteq \cdots \subseteq S_{2}\subseteq S_1$, and their durations $\{t_l\}_{l=1}^L$}
Initialize $t=0$, $l=1$, $S_1=\{i\in\NN:c_i>0\}$\;
\While{$S_l \neq \emptyset$ and $t< T$} {
    $t_l= \min \{\min_{i \in S_l }\frac{c_i}{\phi(i,S_l)}, T-t\}$ \label{line:seq-stockout}\;
    \lFor{$i \in S_l$}{$c_i=c_i-t_l \phi(i,S_l)$}
    $t = t+t_l$, $l = l+1$, $S_l=\{i\in\NN:c_i>0\}$
}
\lIf{$t < T$}{$S_{l}=\emptyset$, $t_{l}= T-t$, $l=l+1$}
$L=l-1$\;
\caption{Sequencing}\label{alg:seq}
\end{algorithm}

The term $\min_{i \in S_l }\frac{c_i}{\phi(i,S_l)}$ in Line~\ref{line:seq-stockout} computes the time until the next product stockout after time $t$. If all inventory is depleted before time $T$, the algorithm sets the final assortment to the empty set for the remaining time period.

\subsection{Sales-based Formulation of Multi-type CDLP} \label{app:multi-sblp}
Under the MNL model, the multi-type CDLP again admits a compact sales-based formulation. Let $\Lambda_k=\sum_{t=1}^{T}\lambda_{k,t}$ denote the expected number of type-$k$ customers over the horizon, let $y_{ki}$ represent the sales of product $i$ to type-$k$ customers, and let $y_{k0}$ denote the number of type-$k$ customers making no purchase. The formulation is as follows:
\begin{equation}
\label{eq:multi-sblp}
\begin{aligned}
  \vspace{-5pt}
\fLPM(\vecc,T)
=
\max_{\mathbf y\ge 0}
\Biggl\{
&\sum_{k \in \MM}\sum_{i\in \NN} r_i y_{ki}
: \quad
\sum_{k \in \MM} y_{ki}\leq c_i, \forall i \in \NN  \\
&y_{k0}+\sum_{i\in \NN} y_{ki}= \Lambda_k , \forall k\in \MM,
\quad y_{ki} \leq v_{ki} y_{k0},  \forall i\in \NN, \forall k\in \MM
\Biggr\}.
\end{aligned}
\end{equation}
The constraints $\sum_{k\in\MM}y_{ki}\leq c_i$ enforce the inventory limits, which are shared across the customer types, and the remaining per-type constraints are the demand-balance and MNL constraints of the single-type SBLP \eqref{sblp} with horizon $\Lambda_k$.

\section{Optimizing the Fluid Relaxations: Algorithms, Lemmas, and Proofs}\label{app:3}
This appendix contains the algorithms and proofs for the first (fluid-optimization) stage: the threshold-based and DP-based algorithms for the single-type CDLP used in \textsc{da}, the proofs of Lemmas~\ref{lem:3-rev} and \ref{lem:sep} and Theorem~\ref{thm:sof}, and the threshold-augmentation algorithms for the multi-type CDLP used in \textsc{dap}.

\subsection{Threshold-Based Algorithm for Maximizing Single-type CDLP}\label{app:single-T}
\medskip

\begingroup
\noindent\textbf{Theorem~\ref{thm:da-fluid-opt} (restated).} {\itshape
For any accuracy level $\epsilon\in (0,1)$, $\max_{\|\vecc\|_1\leq K} \fLP(\vecc,T)$ can be approximated within a factor of $(1 - \frac{1}{e} - \epsilon)$ in $O(\frac{n}{\epsilon^2} \log^2 (\frac{K}{\epsilon})\log K\log n)$ time using {\small THR}. 
}\medskip
\endgroup

\begingroup
Theorem~\ref{thm:da-fluid-opt} in the main text states the cardinality case ($\mathcal F=\{\vecc:\|\vecc\|_1\leq K\}$, Algorithm~\ref{alg:da-thr}); the same argument covers the budget case ($\mathcal F=\{\vecc:\sum_i b_ic_i\leq B\}$, Algorithm~\ref{alg:tau-bud}) with the inner routine replaced by budget-constrained DR-submodular maximization \citep{soma2018maximizing}.
\endgroup

\begingroup
We start with the following proposition showing that $g(\vecc,y_0)$, as defined in Section~\ref{sec:da-fluid-opt}, is a DR-submodular function of the inventory vector when $y_0$ is fixed.
\begin{proposition}[Theorem 4.3 in \cite{bai2022coordinated}]
    For each $i \in \NN $, given $\bar{y}_0\in[0,T]$, and vectors $\vecc,\mathbf{d}\in \mathbb{Z}_+^n$ with $\vecc\geq \mathbf{d}$, the following inequality holds: \(g(\vecc+\vece_i,\bar{y}_0)-g(\vecc,\bar{y}_0)\leq g(\mathbf{d}+\vece_i,\bar{y}_0)-g(\mathbf{d},\bar{y}_0).\)
\end{proposition}
\endgroup
Moreover, $g(\vecc, y_0)$ is monotone in $\vecc$ because increasing the inventory right-hand sides enlarges the feasible region of this maximization LP. 
Leveraging the DR-submodular structure, we develop two variants of our algorithm with a cardinality constraint (Algorithm~\ref{alg:tau}) and with a budget constraint (Algorithm~\ref{alg:tau-bud}). 

\begingroup
{\SetAlgoNoLine%
  \begin{algorithm}[htp]
    \textbf{Input:}  Cardinality $K$, products $\NN$ sorted in descending order of price: $r_1\geq r_2\geq \cdots \geq r_n$, error $\epsilon \in (0,1)$, number of customers $T$\; 
    \textbf{Output: } $\vecc^{\operatorname{ALG}}$\; 
       \For{$i \gets1$ \textbf{to} $\lceil \log_{1+\epsilon }\frac{1}{\epsilon} \rceil+1$} { $y_0^i=\min\{(1+\epsilon)^{i-1}\epsilon T,\, T\}$, 
      $\vecc^i= \text{ Cardinality-Constrained Submodular Max}(K,y_0^i, g, \NN, \epsilon)$\;
    }
    \Return{$\vecc^{\operatorname{ALG}}=\text{Best of } \vecc^1$, $\vecc^2$, $\cdots$, $\vecc^{\lceil \log_{1+\epsilon }\frac{1}{\epsilon} \rceil+1}$, evaluated by $\fLP(\cdot,T)$}
    \caption{Threshold-Based Cardinality-Constrained Optimization for DA}\label{alg:tau}
  \end{algorithm}}

{\SetAlgoNoLine%
  \begin{algorithm}[htp]
    \textbf{Input:}  Budget $B$, weights $(b_i)_{i\in\NN}$, products $\NN$ sorted in descending order of price: $r_1\geq r_2\geq \cdots \geq r_n$, error $\epsilon \in (0,1)$, number of customers $T$\; 
    \textbf{Output: }$\vecc^{\operatorname{ALG}}$\; 
       \For{$i \gets 1$ \textbf{to} $\lceil \log_{1+\epsilon }\frac{1}{\epsilon} \rceil+1$} {  $y_0^i=\min\{(1+\epsilon)^{i-1}\epsilon T,\, T\}$,
      $\vecc^i= \text{Budget-constrained Submodular Max}(B,y_0^i, g, \NN, \epsilon)$\;
    }
    \Return{$\vecc^{\operatorname{ALG}}=\text{Best of } \vecc^1$, $\vecc^2$, $\cdots$, $\vecc^{\lceil \log_{1+\epsilon }\frac{1}{\epsilon} \rceil+1}$, evaluated by $\fLP(\cdot,T)$}
    \caption{Threshold-Based Budget-Constrained Optimization for DA}\label{alg:tau-bud}
  \end{algorithm}}

Here, $\text{Cardinality-Constrained Submodular Max}(K, y_0, g, \NN, \epsilon)$ is the algorithm for maximizing the DR-submodular function $g(\cdot,y_0)$ subject to a cardinality constraint $K$ with accuracy $\epsilon$; it returns an inventory vector $\vecc$ with $g(\vecc,y_0)\geq (1-\frac{1}{e}-\epsilon)\max_{\|\vecc'\|_1\leq K} g(\vecc',y_0)$ \citep{soma2018maximizing}. Similarly, $\text{Budget-constrained Submodular Max}(B, y_0, g, \NN, \epsilon)$ maximizes $g(\cdot,y_0)$ subject to a budget constraint with budget $B$ and weights $(b_i)_{i\in\NN}$, with the same accuracy guarantee. 

For both Algorithms~\ref{alg:tau} and \ref{alg:tau-bud}, we consider a geometric sequence of $y_0$ values from $\epsilon T$ to $T$, with $\lceil \log_{1+\epsilon }\frac{1}{\epsilon} \rceil$ different values; note that $y_0^{\operatorname{OPT}}\in [0,T]$, so this grid covers the entire relevant range up to the factor handled in the proof below. For each $y_0$, we use standard  DR-submodular maximization techniques to find an approximate inventory vector $\vecc^i$ as described above \citep{soma2018maximizing}.
\endgroup

\proof{Proof of Theorem~\ref{thm:da-fluid-opt}.}

Let $y_0^i =\min\{\epsilon T(1+\epsilon)^{i-1},\, T\}$. For any fixed $y_0$, we define $\vecc^{\operatorname{ALG}}(y_0)$ as the inventory vector produced by the DR-submodular maximization algorithm given $y_0$, and $\vecc^{\operatorname{OPT}}(y_0)$ as the optimal inventory vector maximizing $g(\vecc,y_0)$, i.e., $\vecc^{\operatorname{OPT}}(y_0)=\arg \max_{\vecc\in \mathcal{F}} g(\vecc, y_0)$. Let $\vecc^{\operatorname{OPT}}$ denote the optimal inventory for the original problem $\max_{\vecc\in \mathcal{F}} \fLP(\vecc,T)$, and let $y_0^{\operatorname{OPT}}$ denote the optimal value of $y_0$ in $\fLP(\vecc^{\operatorname{OPT}},T)$. Obviously, $\vecc^{\operatorname{OPT}}(y_0^{\operatorname{OPT}})=\vecc^{\operatorname{OPT}}$.

We first consider the case where $y_0^{\operatorname{OPT}}\geq \epsilon T$; the case where $y_0^{\operatorname{OPT}}\in (0, \epsilon T)$ will be addressed later. Since $y_0^{\operatorname{OPT}}\in [\epsilon T,T]$ and the grid covers this range with multiplicative gaps of $1+\epsilon$, there exists an integer $j$ such that $ y_0^{\operatorname{OPT}}(1-\epsilon)\leq y_0^j \leq y_0^{\operatorname{OPT}}$, which implies that there exists $0\leq \epsilon'<\epsilon$ such that $y_0^j=(1-\epsilon') y_0^{\operatorname{OPT}}$. Then we have:
 \begin{align}
  g(\vecc^{\operatorname{ALG}}(y_0^j),y_0^j ) \nonumber 
    &\geq  (1-\frac{1}{e}-\epsilon)g(\vecc^{\operatorname{OPT}}(y_0^j),y_0^j )\nonumber  \\
        & \geq (1-\frac{1}{e}-\epsilon)g(\vecc^{\operatorname{OPT}}(y_0^{\operatorname{OPT}}), y_0^j)\label{eq:app-31-0}  \\
& = (1-\frac{1}{e}-\epsilon)g(\vecc^{\operatorname{OPT}}(y_0^{\operatorname{OPT}}), (1-\epsilon')y_0^{\operatorname{OPT}}).\label{eq:app-31-1}
\end{align}
Here, the first inequality is the $(1-\frac{1}{e}-\epsilon)$ approximation guarantee of the DR-submodular maximization routine \citep{soma2018maximizing}. The second inequality holds since $\vecc^{\operatorname{OPT}}(y_0^j)$ maximizes $g(\vecc,y_0^j)$. Finally, we plug in $y_0^j=(1-\epsilon') y_0^{\operatorname{OPT}}$. 

Next, we analyze how $g(\vecc^{\operatorname{OPT}}(y_0^{\operatorname{OPT}}),y_0)$ changes when $y_0$ decreases from $y_0^{\operatorname{OPT}}$ to $(1-\epsilon') y_0^{\operatorname{OPT}}$. Note that $g(\vecc,y_0)$ corresponds to the optimal value of a knapsack problem where the weight of each product $i$ is $\min\{y_0v_i, c_i\}$ and the capacity is $T-y_0$. Let $\{y_i^{\operatorname{OPT}}\}_{i \in \NN }$ denote the optimal solution of this knapsack problem when evaluated at $g(\vecc^{\operatorname{OPT}}(y_0^{\operatorname{OPT}}), y_0^{\operatorname{OPT}})$.

\begingroup
 Due to the structure of the knapsack problem, the optimal solution includes items in order of decreasing $r_i$ until the capacity is reached. The capacity is reached because an optimal solution of \eqref{sblp} with inventory $\vecc^{\operatorname{OPT}}$ is feasible for this knapsack and has total sales $T-y_0^{\operatorname{OPT}}$, so the total weight $\sum_{i\in\NN}\min\{c_i,y_0^{\operatorname{OPT}}v_i\}$ is at least the capacity. There exists an index $i_{\tau}$ such that $y_i^{\operatorname{OPT}} = \min\{c_i,y_0^{\operatorname{OPT}} v_i\}$ for $i< i_{\tau}$, $y_i^{\operatorname{OPT}} = 0$  for $i> i_{\tau}$, and $y_{i_{\tau}}^{\operatorname{OPT}} = T-y_0^{\operatorname{OPT}}-\sum_{i=1}^{i_{\tau}-1} \min\{c_i,y_0^{\operatorname{OPT}}v_i\}$. So we have 
\endgroup
\begin{equation}\label{eq:app-312}
g(\vecc^{\operatorname{OPT}}(y_0^{\operatorname{OPT}}), y_0^{\operatorname{OPT}})=\sum_{i<i_{\tau}}r_i\min\{v_iy_0^{\operatorname{OPT}}, c_i \}+ r_{i_{\tau}}y_{i_{\tau}}^{\operatorname{OPT}}.
\end{equation}

\begingroup
When $y_0$ decreases to $y_0^{\operatorname{OPT}}(1-\epsilon')$, the weights $v_iy_0$ decrease proportionally and the capacity $T-y_0$ increases. Consider the solution $y_i'= \min \{ y_i^{\operatorname{OPT}}, v_i(1-\epsilon')y_0^{\operatorname{OPT}}, c_i\}$ for $i\leq i_{\tau}$ and $y_i'=0$ for $i>i_{\tau}$. It is feasible for the new knapsack problem: each $y_i'$ is at most the new weight $\min\{c_i, v_i(1-\epsilon')y_0^{\operatorname{OPT}}\}$ by definition, and $y_i'\leq y_i^{\operatorname{OPT}}$ for all $i$, so the total sales fit within the enlarged capacity. Since $\min\{c_i, v_i(1-\epsilon')y_0^{\operatorname{OPT}}\}\geq (1-\epsilon')\min\{c_i, v_iy_0^{\operatorname{OPT}}\}\geq (1-\epsilon')y_i^{\operatorname{OPT}}$ for $i\leq i_{\tau}$, we have $y_i'\geq (1-\epsilon')y_i^{\operatorname{OPT}}$. The decrease in $g$ is therefore bounded by:
\(
\Delta g \leq \sum_{i=1}^{i_{\tau}}r_i(y_i^{\operatorname{OPT}}-y_i')\leq \epsilon' \sum_{i=1}^{i_{\tau}}r_iy_i^{\operatorname{OPT}}.
\)
\endgroup

Hence, we have:
 \begin{align} g(\vecc^{\operatorname{OPT}}(y_0^{\operatorname{OPT}}), (1 - \epsilon') y_0^{\operatorname{OPT}}) &\geq g(\vecc^{\operatorname{OPT}}(y_0^{\operatorname{OPT}}), y_0^{\operatorname{OPT}}) - \Delta g \nonumber \\ &\geq g(\vecc^{\operatorname{OPT}}(y_0^{\operatorname{OPT}}), y_0^{\operatorname{OPT}}) - \epsilon' \sum_{i = 1}^{i_{\tau}} r_i y_i^{\operatorname{OPT}} \nonumber \\ &\geq g(\vecc^{\operatorname{OPT}}, y_0^{\operatorname{OPT}}) - \epsilon' g(\vecc^{\operatorname{OPT}}(y_0^{\operatorname{OPT}}), y_0^{\operatorname{OPT}}) \nonumber \\ &= (1 - \epsilon') g(\vecc^{\operatorname{OPT}}(y_0^{\operatorname{OPT}}), y_0^{\operatorname{OPT}}). \label{eq:1-eps} \end{align}

\begingroup
Combining this inequality with inequality \eqref{eq:app-31-1}, we have $  g(\vecc^{\operatorname{ALG}}(y_0^j),y_0^j )  \geq (1-\frac{1}{e}-\epsilon)(1-\epsilon')g(\vecc^{\operatorname{OPT}},y_0^{\operatorname{OPT}})\geq (1-\frac{1}{e}-\epsilon-\epsilon')g(\vecc^{\operatorname{OPT}},y_0^{\operatorname{OPT}})$.

Moreover, the algorithm returns the candidate with the largest $\fLP$ value among all $y_0^j$ considered, and $\fLP(\vecc,T)\geq g(\vecc,y_0)$ for any fixed $y_0$ (see the definition of $g$ in Section~\ref{sec:da-fluid-opt}). Therefore,
\[ \fLP(\vecc^{\operatorname{ALG}},T) \geq \fLP(\vecc^{\operatorname{ALG}}(y_0^j),T)\geq g(\vecc^{\operatorname{ALG}}(y_0^j),y_0^j ) \geq (1-\frac{1}{e}-\epsilon-\epsilon')g(\vecc^{\operatorname{OPT}},y_0^{\operatorname{OPT}}).\]
Since $\epsilon'<\epsilon $ and $g(\vecc^{\operatorname{OPT}},y_0^{\operatorname{OPT}})=\fLP(\vecc^{\operatorname{OPT}},T)$, we have:
\[
\fLP(\vecc^{\operatorname{ALG}},T)\geq (1-\frac{1}{e}-2\epsilon)\fLP(\vecc^{\operatorname{OPT}},T).
\]

Now we consider the case where $y_0^{\operatorname{OPT}}< \epsilon T$. We show that $g(\vecc, \epsilon T)\geq (1-2\epsilon) g(\vecc, y_0^{\operatorname{OPT}})$ for any $\vecc$. Recall that $g(\vecc,y_0^{\operatorname{OPT}})$ is the optimal objective value of a knapsack problem with weight $\min\{v_iy_0^{\operatorname{OPT}}, c_i\}$ for product $i$ and capacity $T-y_0^{\operatorname{OPT}}$. Moving from $y_0=y_0^{\operatorname{OPT}}$ to $y_0=\epsilon T\geq y_0^{\operatorname{OPT}}$ enlarges every weight, since $\min\{v_i\epsilon T, c_i\}\geq \min\{v_iy_0^{\operatorname{OPT}}, c_i\}$, and reduces the capacity from $T-y_0^{\operatorname{OPT}}$ to $T-\epsilon T$. If the optimal solution at $y_0=y_0^{\operatorname{OPT}}$ does not use the full capacity $T-y_0^{\operatorname{OPT}}$, it either remains feasible at $y_0=\epsilon T$ or can be scaled down by the factor $\frac{T-\epsilon T}{T-y_0^{\operatorname{OPT}}}\geq 1-\epsilon$, so that $g(\vecc,\epsilon T)\geq (1-\epsilon)g(\vecc,y_0^{\operatorname{OPT}})$ and the claim holds. If the capacity is binding, the same scaling argument applies: for an optimal solution $\{y_i\}_{i\in\NN}$ at $y_0=y_0^{\operatorname{OPT}}$, the scaled solution $\{\frac{T-\epsilon T}{T-y_0^{\operatorname{OPT}}}y_i\}_{i\in\NN}$ is feasible at $y_0=\epsilon T$, since each scaled variable is no larger than the original one (hence no larger than the enlarged weight) and the total scaled sales equal $T-\epsilon T$, matching the reduced capacity. Therefore
\begin{align}
\label{2-eps}
g(\vecc, \epsilon T) \geq \frac{T-\epsilon T}{T-y_0^{\operatorname{OPT}}}\, g(\vecc,y_0^{\operatorname{OPT}}) \geq (1-\epsilon)\, g(\vecc, y_0^{\operatorname{OPT}})\geq (1-2\epsilon)\, g(\vecc, y_0^{\operatorname{OPT}}).
\end{align}
Let $\vecc = \vecc^{\operatorname{OPT}}(y_0^{\operatorname{OPT}})=\vecc^{\operatorname{OPT}}$; combining \eqref{2-eps} with inequality \eqref{eq:app-31-0}, we have:
\[
g(\vecc^{\operatorname{ALG}}(\epsilon T),\epsilon T)\geq (1-\frac{1}{e}-\epsilon)g(\vecc^{\operatorname{OPT}},\epsilon T)\geq (1-\frac{1}{e}-3\epsilon) g(\vecc^{\operatorname{OPT}},y_0^{\operatorname{OPT}})=(1-\frac{1}{e}-3\epsilon)\fLP(\vecc^{\operatorname{OPT}},T).
\]
In both cases, $\fLP(\vecc^{\operatorname{ALG}},T)\geq (1-\frac{1}{e}-3\epsilon)\fLP(\vecc^{\operatorname{OPT}},T)$. Replacing $\epsilon$ by $\epsilon/3$ throughout the algorithm, which does not affect the order of the running time, yields the stated guarantee.
\Halmos
\endgroup
\endproof

\begingroup
\paragraph{Runtime Analysis.} In the {\small THR} algorithm (Algorithm~\ref{alg:tau}), we examine $\lceil \log_{1+\epsilon} \frac{1}{\epsilon} \rceil$ values of $y_0$. For each $y_0$, since $g(\vecc,y_0)$ is DR-submodular, we adopt the threshold augmentation algorithm from \cite{soma2018maximizing}. Next we show each call of this routine runs in $O(\frac{n}{\epsilon}\log (\frac{K}{\epsilon})\log K\log n)$ time. 
The threshold augmentation algorithm examines $\lceil \log_{1+\epsilon} \frac{K}{\epsilon} \rceil$ threshold values. For each threshold value $\tau$, the algorithm processes products in descending price order to determine the maximum quantity $k$ of each product that maintains marginal benefit of at least $k \tau$.

Recall that $g(\vecc,y_0)$ corresponds to the optimal value of a knapsack problem where the weight of each product $i$ is $\min\{y_0v_i, c_i\}$ and the capacity is $T-y_0$, and that its optimal solution includes products in order of decreasing $r_i$ until the capacity is reached. Hence, if the total weight is at most $T-y_0$, then $g(\vecc,y_0)=\sum_{j\in\NN}r_j\min\{c_j,y_0v_j\}$; otherwise,
\[
g(\vecc,y_0)=\sum_{j<h}r_j\min\{c_j,y_0v_j\}+r_h\Big(T-y_0-\sum_{j<h}\min\{c_j,y_0v_j\}\Big),
\]
where $h$ is the smallest index at which the cumulative weight reaches $T-y_0$, as in \eqref{eq:app-312}. We maintain the weights $\min\{c_j,y_0v_j\}$ and the revenue contributions $r_j\min\{c_j,y_0v_j\}$ in two Fenwick (binary indexed) trees ordered by descending price; each supports prefix sums and single-coordinate updates and locates the index $h$ in $O(\log n)$ time, so each evaluation of $g$ takes $O(\log n)$ time. The marginal benefit $g(\vecc+k\vece_i,y_0)-g(\vecc,y_0)$ of adding $k$ units of product $i$ is computed by temporarily increasing the weight of product $i$ to $\min\{c_i+k,y_0v_i\}$, updating its revenue contribution accordingly, and re-evaluating $g$; this net value accounts for any sales of lower-priced products displaced by the added units. If the units are accepted, we make this update permanent.

For each product $i$, we determine the maximum quantity $k$, up to the remaining cardinality budget $K-\sum_{j\in\NN}c_j$, whose marginal benefit is at least $k\tau$. Since $g(\cdot,y_0)$ is DR-submodular, the marginal benefit of each additional unit is nonincreasing, so the average marginal benefit $\frac{1}{k}(g(\vecc+k\vece_i,y_0)-g(\vecc,y_0))$ is nonincreasing in $k$; the maximum quantity can therefore be found by binary search over $k$ using $O(\log K)$ evaluations of $g$.

For each threshold value, the algorithm processes $n$ products, and determining the quantity for each product requires $O(\log K)$ evaluations of $g$, each taking $O(\log n)$ time. Therefore, for each candidate value of $y_0$, the total runtime is $O(n\log_{1+\epsilon} (\frac{K}{\epsilon})\cdot\log K\log n) = O(\frac{n}{\epsilon}\log (\frac{K}{\epsilon})\cdot\log K\log n)$, where we use $\log_{1+\epsilon} \frac{K}{\epsilon}=O(\frac{1}{\epsilon}\log \frac{K}{\epsilon})$. Since we examine $\lceil\log_{1+\epsilon} \frac{1}{\epsilon}\rceil+1=O(\frac{1}{\epsilon}\log \frac{1}{\epsilon})$ candidate values of $y_0$, the total running time is $O(\frac{n}{\epsilon^2} \log\frac{1}{\epsilon}\log(\frac{K}{\epsilon})\log K\log n)=O(\frac{n}{\epsilon^2} \log^2(\frac{K}{\epsilon})\log K\log n)$. Finally, the return step evaluates the objective $\fLP(\cdot,T)$ of each of the $O(\frac{1}{\epsilon}\log \frac{1}{\epsilon})$ candidate inventory vectors, which takes $O(n\log n)$ time per candidate using the algorithm of Appendix~\ref{app:fast_cdlp} and is dominated by the bound above. Without loss of generality, we assume $\log K \ge 1$. Then the runtime is $O(\frac{n}{\epsilon^2} \log\frac{1}{\epsilon}\log(\frac{K}{\epsilon})\cdot \log K\log n)=O(\frac{n}{\epsilon^2} \log^2(\frac{K}{\epsilon})\log K\log n)$. 
\endgroup

\hypersetup{bookmarksdepth=section}
\subsection{A Value-Indexed FPTAS for Maximizing the Single-Type CDLP}
\label{app:fptas}

In this section, we present the value-indexed dynamic program used in
Theorem~\ref{thm:dp-fluid-opt}; the algorithm referred to there as {\small DP}
is Algorithm~\ref{alg:fptas}. This algorithm builds upon ideas originally proposed in
an unpublished version of \cite{chen2022approximation}. \begingroup However, a key difference is that now the cardinality consumption is stored as the value of each state rather than as a state coordinate. Consequently,
the number of states does not depend on the numerical value of $K$. \endgroup  

The idea is to first guess the optimal value of $y_0$ and then solve an approximate dynamic program based
on this guess. The structure of the optimal solution $\{y_i^*\}_{i\in\NN}$, which was introduced in Appendix \ref{app:3}, plays
a crucial role in constructing the dynamic program. However, in the original algorithm of \cite{chen2022approximation}, the set of values that are tried is not fully
polynomial but depends on the problem instance. We address this issue by trying $y_0$ values from a geometric series
and provide the full algorithm with proofs for completeness. Since the final published version of \cite{chen2022approximation} does not contain this result, we reproduce all the structural properties and algorithm details here (with permission from the authors) for the sake of
completeness. We first provide the algorithm and then establish the performance guarantee.

\begingroup
\medskip
\noindent\textbf{Theorem~\ref{thm:dp-fluid-opt} (restated).} {\itshape
For any accuracy level $\epsilon\in (0,1)$, under the cardinality constraint
$\mathcal F=\{\vecc:\|\vecc\|_1\leq K\}$, {\small DP} returns an inventory
vector $\widehat{\vecc}\in\mathbb Z_+^n$ with $\|\widehat{\vecc}\|_1\leq K$ and
\(
\fLP(\widehat{\vecc},T)\geq (1-\epsilon)\max_{\vecc\in \mathcal{F}} \fLP(\vecc,T).
\)
Its arithmetic running time is $\tilde{O}\left(\frac{n^4}{\epsilon^4}\right)$,
where $\tilde{O}(\cdot)$ suppresses factors polylogarithmic in the model
parameters and in $\frac{1}{\epsilon}$, and its bit complexity is polynomial in
the binary input size and $\frac{1}{\epsilon}$. 
}\medskip

Let $
\operatorname{OPT}_{\mathrm{LP}}
:=
\max_{\vecc\in\mathbb Z_+^n:\,\|\vecc\|_1\le K}
\fLP(\vecc,T)$.
Throughout we assume $r_i>0$ for every remaining product.

\paragraph{A sales-only reformulation.}
For $u\in[0,T]$, define
\begin{equation}\label{eq:dp-psi}
\Psi(u)
:=
\max_{\mathbf y\in\RR_+^n}
\left\{
\sum_{i\in\NN}r_i y_i:
\sum_{i\in\NN}y_i\le T-u,\quad
 y_i\le uv_i\ \forall i\in\NN,\quad
\sum_{i\in\NN}\lceil y_i\rceil\le K
\right\}.
\end{equation}
Here $u$ is a lower bound on the no-purchase fluid mass. For any target sales
quantity $y_i$, the smallest integral inventory that supports it is
$\lceil y_i\rceil$. The following reformulation allows the dynamic program to
optimize directly over sales rather than inventories.
\vspace{-5pt}
\begin{lemma}[Sales-based reformulation]\label{lem:dp-sales-only}
We have
\(
\operatorname{OPT}_{\mathrm{LP}}
=
\max_{u\in[0,T]}\Psi(u).
\)
\end{lemma}
\endgroup
\vspace{-5pt}
\proof{Proof of Lemma~\ref{lem:dp-sales-only}.}
\begingroup
Let $\vecc^*$ be an optimal inventory vector and let
$(y_0^*,y_1^*,\ldots,y_n^*)$ be an optimal solution to
$\fLP(\vecc^*,T)$. Since $c_i^*$ is integral and $y_i^*\le c_i^*$,
$\lceil y_i^*\rceil\le c_i^*$ for every $i$. Moreover,
$\sum_i y_i^*=T-y_0^*$ and $y_i^*\le y_0^*v_i$. Thus,
$\mathbf y^*$ is feasible for $\Psi(y_0^*)$, and hence
\(
\max_{u\in[0,T]}\Psi(u)
\geq
\Psi(y_0^*)
\geq
\sum_{i\in\NN}r_i y_i^*
=
\operatorname{OPT}_{\mathrm{LP}}.
\)

Conversely, fix $u\in[0,T]$ and a feasible solution $\mathbf y$ to
\eqref{eq:dp-psi}. Set $c_i=\lceil y_i\rceil$ and
$y_0=T-\sum_i y_i$. Then $\|\vecc\|_1\le K$ and $y_0\ge u$. Therefore,
for every $i\in\NN$,
\(
y_i\le uv_i\le y_0v_i,
y_i\le c_i,
\)
and $y_0+\sum_i y_i=T$. Hence $(y_0,\mathbf y)$ is feasible for
$\fLP(\vecc,T)$, which gives $\sum_i r_i y_i\le
\operatorname{OPT}_{\mathrm{LP}}$. Maximizing over $\mathbf y$ and $u$
proves the reverse inequality.
\Halmos
\endgroup
\endproof

\begingroup
Next we quantify the scale of the single-product fluid revenue. 
For each product $i$, consider the feasible inventory $K\vece_i$. Its fluid
value is
\begin{equation}\label{eq:dp-single-value}
L_i
:=
\fLP(K\vece_i,T)
=
r_i\min\left\{K,\frac{Tv_i}{1+v_i}\right\}.
\end{equation}
Define $
L_0:=\max_{i\in\NN}L_i$.
These values can be computed directly in $O(n)$ time.

\begin{lemma}[Single-product scale]\label{lem:dp-scale}
The best single-product value satisfies
$
L_0
\leq
\operatorname{OPT}_{\mathrm{LP}}
\leq
nL_0$.
\end{lemma}
\endgroup

\proof{Proof of Lemma~\ref{lem:dp-scale}.}
\begingroup
The first inequality follows because $K\vece_i$ is feasible for every
$i\in\NN$. For the second, let $(y_0^*,y_1^*,\ldots,y_n^*)$ be an optimal
sales solution associated with an optimal inventory $\vecc^*$. For every
$i\in\NN$, we have $y_i^*\le c_i^*\le K$. In addition,
$y_0^*=T-\sum_j y_j^*\le T-y_i^*$, so the MNL constraint
$y_i^*\le v_i y_0^*$ implies
\(
y_i^*
\leq
v_i(T-y_i^*)\). Thus, \(
 y_i^*\leq\frac{Tv_i}{1+v_i}.
\)
Consequently,
\(
r_i y_i^*
\leq
r_i\min\left\{K,\frac{Tv_i}{1+v_i}\right\}
=L_i
\leq L_0.
\)
Summing over all products gives
$\operatorname{OPT}_{\mathrm{LP}}=\sum_i r_i y_i^*\le nL_0$.
\Halmos
\endgroup
\endproof

\begingroup
When $L_0>0$, define the geometric grid $
\mathcal L
:=
\left\{
2^jL_0:
 j=0,\ldots,\left\lceil\log_2 n\right\rceil
\right\}$.
Lemma~\ref{lem:dp-scale} implies that this grid contains a scale
$\Lambda^*$ satisfying
\begin{equation}\label{eq:dp-good-scale}
\Lambda^*
\leq
\operatorname{OPT}_{\mathrm{LP}}
<
2\Lambda^*.
\end{equation}

\paragraph{The value-indexed dynamic program.}
Fix the internal accuracy parameter $
\eta:=\frac{\epsilon}{3}$,
and define the geometric grid $
\mathcal G_\eta
:=
\left\{
\min\{\eta T(1+\eta)^k,T\}:
 k=0,\ldots,
\left\lceil\log_{1+\eta}\frac{1}{\eta}\right\rceil
\right\}$
for the no-purchase mass. For a scale $\Lambda\in\mathcal L$, set
$
\Delta_\Lambda:=\frac{\eta\Lambda}{n},
Z:=\left\lceil\frac{2n}{\eta}\right\rceil$.
The outer algorithm enumerates the two grids and returns the candidate with
the largest exact CDLP value.

\begin{algorithm}[htp]
\SetAlgoNoLine
\caption{Value-Indexed FPTAS for \textsc{da} Fluid Relaxation}
\label{alg:fptas}
\KwIn{Cardinality $K$, horizon $T$, products $\NN$, accuracy parameter
$\epsilon\in(0,1)$}
\KwOut{Inventory vector $\vecc^{\operatorname{ALG}}$}
Compute $L_i$ from \eqref{eq:dp-single-value} and set $L_0=\max_i L_i$\;
Set $\eta=\epsilon/3$, construct $\mathcal L$ and $\mathcal G_\eta$, and set
$Z=\lceil2n/\eta\rceil$\;
\For{$\Lambda\in\mathcal L$}{
    Set $\Delta=\eta\Lambda/n$\;
    \For{$\widehat y_0\in\mathcal G_\eta$}{
        Set
        $\vecc^{\Lambda,\widehat y_0}
        =\operatorname{ValueDP}(K,T,\widehat y_0,\Delta,Z,\eta,\NN)$\;
    }
}
\Return{the candidate $\vecc^{\Lambda,\widehat y_0}$ with the largest exact
value of $\fLP(\cdot,T)$}\;
\end{algorithm}
\vspace{-5pt}
We next define the inner routine. We first define the discrete grid for the total sales. Fix $\widehat y_0$ and write
\(
C:=T-\widehat y_0,
r_{\max}:=\max_{i\in\NN}r_i,
a:=1+\frac{\eta}{n}.
\)
When $C>0$, define
\begin{equation}\label{eq:dp-p-grid-parameters}
\underline p
:=
\min\left\{C,\frac{\Delta}{r_{\max}}\right\},
\qquad
Q
:=
\left\lceil\log_a\frac{C}{\underline p}\right\rceil,
\end{equation}
and let
\begin{equation}\label{eq:dp-p-grid}
\mathcal P
:=
\{0,C\}
\cup
\{\underline p\,a^j:j=0,\ldots,Q-1\}.
\end{equation}
When $Q=0$, the second set is empty. For $x\in[0,C]$, define
\begin{equation}\label{eq:dp-rounding-map}
P_\eta(x)
:=
\max\{p\in\mathcal P:p\le x\},
\end{equation}
that is, $P_\eta(x)$ rounds $x$ down to the nearest grid point.

For $l\in\{0,\ldots,n\}$, $p\in\mathcal P$, and
$z\in\{0,\ldots,Z\}$, the state value $D_l(p,z)$ is the minimum integral
inventory required after processing products $1,\ldots,l$ and reaching
rounded cumulative sales $p$ and rounded revenue $z\Delta$. An action $q$ for
product $l$ represents the sales quantity $
s_l(q):=\frac{q\Delta}{r_l}$,
which earns revenue $q\Delta$ and requires $\lceil s_l(q)\rceil$ units of
inventory.

\begin{algorithm}[htp]
\SetAlgoNoLine
\caption{Value-Indexed Dynamic Program}
\label{alg:value-dp}
\KwIn{Cardinality $K$, horizon $T$, guess $\widehat y_0$, revenue unit
$\Delta$, maximum revenue index $Z$, accuracy parameter $\eta$, products
$\NN$}
\KwOut{Inventory vector $\vecc$}
Set $C=T-\widehat y_0$; \lIf{$C=0$}{\Return{$\vecc=\mathbf{0}^n$}}Construct the grid $\mathcal P$ and the map $P_\eta(\cdot)$ from
\eqref{eq:dp-p-grid-parameters}--\eqref{eq:dp-rounding-map}\;
Set $D_0(0,0)=0$ and $D_0(p,z)=+\infty$ for every $(p,z)\neq(0,0)$\;
\For{$l\gets1$ \KwTo $n$}{
    Set $D_l(p,z)=+\infty$ for every
    $p\in\mathcal P$ and $z\in\{0,\ldots,Z\}$\;
    \For{$p\in\mathcal P$ and $z\in\{0,\ldots,Z\}$ with
    $D_{l-1}(p,z)\le K$}{
        \For{$q\gets0$ \KwTo $Z-z$}{
            Set $s=q\Delta/r_l$\;
            \If{$s\le\widehat y_0v_l$ and $p+s\le C$}{
                Set $p'=P_\eta(p+s)$ and
                $d'=D_{l-1}(p,z)+\lceil s\rceil$\;
                \If{$d'\le K$ and $d'<D_l(p',z+q)$}{
                    Set $D_l(p',z+q)=d'$ and store the predecessor and action
                    $q$
                }
            }
        }
    }
}
Choose a state $(p^*,z^*)$ with $D_n(p^*,z^*)\le K$ and maximum $z^*$\;
Recover the actions $\widetilde q_1,\ldots,\widetilde q_n$ and set
$\widetilde y_l=\widetilde q_l\Delta/r_l$\;
Set $c_l=\lceil\widetilde y_l\rceil$ for every $l\in\NN$\;
\Return{$\vecc=(c_l)_{l\in\NN}$}\;
\end{algorithm}

Next, we analyze its approximation guarantee. 
\begin{lemma}
\label{lem:dp-y0-grid}
For every $\eta\in(0,1)$,
\(
\max_{u\in\mathcal G_\eta}\Psi(u)
\geq
(1-\eta)\operatorname{OPT}_{\mathrm{LP}}.
\)
\end{lemma}
\endgroup

\proof{Proof of Lemma~\ref{lem:dp-y0-grid}.}
\begingroup
Let $y_0^*$ maximize
$\Psi(u)$, and let $\mathbf y^*$ be an optimal solution to $\Psi(y_0^*)$. Suppose first that $y_0^*\ge\eta T$. There exists
$\widehat y_0\in\mathcal G_\eta$ such that
\(
\widehat y_0
\leq
y_0^*
<
(1+\eta)\widehat y_0.
\)
Set $\theta=\widehat y_0/y_0^*$. Then
$\theta\ge1/(1+\eta)\ge1-\eta$. The scaled sales
$\bar y_i=\theta y_i^*$ satisfy
\(
\bar y_i\le\theta y_0^*v_i=\widehat y_0v_i
\)
and
\(
\sum_i\bar y_i
\leq
\theta(T-y_0^*)
\leq
T-\widehat y_0,
\)
where the last inequality uses $\widehat y_0\le y_0^*$. Since
$\theta\le1$, we also have
$\sum_i\lceil\bar y_i\rceil\le\sum_i\lceil y_i^*\rceil\le K$.
Thus, $\bar{\mathbf y}$ is feasible for $\Psi(\widehat y_0)$ and earns at
least $(1-\eta)\operatorname{OPT}_{\mathrm{LP}}$.

Now suppose that $y_0^*<\eta T$. Take $\widehat y_0=\eta T$ and set
\(
\theta:=\frac{T-\eta T}{T-y_0^*}.
\)
Then $1-\eta\le\theta\le1$. Again setting
$\bar y_i=\theta y_i^*$ gives
\(
\sum_i\bar y_i
\leq
\theta(T-y_0^*)
=T-\eta T
=T-\widehat y_0,
\)
while
$\bar y_i\le y_i^*\le y_0^*v_i\le\widehat y_0v_i$ and
$\sum_i\lceil\bar y_i\rceil\le K$. Hence,
$\bar{\mathbf y}$ is feasible for $\Psi(\widehat y_0)$ and has value at
least $(1-\eta)\operatorname{OPT}_{\mathrm{LP}}$.
\Halmos
\endgroup
\endproof

\begingroup
Next, we establish the guarantee of the output solution for a fixed guess of $y_0$.
\begin{lemma}\label{lem:dp-fixed}
Fix $\widehat y_0\in[0,T]$ and $\Lambda>0$, and set
$\Delta=\eta\Lambda/n$. Assume that the DP's revenue range covers the
optimal value for the fixed guess $\widehat y_0$, namely,
\begin{equation}\label{eq:dp-coverage}
Z\Delta\geq\Psi(\widehat y_0).
\end{equation}
If $\vecc$ is returned by Algorithm~\ref{alg:value-dp}, then
\(
\fLP(\vecc,T)
\geq
\rho\bigl(\Psi(\widehat y_0)-n\Delta\bigr), 
\rho:=\left(1+\frac{\eta}{n}\right)^{-n}.
\)
\end{lemma}
\proof{Proof of Lemma~\ref{eq:dp-coverage}.}
Let $C=T-\widehat y_0$, $a=1+\eta/n$, and $\rho=a^{-n}$. The result is
immediate when $C=0$, because then $\Psi(\widehat y_0)=0$. Hence, suppose
that $C>0$. Let $\mathbf y^*$ be an optimal solution to $\Psi(\widehat y_0)$.
We round $\mathbf y^*$ to obtain $\bar{\mathbf y}$ and show that the
dynamic program can select the actions associated with
$\bar{\mathbf y}$. Because the dynamic program selects a terminal state
with maximum revenue, its recovered solution
$\widetilde{\mathbf y}$ has revenue at least that of
$\bar{\mathbf y}$. Finally, we scale $\widetilde{\mathbf y}$ by $\rho$
to obtain a feasible solution for $\fLP(\vecc,T)$.

For each product $i$, define
\(
\bar q_i
:=
\left\lfloor\frac{r_i y_i^*}{\Delta}\right\rfloor,
\bar y_i
:=
\frac{\bar q_i\Delta}{r_i}.
\)
Because $0\le r_i(y_i^*-\bar y_i)<\Delta$ for every $i$, \(
\sum_i r_i\bar y_i
=
\Delta\sum_i\bar q_i
\geq
\Psi(\widehat y_0)-n\Delta\).
Moreover, $\bar y_i\le y_i^*$ for every $i$, and hence
\begin{equation}\label{eq:dp-rounded-feasibility}
\bar y_i\le\widehat y_0v_i,
\qquad
\sum_i\bar y_i\le C,
\qquad
\sum_i\lceil\bar y_i\rceil\le K.
\end{equation}
The coverage condition \eqref{eq:dp-coverage} also gives
\begin{equation}\label{eq:dp-revenue-index-bound}
\sum_i\bar q_i
\leq
\frac{\Psi(\widehat y_0)}{\Delta}
\leq Z.
\end{equation}
Since $\bar q_i\geq0$, every partial sum satisfies
\[
0\leq\sum_{j=1}^i\bar q_j\leq\sum_{j=1}^n\bar q_j\leq Z.
\]

Set $\bar p_0=0$ and, for $i=1,\ldots,n$, define
\(
\bar p_i=P_\eta(\bar p_{i-1}+\bar y_i).
\)
Since $P_\eta$ rounds downward,
$\bar p_{i-1}\le\sum_{j<i}\bar y_j$, and therefore
\(
\bar p_{i-1}+\bar y_i
\leq
\sum_{j=1}^i\bar y_j
\leq C.
\)
Together with \eqref{eq:dp-rounded-feasibility} and
\eqref{eq:dp-revenue-index-bound}, this shows that the action $\bar q_i$ is
admissible from the state
$\left(\bar p_{i-1},\sum_{j<i}\bar q_j\right)$. Using the initialization
$D_0(0,0)=0$, induction on $i$ yields
\begin{equation}\label{eq:dp-benchmark-state}
D_i\left(\bar p_i,\sum_{j=1}^i\bar q_j\right)
\leq
\sum_{j=1}^i\lceil\bar y_j\rceil
\leq K.
\end{equation}
Let $z^*$ be the revenue index of the terminal state selected by the dynamic
program. Since the dynamic program selects a feasible terminal state with
the largest revenue index, \eqref{eq:dp-benchmark-state} implies
\(
z^*\geq\sum_i\bar q_i.
\)
Let $\widetilde q_1,\ldots,\widetilde q_n$ be the recovered actions and set
$\widetilde y_i=\widetilde q_i\Delta/r_i$. Since
$\sum_i\widetilde q_i=z^*$,
\begin{equation}\label{eq:dp-raw-revenue}
\sum_i r_i\widetilde y_i
=
\Delta z^*
\geq
\Delta\sum_i\bar q_i
\geq
\Psi(\widehat y_0)-n\Delta.
\end{equation}
Furthermore, the stored inventory value and the definition
$c_i=\lceil\widetilde y_i\rceil$ imply that $\sum_i c_i\le K$.

The recovered vector $\widetilde{\mathbf y}$ may not satisfy
$\sum_i\widetilde y_i\le C$, because the state $p$ rounds cumulative sales
downward. The geometric grid satisfies
\begin{equation}\label{eq:dp-rounding-factor}
P_\eta(x)
\geq
\frac{x}{a}
\qquad
\text{for every }x\in[\underline p,C].
\end{equation}
Let $\widetilde p_i$ be the cumulative-sales state after product $i$. If
$\widetilde y_i>0$, then
\[
\widetilde y_i
=
\frac{\widetilde q_i\Delta}{r_i}
\geq
\frac{\Delta}{r_{\max}}
\geq
\underline p,
\]
so \eqref{eq:dp-rounding-factor} applies to the corresponding transition.
If $\widetilde y_i=0$, then $\widetilde p_i=\widetilde p_{i-1}$.
Consequently, induction gives
\(
\widetilde p_i
\geq
a^{-i}\sum_{j=1}^i\widetilde y_j, i=1,\ldots,n.
\)
Since $\widetilde p_n\le C$, it follows that
\begin{equation}\label{eq:dp-scaled-capacity}
\rho\sum_i\widetilde y_i\le C.
\end{equation}

Define $y_i'=\rho\widetilde y_i$ and
$y_0'=T-\sum_i y_i'$. By \eqref{eq:dp-scaled-capacity},
$y_0'\ge\widehat y_0$. In addition,
\(
y_i'\le\widetilde y_i\le c_i,
y_i'\le\widetilde y_i\le\widehat y_0v_i\le y_0'v_i.
\)
Thus, $(y_0',\mathbf y')$ is feasible for $\fLP(\vecc,T)$. Applying
\eqref{eq:dp-raw-revenue}, we obtain
\[
\fLP(\vecc,T)
\geq
\sum_i r_i y_i'
=
\rho\sum_i r_i\widetilde y_i
\geq
\rho\bigl(\Psi(\widehat y_0)-n\Delta\bigr).
\]
\Halmos
\endgroup
\endproof
Finally, we prove Theorem~\ref{thm:dp-fluid-opt}.
\proof{Proof of Theorem~\ref{thm:dp-fluid-opt}.}
\begingroup
We first establish the approximation guarantee. If $L_0=0$,
Lemma~\ref{lem:dp-scale} gives $\operatorname{OPT}_{\mathrm{LP}}=0$, so the
zero inventory is optimal. Suppose $L_0>0$, and let
$\Lambda^*\in\mathcal L$ satisfy \eqref{eq:dp-good-scale}. For this scale,
\(
Z\Delta_{\Lambda^*}
\geq
\frac{2n}{\eta}\cdot\frac{\eta\Lambda^*}{n}
=2\Lambda^*
>
\operatorname{OPT}_{\mathrm{LP}}.
\)
Since $\Psi(u)\le\operatorname{OPT}_{\mathrm{LP}}$ for every $u$ by
Lemma~\ref{lem:dp-sales-only}, condition \eqref{eq:dp-coverage} holds. 

By Lemma~\ref{lem:dp-y0-grid}, there is a candidate
$\widehat y_0\in\mathcal G_\eta$ such that
\(
\Psi(\widehat y_0)
\geq
(1-\eta)\operatorname{OPT}_{\mathrm{LP}}.
\)
Applying Lemma~\ref{lem:dp-fixed} to
$(\Lambda^*,\widehat y_0)$ and using
$n\Delta_{\Lambda^*}=\eta\Lambda^*\le
\eta\operatorname{OPT}_{\mathrm{LP}}$ gives
\begin{align*}
\fLP(\vecc,T) \ge 
\rho\bigl(\Psi(\widehat y_0)-n\Delta_{\Lambda^*}\bigr)\ge 
\rho(1-2\eta)\operatorname{OPT}_{\mathrm{LP}}.
\end{align*}
Moreover,
\(
\rho
=\left(1+\frac{\eta}{n}\right)^{-n}
\geq e^{-\eta}
\geq1-\eta.
\)
Therefore,
\[
\fLP(\vecc,T)
\geq
(1-\eta)(1-2\eta)\operatorname{OPT}_{\mathrm{LP}}
\geq
(1-3\eta)\operatorname{OPT}_{\mathrm{LP}}
=
(1-\epsilon)\operatorname{OPT}_{\mathrm{LP}}.
\]
Algorithm~\ref{alg:fptas} returns the candidate with the largest CDLP
value, so its output satisfies the same guarantee.

We next derive the running time. The scale and no-purchase grids have sizes
\(
|\mathcal L|=O(\log n)\), \(
|\mathcal G_\eta|
=O\left(\frac{1}{\eta}\log\frac{1}{\eta}\right),
\)
and $Z=O(n/\eta)$. The cumulative-sales grid is largest at the smallest
scale $L_0$. Uniformly over all scale and no-purchase guesses,
\begin{align*}
|\mathcal P|
&=O\left(
1+\log_{1+\eta/n}\frac{C}{\underline p}
\right)=O\left(
\frac{n}{\eta}
\log\left(2+\frac{Cr_{\max}}{\Delta}\right)
\right)=O\left(
\frac{n}{\eta}
\log\left(2+\frac{nTr_{\max}}{\eta L_0}\right)
\right).
\end{align*}
The logarithmic factor above is bounded by the encoding length of the input.
Indeed, let $i^*\in\argmax_{i\in\NN}r_i$. Since
$L_0\geq L_{i^*}=r_{\max}\min\{K,Tv_{i^*}/(1+v_{i^*})\}$ and $K\geq1$,
\[
\frac{Tr_{\max}}{L_0}
\leq
\frac{T}{\min\{K,Tv_{i^*}/(1+v_{i^*})\}}
\leq
\max\left\{T,\,1+\frac{1}{v_{i^*}}\right\},
\]
so the factor is $O\left(\log\frac{n}{\eta}+\log T+\log\left(1+\frac{1}{v_{i^*}}\right)\right)$.

For a fixed pair $(\Lambda,\widehat y_0)$,
Algorithm~\ref{alg:value-dp} processes $n$ products,
$O(|\mathcal P|Z)$ states per product, and at most $Z+1$ actions per state.
Evaluating $P_\eta$ by binary search takes $O(\log|\mathcal P|)$ time per
transition. Thus, one call uses
\(
O\bigl(n|\mathcal P|Z^2\log|\mathcal P|\bigr)
\)
arithmetic operations. Multiplying by the two outer-grid sizes yields
\(
O\!\left(
\frac{n^4}{\eta^4}
\log n
\log\frac{1}{\eta}
\log\left(2+\frac{nTr_{\max}}{\eta L_0}\right)
\log|\mathcal P|
\right)
=
\widetilde O\!\left(\frac{n^4}{\epsilon^4}\right).
\)
The exact evaluation of all candidates takes lower-order time by the
$O(n\log n)$ single-type CDLP algorithm in
Theorem~\ref{thm:fast-cdlp}.
\Halmos
\endgroup
\endproof

\begingroup
\subsection{Proof of Lemma~\ref{lem:3-rev}} \label{app:prf-3-rev}
\medskip\noindent\textbf{Lemma~\ref{lem:3-rev} (restated).} {\itshape For any inventory vector $\vecc \succeq 0 $, let $\{y^*_i\}_{i\in\NN\cup\{0\}}$ denote an optimal solution to the sales-based formulation of the single-type CDLP $\fLP(\vecc,T)$. Define $i_{\tau}$ as the highest index $i$ for which $\sum_{k=1}^{i-1}\min\{c_k, v_k y_0^*\}\leq T-y_0^*$ and construct an inventory vector $\CFP\in\RR^n$ by
\[
\cFPi=c_i \ \text{ for } 1\le i<i_\tau,\qquad
\cFPi=0 \ \text{ for } i>i_\tau,\qquad
x_{i_{\tau}}^{\vecc}=T-y_0^*-\sum_{k=1}^{i_{\tau}-1}\min\{c_k, v_k y_0^*\}.
\]
Then $\CFP\preceq\vecc$ and $\fFP(\CFP,T)=\max_{\vecx\preceq\vecc}\fFP(\vecx,T)=\fLP(\vecc,T)$.
}\medskip
\endgroup

To prove this lemma, we use Theorem 2 from \cite{goyal2022dynamic} and establish technical lemmas about the SBLP optimal solution and the consumption of the fluid process. 
\begin{lemma}[Theorem 2 in \cite{goyal2022dynamic}]\label{lem:hatf}  Given any feasible solution $\{y_i\}_{i \in \NN  \cup\{0\}}$ to the SBLP \eqref{sblp} with inventory $\vecc$ and $T$ periods, let $x_i = y_i$ for all $i \in \NN $. Then $\fFP(\vecx,T)=\sum_{i \in \NN }r_iy_i$.
\end{lemma}
This lemma essentially says that if we set the starting inventory of each product as a feasible solution to the SBLP, then the inventory will be depleted within the selling horizon. Therefore, the revenue in \textsc{fp} is equal to the objective function of the SBLP. 

  \begin{lemma}\label{lem:3-equ1} For any inventory vector $\vecc$ and $T$ customers, $\max_{\vecx \preceq \vecc} \fFP(\vecx,T)=\fLP(\vecc,T)$. Moreover, given the optimal solution $\{y^*_i\}_{i \in \NN  \cup\{0\}}$ to the SBLP \eqref{sblp}, let $x^*_i = y^*_i$ for all $i \in \NN $ and let $\vecx^*=(x^*_1,\ldots,x^*_n)$. Then $\fFP(\vecx^*,T)=\max_{\vecx \preceq \vecc} \fFP(\vecx,T)=\fLP(\vecc,T)$.
\end{lemma}

\proof{Proof of Lemma~\ref{lem:3-equ1}.} 
We show that, given any feasible solution to $\fLP(\vecc,T)$, one can construct a feasible solution to $\max_{\vecx\preceq \vecc}\fFP(\vecx,T)$ with the same objective function value, and vice versa.

Let $\{y_i\}_{i \in \NN \cup\{0\}}$ denote a feasible solution to the SBLP. Define $\vecx$ by $x_i=y_i$ for all $i\in \NN$. Since $x_i=y_i\leq c_i$, $\vecx$ is a feasible solution to $\max_{\vecx\preceq \vecc }\fFP(\vecx,T) $. Furthermore, by Lemma~\ref{lem:hatf}, all inventory in $\vecx$ will be depleted; thus, $\fFP(\vecx,T)=\sum_{i\in\NN}{r_iy_i}$. 

\begingroup
Let $\vecx$ be a feasible solution to $\max_{\vecx\preceq \vecc }\fFP(\vecx,T) $, so $x_i\leq c_i$ for all $i \in \NN $. Then let $y_i =\ZFP_i(\vecx,T)$ for all $i \in \NN \cup\{0\}$, where $\ZFP_0(\vecx,T)$ is the number of customers who do not make purchases. Obviously, $\sum_{i \in \NN }r_iy_i = \sum_{i\in\NN}r_i \ZFP_i(\vecx,T)= \fFP(\vecx,T)$. So we just need to show $\{y_i\}_{i\in\NN\cup \{0\}}$ is a feasible solution to $\fLP(\vecc,T)$. We show this by observing that: (i) $\ZFP_i(\vecx,T)\leq x_i\leq c_i$; and (ii) $\ZFP_i(\vecx,T)\leq v_i \ZFP_0(\vecx,T) $ due to the Independence of Irrelevant Alternatives (IIA) property of the MNL model.
\Halmos
\endgroup
\endproof

\begingroup
The following elementary characterization of the fluid process under the MNL model underlies several arguments in this and later appendices; we record it once here and reference it subsequently. It is used in the proofs of Lemma~\ref{lem:3-rev} and the Separability Lemma below, and repeatedly in later appendices (e.g., Appendices~\ref{app:prf-thm-pro} and~\ref{app:fast_cdlp}).
\begin{lemma}[Fluid consumption identity]
Fix an inventory vector $\vecx\in \RR_+^n$ and a horizon $T\geq 0$, and let $\ZFP_0(\vecx,t)$ denote the consumption of the outside option up to time $t$ in \textsc{fp}$(\vecx,T)$. Then $\ZFP_j(\vecx,t)=\min\{v_j\ZFP_0(\vecx,t),x_j\}$ for every $j\in \NN$ and every $t\leq T$, and $\ZFP_0(\vecx,t)$ is the unique solution of
\(
Z_0+\sum_{h\in \NN}\min\{v_hZ_0,x_h\}=t.
\)
Moreover, $\ZFP_0(\vecx,\cdot)$ is continuous and strictly increasing, product $j$ is out of stock at time $t$ if and only if $x_j/v_j\leq \ZFP_0(\vecx,t)$, and consequently products stock out in ascending order of the ratios $x_j/v_j$.
\label{lem:4-clock}
\end{lemma}
\proof{Proof of Lemma~\ref{lem:4-clock}.}
While product $j$ is in stock it is consumed at rate $\phi(j,S(t))=v_j\phi(0,S(t))$, that is, at $v_j$ times the rate at which the outside option is consumed, and after it stocks out its consumption stays at $x_j$; this gives the identity. Since all customers are accounted for, $\ZFP_0(\vecx,t)+\sum_{h\in \NN}\ZFP_h(\vecx,t)=t$, and substituting the identity shows that $\ZFP_0(\vecx,t)$ solves the displayed equation; the solution is unique because the left-hand side is strictly increasing in $Z_0$. Finally, $\phi(0,S(t))\geq 1/(1+\sum_{h\in \NN}v_h)>0$, so $\ZFP_0(\vecx,\cdot)$ is continuous and strictly increasing; product $j$ is out of stock at time $t$ exactly when $v_j\ZFP_0(\vecx,t)\geq x_j$, and therefore product $k$ cannot stock out strictly before product $j$ whenever $x_j/v_j\leq x_k/v_k$.
\Halmos
\endproof
\endgroup

\begingroup
\proof{Proof of Lemma~\ref{lem:3-rev}.}
If $y_0^*=T$, the equality constraint in \eqref{sblp} forces $y_i^*=0$ for every $i\in\NN$, so $\fLP(\vecc,T)=0$; by Lemma~\ref{lem:3-equ1}, $\max_{\vecx\preceq\vecc}\fFP(\vecx,T)=0$ as well, and $\CFP=\mathbf{0}^n$ attains it since $\fFP(\mathbf{0}^n,T)=0$. We therefore assume $y_0^*<T$ in the remainder of the proof.
\endgroup

When $y_0^*$ is fixed, the problem simplifies to a knapsack problem:
\begin{align*}
\max_{\mathbf{y} \in \RR^{n}_+} \left\{ \sum_{i \in \NN } r_i y_i :\sum_{i \in \NN } y_i + y_0^* = T; y_i \leq c_i,  y_i \leq  y_0^* v_i, \forall i \in \NN  \right\}.
\end{align*}
\begingroup
Since the items are indexed in non-increasing order of $r_i$, this knapsack problem is solved by filling the items in index order until the capacity $T-y_0^*$ is exhausted. The capacity is reached, because $\{y_i^*\}$ is itself feasible for this knapsack and has total sales $T-y_0^*$, so the total weight $\sum_{k\in\NN}\min\{c_k,v_ky_0^*\}$ is at least $T-y_0^*$. By the definition of $i_{\tau}$, the cumulative weight $\sum_{k=1}^{i-1}\min\{c_k,v_ky_0^*\}$ is at most $T-y_0^*$ for every $i\leq i_{\tau}$ and exceeds $T-y_0^*$ for every $i>i_{\tau}$, so this solution is $y_i^* = \min\{c_i,y_0^* v_i\}$ for $i< i_{\tau}$, $y_i^* = 0$  for $i> i_{\tau}$, and $y_{i_{\tau}}^* = T-y_0^*-\sum_{i=1}^{i_{\tau}-1} \min\{c_i,y_0^*v_i\}$; as all optimal solutions of the knapsack problem have the same value, we may assume that $\{y_i^*\}$ takes this form. This property of the optimal solution was identified in an unpublished version of \cite{chen2022approximation}. In particular $0\leq y_{i_{\tau}}^*\leq \min\{c_{i_{\tau}},v_{i_{\tau}}y_0^*\}$, so $\CFP\preceq \vecc$. From Lemma~\ref{lem:hatf}, if we construct $\vecx$ by setting $x_i=y_i^*$ for all $i \in \NN $, we have $\fFP(\vecx,T)=\sum_{i \in \NN } r_i y_i^* =\fLP(\vecc,T)$. Therefore, it suffices to show that $\fFP(\vecx,T)=\fFP(\CFP,T)$.
\endgroup

\begingroup
Note that the difference between $\vecx$ and $\CFP$ arises only for $i < i_{\tau}$, where $x_i^{\vecc} = c_i \geq\min\{c_i,y_0^*v_i\}=x_i$. By Lemma~\ref{lem:4-clock}, the consumption of each product $j$ in a fluid process satisfies $\ZFP_j(\vecx,T)=\min\{v_j\ZFP_0(\vecx,T),x_j\}$, and $\ZFP_0(\vecx,T)$ is the unique solution of $Z_0+\sum_{j\in \NN}\min\{v_jZ_0,x_j\}=T$. We verify that $Z_0=y_0^*$ solves this equation for the inventory $\CFP$. For $i<i_{\tau}$, $\min\{v_iy_0^*,\cFPi\}=\min\{v_iy_0^*,c_i\}=y_i^*$. For $i=i_{\tau}$, $x_{i_{\tau}}^{\vecc}=y_{i_{\tau}}^*\leq v_{i_{\tau}}y_0^*$, so $\min\{v_{i_{\tau}}y_0^*,x_{i_{\tau}}^{\vecc}\}=y_{i_{\tau}}^*$. For $i>i_{\tau}$, $\cFPi=0=y_i^*$. Thus $y_0^*+\sum_{i\in \NN}\min\{v_iy_0^*,\cFPi\}=y_0^*+\sum_{i\in \NN}y_i^*=T$. Hence $\ZFP_0(\CFP,T)=y_0^*$, and $\ZFP_i(\CFP,T)=\min\{v_iy_0^*,\cFPi\}=y_i^*=\ZFP_i(\vecx,T)$ for all $i\in \NN$. Consequently, $ \fFP(\CFP,T)=\fFP(\vecx,T)=\fLP(\vecc,T)$.
\Halmos
\endgroup
\endproof

\subsection{Proof of Lemma~\ref{lem:sep} }\label{prf:lem-sep}

\noindent\textbf{Lemma~\ref{lem:sep} (restated).} {\itshape
Fix an inventory vector $\vecx\in\RR_+^n$, a horizon $T$, a product
$i\in\NN$, and an additional inventory level $\delta_i\ge0$. Consider
\textsc{fp}$(\vecx,T)$ and
\textsc{fp}$(\vecx+\delta_i\vece_i,T)$. Define
\(
\alpha_i
=
\ZFP_i(\vecx+\delta_i\vece_i,T)
-
\ZFP_i(\vecx,T),
\)
the increase in the total consumption of product $i$. Then:
\begin{enumerate}[label=(\roman*)]
\item
For product $i$,
\(
\ZFP_i(\vecx+\delta_i\vece_i,T)
=
\alpha_i+\ZFP_i(\vecx,T-\alpha_i),
\)
and, for every $j\neq i$,
\(
\ZFP_j(\vecx+\delta_i\vece_i,T)
=
\ZFP_j(\vecx,T-\alpha_i).
\)
\item
The total revenue satisfies
\(
\fFP(\vecx+\delta_i\vece_i,T)
=
\fFP(\vecx,T-\alpha_i)+\alpha_i r_i.
\)
\end{enumerate}
}\medskip
\vspace{-5pt}
\proof{Proof of Lemma~\ref{lem:sep}.} We first prove (i), from which (ii) will follow. When $\alpha_i=0$, the result holds trivially. So without loss of generality, we consider $\alpha_i>0$.

We first show the first part of (i) holds, i.e., $\ZFP_i(\vecx+\delta_i\vece_i,T)=\alpha_i + \ZFP_i(\vecx, T-\alpha_i)$. Observe that if $\ZFP_i(\vecx,T)=\ZFP_i(\vecx,T-\alpha_i)=x_i$, then it holds since $\ZFP_i(\vecx+\delta_i\vece_i,T)=\ZFP_i(\vecx,T)+\alpha_i$ by the definition of $\alpha_i$. We now show $\ZFP_i(\vecx,T)=\ZFP_i(\vecx,T-\alpha_i)=x_i$, which is equivalent to showing that product $i$ becomes out of stock by $T-\alpha_i$. We prove it by contradiction. First, observe that if $\ZFP_i(\vecx,T)<x_i$, which means $i$ has not stocked out by time $T$, then adding $\delta_i$ units of product $i$ to the starting inventory would not change the consumption in the fluid process, so $\alpha_i=0$, contradicting $\alpha_i>0$. Thus, $i$ is out of stock by $T$ and $\ZFP_i(\vecx,T)=x_i$. Now assume $i$ becomes out of stock at some time between $T-\alpha_i$ and $T$. Let $T-\beta$ denote the time at which this occurs. Then we have $\ZFP_i (\vecx, T-\beta)=x_i$ and $\beta<\alpha_i$. Thus we have:
\begin{align*}
  \ZFP_i(\vecx,T)= \ZFP_i (\vecx, T-\beta)= \ZFP_i(\vecx+\delta_i\vece_i, T-\beta).
\end{align*}
Here the first equality holds since both of them are equal to $x_i$. The second equality holds because adding $\delta_i$ units of product $i$ to the starting inventory would not change the consumption process before time $T-\beta$, by which time the original $x_i$ units have been used up. 

Thus, from the definition of $\alpha_i$, we have:\begin{align*}
\alpha_i=  \ZFP_i(\vecx+\delta_i\vece_i, T)-\ZFP_i(\vecx,T)= \ZFP_i(\vecx+\delta_i\vece_i, T)-  \ZFP_i(\vecx+\delta_i\vece_i, T-\beta).
\end{align*}
Since $\phi(i,S)<1$ for any $S$, $\ZFP_i(\vecx+\delta_i\vece_i, T)-  \ZFP_i(\vecx+\delta_i\vece_i, T-\beta)<\beta$. So we have $\alpha_i<\beta$, which contradicts the assumption that $\beta<\alpha_i$. Therefore, we must have $\ZFP_i(\vecx,T)=\ZFP_i(\vecx,T-\alpha_i)=x_i$.

Next, we show the second part of (i). In the fluid process \textsc{fp}$(\vecx+\delta_i\vece_i ,T)$, we have \(
\sum_{j\in \NN\cup\{0\}}\ZFP_j(\vecx+\delta_i\vece_i,T)=T.
\)
In the truncated fluid process \textsc{fp}$(\vecx,T-\alpha_i)$, we have
\(
\sum_{j\in \NN\cup\{0\}}\ZFP_j(\vecx,T-\alpha_i)=T-\alpha_i.
\)
Thus, we have the following equation:
\begin{equation}
    \label{eq:sep-1}
    \sum_{j\in \NN\cup\{0\}}\ZFP_j(\vecx+\delta_i\vece_i,T)= \sum_{j\in \NN\cup\{0\}}\ZFP_j(\vecx,T-\alpha_i)+\alpha_i.
\end{equation}
Since $\ZFP_i(\vecx,T)=\ZFP_i(\vecx,T-\alpha_i)$ and $\alpha_i =\ZFP_i(\vecx+\delta_i\vece_i,T)-\ZFP_i(\vecx,T) $, we have $\ZFP_i(\vecx+\delta_i\vece_i,T)=\ZFP_i(\vecx,T-\alpha_i)+\alpha_i$. Thus, subtracting the $i$-terms from both sides of \eqref{eq:sep-1} yields:
\begin{equation}
    \label{eq:sep-2}
    \sum_{j\in (\NN\cup\{0\})\setminus \{i\}}\ZFP_j(\vecx+\delta_i\vece_i,T)= \sum_{j\in (\NN\cup\{0\})\setminus \{i\}}\ZFP_j(\vecx,T-\alpha_i).
\end{equation}

\begingroup By Lemma~\ref{lem:4-clock},\endgroup{} the consumptions for $j\in \NN\setminus\{i\}$ can be written as
\[
\ZFP_j(\vecx+\delta_i\vece_i,T)=\min\{v_j\ZFP_0(\vecx+\delta_i\vece_i,T),x_j\}, \quad 
\ZFP_j(\vecx,T-\alpha_i)=\min\{v_j\ZFP_0(\vecx,T-\alpha_i),x_j\}.
\]
Substituting into the equation \eqref{eq:sep-2} gives
\begin{multline}\label{eq:sep-eq}
\sum_{j\in \NN\setminus\{i\}}\min\{v_j\ZFP_0(\vecx+\delta_i\vece_i,T),x_j\}+\ZFP_0(\vecx+\delta_i\vece_i,T)\\
= \sum_{j\in \NN\setminus\{i\}}\min\{v_j\ZFP_0(\vecx,T-\alpha_i),x_j\}+\ZFP_0(\vecx,T-\alpha_i).
\end{multline}
Let $h(Z_0)=Z_0+\sum_{j\in \NN\setminus\{i\}}\min\{v_j Z_0,x_j\}$. Since $h$ is strictly increasing, \eqref{eq:sep-eq} implies
\(
\ZFP_0(\vecx+\delta_i\vece_i,T)=\ZFP_0(\vecx,T-\alpha_i).
\)
Thus for all $j\in \NN\setminus \{i\}$, 
\(
\ZFP_j(\vecx+\delta_i\vece_i,T)=\min\{v_j\ZFP_0(\vecx,T-\alpha_i),x_j\}=\ZFP_j(\vecx,T-\alpha_i).
\)
Thus, we have proved part (i). 

For (ii), we have
\begin{align*}
\fFP(\vecx+\delta_i\vece_i,T)&=\sum_{j\in \NN}r_j\ZFP_j(\vecx+\delta_i\vece_i,T)\\
&=\sum_{j\in \NN\setminus \{i\}}r_j\ZFP_j(\vecx+\delta_i\vece_i,T)+\alpha_i r_i +r_i\ZFP_i(\vecx,T)  \\
&=\sum_{j\in \NN}r_j\ZFP_j(\vecx,T-\alpha_i)+\alpha_i r_i \\
&=\fFP(\vecx, T-\alpha_i)+\alpha_i r_i.
\end{align*}
The first and last equalities follow from the definition of total revenue in the fluid processes.  
The second equality uses the definition of $\alpha_i$, and the third equality follows from our earlier results showing $\ZFP_j(\vecx+\delta_i\vece_i,T)=\ZFP_j(\vecx,T-\alpha_i)$ and $\ZFP_i(\vecx,T)=\ZFP_i(\vecx,T-\alpha_i)$. \hfill\Halmos
\endproof

\subsection{Proof of Theorem~\ref{thm:sof}}\label{app:sof} We begin by establishing several technical lemmas that characterize the structure of optimal inventory
vectors and revenue functions under the fluid process.
 \begin{lemma} 
Let $S_1$ and $S_2$ be two subsets of the product set $\NN$ such that $S_2\subseteq S_1$. Define $A=S_1\setminus S_2$. Then the following holds: (i) If $r_k\geq R(S_1)$ for all $k\in A$, then $R(S_1)\geq R(S_2)$. (ii) If $r_k\geq R(S_2)$ for all $k\in A$, then $R(S_1)\geq R(S_2)$.
\label{lem:subsetie}
\end{lemma}
\proof{Proof of Lemma~\ref{lem:subsetie}.}
Let $V(S)=\sum_{i\in S}v_i +v_0$. Then we have:
$R(S_1) = \frac{R(S_2 )V(S_2)+\sum_{k\in A}r_{k}v_{k}}{V(S_2)+\sum_{k\in A}v_{k}}.
$
Thus, if $r_k\geq  R(S_1)$ for all $k\in A$, $R(S_1)-R(S_2)= \frac{\sum_{k\in A} v_k(r_k-R(S_1)) }{V(S_2)}\geq 0$. Moreover, if $r_k\geq  R(S_2)$ for all $k\in A$, then $R(S_1)-R(S_2)=\frac{\sum_{k\in A}v_k(r_k-R(S_2))}{V(S_2)+\sum_{k\in A}v_k}\geq 0$.
\Halmos
\endproof

\begin{lemma}\label{lem:app-subset}
Let $N(\vecx)$ denote the set of products with positive inventory in vector $\vecx$. Consider two inventory vectors, $\vecx_1,\vecx_2\in \mathbb{R}_+^n$ with $\vecx_1\preceq \vecx_2$. Let $S^1(t)$ and $S^2(t)$ denote the assortment at time $t$ in \textsc{fp}$(\vecx_1,T)$ and \textsc{fp}$(\vecx_2,T)$, respectively. Then, for any $t\leq T$, $S^1(t)\subseteq S^2(t)$.
\end{lemma}
This result holds for any choice model with the substitutability property. 
\proof{Proof of Lemma~\ref{lem:app-subset}.} 

Let $\vecx_1(t)$ and $\vecx_2(t)$ denote the inventory vectors at time $t$ in \textsc{fp}$(\vecx_1,T)$ and \textsc{fp}$(\vecx_2,T)$, respectively. The lemma is equivalent to showing that $\vecx_1(t)\preceq \vecx_2(t)$ for any $t$. 

At $t=0$, we have $\vecx_1(0)\preceq \vecx_2(0)$ by assumption. Assortments change only when a product stocks out, so the two processes together have finitely many stockout times in $[0,T]$, and all consumption rates are constant on each interval between consecutive such times. Suppose $\vecx_1(\tau)\preceq \vecx_2(\tau)$ at the start of such an interval. Then for any product $k\in N(\vecx_1(\tau))$,
 \(\phi(k, N(\vecx_1(\tau)))\geq \phi(k, N(\vecx_2(\tau)))\) due to $N(\vecx_1(\tau))\subseteq N(\vecx_2(\tau))$ and the substitutability property of the choice model. Hence, on this interval, each coordinate of $\vecx_1$ is either zero or decreases at a constant rate at least as large as that of the corresponding coordinate of $\vecx_2$, so $\vecx_1(t)\preceq \vecx_2(t)$ throughout the interval. Induction over the finitely many intervals extends the relationship to all $t\in [0,T]$, proving that $S^1(t) \subseteq S^2(t)$ for any $t$.
\Halmos
\endproof
\begin{lemma}\label{lem:seq-rev1}
Consider any inventory vector $\vecx$, any number $T$ of customers, and the fluid process $\text{\textsc{fp}}(\vecx,T)$. Let $S(t)$ denote the set of products with positive inventory at time $t$, and let $i$ denote a product with $r_i\leq r_k$ for all $k\in N(\vecx)$. For any $t\in [0,T]$:
\begin{enumerate}
\item[\rm(i)] If $r_i\geq R(S(t))$, then $r_i\geq R(S(t))\geq R(S(t'))$ for any $t'\geq t$.
\item[\rm(ii)] If $r_i\leq R(S(t))$, then $R(S(t'))\geq r_i$ for any $t'\leq t$.\end{enumerate}
\end{lemma}
\proof{Proof of Lemma~\ref{lem:seq-rev1}.}
We prove each statement separately:
(i) First, observe that $S(t')\subseteq S(t)$ if $t'\geq t$. Let $A=\{k\in \NN \mid k\in S(t), k\notin S(t')\}$. If $r_i\geq R(S(t))$, then, by the choice of $i$, we have $r_k\geq r_i\geq R(S(t))$ for all $k\in A$. Applying Lemma~\ref{lem:subsetie}, we conclude that $R(S(t'))\leq R(S(t))\leq r_i$ for any $t'\geq t$.

(ii) For any $t'\leq t$, $S(t)\subseteq S(t')$ since the assortment only shrinks over time. Let $A=\{k\in \NN \mid k\in S(t'), k\notin S(t)\}$. Then: \[R(S(t'))=\frac{V(S(t))R(S(t))+\sum_{k\in A}r_kv_k}{V(S(t))+\sum_{k\in A}v_k}.\]
Since $R(S(t))\geq r_i$ by assumption and $r_k\geq r_i$ for all $k$ (by the choice of $i$), every term in the numerator is at least $r_i$ times its corresponding term in the denominator. Therefore, $R(S(t'))\geq r_i$.
\Halmos
\endproof

\begin{lemma}\label{lem:3-yc1}
    Suppose that product $i$ has the highest index among all products in $N(\vecc)\cup \{i\}$. If $\fLP(\vecc^{+i},T) - \fLP(\vecc,T) > 0$, then the following holds: (i) $x_i^{\vecc}=\ZFP_i(\CFP,T)=c_i$, and $x_k^{\vecc}=c_k$ for all products $k\in \NN\setminus \{i\}$. (ii) $\ZFP_i(\vecx^{\vecc^{+i}},T)> c_i$, and $x_k^{\vecc^{+i}}=x_k^{\vecc}=c_k$ for all products $k \in \NN\setminus \{i\}$; 
    (iii) Let $\alpha_i = \ZFP_i(\vecx^{\vecc^{+i}},T)-\ZFP_i(\CFP,T)$ and let $S(t)$ denote the assortment at time $t$ in \textsc{fp}$(\CFP,T)$. Then $r_i-R(S(t))\geq 0$ for all $t\in [T-\alpha_i,T]$. 
   \end{lemma}
This lemma characterizes the revenue-ordered inventories when adding
one unit of a highest-indexed product improves the CDLP value. In the
original solution, the revenue-ordered inventory equals $\vecc$, and
product $i$ is fully consumed. In the augmented solution, every product
other than $i$ retains its original inventory. Part (iii) identifies the
terminal interval on which the additional sales of product $i$ have
nonnegative marginal revenue.

\proof{Proof of Lemma~\ref{lem:3-yc1}.}
 We first prove (i) through contradiction: Let $\{y_k^*\}_{k\in \NN}$ denote an optimal solution to $\fLP$. Since $\fLP(\vecc^{+i},T)-\fLP(\vecc,T)>0$, the constraint $y_i\leq c_i$ must be binding. Otherwise, increasing $c_i$ would not change the optimal solution of the LP, contradicting the improvement in objective value. From the proof of Lemma \ref{lem:3-rev}, we have $\ZFP_i(\vecx^{\vecc},T)=y_i^*=c_i$. Thus, the threshold $i_{\tau}\ge i$, so all products with index smaller than $i_{\tau}$ have full stock in the revenue-ordered inventory vector. Moreover, since $i_{\tau}$ is the highest indexed product with positive inventory in $\vecc$, we have that $x_k^{\vecc}=c_k>0$ for all $k< i_{\tau}$ and $x_k^{\vecc}=c_k=0$ for all $k >i_{\tau}$. 

 Next, we prove (ii) by contradiction. Suppose there exists an optimal solution $\{y_k^* (\vecc^{+i})\}_{k\in \NN\cup \{0\}}$ to $\fLP(\vecc^{+i},T)$ such that $y_i^* (\vecc^{+i})\leq c_i$. 
  Then, this solution is also a feasible solution to $\fLP(\vecc, T)$ since all the constraints are satisfied: For all $k\in \NN\setminus \{i\}$, $y_k^* (\vecc^{+i})\leq  \min \{y_0^* (\vecc^{+i})v_k,c_k\}$; For $i$, we have 
  $y_i^* (\vecc^{+i})\leq   \min \{y_0^* (\vecc^{+i})v_i,c_i\}\leq   \min \{y_0^* (\vecc^{+i})v_i,c_i+1\}$; $\sum_{k\in \NN\cup \{0\}}y_k^* (\vecc^{+i})=T$. 
However, since $\fLP(\vecc^{+i},T) - \fLP(\vecc,T) > 0$, this solution would achieve higher revenue than $\fLP(\vecc,T)$, a contradiction. Therefore, $y_i^* (\vecc^{+i})>c_i\geq 0$. From Lemma \ref{lem:3-rev}, this implies $\ZFP_i(\vecx^{\vecc^{+i}},T)= y_i^*(\vecc^{+i})> c_i$ and thus $x_i^{\vecc^{+i}}\geq \ZFP_i(\vecx^{\vecc^{+i}},T)>c_i$.

Finally, we prove (iii). From the Separability Lemma,
\(\fFP (\vecx^{\vecc^{+i}},T)-\fFP(\CFP,T)
 = \int_{t=T-\alpha_i}^T (r_i-R(S(t)))\;\mathrm{d}t, \) where $S(t)$ is the assortment at time $t$ in the fluid process \textsc{fp}$(\CFP,T)$.
 Then we prove by contradiction.
Suppose there exists $t_0\in (T-\alpha_i,T)$ such that $r_i< R(S(t_0))$. From Lemma \ref{lem:seq-rev1}, $R(S(t))\geq r_i$ for all $t \in [T-\alpha_i, t_0]$. Consider a modified inventory vector $\tilde{\mathbf{x}}^{\vecc^{+i}}$ where $\tilde{x}_i^{\vecc^{+i}}=x_i^{\vecc^{+i}}-(t_0-(T-\alpha_i))$ and 
$\tilde{x}_k^{\vecc^{+i}}=c_k$ for $k\neq i$, then $\ZFP_i(\tilde{\vecx}^{\vecc^{+i}},T)=\ZFP_i(\vecx^{\vecc^{+i}},T)-(t_0-(T-\alpha_i))=\ZFP_i(\CFP,T)+\alpha_i-(t_0-(T-\alpha_i))=\ZFP_i(\CFP,T)-t_0+T$. By the Separability Lemma:
\[\fFP (\tilde{\vecx}^{\vecc^{+i}},T)-\fFP(\CFP,T)
 = \int_{t=t_0}^T (r_i-R(S(t)))\;\mathrm{d}t.\]
Therefore:
 \[\fFP (\vecx^{\vecc^{+i}},T)-\fFP (\tilde{\vecx}^{\vecc^{+i}},T)= \int_{t=T-\alpha_i}^{t_0}(r_i - R(S(t)))\;\mathrm{d}t<0.
\]
This contradicts Lemma~\ref{lem:3-rev}, which states $\fFP(\vecx^{\vecc^{+i}},T)=\max_{\vecx\preceq \vecc^{+i}}\fFP(\vecx,T)$.
\Halmos
\endproof
\begingroup
\begin{lemma}\label{lem:3-ei}
    Suppose $i\neq j$ and product $i$ has the highest index among all products in $N(\vecc^{+j})\cup \{i\}$. If $\fLP(\vecc^{+ij},T)-\fLP(\vecc^{+j},T)>0$, then $\fLP(\vecc^{+i},T)-\fLP(\vecc,T)>0$.
   \end{lemma}

   \proof{Proof of Lemma~\ref{lem:3-ei}.}
Let $\alpha_i = \ZFP_i(\vecx^{\vecc^{+ij}},T)-\ZFP_i(\vecx^{\vecc^{+j}},T)$ denote the additional sales of product $i$ when it is added after product $j$. Let $S^{+j}(t)$ denote the assortment at time $t$ in the fluid process with inventory vector $\vecx^{\vecc^{+j}}$, and let $S(t)$ denote the assortment at time $t$ in the fluid process with inventory vector $\CFP$. From Lemma~\ref{lem:3-yc1}, applied with base inventory $\vecc^{+j}$, we know that $R(S^{+j}(t))\leq r_i$ for $t\in [T-\alpha_i ,T]$, $x_k^{\vecc^{+j}}=c_k$ for $k\in \NN\setminus \{j\}$, and $x_j^{\vecc^{+j}}=c_j+1$. Moreover, since $\vecx^{\vecc^{+ij}}$ and $\vecx^{\vecc^{+j}}$ differ only in product $i$ (Lemma~\ref{lem:3-yc1}(ii)), the Separability Lemma and Lemma~\ref{lem:3-rev} give $0<\fLP(\vecc^{+ij},T)-\fLP(\vecc^{+j},T)=\int_{T-\alpha_i}^T (r_i-R(S^{+j}(t)))\;\mathrm{d}t$, so $\alpha_i>0$ and $R(S^{+j}(t_1))<r_i$ for some $t_1\in[T-\alpha_i,T)$; by Lemma~\ref{lem:seq-rev1}(i), $R(S^{+j}(t))<r_i$ for all $t\in[t_1,T]$. Consider the fluid process with inventory vector $\CFP$. We have $\CFP\preceq \vecx^{\vecc^{+j}}$. Therefore, by Lemma~\ref{lem:app-subset}, we have $S(t)\subseteq S^{+j}(t)$ for all $t$. Applying Lemma~\ref{lem:subsetie}, this implies $R(S(t))\leq R(S^{+j}(t))<r_i$ for all $t\in [t_1,T]$.

Note also that product $i$ is fully consumed strictly before time $T$ in $\fFP(\CFP,T)$: $\ZFP_i(\vecx^{\vecc^{+j}},T)=c_i$ by Lemma~\ref{lem:3-yc1}(i), product $i$ is consumed weakly faster in $\fFP(\CFP,T)$ since $S(t)\subseteq S^{+j}(t)$, and depletion of product $i$ exactly at time $T$ in $\fFP(\vecx^{\vecc^{+j}},T)$ would give $\ZFP_i(\vecx^{\vecc^{+ij}},T)\leq c_i$, contradicting $\alpha_i>0$. Now consider adding a small amount $\epsilon\in(0,T-t_1]$ of product $i$ to the starting inventory of $\fFP(\CFP,T)$. By the Separability Lemma:
\[\fFP(\CFP+\epsilon \vece_i, T)-\fFP(\CFP, T)= \int_{t=T-\epsilon'}^T (r_i-R(S(t)))\;\mathrm{d}t>0,\]
where $\epsilon'$ is the additional amount of product $i$ consumed: $\epsilon'>0$ because product $i$ is fully consumed strictly before $T$; $\epsilon'\leq \epsilon$ since $\ZFP_i(\CFP+\epsilon\vece_i,T)\leq \ZFP_i(\CFP,T)+\epsilon$; and $\epsilon\leq T-t_1$ places $[T-\epsilon',T]$ inside $[t_1,T]$, where the integrand is strictly positive. Since $\epsilon\leq T-t_1\leq\alpha_i\leq x_i^{\vecc^{+ij}}-c_i\leq 1$ and $x_i^{\vecc}\leq c_i$, we have $\CFP+\epsilon \vece_i\preceq \vecc^{+i}$, so it is a feasible solution to $\max_{\vecx\preceq \vecc^{+i}} \fFP(\vecx, T)$ and we can conclude:
\[\fLP(\vecc^{+i},T)=\max_{\vecx\preceq \vecc^{+i}} \fFP(\vecx, T)\geq \fFP(\CFP+\epsilon \vece_i, T)>\fFP(\CFP, T)=\fLP(\vecc,T).\]
\Halmos
\endproof
\endgroup

Combining Lemma~\ref{lem:3-yc1} and \ref{lem:3-ei}, we obtain the following structural result for $\CFP$, $\vecc^{+ij}$, $\vecc^{+i}$, and $\vecc^{+j}$: 
\begin{lemma}
\label{lem:3-yc3}
    Suppose that $i\neq j$ and product $i$ has the highest index among all products in $N(\vecc^{+ij})$. If $\fLP(\vecc^{+ij},T) - \fLP(\vecc^{+j},T) > 0$, then the following holds true: $\vecx^{\vecc^{+j}}=\vecc^{+j}$, $\CFP=\vecc$,
    $x_k^{\vecc^{+i}}=c_k$ for all products $k\in \NN\setminus \{i\}$,  $x_k^{\vecc^{+ij}}=c_k$ for all products $k\in \NN\setminus \{i,j\}$ and $x_j^{\vecc^{+ij}}=c_j+1$. 
\end{lemma}

Next we prove Theorem~\ref{thm:sof}.

\noindent\textbf{Theorem~\ref{thm:sof} (restated).} {\itshape
For any given $T$, sorting products in descending order of price, with ties broken arbitrarily, is a DR-SO for $\fLP(\vecc, T)$. Formally, suppose the products are sorted in non-increasing price order, such that $r_1 \geq r_2 \geq \cdots \geq r_n$. Then for any $\vecc \in \mathbb{Z}_+^n$, and for any $i, j \in \NN$ such that $i \geq j$ and $i \geq k$ for every $k \in \NN$ with $c_k > 0$, we have:
\begin{equation*}
\fLP(\vecc + \vece_j + \vece_i, T) - \fLP(\vecc + \vece_j, T) \leq \fLP(\vecc + \vece_i, T) - \fLP(\vecc, T).
\end{equation*}
}
\proof{Proof of Theorem~\ref{thm:sof}.}
\begingroup
We first dispose of the case $i=j$. In this case the inequality reads $\fLP(\vecc+2\vece_i,T)+\fLP(\vecc,T)\leq 2\fLP(\vecc+\vece_i,T)$, which holds because the optimal value of the sales-based LP is concave in the right-hand side vector: by LP duality, it is the pointwise minimum of affine functions of $\vecc$. In the remainder of the proof we assume $i>j$, so that $i$ and $j$ are distinct products.
First, observe that both sides of the inequality are non-negative since $\fLP(\vecc, T)$ is non-decreasing in $\vecc$ (as $\vecc$ is the right-hand side of an LP). Thus, if $\fLP(\vecc^{+ij},T)-\fLP(\vecc^{+j},T)=0$, the inequality holds trivially. 
Consider the case where $\fLP(\vecc^{+ij},T)-\fLP(\vecc^{+j},T)>0$. Recall that $\fLP(\vecc,T)=\fFP(\CFP,T)$. To establish the desired inequality, it suffices to show:
\[
 \fFP(\CijFPs,T)-\fFP(\CjFPs,T)\leq \fFP(\CiFPs,T)-\fFP(\CFP,T).
\]
\endgroup

Since $\fLP(\vecc^{+ij},T)-\fLP(\vecc^{+j},T)>0$, from Lemma~\ref{lem:3-yc3}, we have $x_k^{\vecc^{+j}}=x_k^{\vecc^{+ij}}=(\vecc^{+j})_k$ for all $k\in \NN\setminus \{i\}$, i.e., all the inventories in $\vecc^{+j}$ are included in $\vecx^{\vecc^{+j}}$ and $\vecx^{\vecc^{+ij}}$. $\vecx^{\vecc^{+ij}}$ and $\vecx^{\vecc^{+j}}$ differ only in the inventory of $i$, and we have $x_i^{\vecc^{+ij}}>(\vecc^{+j})_i=x_i^{\vecc^{+j}}$. Then, by the Separability Lemma, defining $\alpha_1=\ZFP_i(\CijFPs,T)-\ZFP_i(\CjFPs,T)$, we have:
 \begin{align}
 \fFP (\CijFPs,T)-\fFP(\vecx^{\vecc^{+j}},T)&=
 \alpha_1 r_i+ \fFP(\vecc^{+j}, T-\alpha_1)-\fFP (\vecc^{+j},T) \nonumber \\
 &= \int_{t=T-\alpha_1 }^T(r_i-R(S^{+j}(t)))\;\mathrm{d}t. \label{eq:dif1-app}
 \end{align}

Similarly, since $\fLP(\vecc^{+ij},T)-\fLP(\vecc^{+j},T)>0$, we have $x_k^{\vecc^{+i}}=x_k^{\vecc}$ for all $k\in \NN\setminus \{i\}$ from Lemma~\ref{lem:3-yc3}. Let $\alpha_2=\ZFP_i(\CiFPs,T)-\ZFP_i(\CFP,T)$. From the Separability Lemma, we have: 
\begin{equation} \label{eq:dif2-app}
\fFP(\CiFPs,T)-\fFP(\CFP,T)= \int_{t=T-\alpha_2}^T(r_i-R(S(t)))\;\mathrm{d}t.
\end{equation}
So we just need to show \eqref{eq:dif1-app} $\leq$ \eqref{eq:dif2-app}. To establish this inequality, it suffices to show the following: (i) $\alpha_1\leq \alpha_2$; (ii) $ r_i-R(S^{+j}(t))\geq 0$ for $t\in [T-\alpha_1,T]$ and $ r_i-R(S(t))\geq 0$ for $t\in [T-\alpha_2,T ]$; (iii) for $t\in [T-\alpha_1,T ]$, $R(S(t))\leq R(S^{+j}(t))\leq r_i$.

\begingroup
We have shown (ii) in Lemma~\ref{lem:3-yc1}. Now we show (i). Recall that $\fFP(\vecx^{\vecc^{+ij}},T)=\max_{\vecx\preceq \vecc^{+ij}}\fFP(\vecx,T) $ and Lemma~\ref{lem:3-yc3} shows that $\vecx^{\vecc^{+j}}=\vecc^{+j}$, $x_k^{\vecc^{+ij}}=x_k^{\vecc^{+j}}$ for all $k\in \NN\setminus \{i\}$. Thus $\vecx^{\vecc^{+ij}}$ can be written as $\vecc^{+j}+y\vece_i$ 
where $0<y\leq 1$. Let $z$ denote the additional quantity of product $i$ consumed in $\fFP(\vecx^{\vecc^{+ij}},T)$, i.e., $z=\ZFP_i(\CijFPs,T)-\ZFP_i(\vecx^{\vecc^{+j}},T)=\ZFP_i(\CijFPs,T)-(\vecc^{+j})_i$. We have the following from the Separability Lemma:
\begin{equation}\label{eq:sep-alpha}
\fFP(\vecc^{+j}+y \vece_i,T)- \fFP(\vecc^{+j},T)=z r_i+ \fFP(\vecc^{+j},T-z)- \fFP(\vecc^{+j},T)=\int_{T-z}^T (r_i- R(S^{+j}(t))) \;\mathrm{d}t.
\end{equation}
Here $z\leq \ZFP_i(\vecc^{+ij},T)-(\vecc^{+j})_i$ since the maximum consumption of $i$ in $\fFP(\CijFPs,T)$ is $\ZFP_i(\vecc^{+ij},T)$. Conversely, every value $z'\in (0,\ZFP_i(\vecc^{+ij},T)-(\vecc^{+j})_i]$ is attained by adding $y'=z'$ units of product $i$; note that $z'\leq 1$ since $\ZFP_i(\vecc^{+ij},T)\leq (\vecc^{+ij})_i=(\vecc^{+j})_i+1$. Indeed, since $\vecc^{+j}+\vece_i=\vecc^{+ij}$, the fluid process with inventory $\vecc^{+j}+z'\vece_i$ coincides with $\fFP(\vecc^{+ij},T)$ until product $i$ stocks out. The cumulative consumption of product $i$ in $\fFP(\vecc^{+ij},T)$ is continuous in time and reaches $(\vecc^{+j})_i+z'$ by time $T$, so product $i$ stocks out by time $T$ in the former process and $\ZFP_i(\vecc^{+j}+z'\vece_i,T)=(\vecc^{+j})_i+z'$. By \eqref{eq:sep-alpha}, the gain from adding inventory of product $i$ depends on the added amount only through the additional consumption $z$. Since $\CijFPs\in \argmax_{\vecx\preceq \vecc^{+ij}}\fFP(\vecx,T)$, we have:
\begin{align*}
\alpha_1&= \ZFP_i(\CijFPs,T)-\ZFP_i(\vecc^{+j},T)\\
& \in \argmax_{0<z\leq \ZFP_i(\vecc^{+ij},T)-(\vecc^{+j})_i} zr_i+ \fFP(\vecc^{+j},T-z)- \fFP(\vecc^{+j},T) \\
&=\argmax_{0< z\leq \ZFP_i(\vecc^{+ij},T)-(\vecc^{+j})_i} \int_{T-z}^T (r_i- R(S^{+j}(t))) \;\mathrm{d}t.
\end{align*}
Here the membership in the second line holds because every $z$ in the given range is attainable, and the last equality follows from \eqref{eq:sep-alpha}.

Since $\fFP(\CijFPs,T)-\fFP(\vecc^{+j},T)>0$, there exists $0\leq t<T$ such that $r_i-R(S^{+j}(t)) > 0$. Let $\tau_1=\min\{t:r_i- R(S^{+j}(t)) > 0\}$. Then from Lemma~\ref{lem:seq-rev1}, $r_i>R(S^{+j}(t))$ for all $t\geq \tau_1$ and $r_i\leq R(S^{+j}(t))$ for all $t<\tau_1$. Hence $\int_{T-z}^T (r_i- R(S^{+j}(t)))\;\mathrm{d}t$ is strictly increasing in $z$ for $z\leq T-\tau_1$ and nonincreasing for $z\geq T-\tau_1$, so $\min\{T-\tau_1,\ZFP_i(\vecc^{+ij},T)-(\vecc^{+j})_i\}$ is its smallest maximizer. Since every $z$ in the range is attainable, we may take $\CijFPs$ to be the optimal solution with the smallest consumption of product $i$; we adopt the same convention for $\vecx^{\vecc^{+i}}$ below. Therefore, $\alpha_1=\min\{T-\tau_1,\ZFP_i(\vecc^{+ij},T)-(\vecc^{+j})_i\}$. 
\endgroup

Similarly, define
$\tau_2=\min\{t\in[0,T]:r_i>R(S(t))\}$. 
Similarly, $\alpha_2=\min\{T-\tau_2, \ZFP_i(\vecc^{+i},T)-c_i\}$. To show $\alpha_1\leq \alpha_2$, it suffices to prove that $\tau_1\geq \tau_2$ and $\ZFP_i(\vecc^{+ij},T)\leq  \ZFP_i(\vecc^{+i},T)$. The latter inequality is straightforward to show, as the consumption rate of product $i$ in $\fFP(\vecc^{+ij},T)$ is less than or equal to that in $\fFP(\vecc^{+i},T)$ due to the substitutability property. To show $\tau_1\geq \tau_2$, note that $S(\tau_1)\subseteq S^{+j}(\tau_1)$ from Lemma~\ref{lem:app-subset}. Since $ R(S^{+j}(\tau_1))< r_i$ by the definition of $\tau_1$, from Lemma~\ref{lem:subsetie}, we have $R(S(\tau_1))\leq R(S^{+j}(\tau_1))< r_i$. Thus, $\tau_2\leq  \tau_1$.   

Finally, we prove (iii). 
From Lemma~\ref{lem:app-subset}, we have $S(t)\subseteq S^{+j}(t)$ for all $t$. Since $R(S^{+j}(t))< r_i$ for all $t\in [T-\alpha_1,T]$, we have $R(S(t))\leq R(S^{+j}(t))<r_i$ for all $t\in [T-\alpha_1,T]$ from Lemma~\ref{lem:subsetie}.

Combining (i), (ii) and (iii), we have   \[
 \int_{t=T-\alpha_1}^T(r_i-R(S^{+j}(t)))\;\mathrm{d}t\leq  \int_{t=T-\alpha_2 }^T(r_i-R(S(t)))\;\mathrm{d}t,
\]
which implies the desired result.
\Halmos
\endproof

\subsection{Counterexample to the DR-SO Property of Multi-type CDLP}\label{app:ce-SO-multitype}
Here we give a counterexample to show that the descending order of prices is not a DR-SO for the multi-type CDLP objective function, denoted by $\fLPM(\vecc,T)$. Consider a ground set of four products with prices $r_1= 46$, $r_2=46$, $r_3=12$, and $r_4=12$. Suppose there are two types of customers with the following attraction parameters: $v_{1,1}=6$, $v_{1,2}=10$, $v_{1,3}=17$, $v_{1,4}=15$ for type 1, and $v_{2,1}=16$, $v_{2,2}=11$, $v_{2,3}=19$, $v_{2,4}=20$ for type 2. Consider the inventory vector $\vecc=[4,3,7,0]$, and let the expected type arrivals be $\sum_{t=1}^{T}\lambda_{1,t}=2$ and $\sum_{t=1}^{T}\lambda_{2,t}=8$ (e.g., $T=10$ with $\lambda_{1,t}=0.2$ and $\lambda_{2,t}=0.8$ for all $t$). Let $j=1$ and $i=4$. We have: $\fLPM(\vecc,T)\approx 354.559$, $\fLPM(\vecc+\vece_i,T)\approx 354.683$, $\fLPM(\vecc+\vece_j,T)\approx 387.809$, $\fLPM(\vecc+\vece_j+\vece_i,T)\approx 387.944$. Then, $\fLPM(\vecc+\vece_i,T)-\fLPM(\vecc,T)<\frac{225}{1666}=\fLPM(\vecc+\vece_j+\vece_i,T)-\fLPM(\vecc+\vece_j,T)$. Therefore, the descending order of prices is not a DR-SO for $\fLPM$.

\begingroup
\subsection{Proof of Lemma~\ref{lem:alloc-equiv}} \label{app:hat}
\noindent\textbf{Lemma~\ref{lem:alloc-equiv} (restated).} {\itshape
For any inventory $\vecc$ and horizon $T$,
\(
\fLPM(\vecc,T)=\max_{\matc}\ \Bigl\{\hat{f}(\matc,T):\ \textstyle\sum_{k\in\MM}C_{k,i}\le c_i\ \ \forall i\in\NN,\ \ \matc\ge 0\Bigr\}.
\)
}
\proof{Proof of Lemma~\ref{lem:alloc-equiv}.} Let $\{y_{ki}\}$ be an optimal solution to the multi-type CDLP \eqref{eq:multi-sblp} and set $C_{k,i}=y_{ki}$. Then $\sum_{k\in\MM}C_{k,i}=\sum_{k\in\MM}y_{ki}\leq c_i$ for every $i$, so $\matc$ is a feasible allocation; and for each type $k$, the vector $\{y_{ki}\}_{i\in\NN}$ together with the no-purchase value $y_{k0}=\Lambda_k-\sum_{i\in\NN}y_{ki}$ is feasible for the sales-based formulation of the single-type program $\fLP_k(\matc_k,\Lambda_k)$. Hence $\hat{f}(\matc,T)=\sum_{k\in\MM}\fLP_k(\matc_k,\Lambda_k)\geq \sum_{k\in\MM}\sum_{i\in\NN}r_iy_{ki}=\fLPM(\vecc,T)$. Conversely, fix any feasible allocation $\matc$ and, for each $k$, an optimal solution $\{y^k_{i}\}_{i\in\NN\cup\{0\}}$ of the sales-based formulation of $\fLP_k(\matc_k,\Lambda_k)$. The combined solution $y_{ki}=y^k_i$ is feasible for \eqref{eq:multi-sblp}: $\sum_{k\in\MM}y_{ki}\leq \sum_{k\in\MM}C_{k,i}\leq c_i$ for every $i$, and the per-type constraints hold by construction. Hence $\fLPM(\vecc,T)\geq \hat{f}(\matc,T)$ for every feasible $\matc$. Combining the two directions proves the lemma.
\Halmos
\endgroup
\endproof
\subsection{DR-SO Property of \texorpdfstring{$\hat{f}$}{f}} \label{app:subhatf}

Let $\mate_{k,j}\in \RR_+^{mn}$ denote the extended inventory vector where the inventory of product $j$ for type-$k$ customers is $1$ and all other entries are zeros. \begingroup For an extended inventory vector $\matc\in\RR_+^{mn}$, let $N(\matc)=\{i\in\NN:\sum_{k\in\MM}C_{k,i}>0\}$ denote the set of products with positive total allocation. In this subsection, since the horizon \(T\) remains fixed throughout, we suppress it below and write \(\hat f(\matc)\) for \(\hat f(\matc,T)\). 
\endgroup
\begin{lemma}
    For any $\matc\in \mathbb{Z}^{mn}_+$, given any two product indices $i, j \in \NN$ such that $i\geq j$ and $i\geq k$ for every product $k$ with positive inventory in $\matc$, and customer types $w, u \in \MM$, the following inequality holds:\[
    \hat{f}(\matc+\mate_{w,j} +\mate_{u,i} )-\hat{f}(\matc+\mate_{w,j} ) \leq \hat{f}(\matc+\mate_{u,i} )- \hat{f}(\matc).    \vspace{-10pt}
    \]
    \label{lem:sof-mint}
    \end{lemma}
    \vspace{-10pt}
This lemma is a straightforward generalization of Theorem~\ref{thm:sof} by considering different customer types. 
\proof{Proof of Lemma~\ref{lem:sof-mint}.} We analyze two cases based on the relationship between customer types $w$ and $u$. \begingroup
Case 1: When $w\not =u$, we have:
\begin{align*}
&\hat{f}(\matc+\mate_{w,j} +\mate_{u,i} )-\hat{f}(\matc+\mate_{w,j} )\\
&= \fLP_w(\matc_w +\vece_{j},\Lambda_w)-\fLP_w(\matc_w+\vece_{j},\Lambda_w)+\fLP_u(\matc_u +\vece_{i},\Lambda_u)-\fLP_u(\matc_{u},\Lambda_u) \\
& = \hat{f}(\matc+\mate_{u,i} )- \hat{f}(\matc).
\end{align*}
Case 2: When $w=u$, we obtain:
\begin{align*}
\hat{f}(\matc+\mate_{w,j} +\mate_{u,i} )-\hat{f}(\matc+\mate_{w,j} ) & = \fLP_w(\matc_w +\vece_{j}+\vece_{i},\Lambda_w)-\fLP_w(\matc_w+\vece_{j},\Lambda_w) \\
& \leq \hat{f}(\matc+\mate_{w,i} )- \hat{f}(\matc).
\end{align*}
\endgroup
Here the inequality is from the DR-SO property of the single-type CDLP, $\fLP$ (Theorem~\ref{thm:sof}).
\Halmos
\endproof
\begingroup
\begin{remark}[Fractional allocations]
The DR-SO property extends from integral to fractional allocation vectors. 
For the single-type CDLP, homogeneity of the sales-based formulation gives
\(
\fLP(\vecc,\Lambda)
=\delta\,\fLP(\vecc/\delta,\Lambda/\delta)
\) for  \(\delta>0\). 
The proof of Theorem~\ref{thm:sof} applies to every real fluid horizon.
Applying it at horizon $\Lambda/\delta$ and then multiplying by $\delta$
therefore proves the same DR-SO inequality for increments of size $\delta$ on
the grid $\delta\mathbb Z_+^n$.  Summing the single-type inequalities preserves
the property for $\hat f$; increments assigned to different customer types
commute because $\hat f$ is separable across types.  Thus the property holds on
every rational grid.  Finally, the value of a linear program is continuous in
these capacity right-hand sides, so taking rational approximations proves the
inequality for arbitrary fractional allocation vectors.  Subsequent uses of
the DR-SO property of $\hat f$ and of the interleaved-partition inequality \eqref{eq:cor3} invoke this
fractional extension.
\end{remark}

The proofs below apply the interleaved-partition inequality of \cite{sof} to blocks of extended inventory vectors that may contain fractional entries. The following lemma records the form we use.

\begin{lemma}[Interleaved-partition inequality]\label{lem:interleave}
Let $\mathbf{B}_1,\ldots,\mathbf{B}_{\kappa+1}$ and $\mathbf{A}_1,\ldots,\mathbf{A}_{\kappa}$ be extended inventory vectors in $\RR_+^{mn}$, some possibly zero, and consider the sequence $\mathbf{B}_1,\mathbf{A}_1,\mathbf{B}_2,\mathbf{A}_2,\ldots,\mathbf{A}_{\kappa},\mathbf{B}_{\kappa+1}$. Suppose every product appearing in an earlier block of this sequence has index at most every product appearing in a later block. Write $\mathbf{A}(l)=\sum_{l'=1}^{l}\mathbf{A}_{l'}$ with $\mathbf{A}(0)=\mathbf{0}$. Then,
\[
\hat f\Big(\sum_{l\in[\kappa+1]}\mathbf{B}_l+\mathbf{A}(\kappa)\Big)\leq \hat f(\mathbf{A}(\kappa))+\sum_{l\in[\kappa+1]}\hat f(\mathbf{B}_l\mid \mathbf{A}(l-1)),
\]
where $\hat f(\mathbf{B}\mid\mathbf{A})=\hat f(\mathbf{B}+\mathbf{A})-\hat f(\mathbf{A})$ denotes the marginal benefit of adding $\mathbf{B}$ to $\mathbf{A}$.
\end{lemma}
\proof{Proof of Lemma~\ref{lem:interleave}.}
By the continuity of the LP value in the inventory right-hand sides, it suffices to prove the inequality when every entry of every block is rational. Fix $\delta>0$ such that every entry is an integer multiple of $\delta$ and, by the homogeneity relation in the preceding remark, rescale by $1/\delta$ so that all blocks lie in $\mathbb Z_+^{mn}$. Expand the rescaled blocks over a ground set containing one distinct copy of every unit, with the copies ordered by descending product price, ties broken first
by product index, then by the position of their block in the stated
sequence, and finally arbitrarily within a block. Define the set function $F(S)=\hat f(\bm{\chi}(S))$, where $\bm{\chi}(S)\in\mathbb Z_+^{mn}$ counts the copies of each type--product pair in $S$ and $\hat f$ is evaluated at the rescaled per-type horizons $\Lambda_k/\delta$. By Theorem~\ref{thm:sof}, Lemma~\ref{lem:sof-mint}, and the preceding remark, all of which apply at every real horizon, $F$ satisfies the submodular-order condition \eqref{eq:so} with respect to this order, and the blocks correspond to sets of copies interleaved consecutively in the same order. Corollary 3 of \cite{sof} applied to $F$ therefore yields the stated inequality for the rescaled blocks; undoing the rescaling and taking limits of rational approximations completes the proof.
\Halmos
\endproof
\endgroup
\begingroup
We also record a subadditivity property of $\fLPM$ and $\hat{f}$, used repeatedly below.

\begin{lemma}[Subadditivity]\label{lem:subadd}
For any $\vecc,\vecc'\in\RR_+^{n}$ and any $T$, $\fLPM(\vecc+\vecc',T)\leq \fLPM(\vecc,T)+\fLPM(\vecc',T)$. Likewise, for any $\matc,\matc'\in\RR_+^{mn}$, $\hat{f}(\matc+\matc',T)\leq \hat{f}(\matc,T)+\hat{f}(\matc',T)$. In particular, $\fLPM(\vecc,T)\leq \sum_{i\in\NN}c_i\,\fLPM(\vece_i,T)$ for any $\vecc\in\mathbb{Z}_+^n$. All three inequalities are preserved by taking expectations over $T$.
\end{lemma}
\proof{Proof of Lemma~\ref{lem:subadd}.}
Let $\{y_{ki}\}$ be an optimal solution to the multi-type SBLP \eqref{eq:multi-sblp} with inventory $\vecc+\vecc'$. For each $k$ and $i$, split $y_{ki}=a_{ki}+b_{ki}$ with $a_{ki},b_{ki}\geq 0$, $\sum_{k\in\MM}a_{ki}\leq c_i$, and $\sum_{k\in\MM}b_{ki}\leq c'_i$; this is possible because $\sum_{k\in\MM}y_{ki}\leq c_i+c'_i$. For each type $k$, let $y_{k0}$ denote the no-purchase value of the optimal solution and set $a_{k0}=\Lambda_k-\sum_{i\in\NN}a_{ki}\geq y_{k0}$. The vector $\{a_{ki}\}_{i\in\NN}$ with the no-purchase value $a_{k0}$ is feasible for \eqref{eq:multi-sblp} with inventory $\vecc$: the balance constraints hold by construction, and $a_{ki}\leq y_{ki}\leq v_{ki}y_{k0}\leq v_{ki}a_{k0}$ for every $i\in\NN$. Hence $\fLPM(\vecc,T)\geq \sum_{k\in\MM}\sum_{i\in\NN}r_ia_{ki}$, and similarly $\fLPM(\vecc',T)\geq \sum_{k\in\MM}\sum_{i\in\NN}r_ib_{ki}$. Adding the two bounds gives the first inequality. For the second inequality, apply an analogous splitting to an optimal solution of the sales-based formulation of each single-type program $\fLP_k(\matc_k+\matc'_k,\Lambda_k)$, with the no-purchase values set to $\Lambda_k-\sum_{i\in\NN}a_{ki}$ and $\Lambda_k-\sum_{i\in\NN}b_{ki}$, respectively; feasibility for $\fLP_k(\matc_k,\Lambda_k)$ and $\fLP_k(\matc'_k,\Lambda_k)$ follows by the same argument. Summing over $k$ yields the inequality for $\hat{f}$. The third inequality follows by writing $\vecc=\sum_{i\in\NN}\sum_{u=1}^{c_i}\vece_i$ and applying the first inequality repeatedly. The expectation versions follow from linearity of expectation.
\Halmos
\endgroup
\endproof

\subsection{Threshold Augmentation Algorithm for Multi-type CDLP} \label{app:cdlp-m-T-alg}

\subsubsection{Cardinality-Constrained Approximation Algorithm.}
We first present the overall approximation algorithm for the multi-type CDLP under a cardinality constraint. The algorithm uses a threshold-based approach, where we try a geometric sequence of threshold values $\tau$ ranging from $\frac{1}{K}\max_{i \in \NN}\fLPM(\vece_i,T)$ to $\max_{i \in\NN}\fLPM(\mathbf{e}_i,T)$. 

{\SetAlgoNoLine%
  \begin{algorithm}[H]
    \KwIn{Cardinality $K$, error $\epsilon \in (0,1)$, number of customers $T$, products $\NN$ sorted in descending order of price: $r_1\geq r_2\geq \cdots \geq r_n$} \KwOut{Inventory vector $\vecc$, inventory allocation vector $\matc$}
Initialize $\tau = \frac{1}{K}\max_{i \in \NN }\fLPM(\vece_i,T)$\;
       \For{$i \gets 1$ \textbf{to} $\lceil \log_{1+\epsilon }K \rceil$} {
      $\vecc^i, \matc^i= \text{Cardinality-Constrained Threshold Add for DAP}(K,(1+\epsilon)^{i-1}\tau)$\;
    }
    \Return{$\vecc,\matc$=Best of $(\vecc^1,\matc^1)$, $(\vecc^2,\matc^2)$, $\cdots$, $(\vecc^{\lceil \log_{1+\epsilon }K \rceil},\matc^{\lceil \log_{1+\epsilon }K \rceil})$, evaluated by $\hat{f}(\cdot,T)$}
    \caption{Cardinality-Constrained Optimization for DAP}\label{alg:tau-m}
  \end{algorithm}}

The Cardinality-Constrained Threshold Add for DAP algorithm (Algorithm~\ref{alg:cdlp-multi-type}) is defined in Section~\ref{sec:dap-alloc} and returns a feasible inventory vector and inventory allocation vector given a threshold value $\tau$.

\subsubsection{Proof of Theorem~\ref{thm:3-ap-m-ra-T}.} \label{app:prf-3-ap-m-ra-T}

\noindent\textbf{Theorem~\ref{thm:3-ap-m-ra-T} (restated).} {\itshape
For any $\epsilon\in(0,1)$, the problem $\max_{\vecc\in \mathcal{F}}\fLPM(\vecc,T)$ can be approximated within a factor of $\frac{1}{2}-\epsilon$ under a cardinality constraint.
}\medskip

\proof{Proof of Theorem~\ref{thm:3-ap-m-ra-T}.} 

Let $\vecc^{\operatorname{ALG}}$ and $\matc^{\operatorname{ALG}}$ denote the output of Algorithm~\ref{alg:tau-m}. 
 Let $\vecc^{\operatorname{OPT}}$ and $\matc^{\operatorname{OPT}}$ denote the optimal solution to \eqref{OPT-matc} (recall it is $
\max_{\sum_{k\in \MM} \matc_k = \vecc ,\vecc\in \mathcal{F}}  \hat{f}( \matc,T)$). Let $\operatorname{OPT}$ denote the optimal objective value of \eqref{OPT-matc}. Our goal is to show $\fLPM(\vecc^{\operatorname{ALG}}, T) \geq (\frac{1}{2}-\epsilon)\operatorname{OPT}$. Since the LP allocates $\vecc^{\operatorname{ALG}}$ across different customer types optimally, we have $\fLPM(\vecc^{\operatorname{ALG}},T)\geq \hat{f}(\matc^{\operatorname{ALG}},T)$. Thus, proving $\hat{f}(\matc^{\operatorname{ALG}}, T) \geq (\frac{1}{2}-\epsilon)\operatorname{OPT}$ suffices to establish the desired guarantee. For simplicity, 
we omit $T$ in $\hat{f}(\cdot, T)$ throughout the proof.

Let $\tau_i = \tau(1+\epsilon)^{i-1}$ and let $\vecc^i$ and $\matc^i$ denote the output of the Threshold Add algorithm (Algorithm~\ref{alg:cdlp-multi-type}) for $\tau_i$. Since the algorithm outputs the best pair among all $(\vecc^i,\matc^i)$, it suffices to show that some $\tau_i$ satisfies $\hat{f}( \matc^{i})\geq \frac{1-\epsilon}{2}\operatorname{OPT}$.
If $K\tau  \geq \frac{1}{2}\operatorname{OPT}$, then $\hat{f}(\matc^{\lceil \log_{1+\epsilon }K \rceil})\geq \frac{K \tau}{1+\epsilon} \geq \frac{1-\epsilon}{2}\operatorname{OPT}$. We thus focus on the case where $ K\tau<\frac{\operatorname{OPT}}{2}$. Due to subadditivity of $\fLPM$ (Lemma~\ref{lem:subadd}), $\operatorname{OPT}\leq K \max_{i \in \NN }\fLPM(\vece_i,T)=K^2 \tau$, implying $K\tau\geq \frac{\operatorname{OPT}}{K}$.
Therefore, there exists an index $i$ such that $(1-\epsilon)\frac{\operatorname{OPT}}{2K} \leq \tau_i \leq \frac{\operatorname{OPT}}{2K}$. We proceed to show $\hat{f}( \matc^{i})\geq \frac{1-\epsilon}{2}\operatorname{OPT}$ for this $\tau_i$. We distinguish two cases based on whether $\|\matc^{i}\|_1< K$.

If $\|\matc^i\|_1=K$, then $\hat{f}(\matc^i)\geq K\tau_i \geq  \frac{1-\epsilon}{2}\operatorname{OPT} $ from the Threshold Add algorithm. So it remains to discuss the case where $ \|\matc^i\|_1 <K$. For simplicity, we omit the superscript and write $\vecc^i$ as $\vecc$, $\matc^i$ as $\matc$ and $\tau_i$ as $\tau$.
 Let $\kappa$ denote the number of distinct products with positive inventory in $\matc$, and let $j_1<\cdots<j_{\kappa}$ denote the index set of $N(\vecc)$.  When $\kappa\geq1$, define
\[
Q_1=\{1,\ldots,j_1-1\},\qquad
Q_l=\{j_{l-1},\ldots,j_l-1\}\ (2\le l\le\kappa),\qquad
Q_{\kappa+1}=\{j_\kappa,\ldots,n\};
\]
when $\kappa=0$, define $Q_1=\NN$.  Let $\matc_{\cdot j}\in \RR_+^{mn}$ denote the extended inventory vector for product $j$, where the inventory of product $j$ allocated to customer type $k$ is $C_{k,j}$, and all other entries are zeros; $\matc^{\operatorname{OPT}}_{\cdot j}$ is defined analogously.  For each $l$, let
$\matc_l^{\operatorname{OPT}}=\sum_{q\in Q_l}
\matc_{\cdot q}^{\operatorname{OPT}}$.  Then
$\matc^{\operatorname{OPT}}+\matc$ can be expressed as the following
\emph{interleaved partition} (with only $\matc_1^{\operatorname{OPT}}$ when
$\kappa=0$):
\[
\matc^{\operatorname{OPT}}_1,  \matc_{\cdot j_1}, \matc^{\operatorname{OPT}}_{2},  \matc_{\cdot j_2}, \cdots, \matc^{\operatorname{OPT}}_{ \kappa}, \matc_{\cdot j_{\kappa}}, \matc^{\operatorname{OPT}}_{\kappa+1}. 
\]
We define $\matc(l)= \sum_{l'=1}^l \matc_{\cdot j_{l'}}$ and set $\matc(0)=\mathbf{0}$.
Note that some inventory vectors in the partition may be zero. 
Here, the \emph{interleaved partition} is adapted from \cite{sof} with a key difference: while \cite{sof} partitions the set of items, we partition the extended inventory vectors. These extended inventory vectors may include fractional values. The underlying principle is that product $j_l$ is positioned after all products in $N(\matc^{\operatorname{OPT}}_{l})$ and before all products in $N(\matc^{\operatorname{OPT}}_{l+1})$ in the submodular order. Leveraging the DR-SO property of $\hat{f}$, \begingroup we apply Lemma~\ref{lem:interleave} \endgroup to obtain:
\begin{equation}\label{eq:cor3}
    \hat{f}(\matc^{\operatorname{OPT}}+\matc) \leq \hat{f}(\matc(\kappa))+\sum_{l\in [\kappa+1]}\hat{f}(\matc^{\operatorname{OPT}}_l|\matc(l-1)),
\end{equation}
where $\hat{f}(\matc^{\operatorname{OPT}}_l|\matc(l-1))$ denotes the marginal benefit
of adding $\matc^{\operatorname{OPT}}_l$ to $\matc(l-1)$. 
To bound these marginal benefits, we apply Lemma~\ref{lem:3-m-tau}, stated and proved below. This establishes the approximation ratio of $\frac{1}{2}-\epsilon$: 
\begin{align*}
    \operatorname{OPT}& \leq  \hat{f}(\matc^{\operatorname{OPT}}+\matc) \leq \hat{f}(\matc(\kappa) )+\sum_{l\in[\kappa+1]}\hat{f}(\matc^{\operatorname{OPT}}_l|\matc(l-1)) \\
    &=  \hat{f}(\matc )+\sum_{l\in[\kappa+1]}\hat{f}(\matc^{\operatorname{OPT}}_l|\matc(l-1))  \leq \hat{f}(\matc)+\sum_{l\in[\kappa+1]}\|\matc^{\operatorname{OPT}}_l\|_1 \tau \\
  & = \hat{f}(\matc)+\|\matc^{\operatorname{OPT}}\|_1 \tau  \leq  \hat{f}(\matc)+\frac{1}{2}\operatorname{OPT}.
\end{align*}

Here, the first inequality results from the monotonicity of the LP. The second inequality is derived from \eqref{eq:cor3}. The first equality holds since $\matc(\kappa)=\matc $. The third inequality follows from Lemma~\ref{lem:3-m-tau}, which bounds the marginal benefit per product. Finally, the last inequality holds given that $\|\matc^{\operatorname{OPT}}\|_1\leq K $ and $\tau \leq \frac{\operatorname{OPT}}{2K}$. From the above, we conclude $ \hat{f}(\matc)\geq  \frac{1}{2}\operatorname{OPT} $, thus completing our proof.
\Halmos
\endproof

\begin{lemma}\label{lem:3-m-tau}
Use the notation from the proof of Theorem~\ref{thm:3-ap-m-ra-T}, with the horizon argument $T$ and threshold index $i$ suppressed. Suppose that Algorithm~\ref{alg:cdlp-multi-type}, run at threshold $\tau$, returns an allocation $\matc$ satisfying $\|\matc\|_1<K$. For any $l\in [\kappa+1]$, $\hat{f}(\matc^{\operatorname{OPT}}_l|\matc(l-1)) \leq  \|\matc^{\operatorname{OPT}}_l\|_1\tau$.
\end{lemma}
\proof{Proof of Lemma \ref{lem:3-m-tau}.} The proof follows directly from the Threshold Add algorithm and the fractional DR-SO property. Fix $l\in[\kappa+1]$ and
$q\in Q_l$, and let
\(
d_q
=
\|\matc_{\cdot q}^{\operatorname{OPT}}\|_1
=
c_q^{\operatorname{OPT}}
\in\mathbb Z_+.
\)
Because the returned solution is underfilled, the algorithm finishes
processing product $q$ by rejecting one additional unit at
$\matc(l-1)$. This includes the case $q=j_{l-1}$, because
$\matc(l-1)$ already contains all copies of $q$ accepted by the
algorithm. Hence every feasible one-unit allocation $\mathbf X$ of
product $q$ satisfies
\(
\hat f(\mathbf X\mid\matc(l-1))<\tau\). 
If $d_q\geq1$, write
\(
\matc_{\cdot q}^{\operatorname{OPT}}
=
\sum_{s=1}^{d_q}\mathbf X^{q,s}, 
\mathbf X^{q,s}
=
\frac{1}{d_q}\matc_{\cdot q}^{\operatorname{OPT}},
\)
so that every $\mathbf X^{q,s}$ is a one-unit allocation of product $q$.
By fractional DR-SO,
\[
\begin{aligned}
\hat f(
\matc_{\cdot q}^{\operatorname{OPT}}
\mid\matc(l-1))=
\sum_{s=1}^{d_q}
\hat f\!\left(
\mathbf X^{q,s}
\mid
\matc(l-1)+\sum_{h<s}\mathbf X^{q,h}
\right)\leq
\sum_{s=1}^{d_q}
\hat f(\mathbf X^{q,s}\mid\matc(l-1)) <
d_q\tau.
\end{aligned}
\]

Finally, enumerate $N(\matc_l^{\operatorname{OPT}})$ as
$q_1<\cdots<q_r$. Telescoping in product order and applying fractional
DR-SO once more,
\[
\begin{aligned}
\hat f(
\matc_l^{\operatorname{OPT}}
\mid\matc(l-1)) \leq
\sum_{t=1}^{r}
\hat f(
\matc_{\cdot q_t}^{\operatorname{OPT}}
\mid\matc(l-1)) \leq
\sum_{t=1}^{r}
\|\matc_{\cdot q_t}^{\operatorname{OPT}}\|_1\tau=
\|\matc_l^{\operatorname{OPT}}\|_1\tau.
\end{aligned} \Halmos
\]
\endproof

\subsubsection{Budget-Constrained Approximation Algorithms.} \label{app:multi-budget}

For the budget constraint, we modify the Threshold Add algorithm as follows.
\vspace{-10pt}
{\SetAlgoNoLine%
  \begin{algorithm}[htbp]
     \caption{Budget-Constrained Threshold Add for DAP}\label{alg:cdlp-multi-type-budget}
   \textbf{Input:} Budget $B$, weights $(b_i)_{i \in \NN }$, threshold $\tau$, products $\NN$ sorted in descending order of price: $r_1\geq r_2\geq \cdots \geq r_n$, number of customers $T$\; \textbf{Output:} $\vecc,\matc $; Initialize $\vecc=\mathbf{0}^n, \matc=\mathbf{0}^{mn}$\;
       \For{$i \in \NN $ }{ 
         \While{$\sum_{j \in \NN } b_j c_j + b_i \leq B$}{
         Find the optimal allocation of item $i$: $\matz^* = \arg \max_{\matz \in \mathbb{R}^{mn}_+ :\sum_{k \in \MM}Z_{k,i} = 1,  \sum_{k \in \MM}Z_{k,j} = 0, \forall j\not= i} \hat{f}(\matc + \matz)-\hat{f}(\matc)$  \;
    \lIf{$\hat{f}(\matc + \matz^*) - \hat{f}(\matc) < b_i\tau$}{\textbf{break}}
      Add item $i$ to inventory and update allocation: $\matc = \matc + \matz^*$, $c_i = c_i + 1$\;
    }
    }
    \Return{$\vecc,\matc$}\;
  \end{algorithm}}
\vspace{-10pt}
\begingroup
{\SetAlgoNoLine%
  \begin{algorithm}[htbp]
    \KwIn{Budget $B$, weights $(b_i)_{i \in \NN }$, error $\epsilon \in (0,1)$, number of customers $T$, products $\NN$ sorted in descending order of price: $r_1\geq r_2\geq \cdots \geq r_n$}
    \KwOut{$\vecc$, $\matc$}
Exclude from $\NN$ all products $i$ with $b_i>B$\;
Initialize $\tau = \frac{1}{B}\max_{i \in \NN }\fLPM(\vece_i,T)$, $I_{nT}=\max\{\lceil \log_{1+\epsilon }(nT) \rceil,1\}$\;
       \For{$i \gets1$ \textbf{to} $I_{nT}$} {
      $\vecc^i, \matc^i= \text{Budget-Constrained Threshold Add for DAP}(B, (b_i)_{i \in \NN }, (1+\epsilon)^{i-1}\tau)$\;
    }
    \Return{$\vecc,\matc= \text{ Best of }(\vecc^1,\matc^1), (\vecc^2,\matc^2), \cdots, (\vecc^{I_{nT}}, \matc^{I_{nT}})$ and all singletons, evaluated by $\hat{f}(\cdot,T)$}
    \caption{Budget-Constrained Optimization for DAP}\label{alg:tau-m-bg}
  \end{algorithm}}
For each single-item solution $\vece_i$, $i\in\NN$, we solve $\fLPM(\vece_i,T)$ to obtain the optimal allocation of one unit of product $i$.
\endgroup
\begin{theorem}\label{thm:3-ap-m-bud-T}
Algorithm~\ref{alg:tau-m-bg} outputs a $(\frac{1}{3}-\epsilon)$-approximate solution to problem \eqref{OPT-matc} under a budget constraint.
\end{theorem}

\proof{Proof of Theorem~\ref{thm:3-ap-m-bud-T}.} 
Let $\vecc^{\operatorname{OPT}}$ denote the optimal inventory vector for $\max_{\sum_{i \in \NN }b_ic_i\leq B} \fLPM (\vecc,T)$, and let $\matc^{\operatorname{OPT}}$ denote its corresponding optimal allocation that achieves $\max_{\matc\in \mathbb{R}^{mn}_+,\sum_{k\in \MM}\matc_k =\vecc^{\operatorname{OPT}} } \hat{f}(\matc)$. Thus, $\fLPM (\vecc^{\operatorname{OPT}},T)=\hat{f}(\matc^{\operatorname{OPT}})$. For simplicity, let $\operatorname{OPT}$ denote this optimal objective value.

We first exclude any elements $j\in \NN$ with $b_j>B$ as they cannot be part of any feasible solution. Let $\tau_i= \tau(1+\epsilon)^{i-1}$ and let $\matc^i$ denote the allocation output by the Threshold Add algorithm for threshold $\tau_i$. Since our algorithm returns the best solution among all singletons and all $\matc^i$, it suffices to show that either there exists some threshold value $\tau_i$ such that $\hat f(\matc^i,T)\geq(\frac{1}{3}-\epsilon)\operatorname{OPT}$ or there exists $j\in \NN$ such that $\fLPM(\vece_j,T)\geq(\frac{1}{3}-\epsilon)\operatorname{OPT}$.

If $\max_{j\in \NN}\fLPM(\vece_j,T)\geq\frac{\operatorname{OPT}}{2}$, we immediately obtain a $\frac{1}{2}$-approximation through singletons. Thus, we focus on the case where $\max_{j\in \NN}\fLPM(\vece_j,T)=B\tau <\frac{\operatorname{OPT}}{2}$. Additionally, $n T\max_{j\in \NN}\fLPM(\vece_j,T)\geq\operatorname{OPT}$. \begingroup Indeed, summing the type constraints of \eqref{eq:multi-sblp} over $k\in \MM$ shows that the total sales of each product are at most $\sum_{k\in \MM}\sum_{t=1}^{T}\lambda_{k,t}=T$, so replacing each $c_i^{\operatorname{OPT}}$ by $\min\{c_i^{\operatorname{OPT}},T\}$ does not change $\fLPM(\vecc^{\operatorname{OPT}},T)$; Lemma~\ref{lem:subadd} then gives $\operatorname{OPT}\leq \sum_{i\in \NN}\min\{c_i^{\operatorname{OPT}},T\}\,\fLPM(\vece_i,T)\leq nT\max_{j\in \NN}\fLPM(\vece_j,T)$. This implies $nT\tau  \geq\frac{\operatorname{OPT}}{B}$. Therefore, there exists an $i$ such that $(1-\epsilon)\frac{2\operatorname{OPT}}{3B}\leq \tau_i \leq \frac{2\operatorname{OPT}}{3B}$.
\endgroup
For this threshold $\tau_i$, we examine the total weight of selected items $\sum_{j \in \NN} b_j c_j^i$. If this sum is at least $\frac{1}{2}B$, then $\hat{f}(\matc^i) \geq\sum_{j \in \NN} b_j c_j^i \tau_i \geq\frac{1-\epsilon}{3}\operatorname{OPT}$. Therefore, we focus on the case where $\sum_{j \in \NN} b_j c_j^i < \frac{1}{2}B$.

Since \(\fLPM( \vecc^i,T)\geq \hat{f}( \matc^i) \), it suffices to show \(\hat{f}( \matc^i)\geq \frac{1}{3}\operatorname{OPT}\). Similar to the proof of Theorem~\ref{thm:3-ap-m-ra-T}, we consider $\matc^{\operatorname{OPT}}+\matc^i$. The main difference is that now there may exist a product, denoted $q^{\operatorname{E}}$, in $\matc^{\operatorname{OPT}}$ that satisfies the threshold requirement but cannot be added to $\matc^i(k)$ due to the budget constraint for some $k$. Given $\sum_{j\in \NN} b_j c_j^i< \frac{1}{2}B$, we must have $b_{q^{\operatorname{E}}}>\frac{1}{2}B$, implying at most one such item exists. Suppose $j_k\leq q^{\operatorname{E}}<j_{k+1}$ (with $j_{\kappa+1}=n+1$). Then $\matc^{\operatorname{OPT}}+\matc^i$ can be expressed as the following \emph{interleaved partition}:
\[
\matc^{\operatorname{OPT}}_{1},  \matc_{\cdot j_1}^i, \cdots, \matc^{\operatorname{OPT}}_{k},  \matc^i_{\cdot j_k}, \matc_{\cdot q^{\operatorname{E}}}^{\operatorname{OPT}},\cdots, \matc^{\operatorname{OPT}}_{ \kappa}, \matc^i_{\cdot j_{\kappa}}, \matc^{\operatorname{OPT}}_{ \kappa+1}.
\] Here the blocks $\matc^{\operatorname{OPT}}_{l}$ exclude the column of $q^{\operatorname{E}}$. Note that when $q^{\operatorname{E}}$ does not exist, this partition still holds with $\matc_{\cdot q^{\operatorname{E}}}^{\operatorname{OPT}}=\mathbf{0}$. We analyze two cases based on whether $\matc_{\cdot q^{\operatorname{E}}}^{\operatorname{OPT}}=\mathbf{0}$.

Case I: $\matc_{\cdot q^{\operatorname{E}}}^{\operatorname{OPT}}=\mathbf{0}$. Applying Lemma~\ref{lem:3-m-tau-bd}, stated and proved below, gives:
\begin{align*}
\operatorname{OPT}\leq \hat{f}(\matc^{\operatorname{OPT}}+\matc^i)&\leq \hat{f}(\matc^i(\kappa) )+\sum_{l\in[\kappa+1]}\hat{f}(\matc^{\operatorname{OPT}}_l|\matc^i(l-1)) \\
    &=  \hat{f}(\matc^i )+\sum_{l\in[\kappa+1]}\hat{f}(\matc^{\operatorname{OPT}}_l|\matc^i(l-1)) \\
    &\leq \hat{f}(\matc^i )+\sum_{l\in[\kappa+1]}\sum_{q\in N(\matc_l^{\operatorname{OPT}})}\|\matc^{\operatorname{OPT}}_{\cdot q}\|_1b_q \tau_i \\
  & \leq  \hat{f}(\matc^i)+\frac{2}{3}\operatorname{OPT}.
\end{align*}
Here, the first inequality results from the monotonicity of the LP, and the second inequality follows from Lemma~\ref{lem:interleave}. The equality holds as $\matc^i(\kappa)=\matc^i$. The third inequality comes from Lemma~\ref{lem:3-m-tau-bd}. Finally, the last inequality holds given that $\sum_{l\in [\kappa+1]}\sum_{q\in N(\matc^{\operatorname{OPT}}_l)}\|\matc^{\operatorname{OPT}}_{\cdot q}\|_1b_q \leq B$, and $\tau_i \leq \frac{2\operatorname{OPT}}{3B}$. Therefore, we have $\hat{f}(\matc^i)\geq(\frac{1}{3}-\epsilon)\operatorname{OPT}$.

\begingroup
Case II: $\matc_{\cdot q^{\operatorname{E}}}^{\operatorname{OPT}}\neq \mathbf{0}$. The product $q^{\operatorname{E}}$ lies in the support of $\vecc^{\operatorname{OPT}}$ with $j_k\leq q^{\operatorname{E}} < j_{k+1}$; its marginal benefit is at least $b_{q^{\operatorname{E}}}\tau_i$ but it cannot be added due to the budget constraint. Specifically, let $\mathbf{X}_{\cdot q^{\operatorname{E}}}$ denote a one-unit allocation of product $q^{\operatorname{E}}$ with the largest marginal benefit at $\matc^i(k)$. Then $\|\mathbf{X}_{\cdot q^{\operatorname{E}}}\|_1=1$, $\hat{f}(\mathbf{X}_{\cdot q^{\operatorname{E}}}|\matc^i(k))\geq \tau_ib_{q^{\operatorname{E}}}$ by the threshold requirement, and $\sum_{j\in \NN}b_j\|\matc^i_{\cdot j}(k)\|_1 +b_{q^{\operatorname{E}}}>B$. Moreover, each unit of product $j$ accepted by the Threshold Add algorithm has marginal benefit at least $b_j\tau_i$, so $\hat{f}(\matc^i(k))\geq \tau_i\sum_{j\in \NN}b_j\|\matc^i_{\cdot j}(k)\|_1$. \endgroup Thus, $\hat{f}(\matc^i(k)+\mathbf{X}_{\cdot q^{\operatorname{E}}})\geq \tau_i B\geq (1-\epsilon)\frac{2}{3}\operatorname{OPT}$.

Let $s = \argmax_{j \in \NN} \fLPM(\vece_j, T)$ be the index of the best singleton solution, and $\matc_{\cdot s}$ be the corresponding extended inventory vector. 
Let $\matc^{\operatorname{ALG}}$ denote the allocation output by Algorithm~\ref{alg:tau-m-bg}. Since the algorithm outputs the best among $\matc^i$ and all singletons, we have:
\begin{align*}
2\hat{f}(\matc^{\operatorname{ALG}})&\geq \hat{f}(\matc^i)+ \hat{f}(\matc_{\cdot s})\geq \hat{f}(\matc^i)+ \hat{f}(\mathbf{X}_{\cdot q^{\operatorname{E}}})\\ &\geq \hat{f}(\matc^i(k))+ \hat{f}(\mathbf{X}_{\cdot q^{\operatorname{E}}}) \geq \hat{f}(\matc^i(k)+\mathbf{X}_{\cdot q^{\operatorname{E}}})\geq (1-\epsilon)\frac{2}{3}\operatorname{OPT}.
\end{align*}
The first inequality holds because the algorithm outputs the best among singletons and $\matc^i$. The second inequality holds since $s$ is the best singleton solution. The third inequality follows because items are added in submodular order, and the marginal benefit remains positive each time an item is added. The fourth inequality is from the subadditivity property of $\hat{f}$ (Lemma~\ref{lem:subadd}). Finally, the last inequality is from our earlier result that \(\hat{f}(\matc^i(k)+\mathbf{X}_{\cdot q^{\operatorname{E}}})\geq (1-\epsilon)\frac{2}{3}\operatorname{OPT}\). Therefore, we conclude $\hat{f}(\matc^{\operatorname{ALG}})\geq (1-\epsilon)\frac{1}{3}\operatorname{OPT}$.
\Halmos
\endproof

\begin{lemma}\label{lem:3-m-tau-bd}
Use the notation from the proof of Theorem~\ref{thm:3-ap-m-bud-T}, with the horizon argument $T$ suppressed. Suppose that $\sum_{j\in\NN}b_jc_j^i<B/2$ and $\matc_{\cdot q^{\operatorname{E}}}^{\operatorname{OPT}}=\mathbf{0}$. For any $l\in [\kappa+1]$, $\hat{f}(\matc^{\operatorname{OPT}}_l|\matc^i(l-1)) \leq  \sum_{q\in N(\matc^{\operatorname{OPT}}_l)}\|\matc^{\operatorname{OPT}}_{\cdot q}\|_1b_q \tau_i$.
\end{lemma}
\proof{Proof of Lemma~\ref{lem:3-m-tau-bd}.}
The proof follows directly from the Budget-Constrained Threshold Add
algorithm. Fix any product
$q\in N(\matc_l^{\operatorname{OPT}})$ and let
\(
d_q
=
\|\matc_{\cdot q}^{\operatorname{OPT}}\|_1
=
c_q^{\operatorname{OPT}}
\in\mathbb Z_+.
\)
As in the proof of Lemma~\ref{lem:3-m-tau}, the allocation when the
algorithm finishes processing product $q$ is exactly
$\matc^i(l-1)$.

If $d_q\geq2$, then $b_q\leq B/2$ because the optimal solution is budget
feasible. Since the total weight selected by the algorithm is strictly
less than $B/2$, the budget permits one more unit of product $q$, and the
algorithm finishes processing $q$ only because its next one-unit marginal
is smaller than $b_q\tau_i$. If $d_q=1$, the same conclusion holds when
the budget permits another unit. If the budget does not permit another
unit and the one-unit marginal is at least $b_q\tau_i$, then $q$ would be
the exceptional product $q^{\operatorname{E}}$, contradicting the assumed absence of an exceptional product.
Therefore,
\(
\max_{\substack{\matz\in\mathbb R_+^{mn}:\
\sum_{k\in\MM}Z_{k,q}=1,\\
Z_{k,j}=0\ \forall k\in\MM,\ j\neq q}}
\hat f(\matz\mid\matc^i(l-1))
<
b_q\tau_i.
\)

Write
\(
\matc_{\cdot q}^{\operatorname{OPT}}
=
\sum_{s=1}^{d_q}\mathbf X^{q,s},
\mathbf X^{q,s}
=
\frac{1}{d_q}\matc_{\cdot q}^{\operatorname{OPT}},
\)
where each $\mathbf X^{q,s}$ is a one-unit allocation of product $q$.
By the strengthened fractional DR-SO property,
\[
\begin{aligned}
\hat f(
\matc_{\cdot q}^{\operatorname{OPT}}
\mid\matc^i(l-1))
=
\sum_{s=1}^{d_q}
\hat f\!\left(
\mathbf X^{q,s}
\mid
\matc^i(l-1)+\sum_{h<s}\mathbf X^{q,h}
\right)\leq
\sum_{s=1}^{d_q}
\hat f(\mathbf X^{q,s}\mid\matc^i(l-1))<
d_qb_q\tau_i.
\end{aligned}
\]

Finally, since the marginal value of a product is nonincreasing when the
base allocation is augmented by products preceding it in the submodular
order, ordering the products in
$N(\matc_l^{\operatorname{OPT}})$ increasingly and telescoping gives
\(
\hat f(
\matc_l^{\operatorname{OPT}}
\mid\matc^i(l-1))\leq
\sum_{q\in N(\matc_l^{\operatorname{OPT}})}
\hat f(
\matc_{\cdot q}^{\operatorname{OPT}}
\mid\matc^i(l-1))\leq
\sum_{q\in N(\matc_l^{\operatorname{OPT}})}
\|\matc_{\cdot q}^{\operatorname{OPT}}\|_1b_q\tau_i\).
\Halmos
\endproof

\begingroup
\subsubsection{Algorithms for Stochastic $T$.}\label{app:multi-stochT}
Recall that for \textsc{dap} with stochastic $T$, the first-stage fluid relaxation is
\begin{equation}\label{eq:lpm-Stoch-T}
\max_{\vecc\in \mathcal{F}}\EE_T[\fLPM (\vecc,T)].
\end{equation}
Given the total number of arrivals $T$, $\fLPM (\vecc, T)$ allocates inventory according to the expected number of customers of each type. Since the allocation varies with $T$, analyzing how changes in $\vecc$ affect the objective function is more complex. Nevertheless, $\EE_T[\hat{f} (\matc, T)]$ still has the DR-SO property because the property is preserved under addition. Our algorithm for stochastic $T$ follows the same high-level idea as in the deterministic case: given the threshold value $\tau$, we manually allocate fractional inventory across customer types while keeping the inventory vector integer-valued. The key difference is that we define and maintain an inventory allocation vector $\matc^T\in \RR_+^{mn}$ for each possible value of $T$, with $T\in [T_{\operatorname{max}}]$, where $T_{\operatorname{max}}$ is the maximum value in the support of $D_T$ (we will discuss the infinite support in Remark \ref{rem:stoch-T-comp}). If $D_T$ places positive mass at $T=0$, the corresponding realizations generate no customers and $\hat{f}(\matc,0)=0$ for every allocation $\matc$, so they contribute zero to all expectations over $T$; we therefore maintain allocation vectors only for $T\in[T_{\operatorname{max}}]$. For each product $i$, we calculate the optimal allocation for each $T$ and then compute the expected maximum marginal benefit across all values of $T$. Inside an expectation $\EE_T[\cdot]$, a superscript $T$ on an allocation vector refers to the realized value of $T$.

\textbf{Cardinality-Constrained Optimization Algorithm}

\begingroup
{\SetAlgoNoLine%
  \begin{algorithm}[htbp]
    \KwIn{Cardinality $K$, error $\epsilon \in (0,1)$, distribution $D_T$ of the number of customers, products $\NN$ sorted in descending order of price: $r_1\geq r_2\geq \cdots \geq r_n$}
    \KwOut{$\vecc$, $(\matc^T)_{T=1}^{T_{\operatorname{max}}}$}
Initialize $\tau = \frac{1}{K}\max_{i \in \NN }\EE_T[\fLPM(\vece_i,T)]$\;
       \For{$i \gets1$ \textbf{to} $\lceil \log_{1+\epsilon }K \rceil$} {
      $\vecc^i, (\matc^{i,T})_{T=1}^{T_{\operatorname{max}}}= \text{Cardinality-Constrained Threshold Add for DAP with Stochastic $T$}(K,(1+\epsilon)^{i-1}\tau)$\;
    }
  \Return{\( \vecc,(\matc^T)_{T=1}^{T_{\operatorname{max}}} \gets\text{Best of } (\vecc^1, (\matc^{1,T})_{T=1}^{T_{\operatorname{max}}}), \dots, (\vecc^{\lceil \log_{1+\epsilon }K \rceil}, (\matc^{\lceil \log_{1+\epsilon }K \rceil,T})_{T=1}^{T_{\operatorname{max}}}) \)}
    \caption{Cardinality-Constrained Optimization for DAP with Stochastic \texorpdfstring{$T$}{T}}\label{alg:tau-m-DT}
  \end{algorithm}}

{\SetAlgoNoLine%
  \begin{algorithm}[htp]
     \caption{Cardinality-Constrained Threshold Add for DAP with Stochastic \texorpdfstring{$T$}{T}}\label{alg:cdlp-multi-type-DT}
   \textbf{Input: } Cardinality $K$, threshold $\tau$, products $\NN$ sorted in descending order of price: $r_1\geq r_2\geq \cdots \geq r_n$, distribution $D_T$ of the number of customers\; \textbf{Output:} inventory vector $\vecc$, inventory allocation vectors $(\matc^T)_{T=1}^{T_{\operatorname{max}}} $\;  
Initialize $\vecc=\mathbf{0}^n$, $\matc^T=\mathbf{0}^{mn}$ for all $T\in [T_{\operatorname{max}}]$\;
\For{$i \in \NN $ }{
  \While{$\lVert\vecc\rVert_1<K$}{
    Find the optimal allocation of item $i$ for each $T$ at the current
    $\matc^T$: $\matz^{*,T} = \arg \max_{\matz \in \mathbb{R}^{mn}_+ :\sum_{k \in \MM}Z_{k,i} = 1,  \sum_{k \in \MM}Z_{k,j} = 0, \forall j\not= i} \hat{f}(\matc^T + \matz,T)- \hat{f}(\matc^T,T)$ \label{line:m-4-DT}\;
    \lIf{$\EE_T[\hat{f}(\matc^T+\matz^{*,T},T)-\hat{f}(\matc^T,T)] < \tau$}{\textbf{break}}
    Add item $i$ to the inventory and update:
    $\matc^T=\matc^T+\matz^{*,T} $ for all
    $T\in [T_{\operatorname{max}}]$, $c_i=c_i+1 $\label{line:m-5-DT}\;
  }
}
\Return{$\vecc,(\matc^T)_{T=1}^{T_{\operatorname{max}}}$}
  \end{algorithm}}
The set of products in $\matc^T$ remains consistent across different values of $T$, with only the allocation varying. Products are scanned once in decreasing-price order, all accepted copies of one product are contiguous, and the horizon-specific optimal allocations are recomputed after every accepted unit. Next, we establish the performance guarantee, which can be viewed as a stochastic counterpart to Theorem~\ref{thm:3-ap-m-ra-T}.
\endgroup

\begin{theorem}\label{thm:3-ap-m-ra-DT} 
Algorithm~\ref{alg:tau-m-DT} outputs a $(\frac{1}{2}-\epsilon)$-approximate solution to problem \eqref{eq:lpm-Stoch-T} under a cardinality constraint.
\end{theorem}
\proof{Proof of Theorem~\ref{thm:3-ap-m-ra-DT}.}
Let $\vecc^{\operatorname{OPT}}$ denote the optimal solution to \eqref{eq:lpm-Stoch-T}. For any $T\in [T_{\operatorname{max}}]$, let $\matc^{\operatorname{OPT},T}$ denote the optimal inventory allocation for $T$ customers, i.e., the solution to $\max_{\matc\in \RR_+^{mn},\sum_{k\in \MM}\matc_k =\vecc^{\operatorname{OPT}} } \hat{f}(\matc,T)$. Thus, $\EE_T[\fLPM (\vecc^{\operatorname{OPT}},T)]=\EE_T[\hat{f}(\matc^{\operatorname{OPT},T},T)]$. Let $\operatorname{OPT}$ denote this optimal objective value.

Let $\vecc^{\operatorname{ALG}}$ and $(\matc^{\operatorname{ALG},T})_{T=1}^{T_{\operatorname{max}}}$ denote the output of Algorithm~\ref{alg:tau-m-DT}. Since $\fLPM(\vecc^{\operatorname{ALG}},T)\geq \hat{f}(\matc^{\operatorname{ALG},T},T)$ for a given $T$, we have $\EE_T[\fLPM(\vecc^{\operatorname{ALG}},T)]\geq \EE_T[\hat{f}(\matc^{\operatorname{ALG},T},T)]$. We show $\EE_T[\hat{f}(\matc^{\operatorname{ALG},T},T)]\geq (\frac{1}{2}-\epsilon)\operatorname{OPT}$.

Let $\tau_i = \tau (1+\epsilon)^{i-1}$. Let $\vecc^i$ and $(\matc^{i,T})_{T=1}^{T_{\operatorname{max}}}$ denote the output of the Threshold Add algorithm given $\tau_i$. As the algorithm outputs the best solution among all $(\vecc^i,(\matc^{i,T})_{T=1}^{T_{\operatorname{max}}})$, it suffices to show that there exists some threshold value $\tau_i$ such that $\EE_T[\hat{f}(\matc^{i,T},T)]\geq (\frac{1}{2}-\epsilon)\operatorname{OPT}$.

If $ \max_{j\in \NN}\EE_T[\fLPM(\vece_j,T)]\geq \frac{\operatorname{OPT}}{2}$, then $\EE_T[\hat{f}(\matc^{\lceil \log_{1+\epsilon }K \rceil,T},T)]\geq \frac{1-\epsilon}{2}\operatorname{OPT}$. Thus, we only consider the case where $\operatorname{OPT}>2\max_{i \in \NN}\EE_T[\fLPM(\vece_i,T)]$. Moreover, $\operatorname{OPT}\leq K \max_{i \in \NN }\EE_T[\fLPM(\vece_i,T)]=K^2 \tau$ due to subadditivity (Lemma~\ref{lem:subadd}), so we have $K\tau\geq \frac{\operatorname{OPT}}{K}$.
Thus, there exists an $i$ such that $(1-\epsilon)\frac{\operatorname{OPT}}{2K}\leq \tau_i\leq \frac{\operatorname{OPT}}{2K}$.
Then we show $\EE_T[\hat{f}(\matc^{i,T},T)]\geq  \frac{1-\epsilon}{2}\EE_T[\hat{f}(\matc^{\operatorname{OPT},T},T)] $. Note that $\|\matc^{i,T}\|_1$ is the same for all $T$ since only the inventory allocation is different for different $T$, but the total inventory of each product $j$ is the same. We consider two cases based on whether $\|\matc^{i,1}\|_1=K$. If $\|\matc^{i,1}\|_1=K$, then $\EE_T[\hat{f}(\matc^{i,T},T)]\geq K\tau_i \geq \frac{1-\epsilon}{2}\operatorname{OPT}$. So we only need to consider the case where $ \|\matc^{i,1}\|_1 <K$. For simplicity, we omit the superscript and write $\vecc^i$ as $\vecc$, $(\matc^{i,T})_{T=1}^{T_{\operatorname{max}}}$ as $(\matc^{T})_{T=1}^{T_{\operatorname{max}}}$, and $\tau_i$ as $\tau$ from now on. 

Let $\kappa$ be the number of distinct products with positive inventory in $\vecc$, and let $\{j_1, \ldots, j_{\kappa}\}$ denote this set of products, with $j_1< \cdots < j_{\kappa}$. For each $T$, $\matc^{\operatorname{OPT},T}+\matc^T$ can be expressed as the following \emph{interleaved partition}:
\[
\matc^{\operatorname{OPT},T}_{1},  \matc^T_{\cdot j_1}, \matc^{\operatorname{OPT},T}_{2},  \matc^T_{\cdot j_2}, \cdots, \matc^{\operatorname{OPT},T}_{ \kappa}, \matc^T_{\cdot j_{\kappa}}, \matc^{\operatorname{OPT},T}_{ \kappa+1}. 
\]

Here, as in the deterministic case, $\matc^{\operatorname{OPT},T}_{l}$ represents the $l$-th partition of the optimal solution, and $\matc^T_{\cdot j_l}$ represents the allocation of product $j_l$ in our solution. We define the cumulative allocation up to partition $l$ as $\matc^T(l)= \sum_{l'=1}^{l} \matc^T_{\cdot j_{l'}}$.
 Due to the DR-SO property of $\hat{f}$, for any $T$, we have the following result from Lemma~\ref{lem:interleave}:
\begin{equation*}
    \hat{f}(\matc^{\operatorname{OPT},T}+\matc^T,T) \leq \hat{f}(\matc^T(\kappa),T)+\sum_{l\in [\kappa+1]} \left[\hat{f}(\matc^{\operatorname{OPT},T}_{l}+\matc^{T}(l-1),T) -\hat{f}( \matc^{T}(l-1),T)\right]
\end{equation*}

Taking the expectation over $T$, we obtain the following:
\begin{multline}\label{eq:cor3-DT}
\EE_T [\hat{f}(\matc^{\operatorname{OPT},T}+\matc^T,T) ]\leq \EE_T [\hat{f}(\matc^T(\kappa),T)]\\
+\EE_T\!\bigl[\sum_{l\in [\kappa+1]} \bigl(\hat{f}(\matc^{\operatorname{OPT},T}_{l}+\matc^{T}(l-1),T) -\hat{f}( \matc^{T}(l-1),T)\bigr)\bigr]
\end{multline}

To bound the second term in inequality \eqref{eq:cor3-DT}, we apply Lemma~\ref{lem:3-m-tau-stochT}, stated and proved below. This establishes the approximation ratio of $\frac{1}{2}-\epsilon$:

\begin{align*}
   \operatorname{OPT}& \leq   \EE_T[ \hat{f}(\matc^{\operatorname{OPT},T}+\matc^{T},T)]\\
    &\leq \EE_T[\hat{f}(\matc^{T}(\kappa) ,T)]+\sum_{l\in[\kappa+1]}\EE_T[\hat{f}(\matc^{\operatorname{OPT},T}_{l}+\matc^{T}(l-1),T) -\hat{f}( \matc^{T}(l-1),T)] \\
  & =\EE_T[\hat{f}(\matc^{T},T)]+\sum_{l\in[\kappa+1]}\EE_T[\hat{f}(\matc^{\operatorname{OPT},T}_{l}+\matc^{T}(l-1),T) -\hat{f}( \matc^{T}(l-1),T)]\\
  & \leq \EE_T[\hat{f}(\matc^{T},T)]+\sum_{l\in[\kappa+1]}\|\matc^{\operatorname{OPT},1}_l\|_1 \tau \\
  & =\EE_T[ \hat{f}(\matc^{T},T)]+\|\matc^{\operatorname{OPT},1}\|_1 \tau \\
  & \leq \EE_T[ \hat{f}(\matc^{T},T)]+\frac{1}{2}\EE_T[  \hat{f}(\matc^{\operatorname{OPT},T},T)].
\end{align*}
Here the first inequality is from the monotonicity of the LP for each $T$. The subsequent inequality is derived from inequality \eqref{eq:cor3-DT}. The first equality holds as $\matc^{T}(\kappa)=\matc^{T} $ for all $T$. The third inequality is from Lemma~\ref{lem:3-m-tau-stochT}, which is stated and proved below. Finally, the last inequality holds given that $\|\matc^{\operatorname{OPT},1}\|_1\leq K $, after substituting $\tau \leq \frac{\operatorname{OPT}}{2K}$. Thus, we prove $\EE_T[\hat{f}(\matc^{T},T)]\geq (\frac{1}{2}-\epsilon)\operatorname{OPT} $ and complete the proof.
\Halmos
\endproof

\begin{lemma}\label{lem:3-m-tau-stochT}
Use the notation from the proof of Theorem~\ref{thm:3-ap-m-ra-DT}, with the threshold index $i$ suppressed. Suppose that $\|\vecc\|_1<K$. For any $l\in [\kappa+1]$, $\EE_T[\hat{f}(\matc^{\operatorname{OPT},T}_l+\matc^{T}(l-1),T) -\hat{f}( \matc^{T}(l-1),T)]\leq  \|\matc^{\operatorname{OPT},1}_{l}\|_1\tau$.
\end{lemma}
\proof{Proof of Lemma~\ref{lem:3-m-tau-stochT}.}
Let $\Delta_q^T(\mathbf A)$ denote the maximum marginal value of a
one-unit allocation of product $q$ at $\mathbf A$ for horizon $T$, and let
$d_q=\|\matc_{\cdot q}^{\operatorname{OPT},1}\|_1$. Applying the
repeated-copy and product-order telescoping argument from the proof of
Lemma~\ref{lem:3-m-tau} pointwise for every $T$ gives
\[
\hat f(
\matc_l^{\operatorname{OPT},T}
\mid\matc^T(l-1),T)
\leq
\sum_{q\in N(\matc_l^{\operatorname{OPT},1})}
d_q\Delta_q^T(\matc^T(l-1)).
\]
Since the solution is underfilled, the stochastic Threshold Add algorithm
rejects the next unit of every such product at
$\matc^T(l-1)$, and hence
$\EE_T[\Delta_q^T(\matc^T(l-1))]<\tau$. Taking expectations in the
preceding inequality yields
\[
\EE_T\!\left[
\hat f(
\matc_l^{\operatorname{OPT},T}
\mid\matc^T(l-1),T)
\right]
\leq
\sum_{q\in N(\matc_l^{\operatorname{OPT},1})}d_q\tau
=
\|\matc_l^{\operatorname{OPT},1}\|_1\tau.
\]
\Halmos
\endproof
\begingroup
The runtime analysis follows the deterministic-\(T\) case, with the only additional cost coming from evaluating the expectation of the fluid objective over the distribution of \(T\). Let \(T_{\max}\) denote the largest value in the support of \(D_T\). Our algorithm maintains one allocation vector for each integer $T\in[T_{\max}]$, so we state it for distributions with finite support. The resulting runtime is $O\left(
\frac{n}{\epsilon}
T_{\max}(\log K)\log(mn)
(mn+\log K)
\right)$. Note that the algorithm is linear in $T_{\max}$. In fact, we can replace
$T_{\max}$ with the number of points in the support of $D_T$. If $D_T$ is
encoded such that the number of support points is not polynomial in the input
size, then the algorithm may be pseudopolynomial.
\endgroup
\begin{remark}\label{rem:stoch-T-comp}
 A distribution with unbounded support can be truncated first. Let $\bar{R}=\max_{\vecc\in\mathcal{F}}\sum_{i\in\NN}r_ic_i$ and fix $\epsilon'\in(0,1)$. Let $\hat{T}$ be the smallest integer such that $\Prb(T>\hat{T})\,\bar{R}\leq \epsilon' \max_{i:\vece_i\in\mathcal{F}}\EE_T[\fLPM(\vece_i,T)]$; if the right-hand side is zero, every feasible inventory yields zero fluid revenue and the problem is trivial. Set $T'=\min\{T,\hat{T}\}$. Since $\fLPM(\vecc,\cdot)$ is nondecreasing in the horizon and $\fLPM(\vecc,T)\leq\bar{R}$ for every $\vecc\in\mathcal{F}$, we have, for every $\vecc\in\mathcal{F}$,
\(
\EE_T[\fLPM(\vecc,T)]-\EE_T[\fLPM(\vecc,T')]\leq \Prb(T>\hat{T})\,\bar{R}\leq \epsilon'\max_{\vecc'\in\mathcal{F}}\EE_T[\fLPM(\vecc',T)].
\)
Hence running the algorithms with $T'$ in place of $T$ loses at most a factor of $1-\epsilon'$ in the first-stage guarantee, which is absorbed in $\epsilon$. The runtime statements then hold with $T_{\operatorname{max}}=\hat{T}$. The exact value of $\bar{R}$ is not needed: any upper bound on $\bar{R}$ serves the same purpose and only increases $\hat{T}$. 

\end{remark}
\endgroup

\textbf{Budget-Constrained Optimization Algorithm}

For the budget constraint, we modify the Threshold Add algorithm as in Algorithm~\ref{alg:tau-m-DT-budget}.

  {\SetAlgoNoLine%
  \begin{algorithm}[htbp]
    \KwIn{Budget $B$, positive weights $(b_i)_{i \in \NN }$, products $\NN$ sorted in descending order of price: $r_1\geq r_2\geq \cdots \geq r_n$, distribution $D_T$ of the number of customers, error $\epsilon \in (0,1)$}
    \KwOut{$\vecc$, $(\matc^T)_{T=1}^{T_{\operatorname{max}}}$}
Exclude from $\NN$ all products $i$ with $b_i>B$\;
\lIf{$\NN=\emptyset$}{\Return{$\vecc=\mathbf{0}^n$, $\matc^T=\mathbf{0}^{mn}$ for all $T\in[T_{\operatorname{max}}]$}}
Initialize $\tau = \frac{1}{B}\max_{i \in \NN }\EE_T[\fLPM(\vece_i,T)]$, $I_{nT_{\operatorname{max}}}=\max\{\lceil \log_{1+\epsilon }(nT_{\operatorname{max}}) \rceil,1\}$\;
       \For{$i \gets1$ \textbf{to} $I_{nT_{\operatorname{max}}}$} {
      $\vecc^i, (\matc^{i,T})_{T=1}^{T_{\operatorname{max}}}= \text{Budget-Constrained Threshold Add for DAP with Stochastic $T$}(B, (b_i)_{i \in \NN}, (1+\epsilon)^{i-1}\tau)$\;
    }
    \Return{$\vecc,(\matc^T)_{T=1}^{T_{\operatorname{max}}}$=Best of $(\vecc^1,(\matc^{1,T})_{T=1}^{T_{\operatorname{max}}})$, $(\vecc^2,(\matc^{2,T})_{T=1}^{T_{\operatorname{max}}})$, $\cdots$, $(\vecc^{I_{nT_{\operatorname{max}}}},(\matc^{I_{nT_{\operatorname{max}}},T})_{T=1}^{T_{\operatorname{max}}})$ and all singletons, evaluated by $\EE_T[\hat{f}(\cdot,T)]$}
    \caption{Budget-Constrained Optimization for DAP with Stochastic \texorpdfstring{$T$}{T} }\label{alg:tau-m-DT-budget}
  \end{algorithm}}

 {\SetAlgoNoLine%
  \begin{algorithm}[htp]
     \caption{Budget-Constrained Threshold Add for DAP with Stochastic \texorpdfstring{$T$}{T}}\label{alg:cdlp-multi-type-DT-budget}
   \textbf{Input:} Budget $B$, positive weights $(b_i)_{i \in \NN }$, threshold $\tau$, products $\NN$ sorted in descending order of price: $r_1\geq r_2\geq \cdots \geq r_n$, distribution $D_T$ of the number of customers\; \textbf{Output:} $\vecc,(\matc^T)_{T=1}^{T_{\operatorname{max}}} $\;  
Initialize $\vecc=\mathbf{0}^n$, $\matc^T=\mathbf{0}^{mn}$ for all $T\in [T_{\operatorname{max}}]$\;
       \For{$i \in \NN $ }{
         \While{$\sum_{j \in \NN } b_j c_j + b_i \leq B$}{
         Find the optimal allocation of item $i$ for each $T$:\;
         $\matz^{*,T} = \arg \max_{\matz \in \mathbb{R}^{mn}_+ :\sum_{k \in \MM}Z_{k,i} = 1,  \sum_{k \in \MM}Z_{k,j} = 0, \forall j\not= i} \hat{f}(\matc^T + \matz,T)-\hat{f}(\matc^T,T)$\;
         \lIf{$\EE_T[\hat{f}(\matc^T + \matz^{*,T},T) - \hat{f}(\matc^T,T)] < b_i\tau$}{\textbf{break}}
         Add item $i$ to inventory and update allocation:
 $\matc^T=\matc^T+\matz^{*,T} $ for all $T\in [T_{\operatorname{max}}]$, $c_i=c_i+1 $\;
    }
    }
    \Return{$\vecc,(\matc^T)_{T=1}^{T_{\operatorname{max}}}$}\;
  \end{algorithm}}

  For each singleton $\vece_i$, $i\in\NN$, we calculate $\EE_T[\fLPM(\vece_i,T)]$. As discussed for deterministic $T$, computing $\fLPM(\vece_i,T)$ requires $O(nm\log nm)$ time (see Appendix~\ref{app:acc}).   \begin{theorem}\label{thm:3-ap-m-bud-stoch-T}
Algorithm~\ref{alg:tau-m-DT-budget} outputs a $(\frac{1}{3}-\epsilon)$-approximate solution to problem \eqref{eq:lpm-Stoch-T} under a budget constraint.
\end{theorem}
\proof{Proof of Theorem~\ref{thm:3-ap-m-bud-stoch-T}.}
Let $\vecc^{\operatorname{OPT}}$ denote the optimal inventory for $\max_{\sum_{i \in \NN }b_ic_i\leq B} \EE_T[\fLPM(\vecc,T)]$. For any $T\in [T_{\operatorname{max}}]$, let $\matc^{\operatorname{OPT},T}$ denote the optimal solution to $\max_{\matc\in \RR_+^{mn},\sum_{k\in \MM}\matc_k =\vecc^{\operatorname{OPT}} } \hat{f}(\matc,T)$. Thus, $\EE_T[\fLPM (\vecc^{\operatorname{OPT}},T)]=\EE_T[\hat{f}(\matc^{\operatorname{OPT},T},T)]$. Let $\operatorname{OPT}$ denote the optimal objective value for simplicity.
\begingroup
We ignore elements $j\in \NN$ such that $b_j>B$. If every product is excluded, the zero inventory is the only feasible solution and the guarantee holds trivially. Let $\tau_i= \tau(1+\epsilon)^{i-1}$.
\endgroup
 Let $\vecc^i$ and $(\matc^{i,T})_{T=1}^{T_{\operatorname{max}}}$ denote the output of the Threshold Add algorithm given $\tau_i$.
 Since the algorithm outputs the best inventory vector among all singletons and $(\matc^{i,T})_{T=1}^{T_{\operatorname{max}}}$, it suffices to show that there exists some threshold value $\tau_i$ such that $\EE_T[\hat{f}( \matc^{i,T},T)]\geq (\frac{1}{3}-\epsilon)\operatorname{OPT}$ or there exists $j\in \NN$ such that $\EE_T[\fLPM(\vece_j,T)]\geq (\frac{1}{3}-\epsilon)\operatorname{OPT}$.

\begingroup
 If $\frac{\operatorname{OPT}}{2}\leq \max_{j\in \NN}\EE_T[\fLPM(\vece_j,T)]$, then the final solution will give at least $\frac{1}{2}\operatorname{OPT}$ since we consider all singletons. So we need only consider the case where $\max_{i \in \NN }\EE_T[\fLPM(\vece_i,T)]=B\tau <\frac{\operatorname{OPT}}{2}$. We have $nT_{\operatorname{max}}\max_{i \in \NN }\EE_T[\fLPM(\vece_i,T)]\geq \operatorname{OPT}$: for each realization $T$, the truncation argument in the proof of Theorem~\ref{thm:3-ap-m-bud-T} and Lemma~\ref{lem:subadd} give $\fLPM(\vecc^{\operatorname{OPT}},T)\leq \sum_{i\in \NN}\min\{c_i^{\operatorname{OPT}},T\}\,\fLPM(\vece_i,T)\leq T_{\operatorname{max}}\sum_{i\in \NN}\fLPM(\vece_i,T)$, and taking expectations over $T$ and bounding the sum by $n\max_{i \in \NN }\EE_T[\fLPM(\vece_i,T)]$ yields the bound. This implies $nT_{\operatorname{max}}\tau \geq \frac{\operatorname{OPT}}{B}$. Thus, there exists an $i$ such that $(1-\epsilon)\frac{2\operatorname{OPT}}{3B}\leq \tau_i \leq \frac{2\operatorname{OPT}}{3B}$. We consider two cases based on $\sum_{j\in \NN} b_j c_j^i$. If $\sum_{j\in \NN} b_j c^i_j \geq \frac{1}{2}B$, we have $\EE_T[\hat{f}(\matc^{i,T},T)]\geq \frac{1}{2}B\tau_i\geq \frac{1-\epsilon}{3}\operatorname{OPT}$. Thus, we need only examine the case where $\sum_{j\in \NN} b_j c_j^i< \frac{1}{2}B$.
\endgroup

As in Theorem~\ref{thm:3-ap-m-ra-DT}, let $\kappa$ be the number of distinct products with positive inventory in $\vecc^i$, and let $\{j_1, \ldots, j_{\kappa}\}$ denote this set of products, with $j_1< \cdots < j_{\kappa}$. The main difference is that now there may exist a product $q^{\operatorname{E}}$ in $\matc^{\operatorname{OPT},T}$ that satisfies the threshold requirement but cannot be added to $\matc^{i,T}(k)$ due to the budget constraint for some $k$. Since $\sum_{j\in \NN} b_j c_j^i< \frac{1}{2}B$, we must have $b_{q^{\operatorname{E}}}>\frac{1}{2}B$, implying at most one such item exists. Suppose $j_k\leq q^{\operatorname{E}}<j_{k+1}$ (with $j_{\kappa+1}=n+1$). For each $T$, $\matc^{\operatorname{OPT},T}+\matc^{i,T}$ can be expressed as the following \emph{interleaved partition}:
\[
\matc^{\operatorname{OPT},T}_{1},  \matc_{\cdot j_1}^{i,T}, \cdots, \matc^{\operatorname{OPT},T}_{k},  \matc_{\cdot j_k}^{i,T}, \matc^{\operatorname{OPT},T}_{\cdot q^{\operatorname{E}}},\cdots, \matc^{\operatorname{OPT},T}_{ \kappa}, \matc^{i,T}_{\cdot j_{\kappa}}, \matc^{\operatorname{OPT},T}_{ \kappa+1}.
\]

Here the blocks $\matc^{\operatorname{OPT},T}_{l}$ exclude the column of $q^{\operatorname{E}}$. Note that when $q^{\operatorname{E}}$ does not exist, this partition still holds with $\matc^{\operatorname{OPT},T}_{\cdot q^{\operatorname{E}}}=\mathbf{0}$. We consider two subcases based on whether $\matc^{\operatorname{OPT},T}_{\cdot q^{\operatorname{E}}}=\mathbf{0}$. It is important to note that in each set above, while the inventory allocation may differ, the total inventory for each product remains constant across all $T$. 

(i) $\matc^{\operatorname{OPT},T}_{\cdot q^{\operatorname{E}}}=\mathbf{0}$. Applying Lemma~\ref{lem:3-m-tau-bd-2}, stated and proved below, gives:
\begin{align*}
  \operatorname{OPT}&\leq  \EE_T[ \hat{f}(\matc^{\operatorname{OPT},T}+\matc^{i,T},T)]\\
  & \leq \EE_T[\hat{f}(\matc^{i,T}(\kappa),T)]+\sum_{l\in[\kappa+1]}\EE_T[\hat{f}(\matc^{\operatorname{OPT},T}_l+\matc^{i,T}(l-1),T) -\hat{f}( \matc^{i,T}(l-1),T)] \\
    &= \EE_T[\hat{f}(\matc^{i,T},T)]+\sum_{l\in[\kappa+1]}\EE_T[\hat{f}(\matc^{\operatorname{OPT},T}_l+\matc^{i,T}(l-1),T) -\hat{f}( \matc^{i,T}(l-1),T)]  \\
    &\leq \EE_T[\hat{f}(\matc^{i,T},T)]+\sum_{l\in[\kappa+1]}\sum_{q\in N(\matc^{\operatorname{OPT},1}_l)}\|\matc^{\operatorname{OPT},1}_{\cdot q}\|_1 b_q \tau_i  \\
  & \leq  \EE_T[\hat{f}(\matc^{i,T},T)]+\frac{2}{3}\operatorname{OPT}.
\end{align*}

The first inequality results from the monotonicity of the LP. The subsequent inequality follows from Lemma~\ref{lem:interleave}, applied for each $T$ before taking expectations. The first equality holds as $\matc^{i,T}(\kappa)=\matc^{i,T}$ for all $T$. The third inequality comes from Lemma~\ref{lem:3-m-tau-bd-2}. Finally, the last inequality holds given that $\sum_{q\in\NN}\|\matc^{\operatorname{OPT},1}_{\cdot q}\|_1 b_q \leq B $, after substituting $\tau_i \leq \frac{2\operatorname{OPT}}{3B}$. Thus, $\EE_T[\hat{f}(\matc^{i,T},T)]\geq \frac{1}{3}\operatorname{OPT}$.

\begingroup
(ii) $\matc^{\operatorname{OPT},T}_{\cdot q^{\operatorname{E}}}\neq \mathbf{0}$. For each $T$, let $\mathbf{X}^T_{\cdot q^{\operatorname{E}}}$ denote a one-unit allocation of product $q^{\operatorname{E}}$ with the largest marginal benefit at $\matc^{i,T}(k)$. Then $\|\mathbf{X}_{\cdot q^{\operatorname{E}}}^T\|_1=1$ for all $T$, $\EE_T[\hat{f}(\mathbf{X}^T_{\cdot q^{\operatorname{E}}}+\matc^{i,T}(k),T) -\hat{f}( \matc^{i,T}(k),T)] \geq \tau_i b_{q^{\operatorname{E}}}$ by the threshold requirement, and $\sum_{j \in \NN }b_j\|\matc_{\cdot j}^{i,1}(k)\|_1 +b_{q^{\operatorname{E}}}>B$. Moreover, each unit of product $j$ accepted by the Threshold Add algorithm has expected marginal benefit at least $b_j\tau_i$, so $\EE_T[\hat{f}(\matc^{i,T}(k),T)]\geq \tau_i\sum_{j\in \NN}b_j\|\matc^{i,1}_{\cdot j}(k)\|_1$. Thus, $\EE_T[\hat{f}(\matc^{i,T}(k)+\mathbf{X}_{\cdot q^{\operatorname{E}}}^T,T)]\geq \tau_i B\geq (1-\epsilon)\frac{2}{3}\operatorname{OPT}$. According to the Threshold Add algorithm, $\EE_T[\hat{f}(\matc^{i,T},T)]\geq \EE_T[\hat{f}(\matc^{i,T}(k),T)]$. Let $\vecc^{\operatorname{ALG}}$ and $(\matc^{\operatorname{ALG},T})_{T=1}^{T_{\operatorname{max}}}$ denote the output of the algorithm. Since Algorithm~\ref{alg:tau-m-DT-budget} outputs the best of $(\matc^{i,T})_{T=1}^{T_{\operatorname{max}}}$ and all singletons, we have the following:
\begin{align*}
2\EE_T[\hat{f}(\matc^{\operatorname{ALG},T},T)]&\geq \EE_T[\hat{f}(\matc^{i,T},T)]+ \EE_T[\hat{f}(\mathbf{X}_{\cdot q^{\operatorname{E}}}^T,T)]\\
&\geq \EE_T[\hat{f}(\matc^{i,T}(k),T)]+ \EE_T[\hat{f}(\mathbf{X}_{\cdot q^{\operatorname{E}}}^T,T)]\\
&\geq \EE_T[\hat{f}(\matc^{i,T}(k)+\mathbf{X}^T_{\cdot q^{\operatorname{E}}},T)]\\
&\geq (1-\epsilon)\frac{2}{3}\operatorname{OPT}.
\end{align*}

Here, the first inequality holds because the algorithm outputs the best among singletons and $(\matc^{i,T})_{T=1}^{T_{\operatorname{max}}}$, and the best singleton solution is at least as good as $(\mathbf{X}_{\cdot q^{\operatorname{E}}}^T)_{T=1}^{T_{\operatorname{max}}}$. The second inequality holds since $\EE_T[\hat{f}(\matc^{i,T},T)]\geq \EE_T[\hat{f}(\matc^{i,T}(k),T)]$, as shown above. The third inequality follows from the subadditivity of $\hat{f}$ (Lemma~\ref{lem:subadd}). Finally, the last inequality is from our earlier result that $\EE_T[\hat{f}(\matc^{i,T}(k)+\mathbf{X}^T_{\cdot q^{\operatorname{E}}},T)]\geq (1-\epsilon)\frac{2}{3}\operatorname{OPT}$. Therefore, we conclude that $\EE_T[\hat{f}(\matc^{\operatorname{ALG},T},T)]\geq (\frac{1}{3}-\epsilon)\operatorname{OPT}$.
\Halmos 
\endgroup
\endproof

\begin{lemma}\label{lem:3-m-tau-bd-2}
Use the notation from the proof of Theorem~\ref{thm:3-ap-m-bud-stoch-T}. Suppose that $\sum_{j\in\NN}b_jc_j^i<B/2$ and $\matc^{\operatorname{OPT},T}_{\cdot q^{\operatorname{E}}}=\mathbf{0}$ for every $T$. For any $l\in [\kappa+1]$, $\EE_T[\hat{f}(\matc^{\operatorname{OPT},T}_{l}+\matc^{i,T}(l-1),T) -\hat{f}( \matc^{i,T}(l-1),T)] \leq \sum_{q \in N(\matc^{\operatorname{OPT},1}_{l}) } \|\matc^{\operatorname{OPT},1}_{\cdot q}\|_1b_q\tau_i$.
\end{lemma} 
\proof{Proof of Lemma~\ref{lem:3-m-tau-bd-2}.}
Let $\Delta_q^T(\mathbf A)$ denote the maximum marginal value of a
one-unit allocation of product $q$ at $\mathbf A$ for horizon $T$, and let
$d_q=\|\matc_{\cdot q}^{\operatorname{OPT},1}\|_1$. Applying the
repeated-copy and product-order telescoping argument from the proof of
Lemma~\ref{lem:3-m-tau-bd} pointwise for every $T$ gives
\(
\hat f(
\matc_l^{\operatorname{OPT},T}
\mid\matc^{i,T}(l-1),T)
\leq
\sum_{q\in N(\matc_l^{\operatorname{OPT},1})}
d_q\Delta_q^T(\matc^{i,T}(l-1)).
\)
Under these assumptions,
\(
\EE_T[
\Delta_q^T(\matc^{i,T}(l-1))
]
<b_q\tau_i.
\)
If $d_q\geq2$, the budget permits another unit and the inequality follows
from the threshold rejection; if $d_q=1$, its failure would make $q$ the
exceptional product, contrary to the assumptions. Therefore, taking expectations
in the preceding inequality gives
\(
\EE_T\!\left[
\hat f(
\matc_l^{\operatorname{OPT},T}
\mid\matc^{i,T}(l-1),T)
\right]
\leq
\sum_{q\in N(\matc_l^{\operatorname{OPT},1})}
d_qb_q\tau_i.
\)
\Halmos
\endproof

\section{Transformation and Online Assortment Policies: Lemmas and Proofs}\label{app:4}
This appendix contains the lemmas, proofs, and examples supporting the second step: tightness and gap examples, the transformation steps for \textsc{da}, and the analysis of the online policy for \textsc{dap}.

\subsection{\texorpdfstring{$\frac{1}{2}$}{\frac{1}{2}} Is Tight for IFR \texorpdfstring{$T$}{T}}\label{app:geo-12}
Consider a ground set $\NN=\{1\}$ with inventory $x_1=1$, and let $p\in(0,1)$ denote the purchase probability of this product. For any $T$, $\ZDA_1(p,T, x_1)=\EE[\min\{z_T,1\}]=\Prb(\min\{z_T,1\}\geq 1)=1-\Prb(z_T=0)$, where $z_T\sim \text{Binomial} (T,p)$. 
Consider $T\sim \text{Geometric}(p)$; then $\mu_T p=1$. Moreover, $\EE_T[\ZDA_1(p,T, x_1)]=1-\EE_T[\Prb(z_T=0)]$. Then we have:
\begin{align*}
\EE_T[\Prb(z_T=0)]&=\sum_{k=1}^{\infty}(1-p)^{k-1}p(1-p)^k =\frac{p}{1-p}\sum_{k=1}^{\infty}(1-p)^{2k}=\frac{p}{1-p}\cdot\frac{(1-p)^2}{1-(1-p)^2}\\ &=\frac{p(1-p)}{p(2-p)} =\frac{1-p}{2-p}.
\end{align*}

Therefore, $\frac{\EE_T[\ZDA_1(p,T, x_1)]}{\min\{\mu_Tp, 1\}}=1-\frac{1-p}{2-p}=\frac{1}{2-p}$. When $p\rightarrow 0$, this ratio goes to $\frac{1}{2}$, so the constant $\frac{1}{2}$ in Lemma~\ref{lem:4-uni-DT} cannot be improved.

\subsection{Example for Non-IFR \texorpdfstring{$T$}{T}} \label{app:exm-nonIFR}
\begin{example}\label{exm:nonIFR}
Consider a ground set $\NN=\{1\}$ with inventory $x_1=1$, attraction parameters $v_1=M$, $v_0=1$, and price $r_1=a$, where $M>2$ is an integer. Let $T=0$ with probability $1-\frac{1}{M}$, and $T=M$ with probability $\frac{1}{M}$. Then $\mu_T=1$.
The distribution of $T$ does not satisfy the IFR property as $\frac{\Prb(T=2)}{\Prb(T\geq 2)}=0 < \frac{\Prb(T=0)}{\Prb(T\geq 0)}=1-\frac{1}{M}$. We now analyze the revenue. $\fFP(\vecx,\mu_T)=\frac{aM}{M+1}$. $\EE_T[\fDA(\vecx,T)]= \frac{1}{M}\fDA(\vecx,M)=\frac{1}{M}(1-\frac{1}{(M+1)^{M}})a$, which approaches $0$ as $M\rightarrow \infty$. Therefore, $\frac{\EE_T[\fDA(\vecx,T)]}{\fFP(\vecx,\mu_T)}\rightarrow 0$.
\end{example}

\subsection{Proof of Lemma~\ref{lem:4-product}} \label{app:prf-thm-pro}

\medskip\noindent\textbf{Lemma~\ref{lem:4-product} (restated).} {\itshape
For every $\vecx\in \mathbb{Z}_+^{n}$ and $i \in \NN $, we have $\ZDA_i(\vecx,T)\geq (1-\frac{1}{e}) \ZFP_i(\vecx,T)$. 
}\medskip

\begingroup We first establish two auxiliary lemmas; we also make repeated use of the characterization of the fluid process in Lemma~\ref{lem:4-clock} (Appendix~\ref{app:prf-3-rev}). \endgroup{}

\begin{lemma}
Given any inventory $\vecx\in \mathbb{Z}_+^n$ and any number $T$ of customers, let $\ZDA_i(p,T,x_i)$ denote the expected amount of product $i$ consumed in the stochastic consumption process with $T$ customers, inventory $x_i$, and uniform purchase probability $p=\phi(i,N(\vecx))$ for each customer. Then we have $\ZDA_i(\vecx,T)\geq \ZDA_i(p,T,x_i)$.
\label{lem:4-touni}
\end{lemma}

\begingroup
\proof{Proof of Lemma~\ref{lem:4-touni}.}
Let $\vecx_t$ denote the remaining inventory at time $t$. Since inventory decreases over time as products are sold, the set of available products at time $t$ is $N(\vecx_t)\subseteq N(\vecx)$. Then we have $p\leq \phi(i,N(\vecx_t))$ for every $t$ at which product $i$ remains in stock. Therefore, at each time $t$ before product $i$ stocks out, the probability that customer $t$ purchases product $i$ is at least $p$. Consequently, the total expected consumption $\ZDA_i(\vecx,T)$ is greater than or equal to the expected consumption in a process where each customer purchases product $i$ independently with probability $p$, up to the available inventory. Formally, let $X_t$ indicate that customer $t$ purchases product $i$, so that $\Prb(X_t=1\mid \vecx_t)=\phi(i,N(\vecx_t))$ while product $i$ remains in stock. Let $Y_t$ indicate that product $i$ is purchased at time $t$ if the purchase probability is $p$. We can couple $X_t=\1\{U_t\leq \phi(i,N(\vecx_t))\}$ and $Y_t=\1\{U_t\leq p\}$ with independent uniform random variables $(U_t)_{t=1}^T$; then $X_t\geq Y_t$ for every $t$ at which product $i$ remains in stock, and once product $i$ stocks out the consumption already equals $x_i$. Hence, pathwise,
\begin{align*}\ZDA_i(\vecx,T) &=\EE\!\left[\min\left\{\sum_{t=1}^T X_t, x_i\right\}\right] \geq \EE\!\left[\min\left\{\sum_{t=1}^T Y_t, x_i\right\}\right]= \ZDA_i(p,T,x_i). \Halmos
\end{align*}
\endgroup
\endproof
\vspace{-5pt}
\begin{lemma}
Given a number $T$ of customers and an inventory vector $\vecx \in \mathbb{Z}_+^n$, for any product $i \in \NN $ and any set $S \subseteq \NN$ that does not include $i$ \begingroup and satisfies $\sum_{j \in S} x_j \leq T$\endgroup, it holds that $\ZDA_i(\vecx,T) \geq\ZDA_i(\vecx - \sum_{j \in S} x_j \vece_j, T - \sum_{j \in S} x_j)$.
\label{lem:4-sin-subs}
\end{lemma}
We prove this result by induction, relying on the IIA property of the MNL choice model. 
\proof{Proof of Lemma~\ref{lem:4-sin-subs}.} 
We first show that $\ZDA_i(\vecx,T)\geq \ZDA_i(\vecx-x_j\vece_j, T-x_j)$ for any $i,j\in \NN$ with $j\neq i$ and $x_j\leq T$ by using the same ideas as in the simpler case with only two products in Section~\ref{sec:da-transform-alg}.

Since product $j$ has an inventory of $x_j$, at most $x_j$ customers purchase product $j$ before its inventory is depleted in the stochastic process, so at least $T-x_j$ customers do not purchase product $j$. \begingroup We couple the two processes through the embedded subsequence of these customers. Run the process with inventory $\vecx$ and delete every customer who purchases product $j$; because the deleted customers consume only units of product $j$, the remaining inventory of every other product evolves identically along the retained subsequence. Consider a retained customer, and let $A$ denote the set of products other than $j$ that are available upon this customer's arrival. If product $j$ is out of stock, the customer purchases each product $i'\in A$ with probability $\phi(i',A)$. If product $j$ is still in stock, then conditional on not purchasing $j$, the IIA property of the MNL choice model implies that the customer purchases each product $i'\in A$ with probability $\frac{\phi(i',A\cup\{j\})}{1-\phi(j,A\cup \{j\})}=\phi(i',A)$. In both cases, conditional on the history of the retained subsequence, each retained customer purchases exactly as a customer facing availability set $A$ in the process with inventory $\vecx-x_j\vece_j$. The first $T-x_j$ retained customers therefore replicate the stochastic consumption process with inventory $\vecx-x_j\vece_j$ and $T-x_j$ customers, and the consumption of product $i$ in the original process is at least the consumption of product $i$ by these customers. \endgroup Consequently:
\[
\ZDA_i(\vecx,T)\geq \ZDA_i(\vecx-x_j\vece_j,T-x_j).
\]
Repeating this procedure for every product in $S$ ensures that at each step, removing the units of another product and the corresponding number of customers cannot increase the expected consumption of product $i$; the condition $\sum_{j\in S}x_j\leq T$ guarantees that the number of remaining customers stays nonnegative at every step. Since all products in $S$ can eventually be removed in this manner, we obtain $\ZDA_i(\vecx,T) \geq \ZDA_i(\vecx-\sum_{j\in S}x_j \vece_j, T-\sum_{j\in S} x_j) $.
\Halmos
\endproof

\proof{Proof of Lemma~\ref{lem:4-product}.}
First, consider the expected consumption $\ZFP_i(\vecx,T)$ with initial inventory $\vecx$ and $T$ customers, and let $\zeta=\ZFP_0(\vecx,T)$. \begingroup Let
\(
\operatorname{Out}(i)=\Big\{j\in \NN:\ \ZFP_j(\vecx,T)=x_j\ \text{ and }\ \frac{x_j}{v_j}<\frac{x_i}{v_i}\Big\}
\)
denote the set of products that stock out strictly before product $i$ during the process \textsc{fp}$(\vecx,T)$; by Lemma~\ref{lem:4-clock} these are exactly the products that stock out and have a smaller ratio than product $i$, which is every product that stocks out when product $i$ itself does not stock out by time $T$. Write $a=\sum_{j\in \operatorname{Out}(i)}x_j$ and $\vecx'=\vecx-\sum_{j\in \operatorname{Out}(i)}x_j\vece_j$. Every product $j\in \operatorname{Out}(i)$ is fully consumed in \textsc{fp}$(\vecx,T)$, so accounting for all customers gives $a+\zeta\leq T$; hence $a\leq T$, with strict inequality whenever $T>0$, because $\ZFP_0(\vecx,\cdot)$ is strictly increasing (Lemma~\ref{lem:4-clock}) and $\ZFP_0(\vecx,0)=0$, so that $\zeta>0$. Next, removing the products of $\operatorname{Out}(i)$, together with $a$ units of the horizon, leaves the consumption of every remaining product unchanged, that is, $\ZFP_h(\vecx',T-a)=\ZFP_h(\vecx,T)$ for every $h\notin \operatorname{Out}(i)$. This follows by applying the Separability Lemma (Lemma~\ref{lem:sep}) to the products of $\operatorname{Out}(i)$ one at a time: for $j\in \operatorname{Out}(i)$, apply the lemma with base inventory $\vecx-x_j\vece_j$ and increment $\delta_j=x_j$; the additional consumption of product $j$ is $\alpha_j=\ZFP_j(\vecx,T)-\ZFP_j(\vecx-x_j\vece_j,T)=x_j-0=x_j$, so part (i) of that lemma gives $\ZFP_h(\vecx,T)=\ZFP_h(\vecx-x_j\vece_j,T-x_j)$ for every $h\neq j$. In particular, every other product of $\operatorname{Out}(i)$ remains fully consumed in the reduced process, so repeating this argument removes the products of $\operatorname{Out}(i)$ one by one (each intermediate horizon is at least $T-a\geq 0$) and yields the claim; in particular:
\endgroup
\begin{equation}\label{eq:4-uni} 
\ZFP_i(\vecx,T)= \ZFP_i(\vecx-\sum_{j\in \operatorname{Out}(i)}x_j \vece_j, T-\sum_{j\in \operatorname{Out}(i)} x_j).
\end{equation}
 \begingroup Moreover, no product of $N(\vecx')$ stocks out strictly before product $i$ in \textsc{fp}$(\vecx',T-a)$. Indeed, since the consumption of every retained product is unchanged, accounting for all customers gives $\ZFP_0(\vecx',T-a)=(T-a)-\sum_{h\notin \operatorname{Out}(i)}\ZFP_h(\vecx,T)=T-\sum_{h\in \NN}\ZFP_h(\vecx,T)=\zeta$; the retained inventories are also unchanged, so by Lemma~\ref{lem:4-clock} a product $j\notin \operatorname{Out}(i)$ stocks out strictly before product $i$ there if and only if $x_j/v_j\leq \zeta$ and $x_j/v_j<x_i/v_i$, which is precisely the condition for $j\in \operatorname{Out}(i)$ in \textsc{fp}$(\vecx,T)$. The assortment therefore equals $N(\vecx')$ at every time before product $i$ stocks out, and the purchase probability of product $i$ remains constant at $p_i=\phi(i, N(\vecx-\sum_{j\in \operatorname{Out}(i)}x_j \vece_j))$ for as long as product $i$ is in stock. \endgroup{}
 Thus, the expected consumption in this truncated process is:
 \begin{equation} \label{eq:app-4-4}  \ZFP_i(\vecx-\sum_{j\in \operatorname{Out}(i)}x_j \vece_j, T-\sum_{j\in \operatorname{Out}(i)} x_j)= \min\{p_i(T-\sum_{j\in \operatorname{Out}(i)} x_j), x_i\}.
 \end{equation}

 Next, we consider $\ZDA_i(\vecx,T)$. 
Using Lemma~\ref{lem:4-sin-subs} \begingroup(as shown above, $\sum_{j\in \operatorname{Out}(i)}x_j=a\leq T$)\endgroup{}, we have:
\begin{align}
& \ZDA_i(\vecx,T) \geq \ZDA_i(\vecx-\sum_{j\in \operatorname{Out}(i)} x_j\vece_j,T-\sum_{j\in \operatorname{Out}(i)} x_j).
\label{eq:app-4-1}
\end{align}
Then, applying Lemma~\ref{lem:4-touni} to the right-hand side, we have:
\begin{align}
\ZDA_i(\vecx-\sum_{j\in \operatorname{Out}(i)} x_j\vece_j,T-\sum_{j\in \operatorname{Out}(i)} x_j)\geq \ZDA_i(p_i,T-\sum_{j\in \operatorname{Out}(i)} x_j,x_i),
\label{eq:app-4-2}
\end{align}
where $p_i = \phi(i, N(\vecx-\sum_{j\in \operatorname{Out}(i)} x_j\vece_j))$ and $\ZDA_i(p_i,T-\sum_{j\in \operatorname{Out}(i)} x_j,x_i)$ denotes the expected consumption of product $i$ when each customer purchases it independently with probability $p_i$.  

Combining inequalities \eqref{eq:app-4-1} and \eqref{eq:app-4-2}, we have:
\begin{align*}
\ZDA_i(\vecx,T)
&\geq \ZDA_i(p_i,T-\sum_{j\in \operatorname{Out}(i)} x_j,x_i)  \geq (1-\frac{1}{e}) \min\{p_i(T-\sum_{j\in \operatorname{Out}(i)} x_j), x_i\}\\
&= (1-\frac{1}{e})\ZFP_i(\vecx-\sum_{j\in \operatorname{Out}(i)}x_j \vece_j, T-\sum_{j\in \operatorname{Out}(i)}x_j)  = (1-\frac{1}{e}) \ZFP_i(\vecx,T).
\end{align*}
Here the second inequality holds due to the single-product result from Lemma~\ref{lem:4-uni}. The first equality is from equation \eqref{eq:app-4-4}. The last equality follows from equation \eqref{eq:4-uni}.
\Halmos
\endproof

\begingroup
\subsection{Example of Unbounded Stochasticity Gap}\label{app:4-1}

\begin{example}[Unbounded stochasticity gap for a ranking choice model]
\label{exm:fp-da-ranking-gap}
We show that, for a ranking choice model, there is no universal constant
\(\beta>0\) such that
\(
\fDA(\vecx,T)\geq \beta \fFP(\vecx,T)
\)
for all integral inventory vectors \(\vecx\) and deterministic horizons \(T\).

Fix a positive integer \(m\). There are \(m+1\) products: one high-revenue product \(h\)
with revenue \(r_h=1\), and \(m\) low-revenue products
\(\ell_1,\ldots,\ell_m\), each with revenue \(0\). The initial inventory is one unit
of every product:
\(
x_h=1, x_{\ell_j}=1,j=1,\ldots,m.
\)
The selling horizon is \(T=2m\). Let
\(\varepsilon=\frac{1}{2m}\).
Customers follow the following ranking choice model. With probability
\(1-\varepsilon\), a customer draws \(J\) uniformly from \(\{1,\ldots,m\}\) and has
ranking
\(
\ell_J \succ 0 \succ \) all other products. 
Thus this customer buys \(\ell_J\) if it is available and otherwise chooses the outside
option. With probability \(\varepsilon\), a customer draws a uniformly random
permutation \(\sigma\) of \(\{\ell_1,\ldots,\ell_m\}\) and has ranking
\(
\sigma_1 \succ \sigma_2 \succ \cdots \succ \sigma_m \succ h \succ 0.
\)
Thus this customer buys an available low-revenue product if any low-revenue product
remains, and buys \(h\) only after all low-revenue products have stocked out.

We first compute the revenue of the fluid process. As long as all low-revenue products
are available, each low-revenue product is purchased at rate
\(
\frac{1-\varepsilon}{m}+\frac{\varepsilon}{m}=\frac{1}{m}.
\)
Since each low-revenue product starts with one unit of inventory, all \(m\) low-revenue
products deplete simultaneously at time \(m\). From time \(m\) to time \(2m\), the
second class of customers, which arrives at rate \(\varepsilon\), purchases product
\(h\). Since \(x_h=1\) and \(\varepsilon m=\frac{1}{2}\), product \(h\) does not stock out in
the fluid process. Therefore,
\(
\fFP(\vecx,2m)=r_h\,\ZFP_h(\vecx,2m)=\varepsilon m=\frac{1}{2}.
\)

We now upper bound the expected revenue of the stochastic process. Let \(N_H\) be
the number of customers of the second kind, i.e., customers whose ranking places all
low-revenue products above \(h\). Then
\(
N_H\sim \mathrm{Bin}\left(2m,\frac{1}{2m}\right), \EE[N_H]=1.
\)
For each low-revenue product \(\ell_j\), let \(A_j\) be the number of customers of the
first kind whose favorite product is \(\ell_j\). Let
\(
Z=\sum_{j=1}^m \1\{A_j=0\}
\)
be the number of low-revenue products that are never requested by a first-kind
customer.

The high-revenue product \(h\) can be purchased only after all low-revenue products
have stocked out. If a low-revenue product \(\ell_j\) has \(A_j=0\), then the only way
it can stock out is if it is purchased by a second-kind customer. Since each second-kind
customer purchases at most one low-revenue product, the event that all low-revenue
products stock out implies
\(
Z\leq N_H.
\)
Consequently, because \(x_h=1\),
\(
\fDA(\vecx,2m)=r_h\,\ZDA_h(\vecx,2m)
\leq \Prb(Z\leq N_H).
\)

It remains to show that \(\Prb(Z\leq N_H)\to 0\). For each \(j\),
\(
\Prb(A_j=0)
=
\left(1-\frac{1-\varepsilon}{m}\right)^{2m}
=:q_m,
\)
where \(q_m\to e^{-2}\). Hence \(\EE[Z]=m q_m\). Moreover, for \(j\neq k\),
\(
\Prb(A_j=0,A_k=0)
=
\left(1-\frac{2(1-\varepsilon)}{m}\right)^{2m}
\leq
\left(1-\frac{1-\varepsilon}{m}\right)^{4m}
=q_m^2,
\)
so the indicators \(\1\{A_j=0\}\) have nonpositive pairwise covariance.
Therefore,
\(
\operatorname{Var}(Z)\leq m q_m.
\)
By Chebyshev's inequality,
\(
\Prb\left(Z\leq \frac{m q_m}{2}\right)
\leq
\frac{4}{m q_m}.
\)
By Markov's inequality,
\(
\Prb\left(N_H\geq \frac{m q_m}{2}\right)
\leq
\frac{2}{m q_m}.
\)
Thus,
\[
\Prb(Z\leq N_H)
\leq
\Prb\left(Z\leq \frac{m q_m}{2}\right)
+
\Prb\left(N_H\geq \frac{m q_m}{2}\right)
\leq
\frac{6}{m q_m}.
\]
Since \(q_m\to e^{-2}\), this upper bound converges to \(0\). Combining the above,
\[
\frac{\fDA(\vecx,2m)}{\fFP(\vecx,2m)}
\leq
\frac{6/(m q_m)}{\frac{1}{2}}
=
\frac{12}{m q_m}
\to 0.
\]
Therefore, the ratio between the expected stochastic revenue and the fluid-process
revenue can be arbitrarily small, even for a ranking choice model and integral
inventories. The construction is also robust to the low-revenue products carrying a
small positive revenue \(\delta>0\): this adds at most \(m\delta\) to the revenue of
either process, which is negligible for \(\delta=o(1/m)\).
\end{example}
\endgroup
\vspace{-5pt}
\subsection{Proof of Lemma~\ref{lem:4-per-gap}} \label{app:4-prf-T1}
\medskip\noindent\textbf{Lemma~\ref{lem:4-per-gap} (restated).} {\itshape
\begingroup
 Given any inventory vector $\vecc$ and any real horizon $T>0$, let $\CFP$ be the revenue-ordered inventory constructed from $\vecc$, and let $\vecc^{\operatorname{ALG}}\in \mathbb{Z}_+^n$ be the output of Algorithm~\ref{alg:trans-T}. Then $\fFP(\vecc^{\operatorname{ALG}},T) \geq (1-\frac{1}{T+1})\fFP(\CFP,T)$. 
\endgroup
}\medskip

\proof{Proof of Lemma~\ref{lem:4-per-gap}.}
Recall from Lemma~\ref{lem:3-rev} that $\fLP(\vecc,T)=\fFP(\CFP,T)$. Our goal is to prove $\max\{\fFP(\bar{\vecx}^{\vecc}, T), \fFP(\underline{\vecx}^{\vecc},T)\}\geq (1-\frac{1}{T+1})\fFP(\vecx^{\vecc},T)$.

Let $i$ denote the index of the ``bad'' product, where $x_i^{\vecc}$ is fractional.
Define $\eta=\max\{x^{\vecc}_i-\lfloor x^{\vecc}_i \rfloor,0 \}$ as the fractional amount of the ``bad'' product. We only need to consider the case $\eta>0$. From the proof of Lemma~\ref{lem:3-rev}, $\ZFP_i(\vecx^{\vecc},T)=x^{\vecc}_i$. Let $\gamma=\max\{\ZFP_i(\bar{\vecx}^{\vecc},T)-\ZFP_i(\underline{\vecx}^{\vecc},T),0 \}=\max\{\ZFP_i(\bar{\vecx}^{\vecc},T)-\lfloor x^{\vecc}_i \rfloor,0 \}$ denote the fractional amount of the ``bad'' product consumed in \textsc{fp}$(\bar{\vecx}^{\vecc},T)$. Since $\bar{\vecx}^{\vecc}$ only has more inventory of product $i$ compared to $\vecx^{\vecc}$, we have $\eta\leq \gamma \leq 1$. By Lemma~\ref{lem:sep}, we have $\fFP(\vecx^{\vecc},T)=\fFP(\underline{\vecx}^{\vecc},T-\eta)+\eta r_i$, $\fFP(\bar{\vecx}^{\vecc},T)=\fFP(\underline{\vecx}^{\vecc},T-\gamma)+\gamma r_i$.

\begingroup
If $\gamma=\eta$, these two identities give $\fFP(\bar{\vecx}^{\vecc},T)=\fFP(\underline{\vecx}^{\vecc},T-\eta)+\eta r_i=\fFP(\vecx^{\vecc},T)$, so the lemma holds with no loss. In the remainder of the proof we assume $\gamma>\eta$. Note also that $T-\gamma>0$: $\ZFP_0(\bar{\vecx}^{\vecc},\cdot)$ is strictly increasing with $\ZFP_0(\bar{\vecx}^{\vecc},0)=0$ (Lemma~\ref{lem:4-clock}), and accounting for all customers gives $\gamma\leq \ZFP_i(\bar{\vecx}^{\vecc},T)\leq T-\ZFP_0(\bar{\vecx}^{\vecc},T)<T$.
\endgroup

Let $a=\frac{\fFP(\underline{\vecx}^{\vecc},T-\gamma)}{T-\gamma}$ denote the average revenue in $\fFP(\underline{\vecx}^{\vecc},T-\gamma)$. Let $b=\frac{\fFP(\underline{\vecx}^{\vecc},T-\eta)-\fFP(\underline{\vecx}^{\vecc},T-\gamma)}{\gamma-\eta} $ denote the average revenue in the time horizon $[T-\gamma, T-\eta]$ of $\fFP(\underline{\vecx}^{\vecc},T-\eta)$, and let $d= \frac{\fFP(\underline{\vecx}^{\vecc},T)-\fFP(\underline{\vecx}^{\vecc},T-\eta)}{\eta} $ denote the average revenue in the time horizon $[T-\eta, T]$.
Then $\fFP(\vecx^{\vecc},T)=(T-\gamma)a+(\gamma-\eta)b+ \eta r_i$, $\fFP(\bar{\vecx}^{\vecc},T)=(T-\gamma)a+\gamma r_i $, $\fFP(\underline{\vecx}^{\vecc},T)= (T-\gamma)a+(\gamma-\eta)b+\eta d $. We visualize the revenue decomposition in Figure~\ref{fig:rev-decom} for ease of understanding. 
\begin{figure}[htbp]    \centering
\resizebox{0.5\linewidth}{!}{
\begin{tikzpicture}[x=1bp,y=1bp]
  \path[use as bounding box] (0,0) rectangle (703,319);
  \fill[white] (0,0) rectangle (703,319);

  \tikzset{
    process label/.style={
      anchor=east,
      inner sep=0,
      font=\fontsize{30bp}{34bp}\selectfont\bfseries\boldmath
    },
    timeline/.style={
      line width=2.8bp,
      line cap=round,
      -{Latex[length=14bp,width=12bp]}
    },
    tick/.style={line width=2.2bp,line cap=round},
    span/.style={
      line width=2bp,
      line cap=round,
      decorate,
      decoration={brace,amplitude=6bp}
    },
    span label/.style={
      inner sep=0,
      font=\fontsize{23bp}{26bp}\selectfont\bfseries\boldmath
    },
    time label/.style={
      anchor=north,
      inner sep=0,
      font=\fontsize{22bp}{25bp}\selectfont\bfseries\boldmath
    }
  }

  \def\xstart{184}
  \def\xgamma{438}
  \def\xeta{590}
  \def\xend{692}
  \def\xarrow{700}
  \def\ytop{258}
  \def\ymiddle{158}
  \def\ybottom{58}
  \def\ybracetop{266}
  \def\ybracemiddle{166}
  \def\ybracebottom{66}

  \node[process label] at (166,\ytop)
    {$f^{\operatorname{FP}}(\underline{\mathbf{x}}^{\mathbf{c}},T)$};
  \node[process label] at (166,\ymiddle)
    {$f^{\operatorname{FP}}(\mathbf{x}^{\mathbf{c}},T)$};
  \node[process label] at (166,\ybottom)
    {$f^{\operatorname{FP}}(\bar{\mathbf{x}}^{\mathbf{c}},T)$};

  \draw[timeline] (\xstart,\ytop) -- (\xarrow,\ytop);
  \draw[span] (\xstart,\ybracetop) -- (\xgamma,\ybracetop)
    node[midway,above=9bp,span label] {$a$};
  \draw[span] (\xgamma,\ybracetop) -- (\xeta,\ybracetop)
    node[midway,above=9bp,span label] {$b$};
  \draw[span] (\xeta,\ybracetop) -- (\xend,\ybracetop)
    node[midway,above=9bp,span label] {$d$};
  \draw[tick] (\xstart,\ytop) -- (\xstart,\ybracetop);
  \draw[tick] (\xgamma,\ybracetop) -- (\xgamma,248);
  \draw[tick] (\xeta,\ybracetop) -- (\xeta,248);
  \draw[tick] (\xend,\ytop) -- (\xend,\ybracetop);
  \node[time label] at (\xgamma,241) {$T-\gamma$};
  \node[time label] at (\xeta,241) {$T-\eta$};
  \node[time label] at (\xend,241) {$T$};

  \draw[timeline] (\xstart,\ymiddle) -- (\xarrow,\ymiddle);
  \draw[span] (\xstart,\ybracemiddle) -- (\xgamma,\ybracemiddle)
    node[midway,above=9bp,span label] {$a$};
  \draw[span] (\xgamma,\ybracemiddle) -- (\xeta,\ybracemiddle)
    node[midway,above=9bp,span label] {$b$};
  \draw[span] (\xeta,\ybracemiddle) -- (\xend,\ybracemiddle)
    node[midway,above=9bp,span label] {$ r_i$};
  \draw[tick] (\xstart,\ymiddle) -- (\xstart,\ybracemiddle);
  \draw[tick] (\xgamma,\ybracemiddle) -- (\xgamma,148);
  \draw[tick] (\xeta,\ybracemiddle) -- (\xeta,148);
  \draw[tick] (\xend,\ymiddle) -- (\xend,\ybracemiddle);
  \node[time label] at (\xgamma,141) {$T-\gamma$};
  \node[time label] at (\xeta,141) {$T-\eta$};
  \node[time label] at (\xend,141) {$T$};

  \draw[timeline] (\xstart,\ybottom) -- (\xarrow,\ybottom);
  \draw[span] (\xstart,\ybracebottom) -- (\xgamma,\ybracebottom)
    node[midway,above=9bp,span label] {$a$};
  \draw[span] (\xgamma,\ybracebottom) -- (\xend,\ybracebottom)
    node[midway,above=9bp,span label] {$r_i$};
  \draw[tick] (\xstart,\ybottom) -- (\xstart,\ybracebottom);
  \draw[tick] (\xgamma,\ybracebottom) -- (\xgamma,48);
  \draw[tick] (\xend,\ybottom) -- (\xend,\ybracebottom);
  \node[time label] at (\xgamma,41) {$T-\gamma$};
  \node[time label] at (\xend,41) {$T$};
\end{tikzpicture}}
\caption[Total revenue]{Total revenue in $\fFP(\underline{\vecx}^{\vecc},T)$, \(\fFP(\vecx^{\vecc},T)\), and \(\fFP(\bar{\vecx}^{\vecc},T)\).}
    \label{fig:rev-decom}
\end{figure}

\begingroup
Then for any given $T$, the worst-case ratio is bounded below by the optimal value of the following optimization problem, whose constraints are necessary conditions satisfied by every instance.
\begin{align}
& \inf_{a,b,d,\eta, \gamma,r_i}\max\{\frac{\fFP(\bar{\vecx}^{\vecc},T)}{\fFP(\vecx^{\vecc},T)},\frac{\fFP(\underline{\vecx}^{\vecc},T)}{\fFP(\vecx^{\vecc},T)} \} \label{eq:opt-bad}\\
 \text{s.t. }  & \frac{\fFP(\bar{\vecx}^{\vecc},T)}{\fFP(\vecx^{\vecc},T)}= \frac{(T-\gamma)a+\gamma r_i }{(T-\gamma)a+(\gamma-\eta) b+\eta r_i} \nonumber \\
& \quad \frac{\fFP(\underline{\vecx}^{\vecc},T)}{\fFP(\vecx^{\vecc},T)}= \frac{(T-\gamma)a+(\gamma-\eta)b+\eta d }{(T-\gamma)a+(\gamma-\eta) b+\eta r_i} \nonumber\\
&\quad  0 \leq \eta \leq \gamma \leq 1, \gamma \leq T, \quad 0\leq d\leq r_i \leq a , \quad 0\leq d\leq r_i \leq b. \nonumber
\end{align}
Here $d\geq 0$ because $\fFP(\underline{\vecx}^{\vecc},s)$ is non-decreasing in $s$. By Lemma~\ref{lem:3-rev}, $\fFP(\vecx^{\vecc},T)=\max_{\vecx\preceq \vecc} \fFP(\vecx,T)$; comparing $\fFP(\vecx^{\vecc},T)=\fFP(\underline{\vecx}^{\vecc},T-\eta)+\eta r_i$ with $\fFP(\underline{\vecx}^{\vecc},T)=\fFP(\underline{\vecx}^{\vecc},T-\eta)+\eta d$, this optimality gives $d\leq r_i$.

We next show $a\geq r_i$ and $b\geq r_i$. Since
$\vecc$ is integral and $\bar{\vecx}^{\vecc}\preceq\vecc$,
Lemma~\ref{lem:3-rev} implies
\(
0\leq
\fFP(\vecx^{\vecc},T)-\fFP(\bar{\vecx}^{\vecc},T)
=(\gamma-\eta)(b-r_i)\).
Because $\gamma>\eta$, we have $b\geq r_i$. Hence, at some time
$s\in[T-\gamma,T-\eta)$, the instantaneous revenue rate of the fluid
process starting from $\underline{\vecx}^{\vecc}$ is at least $r_i$.
Every product stocked in this process has price at least $r_i$.
Lemma~\ref{lem:seq-rev1}(ii) therefore implies that the revenue rate is
at least $r_i$ at every time in $[0,s]$, and thus throughout
$[0,T-\gamma]$. Averaging over this interval gives $a\geq r_i$. 
Then we discuss two cases based on which of $\fFP(\bar{\vecx}^{\vecc},T)$ and $\fFP(\underline{\vecx}^{\vecc},T)$ is larger.
\endgroup

\textbf{Case 1.} When $(\gamma-\eta)b+\eta d\leq \gamma r_i$, we have $\fFP(\underline{\vecx}^{\vecc},T)\leq \fFP(\bar{\vecx}^{\vecc},T)$, and the objective function becomes $\frac{ta+\gamma r_i}{ta+(\gamma-\eta)b+\eta r_i}$. Let $g(a,b,d, \eta, \gamma,r_i, t)=\frac{ta+\gamma r_i}{ta+(\gamma-\eta)b+\eta r_i}$. Then $\frac{\partial g}{\partial a}=\frac{t(\gamma-\eta)(b-r_i)}{(ta+(\gamma-\eta)b+\eta r_i)^2 }\geq0$. Thus, $\inf_a g(a,b,d,\eta, \gamma,r_i, t)=g(r_i,b,d,\eta, \gamma,r_i,t)$, and the optimization problem \eqref{eq:opt-bad} becomes: 
\begin{align}
   & \inf_{r_i,b,d,\gamma,\eta} \frac{(t+\gamma)r_i}{(\gamma-\eta)b+(t+\eta)r_i} \label{eq:C7-obj} \\
   \text{s.t. }\quad &0\leq d\leq r_i \leq b,\quad 0 \leq \eta \leq \gamma \leq 1,\quad (\gamma-\eta)b+\eta d\leq \gamma r_i. \nonumber
\end{align}

Let $h_1(t,r_i,b,d,\eta, \gamma)$ denote the objective function \eqref{eq:C7-obj}. Then $\frac{\partial h_1}{\partial b}=\frac{-(\gamma-\eta)(t+\gamma)r_i}{[(\gamma-\eta)b+(t+\eta)r_i]^2}\leq 0$. From the constraint \((\gamma-\eta)b+\eta d\leq \gamma r_i\), we have $b^*=\frac{\gamma r_i-\eta d}{\gamma-\eta}$ given $\gamma, r_i, \eta, d$. So the optimization problem becomes: 
\begin{align}
   & \inf_{r_i,d,\gamma,\eta} \frac{(t+\gamma)r_i}{(t+\eta+\gamma)r_i-\eta d} \label{eq:p1} \\
   &\text{s.t. }  0\leq d\leq r_i \leq \frac{\gamma r_i-\eta d}{\gamma-\eta},  0 \leq \eta \leq \gamma \leq 1. \nonumber
\end{align}
Since $d\leq r_i$ and $\eta\leq \gamma$, the first constraint reduces to $0\leq d\leq r_i$. Moreover, $\inf_{0\leq d\leq r_i} {\frac{(t+\gamma)r_i}{(t+\eta+\gamma)r_i-\eta d}}=\frac{(t+\gamma)r_i}{(t+\eta+\gamma)r_i}$. Since $r_i>0$, the objective reduces to $\frac{t+\gamma}{t+\eta+\gamma}$.
Plugging in $t=T-\gamma$ and using $\eta\leq 1$, we have $\frac{T}{T+\eta}\geq \frac{T}{T+1}=1-\frac{1}{T+1}$.

\textbf{Case 2.} When $(\gamma-\eta)b+\eta d\geq \gamma r_i$, $\fFP(\underline{\vecx}^{\vecc},T)\geq \fFP(\bar{\vecx}^{\vecc},T)$, and the objective function becomes $\frac{ta+(\gamma-\eta)b+\eta d}{ta+(\gamma-\eta)b+\eta r_i}$. Let $g(a,b,d, \eta, \gamma,r_i, t)$ denote the objective function. Then $\frac{\partial g}{\partial a}=\frac{t[ta+(\gamma-\eta)b+\eta r_i]-t[ta+(\gamma-\eta) b+\eta d]}{[ta+(\gamma-\eta)b+\eta r_i]^2}\geq 0
$. Thus, $\inf_a g(a,b,d, \eta, \gamma,r_i, t)=g(r_i,b,d, \eta, \gamma,r_i, t)$ and the optimization problem becomes
\begin{align*}
   &\quad \inf_{r_i,b,d,\gamma,\eta} \frac{tr_i+(\gamma-\eta)b+\eta d}{(\gamma-\eta)b+(t+\eta)r_i}\nonumber \\
   &\text{s.t. }  0\leq d\leq r_i \leq b,  0 \leq \eta \leq \gamma \leq 1, (\gamma-\eta)b+\eta d\geq \gamma r_i.
\end{align*}

Let $h_2(t,r_i,b,d,\eta, \gamma)$ denote the objective function. Then $\frac{\partial h_2}{\partial b}=\frac{(\gamma-\eta)[(\gamma-\eta)b+(t+\eta)r_i-tr_i-(\gamma-\eta)b-\eta d]}{[(\gamma-\eta)b+(t+\eta)r_i]^2}\geq 0$. From the last constraint, $b^*=\frac{\gamma r_i-\eta d}{\gamma-\eta}$ for given $\gamma, r_i, \eta, d$. \begingroup Substituting $b^*$ yields exactly the optimization problem \eqref{eq:p1},\endgroup{} so we obtain the same lower bound $1-\frac{1}{T+1}$ in both cases. \Halmos
\endproof

\subsection{Proof of Lemma~\ref{lem:4-Stoch-DT}}\label{app:IFRe-2}
\begingroup
\medskip\noindent\textbf{Lemma~\ref{lem:4-Stoch-DT} (restated).} {\itshape
Let $p\in (0,1)$ be a uniform purchase probability. If $T$ follows an IFR distribution with finite mean $\mu_T$, then for any integer $a\in \mathbb{Z}_+$ with $a<\mu_T$, $\EE_T [ (1-(1-p)^{(T-a)})\1\{T>a\} ] \geq \frac{1}{e}\min\{p(\mu_T-a),1\}$.
}

We first prove the lemma when $(\mu_T-a)p\leq 1$, in which case the right-hand side equals $\frac{1}{e}p(\mu_T-a)$ and the claim is equivalent to the lower bound of $\frac{1}{e}$ for the factor-revealing program \eqref{OPT:mu-c}; the complementary case $(\mu_T-a)p>1$ is treated by a monotonicity argument at the end of this subsection. We begin by reformulating optimization problem \eqref{OPT:mu-c}.
\endgroup

\subsubsection{Reformulation.}\label{app:IFRe-reform}
Given feasible $\mu_T,a$, \eqref{OPT:mu-c} can be reformulated as follows\begingroup, where we defer the division by $\mu_T-a$ to the end of the analysis\endgroup{}:
\begin{align}
  g(\mu_T,a) = & \quad \inf_{p,q_1,\cdots ,q_a} q_1q_2 \cdots q_a h(\mu_T,a,p,q_1,\cdots,q_a) \label{eq:c81-1} \\
   & \text{s.t. } (\mu_T-a)p\leq 1 , \quad  0< p\leq 1\nonumber\\
    &\quad \quad \mu_T = q_1 +q_1q_2+\cdots +q_1q_2\cdots q_a(1+\mu_{T'}) \nonumber, \quad \mu_{T'}\leq \frac{q_a}{1-q_a} \nonumber\\
    & \quad \quad 1\geq q_1\geq q_2\geq q_3\geq \cdots \geq q_a \geq 0\nonumber,
\end{align}
where 
\begin{align}
& h(\mu_T,a,p,q_1,\cdots,q_a) = \inf_{q_{a+1},q_{a+2},\cdots} \frac{ \EE_{T'}[1-(1-p)^{T'}]}{p} \nonumber \\
& \text{s.t.} \quad 1+q_{a+1}+q_{a+1}q_{a+2}+\cdots = (\mu_T-q_1-q_1q_2-\cdots-q_1\cdots q_{a-1})/(q_1\cdots q_a) \label{eq:2-mu} \\
& \quad \quad 1\geq q_a\geq q_{a+1} \geq q_{a+2}\geq \cdots \geq 0. \nonumber
\end{align}

We now analyze the optimal solution of $h(\mu_T,a,p,q_1,\cdots,q_a)$. Let $T'$ have the distribution of $T-a$ conditional on $T\geq a$; then $\Prb(T'\geq k)=\Prb(T-a\geq k|T\geq a)=\Prb(T\geq a+k|T\geq a)=q_{a+1}\cdots q_{a+k} $. The distribution of $T'$ also satisfies the IFR property, as $\Prb(T'\geq k+2|T'\geq k+1)=\frac{q_{a+1}\cdots q_{a+{k+2}} }{q_{a+1}\cdots q_{a+{k+1}}}=q_{a+k+2}\leq q_{a+k+1}=\Prb(T'\geq k+1|T'\geq k)$. Since $q_1,\cdots ,q_a$ are given, $\mu_{T'}=q_{a+1}+q_{a+1}q_{a+2}+q_{a+1}q_{a+2}q_{a+3}+\cdots $ is fixed from constraint \eqref{eq:2-mu}. \begingroup Moreover, since $f(T)=1-(1-p)^T$ is concave for $0<p<1$ ($f'(T)=-(1-p)^T \ln (1-p), f''(T)=-(1-p)^T \ln^2 (1-p)<0$; the boundary case $p=1$ follows by letting $p\uparrow 1$, as both sides of the comparison below are continuous in $p$), from the proposition and lemma below, $\EE_{T'\sim Geometric(\mu_{T'})}[f(T')]\leq \EE_{T'\sim D(\mu_{T'})} [f(T')]$ for all IFR distributions with the same mean $\mu_{T'}$. The cited lemma concerns distributions supported on $\{1,2,\cdots\}$, while $T'$ takes values in $\{0,1,2,\cdots\}$; we therefore apply it to $T'+1$, which is IFR with mean $\mu_{T'}+1$. Translation by one preserves both second-order stochastic dominance and the concavity of $f$, and the extremal geometric distribution on $\{1,2,\cdots\}$ with mean $\mu_{T'}+1$ corresponds, after shifting back, to the geometric distribution on $\{0,1,2,\cdots\}$ with $\Prb(T'\geq k)=b^k$ and $b=\frac{\mu_{T'}}{\mu_{T'}+1}$. So the optimal distribution of $T'$ is a geometric distribution with mean $\mu_{T'}$, i.e., the optimal solution of $h(\mu_T,a,p,q_1,\cdots ,q_a)$ is $q_{a+1}=q_{a+2}=\cdots =b$. Plugging in constraint \eqref{eq:2-mu}, we have $b=\frac{\mu_{T'}}{\mu_{T'}+1}$. Note that $b\leq q_a$ by the constraint $\mu_{T'}\leq \frac{q_a}{1-q_a}$, so this solution is feasible. \endgroup{}
\begin{proposition}[Theorem 1 in \cite{meyer1977second}]If distribution $A$ is second-order stochastically dominant over $B$, and $A$ and $B$ have the same mean, then for any concave function $f:\RR\rightarrow \RR$, $\EE_{x\sim A}[f(x)]\geq \EE_{x\sim B}[f(x)]$. 
\end{proposition}
\begin{lemma}[Lemma 3.4 in \cite{alijani2020predict}]
The geometric distribution with mean $\mu$ is second-order stochastically dominated by any other IFR distribution with the same mean.
\end{lemma}

Then we have:
\begin{align*}
& h(\mu_T,a,p,q_1,\cdots, q_a) = \frac{ \EE_{T'}[1-(1-p)^{T'}]}{p}  \\
 \text{s.t.}& \quad \frac{1}{1-b}= (\mu_T-q_1-q_1q_2-\cdots-q_1\cdots q_{a-1})/(q_1\cdots q_a). \\
\end{align*}
Furthermore, the objective function can be simplified further as $T'$ follows a geometric distribution. 
\begin{align*}
 h(\mu_T,a,p,q_1,\cdots, q_a) &= \frac{ \EE_{T'}[1-(1-p)^{T'}]}{p}    = \sum_{k=0}^{\infty}\Prb(T'=k)[1-(1-p)^k] /p                \\
&  = \sum_{k=0}^{\infty}b^k(1-b)[1-(1-p)^k]/p  = \frac{b}{1-b+bp}.
\end{align*}

Plugging in $h(\mu_T,a,p,q_1,\cdots,q_a)= \frac{b}{1-b+bp}$ to \eqref{eq:c81-1}, we have:
\begin{align}
  g(\mu_T,a) &=  \inf_{p,q_1,\cdots,q_a} q_1q_2 \cdots q_a  \frac{b}{1-b+bp}\nonumber \\
   \text{s.t.} &\quad (\mu_T-a)p\leq 1 \label{con:3-1} \\
    & \quad 0< p\leq 1 \label{con:3-2} \\
    &\quad \mu_T = q_1 +q_1q_2+\cdots +q_1q_2\cdots q_a \frac{1}{1-b} \label{con:3-mu}\\
    &\quad  1\geq q_1\geq q_2\geq q_3\geq \cdots \geq q_a \geq b\geq  0.\nonumber
\end{align}

Since $p$ appears only in the denominator of the objective function, and in constraint \eqref{con:3-1} and \eqref{con:3-2}, we have $p^*= \min\{1, \frac{1}{\mu_T-a}\}$. 

When $\mu_T-a\leq 1$, $p^*=1 $, and the optimization problem $g(\mu_T,a)$ reduces to:
\begin{align}
    & \quad \quad \inf_{q_1,\cdots ,q_a,b}  q_1q_2 \cdots q_ab \label{obj1} \\
    &\text{s.t.} \quad q_1 +q_1q_2+\cdots +q_1q_2\cdots q_a\frac{1}{1-b} =\mu_T \label{eq:con1} \\
    &\quad \quad  1\geq q_1\geq q_2\geq \cdots\geq q_a\geq  b \geq 0. \nonumber
\end{align}

When $\mu_T-a\geq 1$, $p^*=\frac{1}{\mu_T-a}$, and the optimization problem $g(\mu_T,a)$ becomes:
\begin{align}
    & \inf_{q_1,\cdots ,q_a,b}  q_1q_2 \cdots q_a\frac{b}{1-b+\frac{b}{\mu_T-a}} \label{obj2} \\
    &\text{s.t.} \quad q_1 +q_1q_2+\cdots +q_1q_2\cdots q_a\frac{1}{1-b} =\mu_T \nonumber \\
    &\quad \quad  1\geq q_1\geq q_2\geq\cdots\geq q_a\geq b \geq 0. \nonumber
\end{align}
The worst-case original ratio can be expressed as $\inf_{\mu_T,a:a\in \mathbb{Z}_+,\mu_T>a}\frac{g(\mu_T,a)}{\mu_T-a}$. We use $\vecq$ to denote the vector $(q_1, q_2,\cdots ,q_a,b)$ and $D_{\vecq}$ the corresponding probability distribution. Let $\operatorname{OBJ}_1(\vecq)$ and $\operatorname{OBJ}_2(\vecq)$ represent the objective functions in \eqref{obj1} and \eqref{obj2}, respectively.

\subsubsection{Proof of Lemma~\ref{lem:shiftgeo}.} \label{app:lem-4-shiftgeo}
\medskip\noindent\textbf{Lemma~\ref{lem:shiftgeo} (restated).} {\itshape
 Given any feasible $\mu_T$ and $a$ with respect to \eqref{OPT:mu-c}, \begingroup an optimal solution $\vecq^*$ to \eqref{OPT:mu-c} exists and has the following structure\endgroup{}: There exists a constant $b \in (0,1)$ such that either (i) there is a threshold index $1\leq k \leq a$, where $q_i = 1$ for all $i\leq k$ and $q_i = b$ for $ i\geq k+1$, or (ii) $q_i = b$ for all $i\geq 1$.
}\medskip

Note that a standard geometric distribution has $q_i= b$ for all $i$ and some constant $b$. For (i), the distribution of $T$ is characterized by $q_i = 1$ for all $i\leq k$ and $q_i = b$ for $i\geq k+1$; by defining $T'=T-k$, the distribution of $T'$ is a standard geometric distribution. Therefore, we refer to the distributions in the lemma as \emph{shifted geometric distributions}. To prove this lemma, we introduce the concept of a \emph{step}. \begingroup For analytical convenience, we define $q_0=1$ and $q_{a+1}=b$. A \emph{step} occurs at index $k$ (where $1\leq k \leq a+1$) if $q_{k-1}>q_k$. \endgroup{} Thus, a sequence $\vecq$ represents a shifted geometric distribution if and only if it contains exactly one step. Figure~\ref{fig:steps} illustrates examples of distributions with different numbers of steps. 

\begin{figure}[htp]
\caption{Examples of Distributions with Different Numbers of Steps. }
\centering
\begin{minipage}[b]{.31\textwidth}\centering
\includegraphics[width=.8\textwidth]{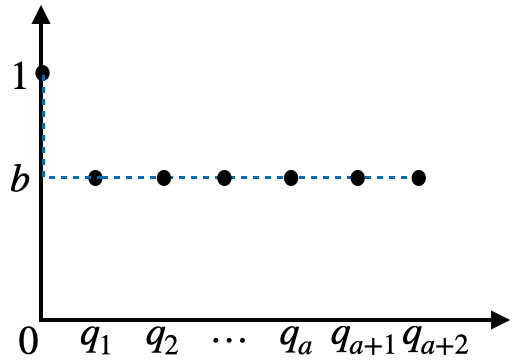}\\
{\small (a) 1 step (Geometric distribution)}
\end{minipage}\hfill
\begin{minipage}[b]{.31\textwidth}\centering
\includegraphics[width=.8\textwidth]{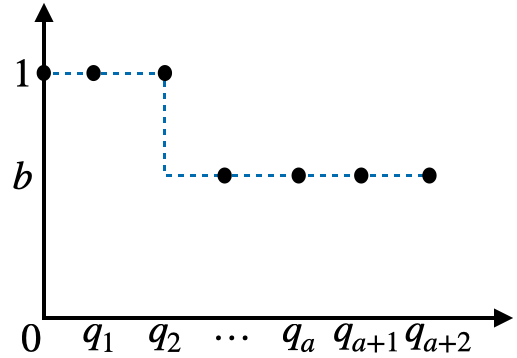}\\
{\small (b) 1 step (Shifted geometric distribution)}
\end{minipage}\hfill
\begin{minipage}[b]{.31\textwidth}\centering
\includegraphics[width=.8\textwidth]{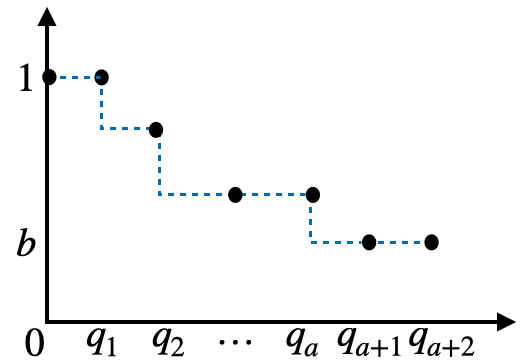}\\
{\small (c) 3 steps}
\end{minipage}
\label{fig:steps}
\end{figure}
We prove that every optimal $\vecq^*$ must have exactly one step by using a step removal construction, which shows that for any feasible distribution with more than one step, we can remove one step to obtain a feasible distribution with a lower objective function value. Consequently, an optimal distribution has at most one step. Figure~\ref{fig:steps-removal} illustrates the step removal approach:  \begingroup the last plateau lying strictly between its neighboring levels is merged either upward into the preceding level or downward into a common level with the tail, and one of the two merges strictly decreases the objective function while removing one \emph{step}. In the following lemma, we establish the step removal approach formally. \endgroup{}
\begin{figure}[htp]
\centering
\includegraphics[width=.5\textwidth]{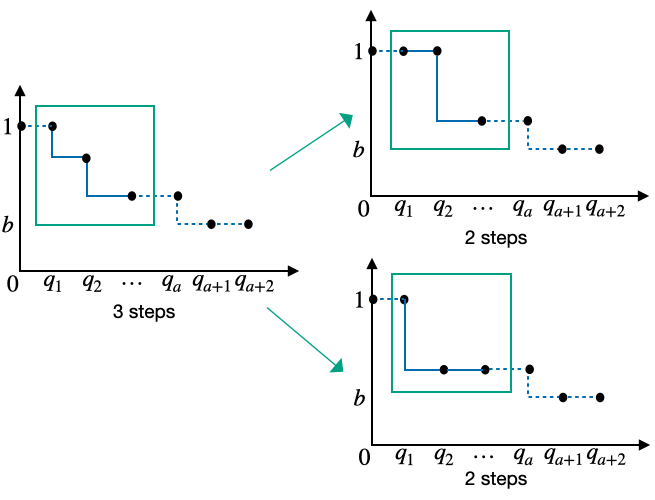}
\caption{Step Removal.}
\label{fig:steps-removal}
\end{figure}

\begingroup
\begin{lemma}\label{lem:4-remstep}
Consider an objective function $\operatorname{OBJ}_i(\vecq)$, where $i \in \{1,2\}$. For any distribution $D_{\vecq}$ that is feasible for the corresponding program \eqref{obj1} or \eqref{obj2} and for which there exist some $k \geq1$ and $m \geq0$ such that
\(
q_{a-k-m} > q_{a-k-m+1} = \cdots = q_{a-m} > q_{a-m+1} = q_{a-m+2} = \cdots = b,
\)
the distributions $D_{\vecq'}$ and $D_{\vecq''}$ constructed below are feasible and satisfy $\min\{\operatorname{OBJ}_i(\vecq'), \operatorname{OBJ}_i(\vecq'')\} < \operatorname{OBJ}_i(\vecq)$:
\begin{itemize}
\item $\vecq'$: fix $q_1,\ldots,q_{a-k-m}$, set $q_{a-k-m+1} = \cdots = q_{a-m} = q_{a-m+1} = \cdots = b$, and solve equation \eqref{eq:con1} to obtain $b$
\item $\vecq''$: fix $q_1,\ldots,q_{a-k-m}$, set $q_{a-k-m+1} = \cdots = q_{a-m} = q_{a-k-m}$ and $q_{a-m+1} = \cdots = q_a = b$, and solve equation \eqref{eq:con1} to obtain $b$
\end{itemize}
\end{lemma}
\endgroup
With this lemma, we can formally state our proof of Lemma~\ref{lem:shiftgeo} in a constructive way.

\proof{Proof of Lemma~\ref{lem:shiftgeo}.} 
    Consider any objective function $\operatorname{OBJ}_i$ where $i \in \{1,2\}$. \begingroup We first show that an optimal solution exists. The feasible set of \eqref{obj1} (equivalently, of \eqref{obj2}) is $\{(q_1,\cdots,q_a,b): 1\geq q_1\geq \cdots \geq q_a\geq b\geq 0 \text{ and } \eqref{eq:con1} \text{ holds}\}$. This set is bounded, and it is also closed: the left-hand side of \eqref{eq:con1} is continuous wherever $b<1$, and no feasible sequence can accumulate at $b=1$, because $q_i\geq b$ for all $i$ implies $q_1\cdots q_a\frac{1}{1-b}\geq \frac{b^a}{1-b}$, which exceeds $\mu_T$ once $b$ is close enough to $1$. Since $\operatorname{OBJ}_i$ is continuous on this set, a minimizer $\vecq^*$ exists. We also note that every feasible point has $b\in (0,1)$: $b=1$ is excluded above, and $b=0$ would reduce \eqref{eq:con1} to $q_1+q_1q_2+\cdots+q_1\cdots q_a=\mu_T$, whose left-hand side is at most $a<\mu_T$.

   Next, suppose that some optimal solution $\vecq^*$ has more than one \emph{step}. Examining the sequence $\{q_0, q_1,\ldots, b\}$ from its last element, let $a-m\geq 1$ be the maximum index where $q_{a-m} > b$ (with $0 \leq m \leq a-1$); such an index exists, for otherwise $q_1 = q_2 = \cdots = b$ and $D_{\vecq^*}$ is a geometric distribution with exactly one step. Then, let $a-k-m$ be the largest index $i$ for which $q_{i}>q_{a-m}$; such an index exists as well (recall $q_0=1$), for otherwise $1=q_0=q_1=\cdots =q_{a-m}>q_{a-m+1}=\cdots = b$ and $D_{\vecq^*}$ is a shifted geometric distribution with exactly one step. Therefore, $\vecq^*$ satisfies the condition of Lemma~\ref{lem:4-remstep}, which produces a feasible distribution with a strictly lower objective value, contradicting the optimality of $\vecq^*$. We conclude that every optimal solution of \eqref{obj1} and \eqref{obj2} has exactly one \emph{step} (a zero-step vector requires $b=1$ and is infeasible), i.e., is a shifted geometric distribution with $b\in (0,1)$. Since the reformulation in Section~\ref{app:IFRe-reform} shows that the optimal value of \eqref{OPT:mu-c} is attained by an optimal solution of \eqref{obj1} or \eqref{obj2}, the lemma follows.
   \Halmos
\endproof
\endgroup

To complete the proof, it remains to prove the key Lemma~\ref{lem:4-remstep}. \begingroup We embed $\vecq$, $\vecq'$, and $\vecq''$ in a one-parameter family of feasible distributions, obtained by varying the common value of the plateau $q_{a-k-m+1}=\cdots=q_{a-m}$ and adjusting $b$ so that constraint \eqref{eq:con1} is preserved exactly. We then show that the logarithm of the objective function is strictly concave along this family, so its minimum over the family is attained only at the two endpoints, which are precisely $\vecq'$ and $\vecq''$.  \endgroup


\proof{Proof of Lemma~\ref{lem:4-remstep}.} \begingroup
Let $\ell:=a-k-m$, recall that $q_0=1$, and write
\(
\bar q:=q_{\ell+1}=\cdots=q_{\ell+k}.
\)
Thus, $\bar q$ is the plateau value in the original vector, and
$q_\ell>\bar q>b$.

Then we move this plateau while adjusting the tail so that the
mean constraint remains satisfied. We may lower the plateau until it
meets the tail, which produces $\vecq'$, or raise it until it meets the
preceding level, which produces $\vecq''$. These two vectors are the
endpoints of a one-parameter family of feasible vectors, and the
original vector is an interior point of this family.
To formalize this idea, we keep $q_1,\ldots,q_\ell$ fixed and define
\(
P:=\prod_{i=1}^{\ell}q_i,
R:=
\frac{\mu_T-\sum_{j=1}^{\ell}\prod_{i=1}^{j}q_i}{P}.
\)

We first construct a feasible family containing the original vector.
For $j\geq 1$, let $S_j(x):=\sum_{r=1}^{j}x^r$, and set $S_0(x):=0$.
If the plateau value is changed to $x$ and the tail is changed to $B$,
then, after summing the geometric tail, constraint~\eqref{eq:con1}
reduces to
\(
S_k(x)+x^k\frac{B}{1-B}=R.
\)
Because the fixed prefix contributes at most $\ell$, $0<P\leq 1$, and
$\mu_T>a=\ell+k+m$, we have
$R>k+m\geq k\geq S_k(x)$ for every $x\in(0,1]$.
Thus, with
\[
W(x):=R-S_k(x),
\qquad
B(x):=\frac{W(x)}{W(x)+x^k},
\]
we have $W(x)>0$ and $B(x)\in(0,1)$. Let $\vecq(x)$ be the vector
obtained by replacing the plateau value by $x$ and all subsequent
coordinates and the tail by $B(x)$. The original constraint gives
$B(\bar q)=b$, and therefore
\(
\vecq(\bar q)=\vecq.
\)

We next identify the values of $x$ for which $\vecq(x)$ is feasible.
Since the constraint holds by construction, it remains only to check
$q_\ell\geq x\geq B(x)$. Moreover,
\begin{align*}
B(x)\leq x
\iff (1-x)W(x)\leq x^{k+1} 
\iff (1-x)R\leq x 
\iff x\geq \frac{R}{R+1},
\end{align*}
where we use $(1-x)S_k(x)=x-x^{k+1}$. Define
\(
q_L:=\frac{R}{R+1}.
\)
Then $B(q_L)=q_L$. Since $B(\bar q)=b<\bar q$ and
$\bar q<q_\ell$, we have
\(
q_L<\bar q<q_\ell.
\)
It follows that $\vecq(x)$ is feasible for every
$x\in[q_L,q_\ell]$. At the lower endpoint, the plateau merges with the
tail, while at the upper endpoint it merges with the preceding level.
Consequently,
\(
\vecq(q_L)=\vecq',
\vecq(q_\ell)=\vecq''.
\)

We now compare the objective functions along this family. For
$p\in(0,1)$, define
\(
\operatorname{OBJ}(\vecq;p)
:=q_1\cdots q_a\frac{b}{1-b+bp}.
\)
Here, $p=1$ gives $\operatorname{OBJ}_1$, and
$p=1/(\mu_T-a)$ gives $\operatorname{OBJ}_2$ in its corresponding
regime. Using
\(
q_1\cdots q_a=Px^kB(x)^m,
\frac{B(x)}{1-B(x)+pB(x)}
=\frac{W(x)}{x^k+pW(x)},
\)
we obtain
\[
L_p(x)
:=\ln \operatorname{OBJ}(\vecq(x);p)
=\ln P+m\ln B(x)
-\ln\left(\frac{1}{W(x)}+\frac{p}{x^k}\right).
\]

It remains to show that $L_p$ is strictly concave. Let
\(
F(x):=\frac{1}{W(x)}, 
G(x):=\frac{p}{x^k}.
\)
Both functions are strictly log-convex because
\(
(\ln F)''
=\frac{S_k''(x)}{W(x)}
+\left(\frac{S_k'(x)}{W(x)}\right)^2>0,
(\ln G)''=\frac{k}{x^2}>0.
\)
If $\lambda:=F/(F+G)$, then
\([\ln(F+G)]''
=\lambda(\ln F)''
+(1-\lambda)(\ln G)''
+\lambda(1-\lambda)
 \big((\ln F)'-(\ln G)'\big)^2
>0.
\)
Hence, $\ln(1/W(x)+p/x^k)$ is strictly convex. On the other hand,
\(
B(x)=\frac{R-S_k(x)}{R-S_{k-1}(x)}.
\)
Suppressing the argument $x$, direct differentiation gives
\(
(\ln B)''
={}-\big[
\frac{S_k''}{R-S_k}
+\big(\frac{S_k'}{R-S_k}\big)^2
\big] +\big[
\frac{S_{k-1}''}{R-S_{k-1}}
+\left(\frac{S_{k-1}'}{R-S_{k-1}}\right)^2
\big]
\leq 0\).
The inequality follows because $S_k'\geq S_{k-1}'$,
$S_k''\geq S_{k-1}''$, and $R-S_k\leq R-S_{k-1}$. Therefore,
\(
L_p''(x)
=m(\ln B(x))''
-\left[\ln\left(\frac{1}{W(x)}+\frac{p}{x^k}\right)\right]''
<0.
\)

Thus, $L_p$ is strictly concave on $[q_L,q_\ell]$, and the original
vector corresponds to the interior point $\bar q$. Hence,
\(
L_p(\bar q)>
\min\{L_p(q_L),L_p(q_\ell)\}.
\)
Since the logarithm is strictly increasing, we conclude that, for
$i\in\{1,2\}$,
\(
\operatorname{OBJ}_i(\vecq)
>
\min\{\operatorname{OBJ}_i(\vecq'),
       \operatorname{OBJ}_i(\vecq'')\}.
\)
\Halmos
\endproof

\endgroup

\subsubsection{Establish the worst-case ratio.}\label{app:4-1e} We show that $\frac{1}{e}$ is a lower bound for the optimal objective value of $\inf_{\mu_T,a:\, a\in \mathbb{Z}_+, \mu_T-a>0, \mu_T>0} \frac{g(\mu_T,a)}{\mu_T-a}$ when limiting to shifted geometric distributions. We discuss two cases based on whether $\mu_T-a\leq 1$.

\noindent \textbf{Case 1: $\mu_T-a\leq 1$.}
 Plugging in $q_i = 1$ for $i\leq m-1$ and $q_i = b$ for $ i\geq m$, or $q_i = b$ for all $i\geq 1$, the factor-revealing program can be reformulated as follows:
\[
\left\{\inf_{m,b, \mu_T,a:\, a\in \mathbb{Z}_+, \mu_T-a>0, \mu_T>0} \frac{b^{a-m+2}}{\mu_T-a} :m-1+\frac{b}{1-b}=\mu_T, \mu_T-a\leq 1, 0< b< 1, m=1,\cdots, a+1\right \}.
\]
From the constraint $m-1+\frac{b}{1-b}=\mu_T$, we have $b=\frac{\mu_T-m+1}{\mu_T-m+2}$. Let $x=a-m+2$, $y=\mu_T-a$, then the optimization problem becomes: 
\(
\left\{ \inf_{x,y,a} \frac{(1-\frac{1}{x+y})^{x}}{y}: 0<y\leq 1, \quad x=1,\cdots, a+1; a\in \mathbb{Z}_+
\right \}
\)
Let $f(x,y)= \frac{(1-\frac{1}{x+y})^{x}}{y}$, then $\frac{\partial f(x,y)}{\partial y}= (1-\frac{1}{x+y})^{x-1}\frac{1}{y^2(x+y)^2}[xy-(x+y)^2(1-\frac{1}{x+y})]$
$=(1-\frac{1}{x+y})^{x-1}\frac{1}{y^2(x+y)^2} [xy-(x+y)^2 +(x+y)]\leq 0$ since $x\geq 1$ and $y\geq 0$. Then $y^*=1$ and the optimization problem becomes: 
\(
\left\{
 \inf_{x=1, \cdots ,a+1; a\in \mathbb{Z}_+} (1-\frac{1}{1+x})^x
\right \}.
\)
 Let $h(x)=\ln (1-\frac{1}{1+x})^x$, then $h'(x)\leq 0$, so the objective function value decreases as $x$ increases. When $x\rightarrow \infty$, the objective function value approaches $\frac{1}{e}$. 
 
\noindent \textbf{Case 2: $\mu_T-a\geq 1$.} In this case, after plugging in $q_i = 1$ for all $i\leq m-1$ and $q_i = b$ for $ i\geq m$, or $q_i = b$ for all $i\geq 1$, we have:
\[g(\mu_T,a)=
 \left\{ \inf_{0\leq b\leq 1;m=1,\cdots, a+1} b^{a-m+2}\frac{\mu_T-a}{(\mu_T-a)(1-b)+b}:m-1+\frac{b}{1-b}=\mu_T
\right \}.  \]

The worst-case ratio is $\inf_{\mu_T>0, a\in \mathbb{Z}_+, \mu_T-a\geq 1} \frac{g(\mu_T,a)}{\mu_T-a}$. Let $x=a-m+1$ and $d=\mu_T-a$. Then the worst-case ratio is: 
\(\left\{ \inf_{0\leq b\leq 1;x=0,\cdots, a;d\geq 1;a\in \mathbb{Z}_+} b^{x+1}\frac{1}{d(1-b)+b}:\frac{b}{1-b}-d=x 
\right \}.\)
From the last constraint, we replace $b$ with $\frac{d+x}{d+x+1}$, which naturally satisfies $0< \frac{d+x}{d+x+1} < 1$. Thus, the factor-revealing program becomes: 
\begin{align*}
& \inf_{x=0,\cdots, a; d\geq 1;a\in \mathbb{Z}_+} (1-\frac{1}{d+x+1})^{x+1}\frac{1}{\frac{d}{d+x+1}+\frac{d+x}{d+x+1}} = \inf_{x=0,\cdots, a;s\geq x+1;a\in \mathbb{Z}_+} (1-\frac{1}{s+1})^{x+1}\frac{s+1}{2s-x}.  \\
\end{align*}
\vspace{-5pt}
In the equality, we let $s=d+x$. Since $d\geq 1$, we have $s\geq x+1$. Let $f(s)= (1-\frac{1}{s+1})^{x+1}\frac{s+1}{2s-x} $, then
\(
f'(s)=\dfrac{x\left(s-x-1\right)\left(1-\frac{1}{s+1}\right)^x}{\left(s+1\right)\left(2s-x\right)^2}\geq 0
\)
as $s\geq x+1$. Thus, for any $x=0,\cdots,a$ and $a\in \mathbb{Z}_+$, the objective value at $s=x+1$ is the smallest among all $s$ satisfying $s\geq x+1$, which is equivalent to $\mu_T-a=1$. Therefore, this is the same as the case when $\mu_T-a\leq 1$, and thus the worst-case ratio should be the same.
Summarizing the two cases, we conclude that the worst-case ratio is $1/e$. 
\begin{remark}
Lemma~\ref{lem:4-Stoch-DT} is stated for every $p\in(0,1)$, and the argument above proves it when
$p(\mu_T-a)\leq 1$. Suppose instead that $p(\mu_T-a)>1$. Since $p\leq 1$, this forces $\mu_T-a>1$,
so $p_0:=\frac{1}{\mu_T-a}\in(0,1)$ and $p_0<p$. On the event $\{T>a\}$ we have $T-a\geq 1$, so
$1-(1-p)^{T-a}$ is nondecreasing in $p$, hence so is
$p\mapsto\EE_T[(1-(1-p)^{T-a})\1\{T>a\}]$. Applying the case already proved to the same IFR
distribution at $p_0$, for which $p_0(\mu_T-a)=1$, gives
\(
\EE_T\!\left[(1-(1-p)^{T-a})\1\{T>a\}\right]
\geq \EE_T\!\left[(1-(1-p_0)^{T-a})\1\{T>a\}\right]
\geq \tfrac{1}{e}\,p_0(\mu_T-a)=\tfrac{1}{e}=\tfrac{1}{e}\min\{p(\mu_T-a),1\},
\)
which is the claim in this case.
\end{remark}
\begingroup
\begin{remark}[Tightness of the $1/e$ constant]\label{rem:ifr-tight}
Let $T_m$ be geometrically distributed on $\{1,2,\ldots\}$ with mean $m$, where $m\geq 2$ is an integer, so that
\(
\Prb(T_m\geq k)
=
\left(1-\frac{1}{m}\right)^{k-1}, k\geq 1.
\)
The distribution of $T_m$ is IFR. Take $a=m-1$ and $p=1$. Then $\mu_{T_m}-a=1$ and
\[
\EE\!\left[
(1-(1-p)^{(T_m-a)_+})\1\{T_m>a\}
\right]
=
\Prb(T_m\geq m)
=
\left(1-\frac{1}{m}\right)^{m-1}
\longrightarrow
\frac{1}{e}
\qquad (m\to\infty).
\]
Hence no universal constant larger than $\frac{1}{e}$ is possible in Lemma~\ref{lem:4-Stoch-DT}.
\end{remark}
\endgroup

\subsection{Proof of Lemma~\ref{lem:IFRe}} \label{app:prf-lem-IFRe}

\medskip\noindent\textbf{Lemma~\ref{lem:IFRe} (restated).} {\itshape
Suppose $T$ follows an IFR distribution. For any inventory vector $\vecx\in \mathbb{Z}_+^n$ and any product $i \in \NN$, $\EE_T[\ZDA_i(\vecx,T)]\geq \frac{1}{e}\ZFP_i(\vecx,\mu_T)$.
}\medskip

\proof{Proof of Lemma~\ref{lem:IFRe}.}
The proof follows a line of reasoning similar to that for Lemma~\ref{lem:4-product}.
\begingroup Let $\operatorname{Out}(i)$ be defined as in the proof of Lemma~\ref{lem:4-product}, with $T$ replaced by $\mu_T$,\endgroup{} let $a=\sum_{j\in \operatorname{Out}(i)} x_j$, and let $p=\phi(i,N(\vecx-\sum_{j\in \operatorname{Out}(i)}x_j \vece_j))$. \begingroup Replacing $T$ by $\mu_T$ in \eqref{eq:4-uni} and \eqref{eq:app-4-4}, we have:\endgroup{}
\begin{equation}\label{eq:4-uni1}
\ZFP_i(\vecx,\mu_T)= \ZFP_i(\vecx-\sum_{j\in \operatorname{Out}(i)}x_j \vece_j, \mu_T-a)= \min \{p(\mu_T-a), x_i \}.
\end{equation}

Next, we consider $\EE_T[\ZDA_i(\vecx,T)]$. For any value of $T>a$, we remove the influence of products in $\operatorname{Out}(i)$ by using Lemma~\ref{lem:4-sin-subs}: $\ZDA_i(\vecx, T)\geq \ZDA_i(\vecx-\sum_{j\in \operatorname{Out}(i)}x_j \vece_j, T-a)$. Using Lemma~\ref{lem:4-touni}, it follows that $ \ZDA_i(\vecx-\sum_{j\in \operatorname{Out}(i)}x_j \vece_j, T-a)\geq \ZDA_i (p,T-a,x_i)$. Recall here that $\ZDA_i(p,T-a, x_i)$ is the expected number of units of $i$ sold in the stochastic consumption process of length $T-a$ with uniform purchase probability $p$ and inventory $x_i$. 
If $x_i=0$, then $\ZFP_i(\vecx,\mu_T)=0$ by \eqref{eq:4-uni1} and the claim is immediate, so we assume $x_i\geq 1$. Moreover, since $\ZDA_i(\vecx,T)\geq 0$, we have $\EE_T[\ZDA_i(\vecx,T)]\geq \EE_T[\ZDA_i(\vecx,T)\1\{T>a\}]$. Thus, we have:
 \begin{align}
    \EE_T[\ZDA_i(\vecx,T)] &\geq \EE_T[\ZDA_i(\vecx,T)\1\{T>a\}]\nonumber\\
    & \geq \EE_T[\ZDA_i(\vecx-\sum_{j\in \operatorname{Out}(i)}x_j \vece_j, (T-a)_+)\1\{T>a\}]\nonumber\\
    & \geq \EE_T[\ZDA_i (p,(T-a)_+,x_i)\1\{T>a\}]\nonumber\\
    &\geq x_i \EE_T[\ZDA_i (\tfrac{p}{x_i},(T-a)_+,1)\1\{T>a\} ]. \label{eq:4-uni2}
\end{align}
 The last inequality is from Lemma~\ref{lem:aux}.
If $\ZFP_i(\vecx,\mu_T)=0$, the claim is immediate because the left-hand side is nonnegative. Otherwise, the separability representation in \eqref{eq:4-uni1} gives $\min\{p(\mu_T-a),x_i\}>0$, which forces $p>0$ and $a<\mu_T$, so Lemma~\ref{lem:4-Stoch-DT} applies to the single-unit purchase probability $\frac{p}{x_i}\in (0,1)$ and yields
\[
\EE_T[\ZDA_i(\tfrac{p}{x_i},(T-a)_+,1)\1\{T>a\}]
\geq \frac{1}{e}\min\Big\{\frac{p}{x_i}(\mu_T-a),\,1\Big\}.
\]
Combining this with inequality \eqref{eq:4-uni2} and equation \eqref{eq:4-uni1}, we have the desired result as follows: 
\[
\EE_T[\ZDA_i(\vecx,T)]
\geq \frac{x_i}{e}\min\Big\{\frac{p}{x_i}(\mu_T-a),1\Big\}
= \frac{1}{e}\min\{p(\mu_T-a),x_i\}
= \frac{1}{e}\ZFP_i(\vecx,\mu_T).  \Halmos
\]
\endproof

\subsection{Proof of Lemma~\ref{lem:tran-d}} \label{app:thm-tran-d}
We first prove the following lemma, which provides a lower bound for the expected consumption in terms of the consumption in a single-unit setting. Note that Proposition 2 in \cite{ma2021dynamic} proves the $1-1/e$ competitive ratio for stationary customer arrivals, while our proof works for nonstationary customer arrivals. 
\begin{lemma}
\label{lem:aux}
Let $Z_i$ denote the expected consumption of product $i$ in the stochastic consumption process with $T$ customers and inventory $c_i\in \mathbb{Z}_{+}$, where each customer $t$ purchases the product independently with a given probability $p_{it}\in(0,1)$. Then we have:
    $Z_i\geq  c_i \EE[\min\{\sum_{t=1}^T\text{Bernoulli}(\frac{p_{it}}{c_i}), 1\}]$.
\end{lemma}

\proof{Proof.}
\begingroup
We first note that $\EE[\min\{\sum_{t=1}^T\text{Bernoulli}(\frac{p_{it}}{c_i}), 1\}]$ is equivalent to the expected consumption in an alternative stochastic process: Each of the $c_i$ units of product $i$ is treated as a distinct item labeled $1,2,\ldots,c_i$. Upon arrival, each customer $t$ is randomly routed to one of these units with equal probability $\frac{1}{c_i}$. If customer $t$ is routed to an available unit, they purchase it with probability $p_{it}$. Thus, the expected consumption of a particular unit of $i$ in this process is $\EE[\min\{\sum_{t=1}^T\text{Bernoulli}(\frac{p_{it}}{c_i}), 1\}]$. Since all units are identical, the total expected consumption in this process is \(c_i\EE[\min\{\sum_{t=1}^T\text{Bernoulli}(\frac{p_{it}}{c_i}), 1\}].\)
\endgroup

We construct a coupling between the two processes by considering the sample path realizations of customer types and choices (without the inventory constraint). More specifically, let $ \boldsymbol{\omega}=(\omega_t)_{t=1}^T$ denote this sample path where $\omega_t= 1$ with probability $p_{it}$ and $\omega_t= 0$ with probability $1-p_{it}$. Note that in both processes, the probability of any realization of $\boldsymbol{\omega}$ is identical. Given any $\boldsymbol{\omega}$, the total consumption in $Z_i$ is $\min\{ \sum_{t=1}^T\omega_t, c_i\}$. However, in the consumption process with routing, a customer may be routed to an out-of-stock unit even when other units are available, resulting in a total consumption that is less than $\min\{ \sum_{t=1}^T\omega_t, c_i\}$. Therefore, for any sample path, the total consumption of product $i$ with routing is less than or equal to that in $Z_i$. Taking the expectation over all sample paths, we conclude $Z_i\geq  c_i \EE[\min\{\sum_{t=1}^T\text{Bernoulli}(\frac{p_{it}}{c_i}), 1\}]$.
\Halmos
\endproof

\medskip\noindent\textbf{Lemma~\ref{lem:tran-d} (restated).} {\itshape
Given any $\vecc\in \mathbb{Z}_+^n$, any number $T$ of customers, and any choice model with the substitutability property, let $\fIND(\vecc,T)$ denote the total expected revenue obtained by the assortment policy above. Then, $\fIND(\vecc,T)\geq (1-\frac{1}{e}) \fLPM(\vecc,T)$.
}\medskip

\proof{Proof of Lemma~\ref{lem:tran-d}.} Let $Z_i$ denote the expected consumption of product $i$ in $\fIND(\vecc,T)$ and $\ZLP_i$ denote the consumption of product $i$ in $\fLPM(\vecc,T)$. It suffices to show $Z_i\geq (1-\frac{1}{e})\ZLP_i$. If $c_i=0$, then $Z_i=\ZLP_i=0$ and the claim is trivial, so we assume $c_i\geq 1$. Our key idea is to first lower bound $Z_i$ using the consumption in a stochastic process with deterministic purchase probabilities derived from the CDLP solution, and then bound the latter through a single-unit reduction.

In $\fIND(\vecc,T)$, at time $t$, if the arriving customer is of type $k$, we sample an assortment $S_t$ from the distribution $\{\Prb(S_t=S)=y_k(S)\}_{S\subseteq \NN}$ and offer the intersection of $S_t$ and the available products $N_t$, so the purchase probability of product $i$ given $N_t$ is $p_{it} = \sum_{k\in \MM, S\subseteq \NN} \lambda_{k,t} y_k(S)\phi_k(i,S\cap N_t)$. Note that $p_{it}$ depends on the realized availability set $N_t$. Let $X_t$ indicate that customer $t$ purchases product $i$, so that $\Prb(X_t=1\mid N_t)=p_{it}$ while product $i$ remains in stock. Thus, $Z_i = \EE[\min\{\sum_{t=1}^T X_t,c_i  \}]$.

Recall that the revenue from $\fLPM(\vecc,T)$ can be interpreted as the total revenue of a fluid process where a mixture of assortments is provided to each customer. If the customer is of type $k$, the mixture includes a $y_k(S)$ portion of assortment $S$. Let $p_{it}'$ represent the expected probability that customer $t$ will purchase product $i$ in $\fLPM(\vecc,T)$. Then, $p_{it}' =\sum_{k\in \MM, S\subseteq \NN} \lambda_{k,t} y_k(S)\phi_k(i,S)$, which does not depend on the sample path. The consumption for product $i$ in this fluid process, denoted by $\ZLP_i$, is $\sum_{t=1}^T p_{it}'=\sum_{t=1}^T \sum_{k\in \MM, S\subseteq \NN} \lambda_{k,t} y_k(S)\phi_k(i,S)$, which is less than or equal to $c_i$ due to the inventory constraint in the CDLP.
Due to the substitutability property, we have $\phi_k(i,S_t\cap N_t)\geq \phi_k(i,S_t)$ whenever product $i$ remains in stock; thus, $p_{it}\geq p_{it}'$ on every sample path before product $i$ stocks out. As in Lemma~\ref{lem:4-touni}, realize $X_t=\1\{U_t\leq p_{it}\}$ and $X_t'=\1\{U_t\leq p_{it}'\}$ with independent uniform random variables $(U_t)_{t=1}^T$; then $X_t\geq X_t'$ while product $i$ remains in stock, and once it stocks out the capped consumption already equals $c_i$, so pathwise: \(
Z_i\geq \EE[\min\{\sum_{t=1}^T X_t',c_i  \}]=\EE[\min\{\sum_{t=1}^T \text{Bernoulli}(p_{it}'),c_i  \}]. \)

The right-hand side is the expected consumption of product $i$ in a stochastic consumption process with $T$ customers, inventory $c_i$, and deterministic purchase probabilities $p_{it}'$. By Lemma~\ref{lem:aux}, it is at least $c_i \EE[\min\{\sum_{t=1}^T\text{Bernoulli}(\frac{p_{it}'}{c_i}), 1\}]$. Here \(\EE[\min\{\sum_{t=1}^T\text{Bernoulli}(\frac{p_{it}'}{c_i}),1\}]\) is the expected consumption of product $i$ in a stochastic process with a single unit in the inventory and purchase probability $\frac{p_{it}'}{c_i}$ for customer $t$. Then we have
\begin{align}
&c_i \EE[\min\{\sum_{t=1}^T\text{Bernoulli}(\frac{p_{it}'}{c_i}), 1\}] = c_i (1-\prod_{t=1}^T (1-\frac{p_{it}'}{c_i}))
\geq c_i (1- (1-\frac{p_{i}'}{c_i})^T)\nonumber \\
&\geq (1-\frac{1}{e})c_i\min\{\frac{p_i'}{c_i} T, 1\} = (1-\frac{1}{e})\min\{\sum_{t=1}^T p_{it}' ,c_i\}.
\label{eq:thm7-2}
\end{align}
 Here the first inequality is from the AM-GM inequality and
$p_i'=\frac{\sum_{t=1}^T p_{it}'}{T}$. The second inequality is from the single-item result of $1-\frac{1}{e}$, see Lemma~\ref{lem:4-uni}. Since $\sum_{t=1}^T p_{it}'=\ZLP_i\leq c_i$, the right-hand side of \eqref{eq:thm7-2} equals $(1-\frac{1}{e})\ZLP_i$.

By combining the two displays above with Lemma~\ref{lem:aux}, we have $Z_i\geq (1-\frac{1}{e})\ZLP_i$.
\Halmos
\endproof

\section{Overall Performance Guarantees: Remaining Proofs}\label{app:5}
\subsection{Performance Guarantees for DA} \label{app:5-da}

\begingroup
\medskip\noindent\textbf{Theorems~\ref{thm:gua-DA-det},  \ref{thm:gua-DA-ifr} and ~\ref{thm:gua-budget-da} (restated)} {\itshape
Consider any accuracy level $\epsilon\in (0,1)$.
When $T$ is deterministic, the threshold-based \textsc{da} algorithm has approximation factor $0.319-\epsilon$ under a cardinality or budget constraint. When $T$ is stochastic and follows a distribution with the IFR property, \textsc{da} can be approximated within a factor of $0.188-\epsilon$ under a cardinality or budget constraint.
}\medskip
\endgroup

\proof{Completion of the proofs of Theorems~\ref{thm:gua-DA-det}, \ref{thm:gua-DA-ifr}, and \ref{thm:gua-budget-da}.}
\begingroup
For deterministic $T$, the approximation guarantees under a budget constraint using the {\small THR} algorithm are the same as those under a cardinality constraint: the budget-constrained first stage achieves the same $(1-\frac{1}{e}-\epsilon)$ factor (Appendix~\ref{app:single-T}), and the transformation step is independent of the form of the inventory constraint (Section~\ref{sec:ext-budget}). For small $T$, the static inventory $\vecc^{\operatorname{static}}$ is redefined as the indicator vector of an optimal budget-feasible static assortment; computing this assortment is a knapsack-type problem that admits an FPTAS \citep{desir2014near}, and the resulting $(1-\epsilon)$ factor is absorbed into $\epsilon$. So we have the same final performance guarantee for deterministic $T$ under a budget constraint.

For stochastic $T$, using the {\small THR} algorithm still yields the same performance guarantees under both cardinality and budget constraints. Combining Theorem~\ref{thm:da-fluid-opt} (for the budget constraint, its budget extension in Appendix~\ref{app:single-T}) and the results in Sections~\ref{sec:da-transform} and \ref{sec:da-stochasticT}, we obtain an inventory vector $\vecc^{\operatorname{ALG}}$ such that $\EE_T[\fDA(\vecc^{\operatorname{ALG}},T)]\geq \frac{1}{e}(1-\frac{1}{e}-\epsilon)(1-\frac{1}{\mu_T+1})\max_{\vecc \in \mathcal{F}}  \fLP(\vecc,\mu_T)$. Therefore, we have $\EE_T[\fDA(\vecc^{\operatorname{ALG}},T)]\geq (\frac{1}{e}(1-\frac{1}{e})(1-\frac{1}{\mu_T+1})-\epsilon) \max_{\vecc \in \mathcal{F}}  \fLP( \vecc, \mu_T)\geq (\frac{1}{e}(1-\frac{1}{e})(1-\frac{1}{\mu_T+1})-\epsilon) \max_{\vecc \in \mathcal{F}}  \EE_T [\fDA (\vecc, T)]$. Here the last inequality is from Proposition~\ref{prop:ub} in Appendix~\ref{app:prop-ub}, which states that $ \fLP( \vecc, \mu_T)$ is an upper bound for $\EE_T[\fDA (\vecc, T)]$.
For the static inventory, since the expected revenue obtained from each customer is at most $\fDA(\vecc^{\operatorname{static}}, 1)$ and the expected number of customers is $\mu_T$, we have $\max_{\vecc \in \mathcal{F}} \EE_T[\fDA (\vecc, T)]\leq \mu_T \fDA(\vecc^{\operatorname{static}}, 1)$. Under the IFR property, $\Prb(T\geq 1)\geq \frac{\mu_T}{\mu_T+1}$. Indeed, the ratios $\frac{\Prb(T\geq k+1)}{\Prb(T\geq k)}=1-\frac{\Prb(T=k)}{\Prb(T\geq k)}$ are non-increasing in $k$ while $\Prb(T\geq k)>0$, so $\Prb(T\geq k)\leq \Prb(T\geq 1)^k$ for every $k\geq 1$. Hence, whenever $\Prb(T\geq 1)<1$, $\mu_T=\sum_{k\geq 1}\Prb(T\geq k)\leq \frac{\Prb(T\geq 1)}{1-\Prb(T\geq 1)}$, which rearranges to $\Prb(T\geq 1)\geq \frac{\mu_T}{\mu_T+1}$; the bound is trivial when $\Prb(T\geq 1)=1$. By the monotonicity of $\fDA(\vecc,t)$ in $t$, $\EE_T[\fDA(\vecc^{\operatorname{static}}, T)]\geq \Prb(T\geq 1)\fDA(\vecc^{\operatorname{static}}, 1)\geq \frac{1}{\mu_T+1}\max_{\vecc \in \mathcal{F}} \EE_T[\fDA (\vecc, T)]$. When $\mu_T< 4.31$, we use $\vecc^{\operatorname{static}}$ as the starting inventory, whereas we use $\vecc^{\operatorname{ALG}}$ when $\mu_T \geq 4.31$. Thus, the worst-case performance guarantee is $\frac{1}{5.31}-\epsilon\geq 0.188-\epsilon$. The same argument applies under a budget constraint, with $\vecc^{\operatorname{static}}$ again redefined as the indicator vector of an optimal budget-feasible static assortment as above.
\Halmos
\endgroup
\endproof
\medskip\noindent\textbf{Theorem~\ref{thm:thr-dp} (restated).} {\itshape Suppose we use the {\small DP} algorithm for the first-stage approximation.
Consider any accuracy level $\epsilon\in (0,1)$.
When $T$ is deterministic, \textsc{da} can be approximated within a factor of $0.474-\epsilon$ under a cardinality constraint. When $T$ is stochastic and follows a distribution with the IFR property, \textsc{da} can be approximated within a factor of $0.268-\epsilon$ under a cardinality constraint.
}\medskip
\proof{Proof of Theorem~\ref{thm:thr-dp}.}
\begingroup
The analysis of the {\small DP} algorithm differs from that of the {\small THR} algorithm only by replacing the \(1-\frac{1}{e}-\epsilon\) factor with \(1-\epsilon\) in the first stage. 
\endgroup

When $T$ is deterministic, combining Theorem~\ref{thm:dp-fluid-opt} and the results in Section~\ref{sec:da-transform}, we obtain an inventory vector $\vecc^{\operatorname{ALG}}$ such that $\fDA(\vecc^{\operatorname{ALG}},T)\geq \left( (1-\frac{1}{e})(1-\frac{1}{T+1}) -\epsilon \right)\max_{\vecc \in \mathcal{F}} \fLP(\vecc,T)$. Since $\fLP( \vecc, T)$ is an upper bound for $\fDA (\vecc, T)$ (Proposition~\ref{prop:ub} in Appendix~\ref{app:prop-ub}), $\vecc^{\operatorname{ALG}}$ provides a $( (1-\frac{1}{e})(1-\frac{1}{T+1}) -\epsilon )$-approximate solution.   
Note that this approximation factor decreases as $T$ decreases. For small values of $T$, we use $\vecc^{\operatorname{static}}$ as the starting inventory, i.e., the indicator vector of an optimal static assortment as defined in Section~\ref{sec:3-bad}.
Since the expected revenue obtained from each customer is at most that of the optimal static assortment, we have $\max_{\vecc \in \mathcal{F}}  \fDA (\vecc, T)\leq T \fDA(\vecc^{\operatorname{static}}, 1)$. Due to the monotonicity of $\fDA(\vecc,t)$ in $t$, it follows that $\fDA(\vecc^{\operatorname{static}}, T)\geq \fDA(\vecc^{\operatorname{static}}, 1) \geq \frac{1}{T}\max_{\vecc \in \mathcal{F}}  \fDA (\vecc, T)$. When $T\leq 2$, we use $\vecc^{\operatorname{static}}$ as the starting inventory, whereas we use $\vecc^{\operatorname{ALG}}$ when $T\geq 3$. Thus, the worst-case performance guarantee is $\min\{\frac{1}{2}-\epsilon, (1-\frac{1}{e})\frac{3}{4}-\epsilon \}\geq 0.474-\epsilon$.

When $T$ is stochastic and we use the {\small DP} algorithm, combining Theorem~\ref{thm:dp-fluid-opt} and the results in Sections~\ref{sec:da-transform} and \ref{sec:da-stochasticT}, we obtain an inventory vector $\vecc^{\operatorname{ALG}}$ such that $\EE_T[\fDA(\vecc^{\operatorname{ALG}},T)]\geq \frac{1}{e}(1-\epsilon)(1-\frac{1}{\mu_T+1})\max_{\vecc \in \mathcal{F}}  \fLP(\vecc,\mu_T)$. Therefore, we have $\EE_T[\fDA(\vecc^{\operatorname{ALG}},T)]\geq (\frac{1}{e}(1-\frac{1}{\mu_T+1})-\epsilon) \max_{\vecc \in \mathcal{F}}  \fLP( \vecc, \mu_T)\geq (\frac{1}{e}(1-\frac{1}{\mu_T+1})-\epsilon) \max_{\vecc \in \mathcal{F}}  \EE_T [\fDA (\vecc, T)]$. Here the last inequality is from Proposition~\ref{prop:ub} in Appendix~\ref{app:prop-ub}, which states that $ \fLP( \vecc, \mu_T)$ is an upper bound for $\EE_T[\fDA (\vecc, T)]$. When $\mu_T< 2.72$, we use $\vecc^{\operatorname{static}}$ as the starting inventory, whereas we use $\vecc^{\operatorname{ALG}}$ when $\mu_T \geq 2.72$: the static factor $\frac{1}{\mu_T+1}$ from the preceding proof exceeds $\frac{1}{3.72}>0.268$ in the former case, and the factor $\frac{1}{e}(1-\frac{1}{\mu_T+1})$ is increasing in $\mu_T$ and satisfies $\frac{1}{e}(1-\frac{1}{3.72})>0.268$ in the latter. Thus, the worst-case performance guarantee is $0.268-\epsilon$.
\Halmos
\endproof

\subsection{Performance Guarantees for DAP}\label{app:5-dap}

\medskip\noindent\textbf{Theorems~\ref{thm:gua-DAP} and \ref{thm:gua-DAP-stoch} (restated; the budget-constraint claims are stated as Theorem~\ref{thm:gua-budget-dap}).} {\itshape
Consider any accuracy level $\epsilon\in (0,1)$. When $T$ is deterministic, \textsc{dap} can be approximated within a factor of $\frac{1}{2}(1-\frac{1}{e})-\epsilon$ under a cardinality constraint, and within a factor of $\frac{1}{3}(1-\frac{1}{e})-\epsilon$ under a budget constraint. When $T$ is stochastic and follows an arbitrary known distribution, \textsc{dap} can be approximated within a factor of $\frac{1}{4}-\epsilon$ under a cardinality constraint and $\frac{1}{6}-\epsilon$ under a budget constraint.
}\medskip
\proof{Proof of Theorems~\ref{thm:gua-DAP}, \ref{thm:gua-DAP-stoch}, and ~\ref{thm:gua-budget-dap}.}
We begin with deterministic $T$ under a cardinality constraint. Let $\vecc^{\operatorname{ALG}}$ denote the output of Algorithm~\ref{alg:tau-m}. Then $\fLPM(\vecc^{\operatorname{ALG}},T)\geq (\frac{1}{2}-\epsilon) \max_{\|\vecc\|_1\leq K}\fLPM( \vecc, T)$ (Theorem~\ref{thm:3-ap-m-ra-T}).
By Lemma~\ref{lem:tran-d}, the expected revenue obtained by the sampling policy (Algorithm~\ref{alg:per-T}), $\fIND(\vecc^{\operatorname{ALG}},T)$, is at least $ (1-\frac{1}{e})\fLPM(\vecc^{\operatorname{ALG}},T)$. Combining these results, and noting that $\fLPM(\vecc,T)$ is an upper bound for the total expected revenue under any online assortment policy \begingroup(Lemma~\ref{lem:multi-ub})\endgroup{}, yields the $(\frac{1}{2}(1-\frac{1}{e})-\epsilon)$ approximation guarantee.

For stochastic $T$ and a cardinality constraint, let $\vecc^{\operatorname{ALG}}$ denote the output of Algorithm~\ref{alg:tau-m-DT}. Then we have $\EE_T[\fLPM (\vecc^{\operatorname{ALG}},T)]\geq (\frac{1}{2}-\epsilon) \max_{\|\vecc\|_1\leq K}  \EE_T[\fLPM (\vecc,T)]$ (Theorem~\ref{thm:3-ap-m-ra-DT}).
By Lemma~\ref{lem:4-DAP-stoc}, the expected revenue obtained by the online greedy policy, $\EE_T[f^{\textsc{Greedy}}(\vecc^{\operatorname{ALG}},T)]$, is at least $\frac{1}{2}\EE_T[\fLPM(\vecc^{\operatorname{ALG}},T)]$.
Multiplying these two approximation ratios, we establish the $(\frac{1}{4}-\epsilon)$ approximation ratio as $\EE_T[\fLPM(\vecc,T)]$ provides an upper bound for the expected revenue of the optimal assortment personalization policy given inventory $\vecc$ \begingroup(Lemma~\ref{lem:multi-ub})\endgroup{}.

When considering the budget constraint, the main difference lies in the first stage. We use Algorithm~\ref{alg:tau-m-bg} for deterministic $T$ and Algorithm~\ref{alg:tau-m-DT-budget} for stochastic $T$. Thus, the first-stage approximation gives a $(\frac{1}{3}-\epsilon)$-approximation guarantee according to Theorems~\ref{thm:3-ap-m-bud-T} and \ref{thm:3-ap-m-bud-stoch-T}. Therefore, the overall performance guarantee is $\frac{1}{3}(1-\frac{1}{e})-\epsilon$ for deterministic $T$ and $\frac{1}{6}-\epsilon$ for stochastic $T$.
\Halmos
\endproof

\begingroup

\section{A Fast Algorithm to Solve the Single-Type CDLP}\label{app:fast_cdlp}
In this section, we introduce a fast algorithm to solve the single-type CDLP, achieving a runtime of $O(n\log n)$. Due to the equivalence between $\fFP(\CFP,T)$ and $\fLP(\vecc,T)$, an optimal solution to the CDLP can be recovered by first computing the revenue-ordered inventory $\CFP$ and then running the Sequencing algorithm (Algorithm~\ref{alg:seq}). By Lemma~\ref{lem:3-rev}, it suffices to identify a threshold index $i_{\tau}$ and a threshold inventory $x_{i_{\tau}}^{\vecc}$ for which the resulting revenue-ordered inventory attains $\fLP(\vecc,T)$.

Recall that we assume products are indexed in descending order of price. Let $\vecc[k]$ denote the inventory vector whose $j$th component equals $c_j$ for $j\leq k$ and $0$ for $j>k$, with $\vecc[0]=\mathbf0^n$. Let $S^k_{\mathrm{term}}$ denote the assortment offered on the final positive-length interval of the fluid process started from $\vecc[k]$. If all inventory depletes strictly before $T$, then $S^k_{\mathrm{term}}=\emptyset$ and $R(\emptyset)=0$. When several revenue-ordered inventories attain $\fLP(\vecc,T)$, we choose one with the smallest support and, among those, the smallest inventory of its threshold product; throughout this appendix we write $\CFP$ for this inventory and $i_{\tau}$ for its threshold index. 

Both this appendix and Appendix~\ref{app:acc} repeatedly analyze the change in fluid revenue from adding inventory of a product whose price is lowest among the stocked products. We will also use the structural properties of the fluid process established in previous parts of the paper. For an inventory vector $\vecx\in\RR_+^n$ and horizon $T>0$, let $q(\vecx)=\ZFP_0(\vecx,T)$ denote the total no-purchase mass of the fluid process \textsc{fp}$(\vecx,T)$. By Lemma~\ref{lem:4-clock}, $\ZFP_j(\vecx,T)=\min\{x_j,v_jq(\vecx)\}$ for every $j\in\NN$, and $q(\vecx)$ is the unique solution of
\begin{equation}\label{eq:fast-fixed-point}
T=q(\vecx)+\sum_{j\in\NN}\min\{x_j,v_jq(\vecx)\}.
\end{equation}
For a product $i$, let $\bar q_i(\vecx,T)$ denote the solution of \eqref{eq:fast-fixed-point} with the term $\min\{x_i,v_iq\}$ replaced by $v_iq$, i.e., the no-purchase mass when product $i$ never depletes; the modified right-hand side is again strictly increasing in $q$, so this solution is also unique. Define
\begin{equation}\label{eq:fast-alpha-bar}
\bar\alpha_i(\vecx,T)=\max\{0,\;v_i\bar q_i(\vecx,T)-x_i\}.
\end{equation}

\begin{lemma}[Terminal-window marginal]\label{lem:fast-terminal-window}
Let $\vecx\in\RR_+^n$, let $T>0$, and let $i$ be a product with $r_i\leq r_j$ for every $j\in N(\vecx)$. Let $S(t)$ denote the assortment at time $t$ in \textsc{fp}$(\vecx,T)$. Then:
\begin{enumerate}[label=(\roman*)]
\item $\ZFP_i(\vecx+\delta\vece_i,T)=\min\{x_i+\delta,\;v_i\bar q_i(\vecx,T)\}$ for every $\delta\geq0$. Let
\(
\bar\alpha_i(\vecx,T)
:=\max\{0,\;v_i\bar q_i(\vecx,T)-x_i\}.
\)
Then
\(
\bar\alpha_i(\vecx,T)
=\max_{\delta\geq0}
\bigl\{
\ZFP_i(\vecx+\delta\vece_i,T)-\ZFP_i(\vecx,T)
\bigr\}.
\)
Moreover, if $\bar\alpha_i(\vecx,T)>0$, then for every
$\alpha\in[0,\bar\alpha_i(\vecx,T)]$, the process
\textsc{fp}$(\vecx+\alpha\vece_i,T)$ consumes all
$x_i+\alpha$ units of product $i$.
\item The set $\{t\in[0,T]:r_i>R(S(t))\}$ is a terminal window of $[0,T]$; let $\ell$ denote its duration. For every $\bar c$ with $0<\bar c\leq\bar\alpha_i(\vecx,T)$, the smallest maximizer of $\alpha\mapsto\fFP(\vecx+\alpha\vece_i,T)$ over $[0,\bar c]$ is $\min\{\ell,\bar c\}$.
\end{enumerate}
\end{lemma}
\proof{Proof.}
(i) Fix $\delta\geq0$ and write $\bar q=\bar q_i(\vecx,T)$. If $x_i+\delta\geq v_i\bar q$, then $\min\{x_i+\delta,v_i\bar q\}=v_i\bar q$, so $\bar q$ solves \eqref{eq:fast-fixed-point} for the inventory $\vecx+\delta\vece_i$; by uniqueness, $q(\vecx+\delta\vece_i)=\bar q$ and $\ZFP_i(\vecx+\delta\vece_i,T)=\min\{x_i+\delta,v_i\bar q\}=v_i\bar q$. If $x_i+\delta<v_i\bar q$, then, since $q(\vecx+\delta\vece_i)\geq\bar q$, we have $v_iq(\vecx+\delta\vece_i)\geq v_i\bar q>x_i+\delta$ and $\ZFP_i(\vecx+\delta\vece_i,T)=x_i+\delta$. In both cases $\ZFP_i(\vecx+\delta\vece_i,T)=\min\{x_i+\delta,v_i\bar q\}$. Applying this identity also at $\delta=0$ gives
\(
\ZFP_i(\vecx+\delta\vece_i,T)-\ZFP_i(\vecx,T)=\min\{x_i+\delta,v_i\bar q\}
  -\min\{x_i,v_i\bar q\}=\min\{\delta,\bar\alpha_i(\vecx,T)\}.
\)
The maximum additional consumption is therefore
$\bar\alpha_i(\vecx,T)$, attained whenever
$\delta\geq\bar\alpha_i(\vecx,T)$.
Finally, if $\bar\alpha_i(\vecx,T)>0$, then
$x_i<v_i\bar q$, and every
$\alpha\in[0,\bar\alpha_i(\vecx,T)]$ satisfies
$x_i+\alpha\leq v_i\bar q$.
The consumption formula therefore implies that all
$x_i+\alpha$ units are consumed.

(ii) By Lemma~\ref{lem:seq-rev1}(i), if $r_i>R(S(t))$, then $R(S(t'))\leq R(S(t))<r_i$ for every $t'\geq t$; hence $\{t:r_i>R(S(t))\}$ is a terminal window. Fix $0\leq\alpha<\alpha'\leq\bar c$; then $x_i<v_i\bar q_i(\vecx,T)$, and (i) gives $\ZFP_i(\vecx+\alpha\vece_i,T)-\ZFP_i(\vecx,T)=\alpha$ and likewise for $\alpha'$. The Separability Lemma therefore gives $\fFP(\vecx+\alpha\vece_i,T)=\fFP(\vecx,T-\alpha)+\alpha r_i$ and the analogous identity for $\alpha'$. Since the fluid process accrues revenue at rate $R(S(t))$, subtracting the two identities gives
\[
\fFP(\vecx+\alpha'\vece_i,T)-\fFP(\vecx+\alpha\vece_i,T)=\int_{T-\alpha'}^{T-\alpha}\bigl(r_i-R(S(t))\bigr)\,\mathrm{d}t,
\]
which is positive when $\alpha'\leq\ell$ and nonpositive when $\alpha\geq\ell$. Thus $\fFP(\vecx+\alpha\vece_i,T)$ is strictly increasing in $\alpha$ on $[0,\min\{\ell,\bar c\}]$ and non-increasing on $[\min\{\ell,\bar c\},\bar c]$, which proves (ii).
\Halmos
\endproof

\begin{lemma}[Threshold-index characterization]\label{lem:5-str}
For $k\leq i_{\tau}$,
\(
r_k>R(S^{k-1}_{\mathrm{term}}),
\)
whereas, for $k>i_{\tau}$,
\(
r_k\leq R(S^{k-1}_{\mathrm{term}}).
\)
Consequently, $i_{\tau}$ is the largest index $k$ such that
$r_k>R(S^{k-1}_{\mathrm{term}})$.
\end{lemma}
\endgroup
\proof{Proof of Lemma~\ref{lem:5-str}.}
\begingroup
First consider $k=i_{\tau}$. By the Separability Lemma, adding the threshold inventory to $\vecc[i_{\tau}-1]$ changes the fluid revenue by
\(
\int_{T-x_{i_{\tau}}^{\vecc}}^T
\bigl(r_{i_{\tau}}-R(S^{i_{\tau}-1}(t))\bigr)\;\mathrm{d}t.
\)
Here the added consumption equals $x_{i_{\tau}}^{\vecc}$ because the process started from $\CFP$ consumes all $x_{i_{\tau}}^{\vecc}$ units of product $i_{\tau}$: if part of that inventory were left unsold, removing it would change neither the assortment sequence nor the fluid revenue, contradicting our choice of the smallest threshold inventory.
Because every product in $\vecc[i_{\tau}-1]$ has price at least $r_{i_{\tau}}$, Lemma~\ref{lem:fast-terminal-window}(ii) shows that $\{t:r_{i_{\tau}}>R(S^{i_{\tau}-1}(t))\}$ is a terminal window of $[0,T]$. Suppose
$r_{i_{\tau}}\leq R(S^{i_{\tau}-1}_{\mathrm{term}})$. Then this window is empty, so the integrand is nonpositive throughout $[0,T]$. Stocking any positive amount of product $i_{\tau}$ therefore cannot increase the fluid revenue. Since $\CFP$ is optimal, removing its positive inventory of product $i_{\tau}$ leaves $\vecc[i_{\tau}-1]$ with the same optimal value, contradicting the choice of $\CFP$ with the smallest support. Hence
$r_{i_{\tau}}>R(S^{i_{\tau}-1}_{\mathrm{term}})$.

Now let $k<i_{\tau}$. Since $\vecc[k-1]\preceq\vecc[i_{\tau}-1]$, Lemma~\ref{lem:app-subset} gives
$S^{k-1}_{\mathrm{term}}\subseteq S^{i_{\tau}-1}_{\mathrm{term}}$. Every product in the difference has price at least $r_{i_{\tau}}>R(S^{i_{\tau}-1}_{\mathrm{term}})$, so Lemma~\ref{lem:subsetie} gives
\(
R(S^{k-1}_{\mathrm{term}})
\leq R(S^{i_{\tau}-1}_{\mathrm{term}})<r_{i_{\tau}}\leq r_k.
\)
This proves the first statement.

For the second statement, suppose first that $i_{\tau}<n$. Optimality of $\CFP$ implies
\begin{equation}\label{eq:fast-next-product}
r_{i_{\tau}+1}\leq R(S^{\CFP}_{\mathrm{term}}),
\end{equation}
where $S^{\CFP}_{\mathrm{term}}$ is the terminal assortment of the process started from $\CFP$. Otherwise, Lemma~\ref{lem:fast-terminal-window}(ii), applied to $\CFP$ and product $i_{\tau}+1$, shows that adding a suitable positive amount of product $i_{\tau}+1$ strictly increases the fluid revenue, contradicting the optimality of $\CFP$.

Fix $k>i_{\tau}$. Since $\CFP\preceq\vecc[k-1]$, Lemma~\ref{lem:app-subset} gives
$S^{\CFP}_{\mathrm{term}}\subseteq S^{k-1}_{\mathrm{term}}$. The revenue
$R(S^{k-1}_{\mathrm{term}})$ is a weighted average of
$R(S^{\CFP}_{\mathrm{term}})$ and the prices of the products in
$S^{k-1}_{\mathrm{term}}\setminus S^{\CFP}_{\mathrm{term}}$. By
\eqref{eq:fast-next-product}, the first quantity is at least
$r_{i_{\tau}+1}\geq r_k$, and every product in the difference has index at most $k-1$ and hence price at least $r_k$. Therefore
$R(S^{k-1}_{\mathrm{term}})\geq r_k$. The case $i_{\tau}=n$ is immediate because there is no $k>i_{\tau}$.
\Halmos
\endgroup
\endproof

\begingroup
It remains to determine the threshold inventory. For a product $i$, run the Sequencing algorithm on $\vecc[i-1]$, and let $S_1,\ldots,S_L$ and $t_1,\ldots,t_L$ be the resulting assortments and durations. Define
\begin{equation}\label{eq:fast-positive-window}
\ell_i=\sum_{q=1}^L t_q\1\{r_i>R(S_q)\}.
\end{equation}
Also define $
z_i=\ZFP_i(\vecc[i],T)$.

\begin{lemma}[Threshold inventory]\label{lem:fast-threshold-level}
For the threshold index $i_{\tau}$, $x_{i_{\tau}}^{\vecc}=\min\{\ell_{i_{\tau}},z_{i_{\tau}}\}$. 
\end{lemma}
\endgroup
\proof{Proof of Lemma~\ref{lem:fast-threshold-level}.}
\begingroup
Fix a product $i$ and apply Lemma~\ref{lem:fast-terminal-window} with $\vecx=\vecc[i-1]$ and horizon $T$; every product in $\vecc[i-1]$ has price at least $r_i$, and $\vecc[i-1]$ has no inventory of product $i$. Since $\vecc[i]=\vecc[i-1]+c_i\vece_i$, part (i) gives
\(
z_i=\ZFP_i(\vecc[i],T)=\min\{c_i,\;v_i\bar q_i(\vecc[i-1],T)\}=\min\{c_i,\;\bar\alpha_i(\vecc[i-1],T)\}>0.\) In particular, every $\alpha\in[0,z_i]$ is fully consumed. If $z_i<c_i$, then $z_i=\bar\alpha_i(\vecc[i-1],T)$, and, by Lemma~\ref{lem:fast-terminal-window}(i), every initial level $\alpha\in[z_i,c_i]$ yields the same no-purchase mass $\bar q_i(\vecc[i-1],T)$, hence the same consumption of every product and the same fluid revenue. Consequently, under our tie-breaking rule, the smallest optimal inventory of product $i$ lies in $[0,z_i]$. By \eqref{eq:fast-positive-window}, the terminal window of Lemma~\ref{lem:fast-terminal-window}(ii) for $\vecx=\vecc[i-1]$ has duration $\ell_i$. Applying part (ii) with $\bar c=z_i$, the smallest maximizer of $\alpha\mapsto\fFP(\vecc[i-1]+\alpha\vece_i,T)$ over $[0,z_i]$ is $\min\{\ell_i,z_i\}$. Taking $i=i_{\tau}$ and applying the tie-breaking rule finishes the proof.
\Halmos
\endgroup
\endproof

\begingroup
By Lemma~\ref{lem:5-str}, a binary search for the largest index $k$ satisfying $r_k>R(S^{k-1}_{\mathrm{term}})$ returns $i_{\tau}$. We use the following $O(n)$ implementation of each Sequencing call. By Lemma~\ref{lem:4-clock}, products deplete in ascending order of the ratios $c_j/v_j$, and, while product $j$ is in stock,
\begin{equation}\label{eq:fast-common-rate}
\frac{\mathrm{d}}{\mathrm{d}t}\frac{c_j(t)}{v_j}
=-\frac{1}{1+\sum_{h\in S(t)}v_h},
\end{equation}
which is the same for every product in stock. Sorting these ratios once allows a prefix sequence to be computed in $O(n)$ time by scanning the common order and maintaining the sums $\sum_{j\in S(t)}v_j$ and $\sum_{j\in S(t)}r_jv_j$.

\begin{algorithm}[H]
\SetAlgoNoLine
\caption{Fast Revenue-Ordered Solution of the Single-Type CDLP}\label{alg_fast}
\KwIn{Products $\NN$ sorted in descending order of price, inventory $\vecc$, attraction weights $\vecv$, and horizon $T$}
\KwOut{Revenue-ordered inventory $\CFP$ and an optimal CDLP solution $y$}
Compute the ascending order of the ratios $c_i/v_i$ and store it as the common ratio order, keeping the product indices in descending price order\;
Initialize $a_l=1$, $a_r=n$\;
\While{$a_l<a_r$}{
  $a_d=a_l+\lceil(a_r-a_l)/2\rceil$\;
  Run Sequencing on $\vecc[a_d-1]$ and let $R_{\mathrm{term}}$ be the revenue of its terminal assortment\;
  \lIf{$r_{a_d}>R_{\mathrm{term}}$}{$a_l=a_d$}
  \lElse{$a_r=a_d-1$}
}
$i_{\tau}=a_l$; run Sequencing on $\vecc[i_{\tau}-1]$ and on $\vecc[i_{\tau}]$, and compute $\ell_{i_{\tau}}$ and $z_{i_{\tau}}$, and $\theta=\min\{\ell_{i_{\tau}},z_{i_{\tau}}\}$\;
Initialize $\CFP=\mathbf0^n$; set $x_i^{\vecc}=c_i$ for $i<i_{\tau}$ and set $x_{i_{\tau}}^{\vecc}=\theta$\;
Scan the common ratio order, retaining only products $i<i_{\tau}$; merge $x_{i_{\tau}}^{\vecc}/v_{i_{\tau}}$ into this order and run Sequencing on $\CFP$, obtaining assortments $S_1,\ldots,S_L$ and durations $t_1,\ldots,t_L$\;
Set $y(S_q)=t_q/T$ for $q\in[L]$; store $y$ as this sparse list of at most $n+1$ positive weights, with $y(S)=0$ implicitly for every other $S\subseteq\NN$\;
\Return{$\CFP,y$}\;
\end{algorithm}

The upper midpoint ensures that the search interval of the threshold binary search strictly decreases in every iteration. Since
$R(S^0_{\mathrm{term}})=0<r_1$, the search inequality holds at $k=1$; if it holds for every product, the search correctly returns $i_{\tau}=n$.

\begin{theorem}\label{thm:fast-cdlp}
Algorithm~\ref{alg_fast} returns an optimal solution to the single-type CDLP \eqref{cdlp} in $O(n\log n)$ time.
\end{theorem}
\endgroup
\proof{Proof.}
\begingroup
Lemmas~\ref{lem:5-str} and~\ref{lem:fast-threshold-level} show that Algorithm~\ref{alg_fast} constructs the revenue-ordered inventory $\CFP$ selected above. Hence,
\(
\fFP(\CFP,T)=\fLP(\vecc,T).
\)
The final Sequencing call returns durations that sum to $T$. Therefore, the weights $y(S_q)=t_q/T$ sum to one, and the fluid consumption of each product is at most $x_i^{\vecc}\leq c_i$. The returned mixture is feasible for \eqref{cdlp}, and its objective value is
\(
T\sum_{S\subseteq\NN}y(S)R(S)
=\sum_{q=1}^L t_qR(S_q)
=\fFP(\CFP,T)
=\fLP(\vecc,T).
\)
Thus, the returned CDLP solution is optimal.

The initial ratio sort takes $O(n\log n)$ time. The threshold binary search makes $O(\log n)$ predicate queries, followed by two further Sequencing calls, each of which takes $O(n)$ time using \eqref{eq:fast-common-rate}. Merging the threshold ratio, running the final sequence, and constructing the sparse CDLP solution take $O(n)$ time. Therefore, the total runtime is $O(n\log n)$ and the space requirement is $O(n)$.
\Halmos
\endgroup
\endproof

\section{Acceleration}\label{app:acc}
In this section, we introduce an acceleration technique for approximately maximizing the multi-type CDLP under a cardinality constraint. \begingroup The accelerated algorithm adds units of a product in batches rather than one unit at a time. We show that the batch size can be found by binary search (Lemma~\ref{lem:bin}), that at most one batch of each product is accepted (Remark~\ref{rem:one-batch}), and that batching preserves the $(\frac{1}{2}-\epsilon)$ guarantee of Theorem~\ref{thm:3-ap-m-ra-T} (Proposition~\ref{prop:acc-dap}). The same marginal-benefit calculation can also be used in the budget-constrained setting; here, we focus on the cardinality-constrained algorithm.

\begin{algorithm}[htp]
\SetAlgoNoLine
\caption{Accelerated Cardinality-Constrained Threshold Add for DAP}
\label{alg:cdlp-multi-type-acc}
\KwIn{Cardinality $K$, threshold $\tau\geq0$, products $\NN$ sorted in descending order of price: $r_1\geq r_2\geq\cdots\geq r_n$, number of customers $T$}
\KwOut{Inventory vector $\vecc$, inventory allocation vector $\matc$}
Initialize $\vecc=\mathbf0^n$, $\matc=\mathbf0^{mn}$\;
\For{$i\in\NN$}{
    Get the largest batch size of product $i$ whose average marginal benefit is at least $\tau$, and its allocation:
    $(z_i,\matz^*)=$ Get Number$(K,\tau,\matc,T,i)$\;
        \lIf{$\|\vecc\|_1\ge K$}{\textbf{break}}
    \lIf{$z_i=0$}{\textbf{continue}}
    Add $z_i$ units of product $i$ and update:
    $\vecc\leftarrow\vecc+z_i\vece_i$, $\matc\leftarrow\matc+\matz^*$\;
}
\Return{$\vecc,\matc$}\;
\end{algorithm}
\endgroup
\vspace{-10pt}
Compared with Threshold Add (Algorithm~\ref{alg:cdlp-multi-type}), the accelerated algorithm uses Get Number to determine the largest batch size $z_i$ whose average marginal benefit is at least $\tau$. 
Products are still scanned once in descending-price order, and all accepted units of one product remain contiguous.

\begingroup
For a current allocation $\matc$ and product $i$, define
\begin{equation}\label{eq:acc-G}
G_i(z\mid\matc)=
\max_{\substack{\matz\geq0:\ \sum_{k\in\MM}Z_{k,i}=z,\\
Z_{k,j}=0\ \forall k\in\MM,\ j\neq i}}
\bigl\{\hat f(\matc+\matz)-\hat f(\matc)\bigr\},
\qquad z\geq0.
\end{equation}
The dependence on the fixed horizon $T$ is suppressed. Because the algorithm maintains $\sum_{k\in\MM}C_{k,i}=c_i$ for every product $i$, we have $\|\matc\|_1=\|\vecc\|_1$.

\begin{algorithm}[htbp]
\SetAlgoNoLine
\caption{Get Number}
\label{alg:cdlp-multi-type-bin}
\KwIn{Cardinality $K$, threshold $\tau$, current inventory allocation vector $\matc$, number of customers $T$, product $i$}
\KwOut{Number of units to be added $z^*$, inventory allocation vector $\matz^*$}
Initialize $z_l=0$, $z_r=K-\|\matc\|_1$\;
\While{$z_l<z_r$}{
    $z_{\mathrm{mid}}=\lceil(z_l+z_r)/2\rceil$\;
    $(\Delta f,\matz)=$ Calculate Benefit$(T,\matc,z_{\mathrm{mid}},i)$\;
    \lIf{$\Delta f\geq z_{\mathrm{mid}}\tau$}{$z_l\leftarrow z_{\mathrm{mid}}$}
    \lElse{$z_r\leftarrow z_{\mathrm{mid}}-1$}
}
$z^*\leftarrow z_l$; $(\Delta f,\matz^*)=$ Calculate Benefit$(T,\matc,z^*,i)$\;
\Return{$z^*,\matz^*$}\;
\end{algorithm}
\endgroup
Get Number uses binary search to find the largest integer $z^*\leq K-\|\matc\|_1$ satisfying $G_i(z^*\mid\matc)\geq z^*\tau$. The upper bound is the remaining cardinality, so no additional demand bound is needed. The rounded-up midpoint ensures termination, as for the threshold search in Algorithm~\ref{alg_fast}. Calculate Benefit returns both $G_i(z\mid\matc)$ and an allocation attaining it; the final allocation is recomputed for $z^*$.
\begingroup
For type $k$, define the single-type marginal curve
\(
h_{ki}(x\mid\matc_k)
=\fLP_k(\matc_k+x\vece_i,\Lambda_k)
-\fLP_k(\matc_k,\Lambda_k), x\geq0.
\)
By separability across customer types,
\begin{equation}\label{eq:acc-resource-allocation}
G_i(z\mid\matc)=
\max_{\substack{x_k\geq0,\ k\in\MM:\ 
\sum_{k\in\MM}x_k=z}}
\sum_{k\in\MM}h_{ki}(x_k\mid\matc_k).
\end{equation}
\begin{lemma}\label{lem:bin}
For every product $i$, $G_i(z\mid\matc)$ is non-decreasing and concave in $z\geq0$. Consequently, the average marginal benefit $G_i(z\mid\matc)/z$ is non-increasing for $z>0$. If $G_i(1\mid\matc)<\tau$, then $G_i(d\mid\matc)<d\tau$ for every integer $d\geq1$.
\end{lemma}
\proof{Proof of Lemma~\ref{lem:bin}.}
For fixed $k$ and $i$, the value of the single-type SBLP is non-decreasing and concave as a function of the capacity $C_{k,i}+x$, as observed for the sales-based LP in the proof of Theorem~\ref{thm:sof}. Hence $h_{ki}(x\mid\matc_k)$ is non-decreasing and concave, with $h_{ki}(0\mid\matc_k)=0$. Convex combinations of feasible allocations for two batch sizes are feasible for the corresponding convex combination of those batch sizes, so concavity of every $h_{ki}$ gives concavity of $G_i$ in~\eqref{eq:acc-resource-allocation}. Monotonicity follows because unused capacity does not reduce an LP value. Finally, for $0<x<y$, concavity and $G_i(0\mid\matc)=0$ imply
\(
G_i(x\mid\matc)
\geq \frac{x}{y}G_i(y\mid\matc),
\)
so $G_i(x\mid\matc)/x\geq G_i(y\mid\matc)/y$. The remaining claims follow immediately.
\Halmos
\endgroup
\endproof

\begingroup
We next describe how Calculate Benefit evaluates~\eqref{eq:acc-resource-allocation}. At the time product $i$ is processed, every product with positive allocation has index at most $i$, and hence price at least $r_i$, so the condition of the lemmas below is satisfied. For each type $k$, let $\vecx^k\preceq\matc_k$ be the generalized revenue-ordered inventory from Lemma~\ref{lem:3-rev} that represents $\fLP_k(\matc_k,\Lambda_k)$. Use Sequencing (Algorithm~\ref{alg:seq}) to obtain the assortment sequence $S_{k1},\ldots,S_{kL_k}$ and durations $t_{k1},\ldots,t_{kL_k}$ of the corresponding fluid process. Thus, $\sum_{q=1}^{L_k}t_{kq}=\Lambda_k$. Applying the terminal-window results of Appendix~\ref{app:fast_cdlp} to the type-$k$ MNL parameters, define
\begin{equation}\label{eq:acc-beta}
\beta_{ki}
=\bar\alpha_i(\vecx^k,\Lambda_k),
\end{equation}
the maximum additional fluid consumption of product $i$ attainable by raising its inventory (Lemma~\ref{lem:fast-terminal-window}(i)). In particular, $0\leq\beta_{ki}\leq\Lambda_k$.

The next lemma extends the structural properties of Lemma~\ref{lem:3-yc1} to real-valued capacities and increments, both of which arise here because the allocation vectors maintained by the algorithm may be fractional. 
The proof of Lemma~\ref{lem:3-yc1} uses neither the integrality of the inventory vector nor the unit size of the increment, so the results hold for arbitrary real inventory vectors. So we get the following corollary directly.
\begin{corollary}\label{lem:acc-yc1}
Let $\vecc\in\RR_+^n$ and $x>0$, and suppose product $i$ has the highest index among all products in $N(\vecc)\cup\{i\}$. If $\fLP(\vecc+x\vece_i,T)-\fLP(\vecc,T)>0$, then $\vecx^{\vecc}=\vecc$, and $\max_{\vecx\preceq\vecc+x\vece_i}\fFP(\vecx,T)$ is attained at $\vecc+\alpha\vece_i$ for some $0<\alpha\leq x$ such that product $i$ is fully consumed in \textsc{fp}$(\vecc+\alpha\vece_i,T)$. Moreover, $r_i-R(S(t))\geq0$ for all $t\in[T-\alpha,T]$, where $S(t)$ denotes the assortment at time $t$ in $\fFP(\vecc,T)$.
\end{corollary}

The terminal-window lemma (Lemma~\ref{lem:fast-terminal-window}) and  Lemma~\ref{lem:acc-yc1} give the following representation.
\begin{lemma}[Marginal curve]\label{lem:acc-curve} 
Suppose product $i$ has the highest index among all products in $N(\matc_k)\cup\{i\}$. Then, for every $x\geq0$,
\begin{equation}\label{eq:acc-single-type-curve}
h_{ki}(x\mid\matc_k)
=\max_{0\leq\alpha\leq\min\{\beta_{ki},x\}}
\int_{\Lambda_k-\alpha}^{\Lambda_k}
\bigl(r_i-R_k(S_k(t))\bigr)\;\mathrm{d}t.
\end{equation}
\end{lemma}
\proof{Proof of Lemma~\ref{lem:acc-curve}.} We first show LHS$\ge$RHS. Fix $\alpha\in[0,\min\{\beta_{ki},x\}]$. Product $i$ has the lowest price among the products stocked in $\vecx^k$, and, since $\alpha\leq\beta_{ki}$, the consumption formula of Lemma~\ref{lem:fast-terminal-window}(i), applied to $\vecx^k$ with horizon $\Lambda_k$, gives $\ZFP_i(\vecx^k+\alpha\vece_i,\Lambda_k)-\ZFP_i(\vecx^k,\Lambda_k)=\alpha$. The Separability Lemma therefore gives $\fFP(\vecx^k+\alpha\vece_i,\Lambda_k)=\fFP(\vecx^k,\Lambda_k-\alpha)+\alpha r_i$. Since the fluid process accrues revenue at rate $R_k(S_k(t))$, this identity, together with $\fFP(\vecx^k,\Lambda_k)=\fLP_k(\matc_k,\Lambda_k)$ and, by Lemma~\ref{lem:3-equ1}, $\fLP_k(\matc_k+x\vece_i,\Lambda_k)\geq\fFP(\vecx^k+\alpha\vece_i,\Lambda_k)$ for $\vecx^k+\alpha\vece_i\preceq\matc_k+x\vece_i$, gives
\[
\fLP_k(\matc_k+x\vece_i,\Lambda_k)\geq \fFP(\vecx^k+\alpha\vece_i,\Lambda_k)=\fLP_k(\matc_k,\Lambda_k)+\int_{\Lambda_k-\alpha}^{\Lambda_k}\bigl(r_i-R_k(S_k(t))\bigr)\;\mathrm{d}t.
\]
Subtracting $\fLP_k(\matc_k,\Lambda_k)$ and maximizing over $\alpha$ gives the lower bound.

Then we show LHS$\leq$RHS. If $h_{ki}(x\mid\matc_k)=0$, the claim holds because the right-hand side of \eqref{eq:acc-single-type-curve} is nonnegative (take $\alpha=0$). If $h_{ki}(x\mid\matc_k)>0$, Lemma~\ref{lem:acc-yc1}, applied with inventory $\matc_k$, increment $x$, and horizon $\Lambda_k$, shows that the generalized revenue-ordered inventory representing $\fLP_k(\matc_k,\Lambda_k)$ is $\matc_k$ itself and that the one representing $\fLP_k(\matc_k+x\vece_i,\Lambda_k)$ may be taken to be $\matc_k+\alpha\vece_i$ for some $0<\alpha\leq x$ with product $i$ fully consumed, so that $\alpha\leq\beta_{ki}$ by the definition \eqref{eq:acc-beta} of $\beta_{ki}$. The Separability Lemma then writes $h_{ki}(x\mid\matc_k)$ as the integral in \eqref{eq:acc-single-type-curve} evaluated at this $\alpha$, which is at most the maximum.
\Halmos
\endproof

For an assortment interval $q$, its overlap with the terminal window $[\Lambda_k-\beta_{ki},\Lambda_k]$ is
\begin{equation}\label{eq:acc-overlap}
t'_{kq}
=\max\!\left\{0,
\min\!\left\{t_{kq},
\beta_{ki}-\sum_{h=q+1}^{L_k}t_{kh}
\right\}\right\}.
\end{equation}
By Lemma~\ref{lem:seq-rev1}(i), $R_k(S_k(t))$ is non-increasing in $t$ on the set
\(
\{t\in[0,\Lambda_k]:R_k(S_k(t))\leq r_i\}.
\) Consequently, the positive slopes of $h_{ki}$ are encountered in non-increasing order as the added quantity $x$ grows: an increment uses the latest part of the fluid process first and then moves backward through the sequence.

Across customer types, problem~\eqref{eq:acc-resource-allocation} is therefore a fractional knapsack problem. Each candidate interval $(k,q)$ has weight $t'_{kq}$ and value density $r_i-R_k(S_{kq})$. Sorting intervals in increasing order of $R_k(S_{kq})$ is the same as sorting them in decreasing order of value density. If two intervals of one type have the same revenue, the later interval is placed first so that the terminal-window precedence is respected. Segments with non-positive density are omitted. If $z$ exceeds the total weight of the positive-density segments, the remaining capacity can be allocated arbitrarily and has zero marginal value, consistent with the maximization in~\eqref{eq:acc-single-type-curve}.

\begin{algorithm}[htbp]
\SetAlgoNoLine
\caption{Calculate Benefit}
\label{alg:calculate-benefit}
\KwIn{Number of customers $T$, current inventory allocation vector $\matc$, batch size $z$, product $i$}
\KwOut{Marginal benefit $y=G_i(z\mid\matc)$, maximizing allocation $\matz^*$}
Initialize $\mathcal D=\emptyset$, $y=0$, $\matz^*=\mathbf0^{mn}$\;
\For{$k\in\MM$}{
    Construct the generalized revenue-ordered inventory $\vecx^k$ representing $\fLP_k(\matc_k,\Lambda_k)$\;
    Use Sequencing to obtain $S_{k1},\ldots,S_{kL_k}$ and $t_{k1},\ldots,t_{kL_k}$ for $\vecx^k$ over $[0,\Lambda_k]$; Compute $\beta_{ki}$ from~\eqref{eq:acc-beta}\;
    \For{$q=1,\ldots,L_k$}{
        Compute $t'_{kq}$ from~\eqref{eq:acc-overlap}\;
        \lIf{$t'_{kq}>0$ and $r_i>R_k(S_{kq})$}{add $(k,q,t'_{kq},r_i-R_k(S_{kq}))$ to $\mathcal D$}
    }
}
Sort $\mathcal D$ in decreasing order of density; Set $s\leftarrow z$\;
\ForEach{$(k,q,t'_{kq},d_{kq})\in\mathcal D$ in the sorted order}{
    \If{$s=0$}{\textbf{break}}
    $a\leftarrow\min\{s,t'_{kq}\}$, $y\leftarrow y+a d_{kq}$, $Z^*_{k,i}\leftarrow Z^*_{k,i}+a$, $s\leftarrow s-a$\;
}
\lIf{$s>0$}{allocate the remaining $s$ units arbitrarily in column $i$ of $\matz^*$}
\Return{$y,\matz^*$}\;
\end{algorithm}
\vspace{-10pt}
\begin{remark}\label{rem:one-batch}
The inner loop of Algorithm~\ref{alg:cdlp-multi-type-acc} calls Get Number once for each product. Suppose Get Number returns $z^*\geq1$ at allocation $\matc$, with maximizing batch allocation $\matz^*$, and suppose the batch does not fill the cardinality bound. Then $z^*<K-\|\matc\|_1$, so $z^*+1$ was rejected by the binary search:
\(
G_i(z^*+1\mid\matc)<(z^*+1)\tau.
\)
Let $\widetilde{\matz}$ be an optimal one-unit allocation at the updated state $\matc+\matz^*$. The combined allocation $\matz^*+\widetilde{\matz}$ is feasible in the definition of $G_i(z^*+1\mid\matc)$, and therefore
\(
G_i(z^*\mid\matc)
+G_i(1\mid\matc+\matz^*)
\leq G_i(z^*+1\mid\matc)
<(z^*+1)\tau.
\)
Since $G_i(z^*\mid\matc)\geq z^*\tau$, it follows that $G_i(1\mid\matc+\matz^*)<\tau$. 
\end{remark}

\begin{proposition}[Guarantee preservation under batching]
\label{prop:acc-dap}
Replacing Algorithm~\ref{alg:cdlp-multi-type} by Algorithm~\ref{alg:cdlp-multi-type-acc} inside Algorithm~\ref{alg:tau-m} preserves the $(\frac{1}{2}-\epsilon)$ approximation guarantee of Theorem~\ref{thm:3-ap-m-ra-T}.
\end{proposition}
\endgroup
\proof{Proof.}
\begingroup
Fix one threshold $\tau$. Whenever Get Number accepts a batch of $z$ units, its allocation satisfies
\[
\hat f(\matc+\matz)-\hat f(\matc)
=G_i(z\mid\matc)\geq z\tau.
\]
Thus, if the accelerated algorithm stocks $K$ units, telescoping over its accepted batches gives $\hat f(\matc)\geq K\tau$, which is exactly the full-cardinality inequality used in the proof of Theorem~\ref{thm:3-ap-m-ra-T}.

Now suppose the ordered scan ends with $\|\matc\|_1<K$. Let $j_1<\cdots<j_\kappa$ be the products to which the algorithm assigns positive inventory, and define $Q_l$, $\matc(l)$, and $\matc_l^{\operatorname{OPT}}$ as in the proof of Theorem~\ref{thm:3-ap-m-ra-T}. For every $q\in Q_l$, the allocation immediately after product $q$ is processed is $\matc(l-1)$; the argument is the one in the proof of Lemma~\ref{lem:3-m-tau}, with the accepted copies of product $j_{l-1}$ replaced by its accepted batch.

If no batch is accepted, batch size one fails the threshold test. If a positive batch is accepted without filling the cardinality bound, Remark \ref{rem:one-batch} implies that the next one-unit marginal at the updated allocation is below \(\tau\). Therefore, $G_q(1\mid\matc(l-1))<\tau$ for every $q\in Q_l$. Then the proof of Lemma~\ref{lem:3-m-tau} applies unchanged and gives $\hat f(\matc_l^{\operatorname{OPT}}\mid\matc(l-1))\leq\|\matc_l^{\operatorname{OPT}}\|_1\tau$ for every $l\in[\kappa+1]$. Then the proof of Theorem~\ref{thm:3-ap-m-ra-T} applies and the accelerated algorithm has the same $(\frac{1}{2}-\epsilon)$ guarantee.
\Halmos
\endgroup
\endproof
\begingroup
\begin{remark}[Stochastic horizons]
\label{rem:acc-dap-stoch}
The same batching argument applies to Algorithm~\ref{alg:tau-m-DT}. For batch size $z$, define the distribution-aware marginal by
\(
G_i^{D_T}(z\mid(\matc^T)_{T=1}^{T_{\operatorname{max}}})
=\EE_T\!\left[G_i(z\mid\matc^T,T)\right],
\)
where the suppressed horizon argument has been restored. The inventory increment $z$ is common across horizons, but its allocation among customer types is optimized separately for every realized $T$, as in Algorithm~\ref{alg:cdlp-multi-type-DT}. Each ratio $G_i(z\mid\matc^T,T)/z$ is non-increasing by Lemma~\ref{lem:bin}, so their nonnegative weighted average is also non-increasing. Accepted batches contribute at least $z\tau$ in expectation. If no batch is accepted, batch size one fails the threshold test. If a positive batch is accepted without filling the cardinality bound, Remark \ref{rem:one-batch} implies that the next one-unit marginal at the updated allocation is below \(\tau\). Then by Lemma~\ref{thm:3-ap-m-ra-DT} and Lemma~\ref{lem:3-m-tau-stochT}, the accelerated stochastic algorithm preserves the $(\frac{1}{2}-\epsilon)$ guarantee of Theorem~\ref{thm:3-ap-m-ra-DT}.
\end{remark}
\endgroup

\textbf{Time Complexity Analysis.} \begingroup
For fixed $\matc$ and $i$, there are at most $L=\sum_{k\in\MM}L_k\leq m(n+1)$ assortment intervals. Under the MNL model, a type-$k$ fluid sequence can be constructed in $O(n\log n)$ time by the ratio-order implementation of Sequencing in Appendix~\ref{app:fast_cdlp}: one sort, then \eqref{eq:fast-common-rate} and $O(n)$ incremental updates. Products with $v_{kj}=0$ are never purchased by type $k$ and contribute nothing to any revenue rate; we adopt the convention $x_j^k/v_{kj}=+\infty$ and omit these products from the type-$k$ sequence. The generalized revenue-ordered inventory can be computed in $O(n\log n)$ time per type by the fast CDLP algorithm of Appendix~\ref{app:fast_cdlp} (Theorem~\ref{thm:fast-cdlp}) followed by one $O(n)$ pass that applies the construction of Lemma~\ref{lem:3-rev} to the optimal SBLP solution it induces, and $\beta_{ki}$ by one additional $O(n)$ scan that solves the fixed point \eqref{eq:fast-fixed-point} with the inventory of product $i$ treated as unbounded, as in \eqref{eq:acc-beta}. Thus, constructing all type-specific sequences and bounds costs $O(mn\log n)$. 
The density list has $O(mn)$ entries, so sorting it costs $O(mn\log(mn))$. After cumulative weights and values are stored, each midpoint benefit query costs $O(\log(mn))$, and the $O(\log K)$ binary-search iterations cost $O(\log K\log(mn))$ in total. Recovering the final allocation takes $O(mn)$. Therefore, one call to Get Number costs
\(
O\!\left(mn\log(mn)+\log K\log(mn)\right).
\)
There is at most one call per product by Remark~\ref{rem:one-batch}, $n$ products, and \(\lceil \log_{1+\epsilon} K\rceil \) threshold values (Algorithm~\ref{alg:tau-m}). The overall running time is therefore
\(
O\!\left(
n \lceil \log_{1+\epsilon }K \rceil
\left[mn\log(mn)+\log K\log(mn)\right]
\right).
\)
\endgroup

\begingroup
\section{Additional Details and Results of Experiments} \label{app:experiments}
This appendix records the omitted implementation details for \textsc{da} and reports the \textsc{dap} results deferred from Section~\ref{sec:exp}.

For \textsc{da} benchmarks, {\small AS} has a publicly available codebase, which we use directly in our experiments. For {\small SG}, we implement the algorithm based on the paper and set the accuracy parameter to $\epsilon=0.01$. We use the same value, $\epsilon=0.01$, for {\small THR}. For {\small DP}, we set $\epsilon=0.1$ to improve computational efficiency. 

\subsection{DAP Experiments}\label{app:dap-experiments}

We next evaluate our algorithms for \textsc{dap} under both deterministic and stochastic horizons. We only describe the procedures that differ from \textsc{da}.

\paragraph{Instance generation. } The arrival probabilities of customer types are stationary over the horizon and are obtained by normalizing five independent $\operatorname{Exp}(1)$ draws.  The size of each type's
consideration set is drawn uniformly over
$[\lceil 0.3n\rceil,\lceil 0.8n\rceil]$, after which the consideration set is sampled uniformly among
subsets of that size. Following \cite{bai2022coordinated}, we reorder these coefficients increasingly across
products of decreasing revenue for the
first two customer types, creating a low-attraction/high-revenue pairing.
For every type, the outside-option weight is calibrated so that the no-purchase
probability is $0.01$ when its entire consideration set is offered.

Since our algorithm for stochastic $T$ works for arbitrary distributions, including non-IFR distributions, we generate the value of $T$ as follows:
for each instance and each value of $T_{\max}$, we independently
generate a four-scenario distribution for the random selling horizon.
Specifically, we draw
\(
H_\ell \stackrel{\mathrm{i.i.d.}}{\sim}
\operatorname{Unif}\{1,\ldots,T_{\max}\},
 \ell=1,\ldots,4,
\)
with replacement. We then independently draw
\(
Z_\ell \stackrel{\mathrm{i.i.d.}}{\sim} N(0,1),
W_\ell=\exp(Z_\ell),
 \ell=1,\ldots,4,
\)
and define the scenario probabilities as
\(
p_\ell
=
\frac{W_\ell}{\sum_{q=1}^{4}W_q}, \ell=1,\ldots,4.
\)
When several scenario draws produce the same horizon value, we aggregate
their probability masses. Consequently, the realized distribution of $T$
has between one and four support points contained in
$\{1,\ldots,T_{\max}\}$. We consider
$T_{\max}\in\{100,200,400,800,1600\}$. For deterministic $T$, we consider $T\in\{100,200,400,800,1600\}$.

\paragraph{Benchmarks.} 
For deterministic $T$, we compare our algorithm ({\small THR}, with $\epsilon=0.05$) with two methods of
\citet{bai2022coordinated}: Greedy Stocking and Randomized Customization
({\small GRA}) and Stocking with Multiple Estimates and Randomized
Customization ({\small MRA}). Let $\tau_k$ be the expected number of type-$k$
arrivals and let $y_{0k}$ denote the no-purchase mass in their exact
linear-programming surrogate $f(\vecc)$. Given $\kappa\in(0,1)$, their
approximate surrogate $f_{\rm app}^{\kappa}(\vecc)$ fixes
$y_{0k}=\kappa\tau_k$ for every type. $f_{\rm app}^{\kappa}(\vecc)$ is submodular in $\vecc$, and {\small GRA} applies greedy stocking
to $f_{\rm app}^{0.5}$. {\small MRA} repeats this greedy construction
for $\kappa\in\{0.1,0.2,\ldots,0.9\}$ and, as specified in \cite{bai2022coordinated}, selects
from the nine resulting inventories the one with the largest value under the
exact surrogate $f(\vecc)$. \cite{bai2022coordinated} report that MRA has the best empirical performance. We use the released implementation from the authors.

For stochastic \(T\), we adapt the algorithms of \citet{bai2022coordinated}, which were originally developed for deterministic \(T\). When determining the stocking level of each product, we replace \(T\) with \(\mu_T\) in the surrogate objective and apply the same optimization algorithms. For online assortment decisions, we consider both the independent-sampling randomized-rounding policy of \citet{bai2022coordinated}, with \(T\) replaced by \(\mu_T\), and the online greedy policy. We use {\small G-R} and {\small G-G} to denote the greedy stocking algorithm paired with the independent-sampling randomized-rounding policy and the online greedy policy, respectively. Similarly, {\small M-R} and {\small M-G} denote the corresponding variants of the multiple-estimates stocking algorithm. Following the same convention, {\small T-G} denotes our threshold-based stocking algorithm paired with the online greedy policy; this is the algorithm of Section~\ref{sec:dap-stochasticT}.

\paragraph{Results.}
Tables~\ref{tab:exp-dap-stochastic-summary} and
\ref{tab:exp-dap-deterministic-p001} demonstrate that our threshold-based
approach provides a strong balance between revenue performance and
computational efficiency. Under stochastic horizons, {\small T-G} achieves the
highest overall normalized revenue, averaging \(99.2\%\) across the 200
instances, and obtains the highest revenue in nine of the ten configurations.
It improves average normalized revenue by \(2.3\) percentage points over
{\small G-G} and by \(4.4\) percentage points over {\small M-G}. Although
{\small T-G} is slower than the simpler greedy benchmarks, its average runtime
remains moderate at \(68.1\) seconds and is approximately \(3.3\) times faster
than the multiple-estimates variants. Its performance also remains stable as
the problem size increases. We also observe that the online greedy policy performs much better than the independent rounding policy for all first-stage algorithms. 

The deterministic-horizon results further highlight the computational
advantage of our method. {\small THR} achieves an average normalized revenue of
\(99.8\%\), matching the best-performing benchmark, {\small MRA}, while
requiring only \(35.3\) seconds on average, compared with \(429.4\) seconds for
{\small MRA}. Thus, {\small THR} is more than \(12\) times faster while
maintaining the same average revenue performance. {\small THR} also
outperforms {\small GRA} on both averages, increasing average normalized
revenue from \(99.4\%\) to \(99.8\%\) while reducing average runtime by
approximately \(34\%\).

\begin{table}[htbp]
\centering
\scriptsize
\setlength{\tabcolsep}{3.2pt}
\renewcommand{\arraystretch}{0.92}
\resizebox{0.8\textwidth}{!}{%
\begin{tabular}{rrrrrrrrrrrr}
\toprule
& &
\multicolumn{5}{c}{Normalized revenue (\%)} &
\multicolumn{5}{c}{Mean runtime (seconds)}\\
\cmidrule(lr){3-7}\cmidrule(lr){8-12}
$n$ & $T_{\max}$ &
G-R & G-G & M-R & M-G & T-G &
G-R & G-G & M-R & M-G & T-G\\
\midrule
\multirow{5}{*}{50}
& 100
& 84.7 & 96.2 & 85.7 & 95.5 & 99.7
& 2.3 & 2.3 & 21.2 & 21.0 & 7.8\\
& 200
& 86.5 & 96.7 & 87.0 & 95.1 & 99.8
& 5.7 & 5.6 & 51.3 & 51.3 & 19.5\\
& 400
& 88.6 & 98.0 & 89.1 & 95.5 & 99.1
& 10.2 & 10.1 & 92.0 & 92.4 & 36.7\\
& 800
& 88.0 & 97.8 & 88.5 & 95.6 & 99.6
& 21.2 & 21.3 & 193.1 & 192.5 & 81.8\\
& 1600
& 89.0 & 97.8 & 89.4 & 95.9 & 99.9
& 42.1 & 42.1 & 383.9 & 382.2 & 178.7\\

\addlinespace

\multirow{5}{*}{100}
& 100
& 86.8 & 97.7 & 87.4 & 95.4 & 97.4
& 6.5 & 6.5 & 61.1 & 61.0 & 15.5\\
& 200
& 86.5 & 96.7 & 86.8 & 94.0 & 99.4
& 13.0 & 13.1 & 121.9 & 121.6 & 27.2\\
& 400
& 88.1 & 97.2 & 88.2 & 95.1 & 99.3
& 23.4 & 23.5 & 220.3 & 220.3 & 46.3\\
& 800
& 86.3 & 94.7 & 86.5 & 91.9 & 98.5
& 45.0 & 44.7 & 420.8 & 420.3 & 95.5\\
& 1600
& 87.4 & 96.3 & 87.6 & 94.1 & 99.1
& 78.5 & 78.5 & 711.9 & 748.6 & 172.4\\

\midrule
\multicolumn{2}{r}{\textbf{Average}}
& 87.2 & 96.9 & 87.6 & 94.8 & \textbf{99.2}
& 24.8 & 24.8 & 227.8 & 231.1 & 68.1\\
\bottomrule
\end{tabular}}
\caption{Stochastic-horizon \textsc{dap} performance.}
\label{tab:exp-dap-stochastic-summary}
\end{table}

\paragraph{Deterministic-horizon results.}
The deterministic $T$ results are summarized in Table~\ref{tab:exp-dap-deterministic-p001}.
\begin{table}[htbp]
\centering
\scriptsize
\setlength{\tabcolsep}{3.2pt}
\renewcommand{\arraystretch}{0.92}
\resizebox{0.63\textwidth}{!}{%
\begin{tabular}{rrrrrrrr}
\toprule
& & \multicolumn{3}{c}{Normalized revenue (\%)} &
\multicolumn{3}{c}{Mean runtime (seconds)}\\
\cmidrule(lr){3-5}\cmidrule(lr){6-8}
$n$ & $T$ &
THR & GRA & MRA &
THR & GRA & MRA\\
\midrule
\multirow{5}{*}{50}
& 100
& 99.8 & 99.2 & 99.7
& 2.0 & 4.5 & 41.4\\
& 200
& 99.8 & 99.2 & 99.7
& 5.0 & 10.0 & 91.1\\
& 400
& 100.0 & 99.3 & 99.8
& 11.0 & 20.0 & 182.8\\
& 800
& 100.0 & 99.5 & 99.8
& 28.0 & 38.6 & 349.9\\
& 1600
& 100.0 & 99.7 & 99.9
& 89.0 & 76.7 & 700.3\\

\addlinespace

\multirow{5}{*}{100}
& 100
& 99.8 & 98.8 & 99.7
& 5.3 & 12.6 & 118.8\\
& 200
& 99.8 & 99.3 & 99.7
& 10.8 & 24.0 & 229.6\\
& 400
& 99.6 & 99.5 & 99.8
& 22.7 & 46.1 & 438.8\\
& 800
& 99.2 & 99.6 & 99.9
& 56.8 & 102.4 & 958.4\\
& 1600
& 99.9 & 99.6 & 99.8
& 122.6 & 198.6 & 1{,}183.3\\

\midrule
\multicolumn{2}{r}{\textbf{Average}}
& \textbf{99.8} & 99.4 & \textbf{99.8}
& 35.3 & 53.3 & 429.4\\
\bottomrule
\end{tabular}}
\caption{Deterministic-horizon \textsc{dap} performance.}
\label{tab:exp-dap-deterministic-p001}
\end{table}
\endgroup

\end{APPENDICES}
\end{document}